  \documentstyle[12pt]{article}
\newtheorem{ournothing}{}[subsection]

\newtheorem{ourlemma}[ournothing]{Lemma}
\newtheorem{ourtheorem}[ournothing]{Theorem}

\newtheorem{ourcorollary}[ournothing]{Corollary}
\newtheorem{ourproposition}[ournothing]{Proposition}
\newtheorem{ourconjecture}[ournothing]{Conjecture}

\newcommand{\numero}[1]{
\addtocounter{section}{1}
\newpage
\begin{center}{\bf \thesection .\
#1\vspace{-.1in}}\end{center}
\setcounter{subsection}{0}
\setcounter{subsubsection}{0}
\setcounter{ournothing}{0}\indent}

\newcommand{\oldsubnumero}[1]{
\begin{center}
{\em #1}
\end{center}
}

\newcommand{\subnumero}[1]{
\addtocounter{subsection}{1}
\bigskip

\noindent{\bf \thesubsection \ }
{\sc #1}
\setcounter{subsubsection}{0}
\setcounter{ournothing}{0}
\newline
\indent}

\newenvironment{parag}{
\addtocounter{subsubsection}{1}
\begin{ournothing}\rm }{
\end{ournothing}}

\newcommand{\newparag}[1]{
\end{parag}
\begin{parag}
\label{#1}
}

\newcommand{\eop}{\hfill $/$\hspace*{-.1cm}$/$\hspace*{-.1cm}$/$\vspace{.1in}}

\newenvironment{lemma}{
\addtocounter{subsubsection}{1}
\begin{ourlemma}}{
\end{ourlemma}}

\newenvironment{corollary}{
\addtocounter{subsubsection}{1}
\begin{ourcorollary}}{
\end{ourcorollary}}

\newenvironment{theorem}{
\addtocounter{subsubsection}{1}
\begin{ourtheorem}}{
\end{ourtheorem}}

\newenvironment{proposition}{
\addtocounter{subsubsection}{1}
\begin{ourproposition}}{
\end{ourproposition}}

\newenvironment{conjecture}{
\addtocounter{subsubsection}{1}
\begin{ourconjecture}}{
\end{ourconjecture}}

\newenvironment{definition}{
\addtocounter{subsubsection}{1}
\begin{ournothing}\hspace*{-.2cm}---{\bf Definition:}\, \rm }{
\end{ournothing}}

\newenvironment{problem}{
\addtocounter{subsubsection}{1}
\begin{ournothing}\hspace*{-.2cm}---{\bf Problem:}\, \rm }{
\end{ournothing}}

\newenvironment{exercise}{
\addtocounter{subsubsection}{1}
\begin{ournothing}\hspace*{-.2cm}---{\bf Exercise:}\, \rm }{
\end{ournothing}}

\newenvironment{example}{
\addtocounter{subsubsection}{1}
\begin{ournothing}\hspace*{-.2cm}---{\bf Example:}\, \rm }{
\end{ournothing}}

\newenvironment{remark}{
\addtocounter{subsubsection}{1}
\begin{ournothing}\hspace*{-.2cm}---{\bf Remark:}\,  \rm }{
\end{ournothing}}

\newenvironment{hypothesis}{
\addtocounter{subsubsection}{1}
\begin{ournothing}\hspace*{-.2cm}---{\bf Hypothesis:}\,  \rm }{
\end{ournothing}}

\newenvironment{hypotheses}{
\addtocounter{subsubsection}{1}
\begin{ournothing}\hspace*{-.2cm}---{\bf Hypotheses:}\,  \rm }{
\end{ournothing}}

\newcommand{\Gm}{{\bf G}_m}
\newcommand{\af}{{\bf A}^1}

\newcommand{\cc}{{\bf C}}
\newcommand{\rr}{{\bf R}}
\newcommand{\qq}{{\bf Q}}

\newcommand{\zz}{{\bf Z}}
\newcommand{\pp}{{\bf P}}

\newcommand{\Ee}{{\cal E}}
\newcommand{\Ff}{{\cal F}}
\newcommand{\Gg}{{\cal G}}

\newcommand{\Oo}{{\cal O}}

\newcommand{\Bb}{{\cal B}}
\newcommand{\Mm}{{\cal M}}

\newcommand{\Uu}{{\cal U}}
\newcommand{\Ll}{{\cal L}}

\newcommand{\Pp}{{\cal P}}

\newcommand{\tworightarrows}{\stackrel{\displaystyle \rightarrow}{\rightarrow}}

\newcommand{\ngr}{n^{\rm gr}} 

\begin{document}

\section*{Algebraic aspects of higher nonabelian Hodge theory}

\noindent
Carlos Simpson\newline

\begin{center}
{\bf Table of Contents}
\end{center}

\noindent
{\bf 1. Introduction}---p. \pageref{intropage}
\newline
{\small
Principal questions;  plan of paper.
}

\noindent
{\bf 2. Varieties with abelian fundamental group}---p. \pageref{abelianpage}
\newline
{\small
The large size of $\pi _3$ in the presence of a free abelian $\pi _1$; existence
of this problem for algebraic varieties; a nonabelian cohomological alternative.
}

\noindent
{\bf 3. Nonabelian cohomology}---p. \pageref{nonabelianpage}
\newline
{\small
Motivation, historical remarks; $n$-stacks; nonabelian cohomology; the main type
of problem; examples.  }

\noindent
{\bf 4. Zoology}---p. \pageref{zoologypage}
\newline
{\small
Realms; index of realms; closure properties; nonabelian cohomology of a constant
stack; directories of Serre classes; the Postnikov-type argument.
}

\noindent
{\bf 5. Cartesian families and base change}---p. \pageref{cartesianpage}
\newline
{\small
A fibrant replacement for $nCAT$; relative nonabelian cohomology and
base-change;
cartesian families; construction of the arrow family.
}

\noindent
{\bf 6. Very presentable $n$-stacks}---p. \pageref{presentablepage}
\newline
{\small
(Very) presentable $n$-stacks; vector sheaves, the Breen calculations;
representability of shape; very presentable suspension.
}

\noindent
{\bf 7. Geometric $n$-stacks}---p. \pageref{geometricpage}
\newline
{\small
Abstract and concrete approaches to geometric $n$-stacks; doing geometry
therein; a criterion for geometricity; a semicontinuity-type result;
geometricity of representing objects.
}

\noindent
{\bf 8. Formal groupoids of smooth type}---p. \pageref{formalpage}
\newline
{\small
Review of formal categories of smooth type; calculation of cohomology via the de
Rham complex; a construction using blow-up; the Hodge filtration; morphisms of
smooth type.}

\noindent
{\bf 9. Formal categories related to Hodge theory}---p. \pageref{hodgepage}
\newline
{\small
Examples of formal categories giving de Rham theory, Dolbeault theory,
the Gauss-Manin connection, logarithmic de Rham theory, the Hodge filtration,
Griffiths transversality, regularity.
}

\noindent
{\bf 10. Presentability and geometricity results}---p. \pageref{resultspage}
\newline
{\small
The main results concerning the properties of the nonabelian cohomology
$\underline{Hom}(\Ff /\Ee , T)$ for $\Ff \rightarrow \Ee$ a projective morphism
of smooth type between formal categories, and $T$ a coefficient $n$-stack
satisfying various properties such as very-presentability or geometricity.
}

\numero{Introduction}
\label{intropage}

Hodge theory involves a large number of analytical aspects, but also a certain
number of algebraic aspects pertaining notably to the algebraic de Rham
cohomology \cite{Grothendieck}(ii), \cite{Hartshorne}. Among these latter are:
the Hodge filtration; the Gauss-Manin connection in the case of a family of
varieties; Griffiths transversality for the Hodge filtration with respect to the
Gauss-Manin connection; and regular singularities of the Gauss-Manin connection
in case of a degenerating family.

The {\em degree $1$ nonabelian de Rham cohomology} (of a smooth projective
variety $X$ in characteristic zero) is the moduli stack of  principal
$G$-bundles
with integrable connection which we can denote by
$$
\underline{H}^1(X_{DR}, G)= \underline{Hom}( X_{DR}, K(G,1)).
$$
By the {\em higher nonabelian de Rham cohomology} of $X$ we mean $n$-stacks
of the form
$$
\underline{Hom}( X_{DR}, T)
$$
for coefficient $n$-stacks $T$ with $\pi _0(T)=\ast$, $\pi _1(T)=G$ an affine
algebraic group, and $\pi _i(T)$ represented by finite dimensional vector spaces
for $i\geq 2$.

The purpose of this paper is to treat the above-mentionned algebraic aspects of
Hodge theory, for the higher nonabelian de Rham cohomology, generalizing from
the degree $1$ case which was treated in \cite{Simpson}(vi). The main
observation
is that all of the algebraic aspects cited above can be interpreted in terms of
certain ``formal categories'' \cite{Grothendieck} \cite{Berthelot}
\cite{Illusie},
starting from the formal category for $X_{DR}$ whose object-object is $X$ and
whose morphism-object is the formal completion of the diagonal in $X\times
X$. The
formal categories (and morphisms of such) necessary for the Hodge
filtration, the
Gauss-Manin connection, Griffiths transversality, and regular singularities,
are reviewed in \S 9 below.

Here are two examples:
first, if $X$ is a smooth projective variety then the Hodge filtration
comes from
a morphism of formal categories
$$
X_{Hod} \rightarrow \af
$$
where the preimage of $\Gm \subset \af$ is isomorphic to $X_{DR}\times \Gm$,
and where the fiber over $0\in \af$ is
$$
X_{Dol} := K(\widehat{TX}/X , 1)
$$
(i.e. the relative classifying stack for the formal completion of $TX$  along
the zero-section).  Second, if $X\rightarrow S$ is a smooth projective
morphism then the Gauss-Manin connection comes from the morphism of formal
categories
$$
X_{DR} \rightarrow S_{DR}.
$$

If $\Ff \rightarrow \Ee$ is one of the morphisms of formal categories which
enters into the above aspects of Hodge theory, and if $T$ is a coefficient
$n$-stack, then we obtain the relative nonabelian cohomology
$n$-stack
$$
\underline{Hom}(\Ff /\Ee , T) \rightarrow \Ee .
$$
Much of the paper is concerned with the following basic questions:
\newline
(A) \, how to interpret the above, for example as giving a morphism $\Ee
\rightarrow n\underline{STACK}$;
\newline
(B) \, how to interpret the above as a functor in the variable $T$, giving a
notion of ``shape'' denoted
$$
{\bf Shape}(\Ff /\Ee )  : \Ee \rightarrow \underline{Hom}(n\underline{STACK},
n\underline{STACK});
$$
and
\newline
(C) \, what are the properties of $\underline{Hom}(\Ff /\Ee , T)$ as a function
of the properties of $T$?

Note that in (B) the stack $T$ is naturally brought to vary in a family
parametrized  by a base scheme, and question (C) concerns mostly the
properties of the resulting family of nonabelian cohomology stacks, as a function
of the properties of the family of stacks $T$. Our results in this direction are
given in the last chapter \S 10.

Here is a simplified version of the results of \S 10.
Leaving aside the question of variation of $T$ in a family, if
$T$ is a coefficient $n$-stack chosen as described at the start (connected, with
$\pi _1$ an affine algebraic group and $\pi _i$ vector spaces for $i\geq 2$),
and if we base-change by the object-scheme $Z\rightarrow \Ee$ of the base formal
category, then the resulting nonabelian cohomology $n$-stack
$\underline{Hom}(\Ff \times _{\Ee} Z /Z, T) \rightarrow Z$ is a {\em
geometric very presentable $n$-stack}. The descent of this $n$-stack down
to $\Ee$
represents the ``action'' of the morphisms of $\Ee$ on it.

Applying these results to $X_{Hod}\rightarrow \af$ we obtain the {\em Hodge
filtration} on the nonabelian de Rham cohomology,
a $\Gm$-equivariant morphism
$$
\underline{Hom}(X_{Hod}/\af , T) \rightarrow \af .
$$
Considered as a functor in the variable $T$, this gives the {\em Hodge
filtration on the de Rham shape}
$$
{\bf Shape}(X_{Hod}/\af ): \af / \Gm \rightarrow
\underline{Hom}(n\underline{STACK}, n\underline{STACK}).
$$
If $X$ is simply connected, then the shape (restricted to ``$1$-connected very
presentable $n$-stacks'' cf \S 6) is representable by a $1$-connected very
presentable geometric $n$-stack
$$
{\bf rep}(X_{Hod}/\af )\rightarrow \af ,
$$
which together with its $\Gm$-equivariance, is the Hodge
filtration on the full complexified homotopy type of $X$.

Using the example $X_{DR}\rightarrow S_{DR}$ we obtain in a similar way, the
Gauss-Manin connection on the de Rham shape, and in the simply connected case,
the Gauss-Manin connection on the de Rham homotopy type of the family $X/S$
(this adds the higher levels of homotopy coherence to the result of Navarro
Aznar
\cite{Navarro-Aznar}).

Using other examples of morphisms of formal categories (see \S 9) we obtain
the other algebraic aspects of Hodge theory such as
Griffiths transversality, the Kodaira-Spencer classes for Dolbeault
cohomology, regular singularity of the Gauss-Manin connection (again extending
\cite{Navarro-Aznar}), logarithmic de Rham cohomology, etc.

\bigskip

The plan of the paper is as follows. In \S 2, we start
with a digression on the example of the higher homotopy of varieties with
abelian
fundamental group. This digression is intended to show the need for a higher
nonabelian cohomological formalism, by showing that direct consideration of the
higher homotopy groups, in the case of non-simply connected varieties, seems
problematic.

In \S 3, we define higher nonabelian cohomology, and give some examples. We
treat question (B).

In \S 4 entitled ``zoology'', we consider some questions related to the
notion of
a ``property'' of an $n$-stack (the notion of ``property''  comes into play
in question (C) in the above list). Notably, we define the notion of {\em realm}
as a saturated full sub-$n+1$-stack $\underline{R} \subset n\underline{STACK}$.
A realm should be thought of as the sub-$n+1$-stack of $n$-stacks having a given
property. We discuss various things about realms, and at the end of \S 4 we give
the main ``Postnikov reduction'' argument which will be used in \S 10.

Possibly useful in \S 4 is the ``index'' \ref{index1} of the various realms
which will come into play later on.

In \S 5, we discuss base-change and the notion of ``cartesian family''. This is
intended to bring some clarification to question (A) in the above list.
Unfortunately, it requires development of some new machinery for $n$-categories.
This new machinery should be useful in many contexts: we construct a canonical
fibrant replacement denoted $nFAM$ for the $n+1$-category $nCAT$ of all
$n$-categories which had been constructed in \cite{Simpson}(x). On the
downside, the discussion is rather technical and also incomplete, and the reader
will probably complain that we have shedded darkness rather than light on
question (A).

In \S 6, we begin our treatment of concrete issues surrounding the above
questions, by introducing the first main realm: that of {\em very presentable
$n$-stacks}. Our discussion here is really just a review of things which occur
in previous preprints such as \cite{Simpson}(vii) or (xii). In turn, most
of what we review from these preprints comes from Auslander
\cite{Auslander}, Hirschowitz \cite{Hirschowitz}, and Breen \cite{Breen}.

In \S 7 we continue as in the previous section, this time discussing the second
main realm: that of {\em geometric $n$-stacks}. Again, this starts with a review
of the definition from \cite{Simpson}(ix). We go on, however,
to  discuss some new things such as the type of geometry which one could expect
to do with geometric $n$-stacks, and at the end a criterion for being a
geometric $n$-stack.

In \S 8, we discuss formal categories of smooth type. The main work here is to
explain why the cohomology of a formal category is calculated by a de Rham
complex. This is well-known since Berthelot \cite{Berthelot} and Illusie
\cite{Illusie}, but we give an argument anyway.

In \S 9, we describe the formal categories and morphisms of formal categories
which enter into the aspects of Hodge theory listed at the beginning of the
introduction.

In \S 10 we give our main results concerning question (C) in the above list.
These results say that for appropriate types of coefficient $n$-stacks $T$,
the resulting nonabelian cohomology of a formal category, is very presentable
(cf \S 6) and sometimes geometric (cf \S 7).

{\em Acknowledgements:} I would like to thank L. Katzarkov for organizing the
conference and lectures on which this paper is (very loosely) based. I would
also like to thank him and T. Pantev for many useful questions, suggestions and
remarks; and to thank A. Hirschowitz for continuing discussions about
$n$-categories and $n$-stacks. From a long time ago, I would like to thank again
P. Deligne for communicating his notion of $\lambda$-connection which led to the
Hodge filtration for degree $1$ nonabelian cohomology; to thank J. Propp for
one day mentionning the words ``shape theory''; and to thank V. Navarro
Aznar for
discussions about rational homotopy theory and his Gauss-Manin connection.

\numero{Varieties with abelian fundamental group}
\label{abelianpage}

Rational homotopy theory for a space $X$ has classically  treated either the
pro-nilpotent completion
$$
\pi _1^{\rm nil}(X)  := {\rm prolim} \pi _1(X) /\Gamma _r \pi _1(X),
$$
or else the higher homotopy groups $\pi _i(X)$ in the case where $X$ is simply
connected. These can be combined together in case the space $X$ is {\em
nilpotent}, that is if the action of $\pi _1$ on the $\pi _i$ is nilpotent.
However, this latter condition doesn't usually hold.

To illustrate what can happen in the non-nilpotent case, we look at
the case where $\pi _1$ is an abelian group.  The action of $\pi _1$ on the $\pi
_i$ need not be nilpotent, and classical rational homotopy theory isn't really
adequate to deal with the situation. One of the reasons for this is that the
higher $\pi _i$ get very big. We illustrate this problem by looking carefully at
$\pi _2$ and $\pi _3$ in the case where $\pi _1(X)=A$ is a {\em free abelian
group of positive rank}.

This section is only intended to point out the type of problem which leads to a
need for higher nonabelian cohomology. We don't necessarily solve this problem
in the rest of the paper, though. The first two subsections are
mathematically independant from the rest of the paper, whereas in \S 2.3 we give
a brief introduction to ``nonabelian cohomology'' with a few examples.

Everything
in \S\S 2.1, 2.2 is well-known in one way or another, and we don't endeavour to
give references. Many of the necessary references can be found in the
bibliography of \cite{Simpson}(ii).

\subnumero{On the size of $\pi _3$ in the presence of a free abelian $\pi _1$}

Assume that $X$ is a CW complex with finitely many cells, with $\pi _1(X)=A$ a
free abelian group of nonzero (finite) rank.

For simplicitly,
tensor everything with $\qq$. Recall that $\pi _1(X)$ acts on $\pi _i(X)$, on
$H^i(\tilde{X} ,\qq )$, etc. This gives an action of the (commutative, in this
case) algebra $\qq [A]$ on $H_i(\tilde{X}, \qq )$ and on $\pi _i(X)\otimes
\qq $.

If $X$ is a CW complex with finitely many cells, then the chain complex
$C_{\cdot}(\tilde{X}, \qq )$ (with respect to the lifted cell decomposition of
$\tilde{X}$) is a complex of finitely generated free $\qq [A]$-modules. Since
$\qq [A]$ is noetherian, this implies that the cohomology of this complex
consists of $\qq [A]$-modules of finite type. Thus $H_i(\tilde{X}, \qq )$ is a
finitely generated $\qq [A]$-module.

Let $M_B(X, \Gm ):= Hom (\pi _1(X), \Gm )$. It is a product of tori defined
over $\qq$, and we have
$$
M_B(X, \Gm )= Spec (\qq [A]).
$$

If $\rho : \pi _1(X)\rightarrow \Gm $ is a rank one representation (defined over
$\cc$, say), let $L_{\rho}$ denote the corresponding local system of rank one
$\cc$-vector spaces on $X$, and let $V_{\rho}$ denote the corresponding $\qq
[A]$-module (in fact it is a $\cc [A]$-module).  Then
$$
H_i(X, L_{\rho})= H_i (C_{\cdot}(\tilde{X}, \qq )\otimes _{\qq [A]} V_{\rho}).
$$
There is a connected dense Zariski open subset $U\subset M_B(X,\Gm )$ such that
for $\rho \in U$, we have
$$
H_i(X, L_{\rho})= H_i(\tilde{X}, \qq )\otimes _{\qq [A]} V_{\rho}.
$$
Note that $U$
is connected because of the fact that we have assumed that $A$ is without
torsion.

We may assume that the function $h_i(X, L_{\rho })$ is
constant on $U$, and the coherent sheaf associated to the $\qq [A]$-module
$ H_i(\tilde{X}, \qq )$ is locally free over $U$ of rank equal to the value of
this function.
Thus we may speak of the ``generic dimension'' which we denote $h_i(X, L_{\rm
gen})$. The ring $\qq [A]$ is integral and we have
$$
h_i(X, L_{\rm gen})= dim _{Frac (\qq [A])} H_i(\tilde{X}, \qq )\otimes
_{\qq [A]}
Frac (\qq [A]).
$$

Cohomology with local coefficients is dual to homology with local coefficients:
$$
H^i(X, L_{\rho})= H_i(X, L_{\rho}^{\ast})^{\ast}
$$
(here we are taking duals of complex vector spaces). The dual of a generic
representation $\rho$ is again a generic representation. Thus we may write
$$
h^i(X, L_{\rm gen})= h_i(X, L_{\rm gen}).
$$
Let $b$ be the largest integer such that there is an injection of $\qq
[A]$-modules
$$
\qq [A]^b \hookrightarrow H_i (\tilde{X}, \qq ).
$$
Then $b= h_i(X, L_{\rm gen})$. Fix such an injection.

We now look at $\pi _i(X)$, $i=2,3$. Note that $\pi _2(X)= H_2(\tilde{X}, \zz )$
so
$$
\pi _2(X)\otimes \qq = H_2(\tilde{X}, \qq )
$$
is a $\qq [A]$-module of finite type ({\em Note:} in a more general
situation this
statement would require that the group ring $\qq [\pi _1]$ be noetherian).
On the other hand, the coproduct in homology is the first morphism
$\xi$  in an exact sequence $(\ast )$
$$
H_4(\tilde{X} , \qq )\stackrel{\xi}{\rightarrow} H_2(\tilde{X} , \qq ) \otimes
_{\qq } H_2(\tilde{X} , \qq ) \rightarrow \pi _3(X)\otimes \qq \rightarrow
H_3(X,
\qq ). $$
Existence and functoriality of this exact sequence for finite subcomplexes of
$\tilde{X}$ implies existence of the exact sequence for an infinite complex such
as $\tilde{X}$. Note that this exact sequence is compatible with the action of
$\qq [A]$, with  diagonal action of $A$ on the tensor product.

\begin{hypothesis}
\label{hfourgen}
Suppose that $H_4(X, L_{\rm gen})=0$. This means
that $H_4(\tilde{X}, \qq )$ is a torsion $\qq [A]$-module.
\end{hypothesis}

Recall from above
that we have fixed an injection from a free $\qq [A]$-module of rank
$b=h_2(X, L_{\rm gen})$ into $H_2(\tilde{X},\qq )$. This gives an injection
$$
i: \qq [A] ^b \otimes _{\qq } \qq [A]^b \hookrightarrow H_2(\tilde{X} , \qq
)\otimes _{\qq} H_2(\tilde{X}, \qq ),
$$
compatible with the diagonal action of $A$. Now
$\qq [A] ^b \otimes _{\qq } \qq
[A]^b$ with its diagonal action of $A$ is a free $\qq [A]$-module (of infinite
rank), in particular it is torsion-free so the image of $i$ doesn't meet the
image of $\xi$ (here we use our hypothesis \ref{hfourgen}). Therefore we
obtain an injection of $\qq [A]$-modules
$$
\qq [A] ^b \otimes _{\qq } \qq [A]^b \hookrightarrow \pi _3(X)\otimes \qq .
$$
Assuming that $b>0$ and also of course that $A$ is an abelian group of nonzero
rank, we find that $\pi _3(X)\otimes \qq $ contains a free $\qq [A]$-module of
infinite rank. In particular, it cannot be of finite type over $\qq [A]$. This
could pose a problem in an algebraic-geometric approach where one can only
adequately deal with things of finite type.

Aside from this infinite-dimensionality, (or partly, because of it) there is
another type of problem which occurs.
The exact sequence $(\ast )$ sets up $\pi _3(X)\otimes \qq $ as a curious
mixture of a $\qq [A]$-module of finite type (the image in $H_3(\tilde{X}, \qq
)$) with something which is itself a $\qq [A \times A]= \qq [A] \otimes
_{\qq}\qq
[A]$-module of finite type divided out by a $\qq [A]$-module of finite type.
Unfortunately, there doesn't seem to be any way to provide $\pi _3(X)\otimes \qq
$ with a natural $\qq [A\times A]$-structure making it into something of finite
type. On the other hand, it is probably possible, but complicated, to develop a
theory taking into account this mixture of different types of structures. This
behavior continues in $\pi _i(X)$ for larger values of $i$, where one meets (via
the Curtis spectral sequence for example) pieces which naturally should be seen
as $\qq [A \times \ldots \times A]$-modules of finite type for various different
numbers of factors.

This mixing of actions of various $A\times \ldots \times A$ poses a problem if
one wants to take completions or other types of tensor products over $\qq [A]$.
If we have a $\qq [A]$-algebra $B$ (e.g. the completion at some ideal) and
try to
tensor everything with $B$ over $\qq [A]$, then what becomes of the factor
$H_2(\tilde{X}, \qq )\otimes _{\qq } H_2(\tilde{X}, \qq )$? It could either be
tensored once with $B$ over $\qq [A]$, but then it wouldn't really represent the
``completion''---for example our term $\qq [A]^b\otimes _{\qq } \qq [A]^b$
becomes something approximately (although not exactly) like $\qq [A] ^b \otimes
_{\qq} B^b$; or else we might get a term of the form  $$ (H_2(\tilde{X}, \qq
)\otimes _{\qq [A]}B)\otimes _{\qq} (H_2(\tilde{X}, \qq )\otimes _{\qq [A]}B)
$$
in the exact sequence $(\ast )$,
but it seems difficult to describe the operation on $\pi _3(X)\otimes \qq $
which corresponds to this latter case.

In the comparison between Betti and de Rham cohomology, we eventually need to do
a tensor product as above with $B$ equal to the ring of holomorphic functions on
$M_B(X, \Gm )(\cc )$. Thus, the problem of what this ``tensoring''
should mean on $\pi _3$ seems to block any program of trying to compare
Betti and
algebraic de Rham theories.

This type of problem---as it occurs when one wants to do $\ell$-adic
completions---was pointed out by  Artin and Mazur \cite{ArtinMazur}, and also
Bousfield \cite{Bousfield}.  It can occur, for example, when $X$ is obtained by
attaching a bouquet of $2$-spheres to a torus.

One way of getting around this problem would be
to treat the rational homotopy theory of the universal covering $\tilde{X}$
equivariantly with respect to the action of $\pi _1$ (see for example Brown and
Szczarba \cite{BrownSzczarba}, also Gomez-Tato \cite{Gomez-Tato} and
others). This
assumes, of course, that one can get a handle on $\pi _1(X)$ to start with. For
this reason, the approach of Brown and Szczarba seems somewhat unsuited for
algebraic geometry e.g. for an algebraic de Rham version of the theory.

Instead of tackling the problem directly, we shall get around it by
looking at ``nonabelian cohomology'' and shape. This is one of the motivations
for what we do in subsequent chapters; we'll give a brief look in \S 2.3 below.

\subnumero{Smooth projective varieties with abelian $\pi _1$}

This subsection is devoted to pointing out that the problem alluded
to above effectively does occur, and is even in some sense the ``generic''
situation, for smooth projective varieties with abelian fundamental group.
Again, it should be stressed that all of the results herein are well-known by
various authors and we don't give references.

Suppose $X$ is a smooth projective variety over $\cc$. Then we have several
other ways of computing cohomology with local system coefficients. Let
$M_{DR}(X,
\Gm )$ (resp. $M_{Dol}(X, \Gm )$) be the moduli space of $(\Ll , \nabla )$ where
$\Ll \in Pic ^{\tau } (X)$ and $\nabla $ is an integrable connection (resp. of
$(\Ee  , \theta )$ where $\Ee \in Pic ^{\tau} (X)$ and $\theta \in H^0(X, \Omega
^1_X)$). We have the {\em de Rham cohomology}
$$
H^i_{DR}(X, (\Ll , \nabla )):= {\bf H}^i(\Ll \stackrel{\nabla}{\rightarrow}
\Omega ^1_X \otimes \Ll \stackrel{\nabla}{\rightarrow} \ldots )
$$
(resp. the {\em Dolbeault cohomology}
$$
H^i_{Dol}(X, (\Ee , \theta )):= {\bf H}^i(\Ee \stackrel{\theta }{\rightarrow}
\Omega ^1_X \otimes \Ee \stackrel{\theta}{\rightarrow} \ldots ) \; )
$$
and these calculate the Betti cohomology. If $\rho \in M_{B}(X, \Gm )$,  let
$(\Ll , \nabla )_{\rho}$ be the corresponding line bundle with integrable
connection, and let $(\Ee , \theta )_{\rho}$ be the corresponding Higgs bundle.
Then
$$
H^i(X, L_{\rho})\cong H^i_{DR}(X, (\Ll , \nabla )_{\rho}) \cong H^i_{Dol}(X,
(\Ee , \theta )).
$$

Under the homeomorphism $M_{B}(X, \Gm )\cong M_{Dol}(X, \Gm )$, the open set
$U$ (which we now call $U_B$) where the cohomology dimension equals the generic
dimension, corresponds to a Zariski open set $U_{Dol} \subset M_{Dol}(X, \Gm )$
where the Dolbeault cohomology dimension equals the generic dimension.
Note that the set of Higgs bundles of the form $(\Oo , \alpha )$ goes to a
totally real subset of $M_B(X, \Gm )$ and in particular it must intersect $U_B$.
Therefore, $U_{Dol}$ contains a point of the form $(\Oo , \alpha )$ with
$\alpha$ a generic section of $\Omega ^1_X$.

We  now {\em assume that the
Albanese map $X\rightarrow Alb(X)$ is finite}. This implies that a generic
section $\alpha \in H^0(X, \Omega ^1_X)$ has only isolated zeros on $X$; for
any curve in the zero set of $\alpha$ maps to a positive-dimensional abelian
subvariety of $Alb(X)$ on which $\alpha$ restricts to zero, and  a generic
$\alpha$ will never do this, so no such curve can exist.

The fact that $\alpha$ has isolated zeros implies that the complex
$$
\Oo _X \stackrel{\wedge \alpha}{\rightarrow} \Omega ^1_X
\stackrel{\wedge \alpha}{\rightarrow} \ldots
\stackrel{\wedge \alpha}{\rightarrow} \Omega ^n_X
$$
is exact except at the last place, where it has a cokernel which we denote by
$S(\alpha )$ which is a skyscraper sheaf supported on the zero scheme of
$\alpha$. If $\alpha$ has at least one zero then the skyscraper sheaf is
nontrivial. We get
$$
H^i_{Dol}(X, (\Oo , \alpha )) = 0,\;\;\;\; i<n,
$$
and
$$
H^n_{Dol}(X, (\Oo , \alpha )) = H^0(S(\alpha )).
$$
This gives the following statement.

\begin{lemma}
Suppose that $X\rightarrow Alb(X)$ is a finite morphism, and let $n:= dim (X)$.
Then
$$
h^i(X, L_{\rm gen})=0,\;\;\;\; i < n.
$$
If, furthermore, every generic section of $H^0(X, \Omega ^1_X)$ vanishes
somewhere
on $X$, then $h^n(X, L_{\rm gen})>0$.
\end{lemma}
\eop

\begin{remark} Under the condition of the first part of the lemma, the
topological Euler characteristic of $X$ (which doesn't depend on the choice of
rank one local system) is equal to $h^n(X, L_{\rm gen})$. Thus, if this Euler
characteristic is different from zero then the latter generic dimension doesn't
vanish.
\end{remark}

We now look at the case of a surface $X$ and plug this back into the previous
discussion.

\begin{corollary}
Suppose $X$ is a smooth projective surface with $\pi _1(X)=A$ free of
finite rank,
such that the Albanese map of $X$ is finite, and such that
the topological Euler characteristic of $X$ is nonzero.
Then $b= h^2(X, L_{\rm gen})> 0$, and $h^4(X, L_{\rm gen})=0$, so by the
previous discussion, $\pi _3(X)\otimes \qq$ contains a free $\qq [A]$-module of
infinite rank.
\end{corollary}
\eop

The hypotheses of the corollary will hold, for
example, if $X$ is a complete intersection of hyperplane sections in $Alb(X)$.

The previous corollary is the sense in which we say that the ``generic''
behavior for a smooth projective surface with nontrivial free abelian $\pi _1$,
is to have an infinite-rank free $\qq [\pi _1]$-module contained in $\pi _3$.

\subnumero{A different way of proceeding}

From the above examples, it looks like it might not be so easy to look directly
at the higher homotopy groups or more generally the higher homotopy type, of a
non-simply connected variety. Instead, we suggest to take a different approach
based on {\em nonabelian cohomology}. Before getting to the actual definition
and technical discussion in the succeeding chapters, we give a brief
introduction to what this looks like in the context of abelian fundamental
groups.

Let $X_B$ denote the topological space $X^{\rm top}$.
We will be looking at a nonabelian cohomology object denoted
$\underline{Hom}(X_B, T)$. The ``coefficients'' are an {\em $n$-stack} $T$.
Without yet giving the definition, roughly speaking this is something whose
homotopy groups (occuring in degrees $0\leq i \leq  n$) will be sheaves of
groups
on the site $Sch /\cc$. For the present discussion we restrict to the case $\pi
_0(T)=\ast$ and where the $\pi _i(T)$ are sheaves represented by group schemes.
We then obtain concrete cohomological data for $X_B$ with coefficients in
$T$, by
taking the sheaf of sets $\pi _0\underline{Hom}(X_B, T)$.

Part of a map $X_B\rightarrow T$ involves a representation of $\pi _1(X_B)$ in
the group $\pi _1(T)$. If, as in the present chapter, $\pi _1(X_B) = A$ is
abelian, then for all practical purposes we might as well assume that $\pi
_1(T)$
is abelian. Assume furthermore that it is a reductive group, and even
indecomposable: this says $\pi _1(T)=\Gm$ is the multiplicative group scheme.
For the rest of our discussion in this subsection, make this hypothesis.

Now assume furthermore that $\pi _i(T)$ are finite dimensional complex vector
spaces (with $\Gm$ action, representing the action of $\pi _1$ on the $\pi _i$
for $T$). To be precise this means that the sheaves $\pi _i(T)$ on $Sch
/\cc$ are
represented by finite dimensional complex vector spaces.

\begin{parag}
\label{caseA}
In the simplest case (which really has $n=1$)
$$
\pi _1(T) = \Gm ,\;\;\;\; \pi _i (T) = 0, \;\; i\geq 2.
$$
Then,
$$
\pi _0\underline{Hom}(X_B, T)= M_B(X, \Gm ) = Spec (\cc [A])
$$
is the moduli space of rank
one local systems on $X_B$ which occurs above.

\newparag{caseB}
Take the next nontrivial case:
$$
\pi _1(T) = \Gm ,\;\;\;\; \pi _n (T) = V, \;\;\;  \pi _i (T) = 0,\;\; 2 \leq i
\neq n.
$$
Here $V$ is the standard $1$-dimensional representation of $\Gm$. In this case,
we have a map
$$
\pi _0\underline{Hom}(X_B, T)\rightarrow  M_B(X, \Gm ),
$$
and the fiber over a point $\rho$ is the cohomology vector space $H^n(X_B,
L_{\rho})$. There may be local problems around values of $\rho$ for which the
dimensions jump, making it so that $\pi _0\underline{Hom}(X_B, T)$ is not a
scheme; however, the inverse image of the open subset where the cohomology
dimension attains its generic value is indeed (a sheaf represented by) a scheme,
and it is just the tautological vector bundle over this open set whose fiber is
the cohomology.

We can write heuristically
$$
\pi _0\underline{Hom}(X_B, T) = \{ (\rho , \eta ): \;\; \rho \in M_B(X, \Gm
),\;\; \eta \in H^n(X_B, L_{\rho }) \} .
$$

The reader who at this point is interested in just what type of object
$\pi _0\underline{Hom}(X_B, T)$ is, is refered to the body of the paper:
this $0$-stack is ``very presentable'' in the terminology of \S 6,
and more generally the $n$-stack $\underline{Hom}(X_B, T)$ is very presentable.
By keeping the full $n$-stack here  we also gain the additional property that
it is a {\em geometric} $n$-stack, which is a higher analogue of the notion of
Artin algebraic stack---see \S 7 below.

\newparag{caseC}
Finally, one can incorporate nontrivial Whitehead products into the coefficient
stack $T$, which basically means that we take into account cup products and the
like. For example, fix $m \leq n/2$ and define an $n$-stack $T$ with
$$
\pi _1(T) = \Gm ,\;\;\;\; \pi _m (T) = V, \;\;\;  \pi _{2m-1} (T) =
V^{\otimes 2},
$$
and the remaining homotopy groups vanishing. Here we assume that there is a
nontrivial Whitehead product $\pi _m \times \pi _m \rightarrow \pi _{2m-1}$
(in order for such to exist, the $\Gm$ action on $\pi _{2m-1}$ has to be the
tensor square of that on $\pi _m$).

Let $T'$ be the $n$-stack of the previous
example with $\pi _1(T')=\Gm$ and $\pi _m(T')=V$ but all others vanishing.
For cohomology with coefficients in $T'$, the previous discussion
applies; in particular we denote the points of
$\pi _0\underline{Hom}(X_B, T')$ by $(\rho , \eta )$ following the expression
given above.

We have a morphism $T\rightarrow T'$ (``truncation'') which induces a
sequence of morphisms
$$
\pi _0\underline{Hom}(X_B, T)\rightarrow
\pi _0\underline{Hom}(X_B, T')\rightarrow  M_B(X, \Gm ).
$$
Roughly speaking this morphism maps
$\pi _0\underline{Hom}(X_B, T)$ to the subset of points $(\rho , \eta )$
such that $\eta \cup \eta =0$ in $H^{2m}(X_B , L_{\rho}^{\otimes 2})$ and a
little bit more precisely (but still somewhat heuristically) we can write
$$
\pi _0\underline{Hom}(X_B, T) =
$$
$$
\{ (\rho , \eta , \varphi ): \; \rho \in M_B(X,
\Gm ),\; \eta \in H^n(X_B, L_{\rho }), \varphi \in C^{2m-1}(X_B,
L_{\rho}^{\otimes 2})/B^{2m-1}, \; d\varphi = \eta \cup \eta  \} .
$$
\end{parag}

One can construct more complicated examples but the above show essentially what
is going on. It seems intuitively clear that the nonabelian cohomology we are
looking at above contains all of the ``rational'' homotopic information, at
least
in the example where $\pi _1$ is abelian. However, it is not immediately obvious
exactly how this works; we formulate that as a problem for further study.

\begin{problem} How can we get back the homotopical information e.g. of $\pi
_2(\tilde{X})$ and $\pi _3(\tilde{X})$ from the nonabelian cohomology, notably
in the example considered in this chapter where $\pi _1$ is abelian?
\end{problem}

Part of the goal of this paper is to explain the ``de Rham'' and ``Dolbeault''
counterparts of the nonabelian Betti cohomology illustrated in the above
examples. We will obtain nonabelian cohomology $n$-stacks denoted
$\underline{Hom}(X_{DR}, T)$ and $\underline{Hom}(X_{Dol}, T)$. We will also
obtain a ``Hodge filtration'' which is a family over $\af$ linking these two.
If the variety $X$ varies in a family then we obtain a ``Gauss-Manin
connection'' on the family of nonabelian de Rham cohomologies, and together with
the Hodge filtration this satisfies an analogue of Griffiths transversality; and
at singularities of the family, it satisfies a regular singularity condition. In
the present paper we restrict to looking at the above algebraic-type aspects of
Hodge theory.

\numero{Nonabelian cohomology}
\label{nonabelianpage}

\subnumero{Motivation and historical remarks}

I'll begin by presenting my own reasons for approaching this theory.
It starts with the search for ``higher'' nonabelian cohomological invariants for
which Hodge theory would make sense. The first observation is that Hodge
theory is essentially about complexified invariants. For example, the Hodge
decomposition is defined on the cohomology with complex coefficients $H^i(X, \cc
)$. Similarly, the degree one nonabelian Hodge fitration
is defined on the space of representations in a complex algebraic group, such as
$H^1(X, GL_n(\cc ))$. This observation led to the idea in the last paragraph of
\cite{Simpson}(i), that one should look for other functors
$$
F: Top \rightarrow Sch
$$
from spaces to schemes, and try to obtain a Hodge theory for the complexified
versions $F(X)_{\cc}$. We can also operate the complexification beforehand, and
look for functors
$$
F: Top \rightarrow Sch /\cc .
$$

Thinking of a scheme as representing a sheaf of sets on the site $Sch /\cc $,
we can more generally look for functors from $Top$ to the category of sheaves
of sets on $Sch/\cc$. Thus, for each space $X$ and scheme $Y$ we would like to
have a set $F(X)(Y)$. This set should behave somewhat like a cohomology
theory in
$X$.

The next step is to apply the ``Brown representability'' idea, which says that a
functor which behaves like a cohomology theory in the variable space $X$, should
really be the set of homotopy classes of maps from $X$ to a space $T$. In our
case, the functor $X\mapsto F(X)(Y)$ in question depends on the scheme $Y$,
so $T$
should depend on $Y$.

From this (somewhat heuristic) reasoning, we finally come to the following
situation: we would like to have a family of spaces $T_Y$ indexed
contravariantly by schemes $Y$, and then we would like to set
$$
F(X)(Y):= Hom (X,T_Y).
$$

The previously known complexified cohomological invariants do indeed fit into
this picture: if we set $T_Y:= K(\Oo (Y), i)$ then
$$
Hom (X,T_Y)= H^i(X, \cc )
\otimes _{\cc} \Oo (Y)
$$
so the functor
$Y\mapsto Hom(X, T_Y)$ is represented by the complex vector space $H^i(X,
\cc )$.

Similarly, if we set $T_Y:= K(GL_n(\Oo (Y)) ,1)$ with basepoint $0$, then for a
pointed space $(X,x)$ we have
$$
Hom ((X,x), (T_Y, 0))= Hom (\pi _1(X, x), GL_n(\Oo (Y))
$$
and again, as a functor of $Y$, this is represented by the complex
representation scheme $R(X,x, GL_n(\cc ))$.

As one can already see from the need to go to pointed spaces in the previous
example, it is not actually a good idea to take homotopy classes of maps
separately over each scheme $Y$. Rather, one should take {\em homotopy coherent
maps} from the constant presheaf of spaces $\underline{X}$ with values $X$, to
$T$. One obtains in this way a functor of schemes $Y$, by looking at the
over-sites $Sch /Y$; more precisely it is a functor from schemes $Y$ to spaces,
which we can denote
$$
\underline{Hom}(\underline{X}, T):= \left( Y\mapsto Hom
_{Sch/Y}(\underline{X}|_{Sch /Y}, T|_{Sch/Y}) \right) .
$$
Here one should consider, for example, $T$ as being a homotopy-coherent family
of spaces indexed contravariantly by $Sch /\cc$, and the $Hom$ on the right is
the space of homotopy-coherent maps. The theory of homotopy-coherent maps,
due to
Leitch \cite{Leitch}, Vogt \cite{Vogt} and Mather \cite{Mather} in the
1970's, was
reinvented for the purposes of the above discussion in \cite{Simpson}(iv),
and the
application to nonabelian de Rham cohomology was written up in
\cite{Simpson}(v).
These latter use the internal $\underline{Hom}$ between homotopy-coherent
diagrams or presheaves of spaces, as ``nonabelian cohomology''.

Now $\underline{Hom}(\underline{X}, T)$ is
a functor from schemes to spaces. If one wants to get a functor from schemes to
sets, compose with $\pi _0$ (and then, sheafify, if you want to get a sheaf of
sets).

Before getting to the optimal technical way to approach the above theory (via
Jardine's theory of simplicial presheaves), we should back up a bit and talk
about categorical shape theory and $n$-stacks.

The ``arrow family'' or $Hom$ functor plays a central role in the categorical
shape theory developed by Armin, Frei \cite{Frei}, Cordier and
Porter \cite{Cordier} \cite{CordierPorter} \cite{Porter}, and others. This was
early on combined with ``enriched category theory'' in a way which is very close
to the use of internal $Hom$ that we have in mind (see \cite{Porter} for
example). A recent paper which also looks relevant is Carboni {\em et al}
\cite{CKVW}. The notion of topological shape is closely related to \v{C}ech
cohomology, and can be thought of as a ``nonabelian'' version of \v{C}ech
cohomology.  We will often employ a shape-theoretic terminology.

Then comes the idea of looking at ``$n$-stacks''. For $n=1$ this theory was
developped by Grothendieck and Giraud in the early 1960's, and it led to
Giraud's
book on nonabelian cohomology.  The generalization to higher $n$, while
present in
spirit throughout Grothendieck's work, came into focus in ``Pursuing
Stacks''. In
particular, the use of an internal mapping stack $\underline{Hom}(X,T)$ between
two $n$-stacks, as ``nonabelian cohomology of $X$ with coefficients in $T$'',
was explicitly mentionned there. Grothendieck writes
(\cite{Grothendieck}(iv), letter dated 27.2.83):

{\small ``Thus $n$-stacks, relativized  over a topos to `$n$-stacks over $X$',
are viewed primarily as the natural `coefficients' in order to do
[non-commutative] (co)homological algebra of dimension $\leq n$ over $X$.''}

My own use of such an internal $\underline{Hom}$ was directly inspired by this
sentence.

Grothendieck discusses the relationship between  nonabelian cohomology with
finite locally constant coefficient $n$-stacks, and the Artin-Mazur etale
prohomotopy theory.  The only thing lacking in Grothendieck's
discussion was the definition of $n$-stack, a gap which was closed by the
letter of Joyal \cite{Joyal}.

Joyal's letter takes us back to the work of K. Brown \cite{KBrown} and Heller
\cite{Heller}, closely related to Quillen's \cite{Quillen}(i) and also to
several early works of Breen \cite{Breen}. These in turn eventually lead back to
Eilenberg and MacLane, and the representation of cohomology by Eilenberg-MacLane
spaces $K(\pi , n)$. The idea is to relativize over a category, site, or just a
base topological space, this representation of cohomology by spaces.

One historical point which is not often recognized is that the notion of
``topological space relative to a topos'' was discovered quite early on by
Grothendieck and his entourage, in the notion of ``topos bifibered in toposes
over a base topos'' which appears in SGA 4 \cite{SGA4}.  If one considers the
topos fibers as being spaces, then one obtains directly the notion of family of
spaces indexed by a base topos; and it would have been possible to start the
theory of $n$-stacks of groupoids directly from there. However, this was not
pursued at the time, nor was it taken up later in ``Pursuing Stacks''.

It is not clear to what
extent Grothendieck was aware of the works of Breen and Brown refered to two
paragraphs ago, and indeed after recieving Joyal's letter, he seems to have more
or less (but not completely) ignored this idea, instead changing direction (a
couple of times) in the middle of ``Pursuing stacks'' before eventually dropping
the matter entirely to look at Teichmuller spaces.

In particular, Grothendieck should have pointed out that the theory
of simplicial sheaves explained by Joyal, in fact responds quite well to what he
was looking for in a theory of $n$-stacks of groupoids; and thus enables one to
define ``nonabelian cohomology''.

Joyal's theory was taken up by Jardine, who noticed the
essential observation that the sheaf condition levelwise for simplicial
presheaves, is basically irrelevant if not counterproductive. (For example
the problem encountered by K. Brown \cite{KBrown}, which meant that he couldn't
define a closed model structure, was related to this sheaf condition). This
observation led Jardine to look at {\em simplicial presheaves} over a site
$\Gg$,
and to define a closed model structure \cite{Jardine}(i). The topology
on the site comes into the definition of {\em Illusie weak equivalence}
\cite{Illusie} where one takes the sheaf associated to the presheaf of homotopy
groups. (In the absence of a topology, this closed model structure was due to
Heller \cite{Heller} and another closed model structure had been given by
Bousfield and Kan \cite{BousfieldKan}.)

Thomason, in his ENS paper \cite{Thomason}, was the first (or maybe the second
after Brown and Gersten \cite{BrownGersten}) to actually use nonabelian
cohomology
in the sense we are talking about here in a serious way to prove results in
$K$-theory. He developped many techniques such as the notion of ``Godement
resolution'', etc. In a paragraph which he says was asked for by the referee, he
develops a ``Leray theory'' calculating the cohomology of $X$ with coefficients
in $T$ by taking the cohomology of $Y$ with coefficients in $H(X/Y,
T)\rightarrow
Y$. There are two {\em caveats}: his theory was developed exclusively in the
``stable'' case, i.e. in the context of presheaves of spectra; and he doesn't
seem to mention the internal $\underline{Hom}$, sticking to the external
$Hom(X,T)$ which is a space rather than a stack.

While Jardine
doesn't seem to have explicitly stated it in \cite{Jardine}, it is immediate to
obtain the internal $\underline{Hom}(X,T')$ between two simplicial presheaves, by
setting
$$
\underline{Hom}(X,T')(Z):= Hom _{\rm spl}(X|_{\Gg /Z}, T'|_{\Gg /Z})
$$
using his simplicial model category structure \cite{Jardine}(i).

An internal mapping
presheaf of spectra between two presheaves of spectra was used (for
generalized etale cohomology) by R. Joshua in \cite{Joshua}.
This concerns only the stable case (although, again, it is immediate
to extend this definition to the unstable case using Jardine's closed model
structure). Joshua discusses many aspects of ``generalized cohomology'' which
are precursors to things that we shall discuss below, and his paper is
a point of origin for much of the theory.

More recently, the internal $\underline{Hom}$ between simplicial presheaves has
entered into Morel's and Voevodsky's theory of motivic cohomology in
\cite{Voevodsky} and \cite{MorelVoevodsky}, see also Kahn
\cite{Kahn}.

\subnumero{$n$-stacks and $n$-categories}

We will employ Grothendieck's notation ``$n$-stack'' since it is
convenient, historically early (even though it wasn't well defined
when Grothendieck invented it), 	and compact. See, however, the {\em
Alternative} below.

By {\em $n$-category} we shall mean weak $n$-category or $n$-nerve in the
sense of
Tamsamani \cite{Tamsamani}. He defines a groupoid condition, and
{\em $n$-groupoid} shall mean $n$-category satisfying this condition. We
sometimes denote by ``$\ngr$-category'' an $n$-groupoid, and by ``$\ngr$-stack''
an $n$-stack of $n$-groupoids.

There are currently available several other definitions of weak $n$-category,
such as Baez-Dolan \cite{BaezDolan} or Batanin \cite{Batanin}(iv). It should in
principle be possible to develop the present theory with those definitions; or
alternatively to wait and transfer the theory into those frameworks once the
appropriate comparison results are known. We don't delve into these questions
here.

The reader
is referred to \cite{Tamsamani}, \cite{Simpson} for the definitions and basic
properties of $n$-categories. In particular, there are introductory sections
(basically the same) in \cite{Simpson}(xi) and (xiii), and it
would be redundant to repeat that a third time here.

Instead, we just recall the
basic  notations: if $A$ is an $n$-category, then its set of objects is denoted
$A_0$. Its $n-1$-category of morphisms is denoted $A_{1/}$ and more generally
the $n-1$-category of composable $p$-uples of morphisms is denoted $A_{p/}$. By
definition these fit together into a simplicial object indexed by $p$
(this simplicial structure is Segal's way of encoding the composition law plus
higher-order homotopy coherencies). If $x,y\in A_0$ are two objects, the
$n-1$-category of morphisms from $x$ to $y$ is denoted $A_{1/}(x,y)$ and more
generally if $x_0,\ldots , x_p$ are objects then the $n-1$-category of $p$-uples
of composable morphisms having sources/targets the $x_i$ in order, is denoted by
$$
A_{p/}(x_0, \ldots , x_p).
$$

In \cite{Simpson}(x) was defined the {\em $n+1$-category of $n$-categories}
denoted $nCAT$.
The basic technique was to define a closed model category of {\em $n$-precats}
denoted $nPC$, and to note that this closed model category admits a
homotopically
correct internal $\underline{Hom}$. Note that we have the inclusions
$$
(\mbox{fibrant} \;\; n-\mbox{precats})\subset (n-\mbox{categories})
\subset nPC,
$$
and any object is equivalent to a fibrant one.
Now $nCAT$ is
defined as the $n+1$-category having for objects the fibrant $n$-precats,
and with
$$
nCAT_{p/}(U_0,\ldots , U_p):= \underline{Hom}(U_0,U_1)\times
\ldots \times \underline{Hom}(U_{p-1}, U_p).
$$

We recall the notation $\tau _{\leq k}$ for the {\em truncation} of an
$n$-category to a $k$-category. In the case of $\ngr$-categories this
corresponds exactly to the Postnikov truncation which kills off the homotopy
groups in degrees $i>k$. According to context, it will be useful either to
consider this as a $k$-category or else to consider it as an $n$-category. In
the few places where needed for clarity, we denote the operation of considering
a $k$-category as an $n$-category ($n\geq k$) by $Ind_k^n$.

If $A$ is an $n$-category then we obtain its {\em opposite} denoted $A^o$ as
follows. The set of objects is the same, and we define
$$
A^o_{p/}(x_0,\ldots , x_p):= A_{p/}(x_p, \ldots , x_0).
$$
The maps giving the simplicial structure are defined in the obvious way using
those of $A$ (and in fact this operation comes from an involution of the
category $\Delta$ reversing the ordering of the finite ordered sets). The
operation $A\mapsto A^o$ reverses the direction of the $1$-morphisms but doesn't
reverse the directions of the $i$-morphisms for $i\geq 2$. One could of course
imagine having a number of different ``opposite'' constructions, reversing
directions of arrows at various different levels. This isn't useful in the
present paper so we don't establish a notation for it.

Throughout the paper, $\Gg$ will denote a Grothendieck site. Most often this
will be the site of noetherian schemes over a base field $k$ of characteristic
zero. We use the theory of {\em $n$-stacks} over $\Gg$ developped in
\cite{HirschowitzSimpson}. Briefly, one defines a closed model category $nPS(\Gg
)$ whose objects are called {\em $n$-prestacks over $\Gg$}. These are just
presheaves of $n$-precats. An $n$-prestack $A$ is said to be an {\em $n$-stack}
if the $A(X)$ are $n$-categories (rather than just $n$-precats) and if
the morphism to a fibrant replacement $A\rightarrow A'$ induces
object-by-object equivalences $A(X)\cong A'(X)$. This condition is equivalent
to a descent condition analogous to the standard one for $1$-stacks.
Furthermore, the fact that we have taken strict presheaves is not a restriction
because any weak functor from $\Gg$ to $nCAT$ can be ``strictified''. Thus, an
alternative notion of $n$-stack is obtained by looking at the weak functors
$\Gg \rightarrow nCAT$ satisfying the analogue of the usual descent condition.
This equivalence is proven in \cite{HirschowitzSimpson} Th. 12.1.
One defines the $n+1$-category $nSTACK(\Gg )$ as having objects the fibrant
$n$-prestacks over $\Gg$, and having
$$
nSTACK(\Gg )_{p/}(U_0,\ldots , U_p):= \underline{Hom}(U_0,U_1)\times
\ldots \times \underline{Hom}(U_{p-1}, U_p).
$$
We obtain an $n+1$-prestack denoted $n\underline{STACK}$ by setting
$$
n\underline{STACK} (X):= nSTACK (\Gg /X),
$$
and \cite{HirschowitzSimpson} Th. 20.5 states that this is an $n+1$-stack.
We call it the {\em $n+1$-stack of $n$-stacks}.

In the body of the present paper, we will often use the notation
$\underline{Hom}(A,B)$ for $n$-stacks $A$ and $B$. We make the convention that
$B$ should be replaced by an equivalent fibrant
object in the model category of $n$-prestacks. Thus what we shall often write as
$\underline{Hom}(A,B)$ really means (in the notation of
\cite{HirschowitzSimpson} for example) $\underline{Hom}(A,B')$ where $B\mapsto
B'$ is the fibrant-replacement functor.

We often want to restrict our attention to $n$-groupoids. An {\em $n$-stack of
groupoids} is an $n$-stack $A$ such that the $A(X)$ are $n$-groupoids. For
brevity (and following---up to a slight typographic repositioning---a convention
from \cite{Grothendieck}(iv)) we call this an {\em $\ngr$-stack}.
We denote by $\ngr \underline{STACK}$ the full substack of $n\underline{STACK}$
consisting of $\ngr$-stacks. Similarly, $\ngr CAT$ is the full subcategory of
$nCAT$ consisting of $n$-groupoids.

The reader who can't resist trying to take $n=\infty$ in the above discussion,
is referred to \cite{HirschowitzSimpson} where we employ a notion of {\em Segal
$n$-category} (and consequently, {\em Segal $n$-stack} etc.) which corresponds
to the notion of $\infty$-category where the morphisms are invertible in degrees
$>n$. Some of what we shall say below is applicable in the case $n=\infty$ but
other parts would require additional hypotheses so the reader is advised to
tread with care. In the end, since we don't look at ``phantom maps''
\cite{Handbook}, we are only interested in homotopy in finite degrees so it is
sufficient to take a (sometimes large and unspecified) finite value of $n$.

\oldsubnumero{Alternative}

The alternative to the above notation and way of proceeding is the notion of
``simplicial presheaf''. This has priority in a historical sense and many
readers will be more familiar with it.

We indicate here the changes which should be made in order to read the present
paper from a ``simplicial presheaves'' perspective. See the first section of
\cite{HirschowitzSimpson} for a more precise version of the present remarks. The
main remark is that by \cite{Tamsamani}, an $n$-groupoid (or $\ngr$-category) is
the same thing as an $n$-truncated simplicial set i.e. a simplicial set whose
homotopy groups vanish in degrees $i>n$. Thus, instead of an $\ngr$-stack i.e.
presheaf of $\ngr$-categories, one can look at a simplicial presheaf such that
the values are $n$-truncated. Jardine's closed model structure on the category
of simplicial presheaves \cite{Jardine}(i) serves to define the homotopy
theory in question and replaces the closed model structure refered to above (and
indeed the closed model structure refered to above was just a modification of
Jardine's).  Concretely, if $Y$ and $T$ are simplicial presheaves then to
calculate the nonabelian cohomology, suppose that $T$ is fibrant (otherwise
replace it by an equivalent fibrant object) and take the internal
$\underline{Hom}(Y,T)$ in the category of simplicial presheaves. This again is a
simplicial presheaf. If $T$ is $n$-truncated then $\underline{Hom}(Y,T)$ will
also be $n$-truncated. (Note also in this connection that $n$-truncating $Y$ by
adding cells of dimension $\geq n+2$ doesn't affect $\underline{Hom}(Y,T)$
up to homotopy, if $T$ is $n$-truncated; so in this situation we are not
necessarily obliged to assume that $Y$ itself is $n$-truncated.)

The translation becomes a bit more complicated when it comes to
$n+1$-categories such as $\ngr CAT$ and $n+1$-stacks such as $\ngr
\underline{STACK}$. These are $n+1$-categories or $n+1$-stacks which are {\em
$1$-groupic}, i.e. where the $i$-morphisms are invertible for $i>1$. However,
there are $1$-morphisms which are not invertible, so these objects don't
correspond to simplicial sets. Rather, they correspond to {\em simplicial
categories}. Thus $\ngr CAT$ should be viewed as the simplicial category
of fibrant $n$-truncated simplicial sets. The $n+1$-category of global sections
$$
\ngr  STACK (\Gg ):= \Gamma (\Gg ,\ngr \underline{STACK})
$$
should be viewed as the simplicial category of fibrant $n$-truncated simplicial
presheaves. For variable $X\in \Gg$, the simplicial categories
$\ngr  STACK (\Gg /X)$ fit together into a presheaf of simplicial categories,
and this is what we call $\ngr \underline{STACK}$.

The only difficulty in this point of view is when we wish to speak of an
internal $\underline{Hom}(\underline{R}, \underline{A})$ between two
$n+1$-stacks
$$
\underline{R}, \, \underline{A} \subset \ngr \underline{STACK}
$$
which are viewed as presheaves of simplicial categories. For this we need a
closed model structure on the presheaves of simplicial categories. One should
be able to obtain this using the closed model structure on simplicial categories
of Dwyer-Hirschhorn-Kan \cite{DHK} (in this structure the domain $\underline{R}$
would have to be replaced by a cofibrant object).

An alternative to \cite{DHK}
is to use the notion of {\em Segal category} which is the Segal-delooping
machine weakened version of the notion of simplicial category. The definition
of Segal category first appears in Dwyer-Kan-Smith \cite{DwyerKanSmith}, and
they prove that Segal categories are equivalent to strict simplicial
categories.

In \cite{HirschowitzSimpson} we obtain a closed model structure for the
presheaves
of Segal categories (which are called ``Segal $1$-prestacks'' there) and
one can use that to calculate the  $\underline{Hom}(\underline{R},
\underline{A})$. That construction can be read without too much reference to the
notion of $n$-category for general $n$.

The notion of Segal category in the paper \cite{DwyerKanSmith} seems to be the
first place where a notion of ``$A_{\infty}$-category'' (i.e. weak simplicial
category) appears.  A subsequent appearence of this notion in a somewhat more
general form is in Batanin \cite{Batanin}(ii).
It might also be possible to use
\cite{Batanin}(ii) to do the theory of internal $\underline{Hom}$ between
presheaves of simplicial categories.

\subnumero{Nonabelian cohomology---basic definitions}

If $X$ and $T$ are $n$-stacks (over our site $\Gg$ or eventually over some $\Gg
/Z$) then we obtain the $n$-stack of {\em nonabelian cohomology of $X$ with
coefficients in $T$} denoted
$$
\underline{Hom}(X,T)
$$
which is, to be more precise in the framework of e.g. \cite{Simpson}(x),
\cite{HirschowitzSimpson}, the internal $\underline{Hom}(X,T')$ to a fibrant
replacement $T'$ of $T$. Nonabelian cohomology is a morphism of $n+1$-stacks
$$
\begin{array}{ccc}
n\underline{STACK} ^o \times n\underline{STACK} &\rightarrow
&n\underline{STACK}\\
X\;\;\; , \;\;\; T & \mapsto & \underline{Hom}(X,T).
\end{array}
$$
This morphism is basically the ``arrow family'' for the $n+1$-stack
$n\underline{STACK}$.

More generally for any $n+1$-stack $\underline{A}$ one
should be able to define a weak morphism (i.e. morphism to the fibrant
replacement denoted by $( - )'$)
$$
{\bf Arr}_{\underline{A}}:\underline{A}^o \times \underline{A} \rightarrow
n\underline{STACK}'.
$$
This was done for simplicial categories by Cordier and Porter in
\cite{CordierPorter}(vi).
The general situation is complicated by the fact that the morphism in
question will be defined only as a morphism into the fibrant replacement
$n\underline{STACK}'$---so in order to define the arrow family one must
therefore
have a good handle on a nice fibrant replacement. We treat that question in
Chapter 5.

\begin{parag}
\label{arrownstack1}
In the particular case $\underline{A} = n\underline{STACK}$ (or for any
full substack thereof) everything is made easier by the fact that
$n\underline{STACK}$ is constructed as a presheaf of enriched categories.
More precisely as was recalled above, for any
$Y\in \Gg$ we have
$$
n\underline{STACK} (Y):= nSTACK (\Gg /Y),
$$
and in turn the $n+1$-category $nSTACK (\Gg /Y)$ comes from the category of
fibrant $n$-prestacks over $Y$, enriched by the internal $\underline{Hom}$
$n$-stacks. Further helping is the fact that these internal $\underline{Hom}$
are themselves fibrant. Thus, the objects of $nSTACK (\Gg /Y)$ are the fibrant
$n$-prestacks $A$ over $Y$, and for any $p+1$-uple of objects $A_0,\ldots , A_p$
we  have
$$
nSTACK (\Gg /Y) _{p/} (A_0,\ldots , A_p):= Hom (A_0,A_1)\times \ldots \times Hom
(A_{p-1}, A_p).
$$
Here $Hom (A,B)$ denotes the $n$-precat of global sections of the internal
$\underline{Hom}(A,B)$.

\newparag{arrownstack2}
We can now complete the construction of the morphism
$$
{\bf
Arr}_{n\underline{STACK}}:
n\underline{STACK}^o\times n\underline{STACK} \rightarrow n\underline{STACK}.
$$
Fix an object $X\in \Gg$ and we construct the morphism between the
$n+1$-categories of sections over this object. On the level of objects,
if $A$ and
$B$ are objects of $n\underline{STACK}(X)$, this means that they are fibrant
$n$-prestacks over $X$. Put
$$
{\bf
Arr}_{n\underline{STACK}}(A,B):= \underline{Hom}(A,B).
$$
Suppose now that $A_0,\ldots , A_p$ and $B_0,\ldots , B_p$ are two sequences of
$p+1$ objects of the $n$-category $n\underline{STACK}(X)$. For brevity of
notation
set
$$
H_i:= \underline{Hom}(A_i, B_i).
$$
The composition map
$$
\underline{Hom}(A_1,A_0) \times  \underline{Hom}(B_0, B_1)
\times H_0 \rightarrow H_1
$$
yields, in view of the definition of internal $\underline{Hom}$, a map
$$
\underline{Hom}(A_1,A_0) \times  \underline{Hom}(B_0, B_1)
\rightarrow \underline{Hom}(H_0,H_1).
$$
Putting together several copies of the same type of map gives a map
$$
[n\underline{STACK}^o(X)\times n\underline{STACK}(X)]_{p/}((A_0, B_0), \ldots
,(A_p, B_p)) =
$$
$$
\underline{Hom}(A_1, A_0)\times
\underline{Hom}(B_0,B_1)\times  \ldots \times \underline{Hom}(A_p,
A_{p-1})\times \underline{Hom}(B_{p-1}, B_p)
$$
$$
\rightarrow
\underline{Hom}(H_0,H_1) \times \ldots \times \underline{Hom}(H_{p-1}, H_p)
$$
$$
= n\underline{STACK}(X)_{p/}(H_0,\ldots , H_p).
$$
The fact that the composition morphisms between internal $\underline{Hom}$ are
strictly associative implies that the above morphism respects the morphisms in
the simplicial structure, so it defines a morphism of $n+1$-precats
$$
n\underline{STACK}^o(X)\times n\underline{STACK}(X)
\rightarrow n\underline{STACK}(X).
$$
This is functorial in $X$, so it defines a morphism of $n+1$-prestacks
$$
{\bf
Arr}_{n\underline{STACK}}:
n\underline{STACK}^o\times n\underline{STACK} \rightarrow n\underline{STACK}
$$
as desired.

\end{parag}

\subnumero{A general type of question}

One of the main problems which arises in the context of nonabelian cohomology is
the following, which was stated as ``Question (C)'' in the introduction.

\begin{problem}
\label{zoo}
What properties does the nonabelian cohomology $\underline{Hom}(X,T)$ have, as a
function of the properties of $X$ and of $T$?
\end{problem}

In order to formulate this problem in a reasonable way, we will introduce (in
the next chapter) the notion of {\em realm}: this will be a full saturated
sub-$n+1$-stack $\underline{R} \subset n\underline{STACK}$. It corresponds
nicely to the notion of ``property'', putting
$$
\underline{R}(X)\subset n\underline{STACK}(X)= nSTACK (\Gg /X)
$$
equal to the full sub-$n+1$-category of those $n$-stacks on $\Gg /X$ which
locally have the property in question.

We can think of the problem of classifying interesting realms $\underline{R}
\subset n\underline{STACK}$, as a sort of ``zoology'' of $n$-stacks.  We will
begin to look at this in the next section.

In terms of realms, we can rephrase the above problem in the following way:

\begin{problem}
\label{zoo2}
Find triples of realms $\underline{P}$, $\underline{R}$, $\underline{A}$
such that the arrow family restricts to a morphism
$$
{\bf Arr}: \underline{P}^o \times \underline{R} \rightarrow \underline{A}.
$$
In other words we would like $\underline{Hom}(X,T)$ to be in $\underline{A}$
whenever $X$ is in $\underline{P}$ and $T$ is in $\underline{R}$. We call a
triple of realms satisfying this condition {\em well-chosen}.
\end{problem}

In a well-chosen triple $(\underline{P}, \underline{R}, \underline{A})$
we call $\underline{P}$ the {\em realm of domains}, we call $\underline{R}$ the
{\em realm of coefficients}, and we call $\underline{A}$ the {\em realm of
answers} (all of these terminologies are for obvious reasons).

We say that a realm $\underline{R}$ is {\em cohomologically self-contained}
if the
triple $(\underline{R}, \underline{R}, \underline{R})$ is well-chosen. This
property seems to be fairly rare, but it is interesting to try to get as close
to it as possible.

\oldsubnumero{Shape}

The notion of nonabelian cohomology, and in particular the main problem we have
formulated above, is closely related to  (and motivated by) {\em shape theory}.
See Borsuk \cite{Borsuk}, Cordier and Porter \cite{Cordier} \cite{CordierPorter}
\cite{Porter}, Deleanu-Hilton \cite{DeleanuHilton}, Frei \cite{Frei}, Mardesic
and Segal \cite{MardesicSegal} \cite{Siegel}, Batanin \cite{Batanin} etc.
Start with the following remark.

\begin{parag}
\label{shape1} 
Using the definition of internal
$\underline{Hom}$ for $n+1$-stacks, the arrow family ${\bf
Arr}_{n\underline{STACK}}$ constructed in \ref{arrownstack2} can be viewed as a
morphism
$$
{\bf Shape}: n\underline{STACK}^o \rightarrow
\underline{Hom}(n\underline{STACK}, n\underline{STACK}).
$$
Thus if $\Ff$ is an $n$-stack, we denote by ${\bf Shape}(\Ff )$ the object of
$\underline{Hom}(n\underline{STACK},
n\underline{STACK})$ given by $T\mapsto \underline{Hom}(\Ff , T)$.
\end{parag}

\begin{parag}
\label{shape2}
Now suppose that $\underline{R}\subset n\underline{STACK}$ is a realm. If
$\Ff$ is an $n$-stack then we denote by
$$
{\bf Shape}_{\underline{R}}(\Ff )
$$
the restriction of $Shape(\Ff )$ to an object in $\underline{Hom}(\underline{R},
n\underline{STACK})$. We call this functor the {\em $\underline{R}$-shape of
$\Ff$}. It is the nonabelian cohomology functor $T\mapsto
\underline{Hom}(\Ff ,T)$ but only for coefficient $n$-stacks $T$ in the realm
$\underline{R}$.

\newparag{shape3}
Suppose now that $\underline{A}$ is another realm. We say
that {\em the $\underline{R}$-shape of $\Ff$ takes values in
$\underline{A}$} if  ${\bf Shape} _{\underline{R}}(\Ff )$ actually lies in
$$
\underline{Hom}(\underline{R},
\underline{A})\subset
\underline{Hom}(\underline{R},
n\underline{STACK}).
$$
If this is the case, then we denote the resulting object by
$$
{\bf Shape} _{\underline{R}}^{\underline{A}}(\Ff )\in
\underline{Hom}(\underline{R}, \underline{A}).
$$
For short, we can say that ``${\bf Shape}
_{\underline{R}}^{\underline{A}}(\Ff )$
exists'' instead of saying that
``${\bf Shape} _{\underline{R}}(\Ff )$ takes values in $\underline{A}$'', and we
sometimes even write this condition as
$$
\underline{R} \stackrel{{\bf Shape}(\Ff )}{\longrightarrow} \underline{A}.
$$

\newparag{shape4}
The class of $n$-stacks
$\Ff$ such that ${\bf Shape} _{\underline{R}}^{\underline{A}}(\Ff
)$ exists, forms  a saturated full sub-$n+1$-stack of $n\underline{STACK}$, i.e.
a realm. It is the largest realm $\underline{P}$ such that the triple
$(\underline{P}, \underline{R}, \underline{A})$ is well-chosen, and we have
$$
{\bf Shape} _{\underline{R}}^{\underline{A}} : \underline{P}^o \rightarrow
\underline{Hom}(\underline{R}, \underline{A}).
$$
This morphism (which of course also exists for any smaller realm $\underline{P}$
making the triple well-chosen) is the ``shape functor'', or ``nonabelian
cohomology functor'' for  $(\underline{P}, \underline{R}, \underline{A})$.

The full subcategory which is
the image of this morphism is known in the above references for shape theory, as
the ``shape category''. The realm of coefficients $\underline{R}$ is sometimes
known in the above references as the ``category of models''.
\end{parag}

We can again reformulate our main problem:

\begin{problem}
\label{zoo3}
Find triples of realms $\underline{P}$, $\underline{R}$, $\underline{A}$
such that $\underline{P}$ is contained in the realm of $n$-stacks $\Ff$ such
that ${\bf Shape} _{\underline{R}}^{\underline{A}}(\Ff
)$ exists.
\end{problem}

The possibility of having interesting problems of the above form \ref{zoo}
or equivalently \ref{zoo2}, \ref{zoo3} is closely related to the fact that we
look at nonabelian cohomology in a world of ``stacks'' i.e. relative to a site
$\Gg$. If $\Gg$ is the punctual site $\ast$ so we are really just looking at
mapping spaces between two spaces, the notion of ``property'' of a space is
relatively quickly exhausted (for example, a space with vanishing higher
homotopy
groups is just a set, classified by its cardinality), and the above
questions are
not necessarily the most important ones. For a general site $\Gg$, even a
$0$-stack  (i.e.  sheaf of sets) can have interesting properties---for
example if
it is represented by an object of the site, one recovers all of the standard
questions that can be asked about, say, schemes. Thus the above formulation of
three principal questions about nonabelian cohomology comes into play when we
look at an interesting site $\Gg$.

\subnumero{The relative case}

Suppose $X\rightarrow Y$ is a morphism of $n$-stacks, and $T$ is an $n$-stack.
Sometimes we can then define a relative nonabelian cohomology $n$-stack
$$
\underline{Hom}(X/Y,T)\rightarrow Y.
$$
Recall from \cite{HirschowitzSimpson} \S 11 that a morphism $X\rightarrow Y$ of
$n$-stacks (realized by a fibrant morphism of $n$-prestacks) is called {\em
compatible with change of base} (ccb) if, for every diagram of $n$-prestacks
$$
B'\stackrel{a}{\rightarrow} B'' \rightarrow Y
$$
where $a$ is a weak equivalence, the morphism
$$
B'\times _YX \rightarrow B'' \times _YX
$$
is a weak equivalence. This condition is independant of the realization of
our morphism as a fibrant morphism of prestacks.

From  \cite{HirschowitzSimpson}
Proposition 11.13,  if $X\rightarrow Y$ is a fibrant morphism of $n$-prestacks
which is ccb, and if $T$ is any fibrant $n$-prestack, then the functor
$$
H(E):= \{ E\rightarrow Y, E\times _YE\rightarrow T \}
$$
is representable by an object
$$
\underline{Hom}(X/Y, T)\rightarrow Y
$$
and furthermore this object is homotopically well-defined as a function of $X/Y$
and $T$. We call it the {\em relative internal $\underline{Hom}$}, or also the
{\em relative nonabelian cohomology of $X/Y$ with coefficients in $T$}.
We refer the reader to  \cite{HirschowitzSimpson}
Proposition 11.13 and Lemma 11.14 for more precise statements of the
properties of
this construction.

The main property (Lemma 11.14), a version of Leray theory, is that the
$n$-stack
of sections of $\underline{Hom}(X/Y,T)$  over $Y$, is equivalent to
$\underline{Hom}(X,T)$.

It isn't hard to see that there are morphisms $X\rightarrow Y$ which are not
ccb, and for which a homotopically well-defined relative nonabelian cohomology
cannot exist. On the other hand, the morphisms we are interested in will all be
ccb, thanks to the result of \cite{HirschowitzSimpson} Lemma 11.15 which says
that if $Y$ is an $n$-stack of groupoids, then any morphism $X\rightarrow Y$ is
ccb.

Thus the relative nonabelian cohomology $\underline{Hom}(X/Y,T)$ will always
exist if the base $Y$ is an $\ngr$-stack. This is the only case in which we
shall use it.

In \S 5 below we will explore a somewhat different approach to defining relative
nonabelian cohomology. This approach will allow us to obtain a relative version
of the shape functor: if $\Ff \rightarrow \Ee$ is an appropriate type of
morphism of $n$-stacks (or more generally a ``cartesian family'' of $n$-stacks
parametrized by a base $n+1$-stack $\Ee$) and if $\underline{R}$ is a realm,
then we will obtain the relative shape functor
$$
{\bf Shape}_{\underline{R}}(\Ff /\Ee ): \Ee \rightarrow
\underline{Hom}(\underline{R}, n\underline{STACK}).
$$
Again, the main case in which that approach works (i.e. when a morphism of
$n$-stacks can be considered as a cartesian family) is when the base stack is a
stack of $n$-groupoids (cf Proposition \ref{correspondence} below). It would be
interesting to analyse more carefully the relationship between the condition ccb
and the notion of ``cartesian family'' of \S 5 below.

\subnumero{Example: Artin-Mazur shape}

We now give some examples of nonabelian cohomology. The first and most basic one
corresponds to the Artin-Mazur etale homotopy
theory. This reproduces a discussion in ``Pursuing stacks''
\cite{Grothendieck}(iv). It also seems to have been one of the main motivations
for example for the work of Cordier and Porter \cite{Cordier}
\cite{CordierPorter}
\cite{Porter}. We require the homotopy group sheaves of $T$ to be finite locally
constant etale sheaves over the base $Z$. To be more precise, $T\in
\underline{AM}(Z)$ if and only if $\pi _0(T/Z)$ is a finite locally constant
etale sheaf on $Sch /Z$, and if for any $Z'\rightarrow Z$ and basepoint $t\in
T(Z')$, the $\pi _i(T|_{Z'}, t)$ are finite etale sheaves over $Sch /Z'$. This
condition defines the {\em Artin-Mazur realm}
 $\underline{AM}$ (cf \S 4 for the definition of ``realm''), and the
corresponding
shape theory is essentially the same thing as Artin-Mazur etale homotopy.

Suppose $X$ is a scheme. It represents a $0$-stack on the big etale site $\Gg$.
We obtain the {\em Artin-Mazur shape}
$$
{\bf Shape}_{\underline{AM}}(X)\in \underline{Hom}(\underline{AM},
n\underline{STACK}),
$$
whose value on $Z\in \Gg$ is the morphism
$$
T\mapsto \underline{Hom}_{\Gg /Z}(X\times Z, T).
$$
For most purposes it suffices to look only on the values over $Z=Spec(k)$
(the final object of the site), in which case we can write the shape as being
the morphism
$$
T\mapsto \underline{Hom}(X,T).
$$
For example, if $T=K(A,n)$ for $A$ a finite group, then
$\underline{Hom}(X,T)$ is an $n$-stack with $\pi _i = H^{n-i}(X_{\rm et}, A)$.
Any $T\in \underline{AM}(Spec(k))$ is obtained as a Postnikov tower whose
stages are of the form $K(A,n)$ (at least if $k=\overline{k}$; otherwise there
could be Galois twisting).

We can conclude this (overly brief) discussion by saying that
${\bf Shape}_{\underline{AM}}(X)$ contains all of the homotopical information
about $X$ which can be ``seen'' by etale cohomology with finite coefficients.

One of the important aspects of the theory of Artin-Mazur is the
``pro-representability'' of the shape in the simply connected case. This can
serve as an avatar for some of what we do in \S 6 and \S 10 below.

A generalized-cohomology (i.e. stable) version of the nonabelian cohomology
which enters in here is used by Thomason \cite{Thomason}, Jardine \cite{Jardine}
and others, notably Joshua \cite{Joshua}.

\subnumero{Example: nonabelian de Rham cohomology}

Here is the example that we are interested in for nonabelian Hodge theory:
the {\em  nonabelian de Rham cohomology} of a smooth projective variety $X$.
We have to define the type of domain and coefficient stacks.

Work on the site $\Gg$ of noetherian schemes over $\cc$.
The ``domain'' is the sheaf of sets $X_{DR}$ defined by
$$
X_{DR}(Y):= X(Y^{\rm red}).
$$
Heuristically, we look at points $Y\rightarrow X$ but glue together
infinitesimally near points, i.e. we identify two points which agree on the
reduced subscheme $Y^{\rm red}$. This sheaf is closely related to the
``crystalline site'' in a way which we leave to the reader's imagination.
It also has an interpretation as a ``formal category'' cf Grothendieck
\cite{Grothendieck}(iii), Berthelot \cite{Berthelot}, Illusie \cite{Illusie}.
It is this interpretation which we shall use in chapters 8-10 below.

This definition of $X_{DR}$ gives the correct answer for the usual cohomology,
for example:

\begin{parag}
\label{degree1}
$H^1(X_{DR}, GL_r (\Oo ))$ is the moduli stack of rank $r$
vector bundles with integrable connection on $X$ (the proof is relatively easy,
see \cite{Simpson}(v)); and

\newparag{degree2}
  $H^i(X_{DR}, \Oo )$ is the sheaf represented by
the complex vector space
$$
H^i_{DR} (X,\cc ):= {\bf H}^i(X, (\Omega ^{\cdot}_X, d))
$$
(i.e. the algebraic de Rham cohomology cf \cite{Grothendieck}
\cite{Hartshorne}).
For the proof of this second part, see Berthelot \cite{Berthelot}, Illusie
\cite{Illusie}, C. Teleman \cite{CTeleman}, or \cite{Simpson}(v). A proof
is recapitulated in section 8 below.

For the ``coefficients'', we first investigate which $n$-stacks give rise to the
above two cases.

\newparag{degree1bis}
 The $1$-stack denoted either $B\, GL_r(\Oo )$
or $K(GL_r(\Oo ), 1)$ is obtained by taking the stack associated to the
$1$-prestack (presheaf of $1$-truncated spaces)
$$
K^{\rm pre}(GL_r(\Oo ), 1)(Y):= K(GL_r(\Oo (Y)), 1).
$$
We have that
$$
\underline{Hom}(X_{DR}, K(GL_r(\Oo ), 1)),
$$
which is a $1$-stack, is the moduli stack $H^1(X_{DR}, GL_r(\Oo ))$ of rank $r$
vector bundles with integrable connection.

\newparag{degree2bis}
The $n$-stack $K(\Oo , n)$ is defined as the $n$-stack
associated to the presheaf of $n$-truncated spaces
$$
K^{\rm pre}(\Oo , n)(Y):= K(\Oo (Y), n).
$$
We have that
$$
\underline{Hom}(X_{DR}, K(\Oo , n))
$$
is an $n$-stack whose $\pi _0$ is given by the formula
$$
\pi _0\underline{Hom}(X_{DR}, K(\Oo , n)) = H^n(X_{DR}, \Oo ) = H^n_{DR}(X,
\cc )
$$
i.e. it is the algebraic de Rham cohomology. More generally
$$
\pi _i\underline{Hom}(X_{DR}, K(\Oo , n)) = H^{n-i}(X_{DR}, \Oo )
$$
and one can even write (non-canonically)
$$
\underline{Hom}(X_{DR}, K(\Oo , n))= \prod _{i=0}^n K(H^{n-i}(X_{DR}, \Oo ), i).
$$

\newparag{degree3}
Next, our idea for nonabelian de Rham cohomology is to make a definition of
coefficient stack which allows a mixing-up of the above examples. We say that an
$\ngr$-stack $T$ is {\em connected very presentable} if:
\newline
---$\pi _0(T) = \ast$, which implies  that we can choose a basepoint $t\in
T(Spec
(\cc )$ unique up to homotopy;
\newline
---$\pi _1(T,t)$ is a sheaf of groups on the site $\Gg$ represented by an
affine algebraic group of finite type over $\cc$; and
\newline
---for $i\geq 2$, $\pi _i(T, t)$ is a sheaf of abelian groups over $\Gg$
represented by a finite dimensional vector space, i.e. it is isomorphic to $\Oo
^{\oplus a}$ for a finite $a$.

We now obtain the {\em nonabelian de Rham cohomology of $X$ with coefficients in
$T$} which is an $n$-stack on $\Gg$
$$
\underline{Hom} (X_{DR}, T).
$$
One could get a sheaf of sets by looking at $\pi_0\underline{Hom} (X_{DR}, T)$
(this is a bit like looking at the coarse moduli problem for example in the case
of $H^1(X_{DR}, GL_r(\Oo ))$). However, this destroys a nice property: the
$n$-stack $\underline{Hom} (X_{DR}, T)$ is {\em geometric}, which is the
analogue for $n$-stacks of Artin's notion of ``algebraic $1$-stack''. We will
discuss geometric $n$-stacks in \S 7 below, and we prove the geometricity in
the present example, in \S 10.

\end{parag}

\subnumero{Example: constructible stacks}

We can tweak slightly the example of Artin-Mazur shape which occurs above, by
looking at nonabelian cohomology with ``constructible'' coefficients.

Suppose that $\Gg = X^{\rm et}$ is the small etale site of a scheme $X$.
If $T$ is an $n$-stack over $X^{\rm et}$ we say that $T$ is {\em finite
constructible} if:
\newline
---$\pi _0(T)$ is a finite constructible sheaf of sets
(recall that ``constructible'' means that it becomes locally constant over a
stratification of $X$ by locally closed subvarieties; ``finite'' means that the
values are finite sets); if
\newline
---for every section $t\in T(Y)$ for $Y\rightarrow X$ etale, the sheaf of groups
$\pi _1(T|_{Y^{\rm et}}, t)$ is finite constructible on $Y$; and
\newline
---for $i\geq 2$ and for every section $t\in T(Y)$ with $Y\rightarrow X$ etale,
the sheaf of abelian groups $\pi _i(T|_{Y^{\rm et}}, t)$ is finite
constructible.

We obtain a ``constructible nonabelian etale cohomology''
$H(X^{\rm et}, T):= \underline{Hom}(\ast _{X^{\em et}}, T)$
where $\ast _{X^{\em et}}$ is the sheaf of sets whose values are the one-point
set.

More generally, we could look at the nonabelian cohomology of any $n$-stack
$\Ff$ on $X^{\rm et}$, with coefficients in a finite constructible $n$-stack
$T$. This would be the $n$-stack $\underline{Hom}(\Ff , T)$. This might be
interesting in the following cases, among others: when $\Ff$ is some sheaf of
sets; when $\Ff = K(A, n)$ for a sheaf of groups $A$; or when $\Ff$ itself is a
finite constructible $n$-stack. These cases come up in classifying the
Postnikov invariants for a finite constructible $n$-stack $T$, and we can
formulate the following problem.

\begin{problem}
Suppose $A$ and $G$ are finite constructible sheaves of groups. Calculate the
internal cohomology (which is a sheaf of groups or sets on $X^{\rm et}$)
$$
\underline{H}^i(K(A, n), G).
$$
Here, $G$ (resp. $A$) is required to be abelian if $i\geq 2$ (resp. $n\geq 2$).
Similarly, calculate $\underline{H}^i(F, G)$ for a finite
constructible sheaf of sets $F$.

Also calculate the ``external'' cohomologies $H^i(K(A, n), G)$ or
$H^i(F, G)$ (i.e. the global cohomologies over the site).
\end{problem}

This type of problem (for various $A$ and $G$) was formulated and extensively
investigated by Breen \cite{Breen}.

R. Joshua's paper \cite{Joshua} formulates the notion of nonabelian
cohomology on
the small etale site $X^{\rm et}$. While he doesn't explicitly discuss finite
constructible stacks (at least in part I of \cite{Joshua}),
it is clear that this is what he has in mind because of the following sentence
in the introduction of \cite{Joshua}:

{\small ``One of the goals of the present work is to set up a broad
framework for
defining a sheaf-theoretic version of generalised intersection cohomology and
this will be continued in [part II].''}

Joshua's paper is concieved in the ``stable'' case. To simplify, we can think of
this as meaning that he is interested in the problem formulated above, in the
``stable range'' $i<2n-1$.

\oldsubnumero{Constructible stacks on the small etale analytic site}

Suppose $X$ is a complex analytic space. We denote the small analytic site
of $X$ again by $X^{\rm et}$ and we can make exactly the same definitions as
above. In this case, it becomes reasonable to consider arbitrary rather than
just finite sets or groups in our constructible sheaves. One says that a
sheaf (of sets, groups or abelian groups) is {\em constructible}
if it is locally constant when restricted to the strata of a stratification of
$X$ by locally closed analytic subspaces. (One could also fix a stratification
and look at constructible sheaves with respect to the given stratification; this
distinction isn't important here, however it would allow us to extend the
entire present discussion to arbitrary stratified topological spaces.)

With the possibility of infinite groups or sets as values, it becomes reasonable
to look at $n$-stacks whose values are not necessarily $n$-groupoids. Thus we
can make the following definition (which more or less occurs already in
``Pursuing stacks'', or at least it is strongly suggested in certain passages
there). An $n$-stack $T$ on $X^{\rm et}$ is {\em constructible} if there is a
stratification of $X^{\rm et}$ by locally closed analytic subspaces such that in
a neighborhood of any point $P$, the stratified space together with the
$n$-stack
$T$ are trivialized in the direction of the stratum containing $P$ (this means
that in a neighborhood of $P$, the stratified space may be written  as a product
of a stratified space with an isolated point as the closed stratum, with a
smooth
manifold; and that the $n$-stack $T$ is pulled back from the first factor).

It is easy to see that an $\ngr$-stack on $X^{\rm et}$ is  constructible if and
only if:  \newline
---$\pi _0(T)$ is a  constructible sheaf of sets;
\newline
---for every section $t\in T(Y)$ for $Y\rightarrow X$ etale, the sheaf of groups
$\pi _1(T|_{Y^{\rm et}}, t)$ is constructible on $Y$; and
\newline
---for $i\geq 2$ and for every section $t\in T(Y)$ with $Y\rightarrow X$ etale,
the sheaf of abelian groups $\pi _i(T|_{Y^{\rm et}}, t)$ is constructible.

One easy statement is the following. We can state it even for
$n$-stacks whose values are not necessarily $n$-groupoids.

\begin{lemma}
\label{loctriv}
Suppose $X$ is a complex analytic space and let $X^{\rm et}$ denote the small
etale analytic site.
If $U$ and $T$ are constructible (resp. finite constructible) $n$-stacks on
$X^{\rm et}$ then the nonabelian cohomology stack
$$
\underline{Hom}(U,T)
$$
is a constructible (resp. finite constructible) $n$-stack on  $X^{\rm et}$.
\end{lemma}
{\em Proof:}
This follows from the local topological triviality of $X^{\rm et}$ (with its
stratification) along the strata.
\eop

An equivalent formulation, in the case of $\ngr$-stacks, is the following
statement (which could alternatively be proved first and then it implies
the above lemma for the case of groupoids, by the usual Postnikov and Leray
reductions).

\begin{corollary}
\label{cohoisconstruc}
Suppose $X$ is a complex analytic space.
Suppose $A$ and $G$ are constructible (resp. finite constructible) sheaves of
groups on the etale analytic site $X^{\rm et}$ (abelian if respectively
$i\geq 2$
or $n\geq 2$). Then the cohomology sheaf
$$
\underline{H}^i(K(A, n), G).
$$
is constructible (resp. finite constructible). Similarly if $\Ff$ is a
constructible sheaf of sets then
$\underline{H}^i(\Ff , G)$ is constructible, and finally if $\Ee $ is another
sheaf of sets then $\underline{Hom}(\Ff , \Ee )$ is constructible.
 \end{corollary}
\eop

In terms of cohomology, we can formulate the same problem as previously:

\begin{problem}
Suppose $X$ is a complex analytic space.
Suppose $A$ and $G$ are constructible sheaves of groups on the etale analytic
site $X^{\rm et}$. Calculate the cohomology sheaf (which is a sheaf of groups or
sets on $X^{\rm et}$)
$$
\underline{H}^i(K(A, n), G).
$$
Here, $G$ (resp. $A$) is required to be abelian if $i\geq 2$ (resp. $n\geq 2$).
Similarly, calculate $\underline{H}^i(F, G)$ for a
constructible sheaf of sets $F$.

Also calculate the ``external'' cohomologies $H^i(K(A, n), G)$ or
$H^i(F, G)$ (i.e. the global cohomologies over the site).
\end{problem}

In the above problem, one particular case which might be of interest is when
$A$ and $G$ are constructible sheaves of $\qq$-vector spaces.

\subnumero{Example: nonabelian coherent sheaf cohomology}

One of the main features of coherent sheaf cohomology is that one looks at
coherent sheaves with their $\Oo$-module structures. This isn't easy to
integrate into a nonabelian point of view (one possible path would
be to look in general at module-spectra and algebra-spectra over
ring-spectra in a
topos, see \cite{Joshua} for example; but we don't get into that here). Instead,
we make use of the following observation, which in a certain sense dates back to
Breen \cite{Breen}: on the big site of all noetherian schemes
in characteristic zero, the morphisms of sheaves of groups $\Oo \rightarrow \Oo$
are automatically morphisms of sheaves of $\Oo$-modules, and this property
extends to the higher $Ext$ groups too.  One can note that this is
certainly not true in characteristic $p>0$ where the Frobenius morphism is new,
so our discussion here will not work in its present form in characteristic
$p>0$.
(Nevertheless, the definitions we give might be interesting to look at in that
case too---it is basically what Breen looked at in \cite{Breen}(i).)

If $X$ is a scheme in characteristic zero, a coherent sheaf $\Ff$ on $X$ may be
thought of as a sheaf of groups on the big site $Sch /X$ (in the etale topology,
say). For this, use the formula
$$
\Ff (Y):= p^{\ast} (\Ff )(Y)
$$
for
$p:Y\rightarrow X$ the projection, and where $p^{\ast}$ denotes the pullback of
coherent sheaves.

We can now say that a $1$-connected $n$-stack $T$ on $Sch/X$ is {\em coherent}
if the $\pi _i(T)$ are coherent sheaves. We obtain the {\em nonabelian coherent
sheaf cohomology} $\underline{Hom}(U, T)$ of any $n$-stack $U$ on $Sch /X$ with
coefficients in a simply connected coherent $n$-stack $T$. This might be
particularly interesting to look at if $U$ is itself a coherent $n$-stack.
This latter type of problem comes up for example in the ``Dolbeault
cohomology'' of \cite{Simpson}(xii) where one looks at the cohomology of
$U= K(TX/X, 1)$. It would probably be interesting to look at this for more
complicated $U$.

We can of course formulate the abelian cohomology problems which arise in this
context: essentially the problem is to calculate $\underline{H}^i(K(\Ee ,
m), \Ff
)$ for coherent sheaves $\Ee$ and $\Ff$. Up to a problem of spectral
sequences, we
can calculate these. In order to present that result, one needs the more general
notion of {\em vector sheaf} which will be explained in \S 6 below. The
basic calculation is then Theorem \ref{bc}, but to obtain the calculation for
coherent sheaves $\Ee$ one needs to consider some spectral sequences as
explained in Remark \ref{specseqcalx}.

For the moment, we describe the calculation in a more restrained case: we say
that a $1$-connected $n$-stack $T$ is {\em locally free} if the $\pi _i(T)$ are
locally free sheaves of finite rank on $Sch /X$, in other words they are locally
isomorphic to $\Oo ^a$ for some $a$. For this case the resulting calculation of
cohomology is easy to state: we call it the ``Eilenberg-MacLane-Breen
calculations''. It is the analogue for coherent sheaf cohomology of the
Eilenberg-MacLane calculations, and a more complicated situation was treated by
Breen in characteristic $p>0$ \cite{Breen}(i) (he mentions the characteristic
zero case in \cite{Breen}(ii)).

\begin{proposition}
Suppose $\Ee \cong \Oo ^a$. If $m$ is even then
$$
\underline{H}^i(K(\Ee , m), \Oo )= Sym ^{\frac{i}{m}}(\Ee ^{\ast}),
$$
whereas if $m$ is odd then
$$
\underline{H}^i(K(\Ee , m), \Oo )= \bigwedge ^{\frac{i}{m}}(\Ee ^{\ast}).
$$
Here the convention is that the symmetric or exterior powers are said to be zero
if the exponent is not an integer.
\end{proposition}
{\em Proof:} See \cite{Breen}(i), (ii); or \cite{Simpson}(xii).
\eop

Using this proposition, we can make a number of constructions and calculate
examples (up to some spectral sequences and extension problems at least). For
example, for $m$ even, we define the {\em complexified $m$-sphere} to be the
$2m-1$-stack on $Sch /Spec (\cc )$ fitting into the following diagram:
$$
\begin{array}{ccccc}
K(\Oo , 2m-1) & \rightarrow & S^m_{\cc } & \rightarrow & \ast \\
&&\downarrow && \downarrow \\
&&K(\Oo , m) & \rightarrow & K(\Oo , 2m).
\end{array}
$$
Here the right-hand square is cartesian, whereas the upside-down ``L'' on the
left is a fibration sequence. The map on the bottom is the map corresponding to
the cohomology operation ``cup-product'': its classifying element is the
generator of
$$
H^{2m}(K(\Oo , m), \Oo )\cong \Oo ,
$$
this isomorphism being given by the preceding proposition.

For $m$ odd we simply
put $S^m_{\cc} := K(\Oo , m)$.

If $X$ is a scheme, one then obtains the nonabelian cohomology $n$-stack
$\underline{Hom}(X, S^m_{\cc })$. This was investigated in some depth in
\cite{Simpson}(xii).

\numero{Zoology}
\label{zoologypage}

In working with $n$-stacks as the natural coefficients for nonabelian
cohomology, one relatively quickly runs into a problem of ``zoology'':
there is a wide array of possible properties which one can require of an
$n$-stack, and this creates a need for some organizational principles.  One
might say that we would like to ``classify'' $n$-stacks, but in fact it would
be more to the point to say that we would like to classify the properties that
one can, or wants to, require of $n$-stacks.
To this end, we start by introducing the notion of ``realm''.

\subnumero{Preliminaries}

If $A$ is an
$n+1$-category, recall that a {\em full sub-$n+1$-category} $B\subset A$ is a
subobject (say, in the category $(n+1)PC$ of $n+1$-precats) which has the
property that for any sequence of objects $x_0,\ldots , x_p\in B_0$,
we have
$$
B_{p/}(x_0, \ldots , x_p) = A_{p/}(x_0, \ldots , x_p).
$$
Homotopically speaking, this amounts to saying that we have a functor
$i:B\rightarrow A$ of $n+1$-categories, which is fully faithful in that
$$
i: B_{1/}(x,y)\rightarrow A_{1/}(i(x), i(y))
$$
is an equivalence of $n$-categories for any pair of objects $(x,y)$. This
second version, apparently somewhat larger than the previous version of the
definition, is equivalent, namely if $i: B\rightarrow A$ is a fully faithful
functor then the full sub-$n+1$-category $B'$ of $A$ whose objects are the image
of $B_0$, is a full sub-$n+1$-category in the first sense and $i$ induces an
equivalence between $B$ and $B'$. We shall henceforth use these two notions
interchangeably.

We say that a full sub-$n+1$-category $B\subset A$ is {\em saturated} if it
satisfies the following condition: that if $x\in B_0$ and if $y\in A_0$ is an
object which is equivalent (in $A$) to $x$, then $y\in B_0$. A saturated full
sub-$n+1$-category $B\subset A$ is completely determined by the subset
$$
\tau _{\leq 0}(B)\subset \tau _{\leq 0}(A)
$$
and conversely any such subset determines a saturated full sub-$n+1$-category.

Now, we say that an inclusion of $n+1$-stacks $B\subset A$ is {\em full} (resp.
{\em saturated}) if for every $X\in \Gg$, the inclusion of $n+1$-categories
$B(X)
\subset A(X)$ is full (resp. saturated).

We note  that if
$A$ is an $n+1$-stack and if
$B'\subset A$ is a saturated full substack with respect to the coarse topology,
then the $\Gg$-stack $B$ associated to $B'$ is again a saturated full substack
of $A$.

\subnumero{Realms}

We now give the definition of ``realm'':

\begin{definition}
\label{realm}
Assume that the site $\Gg$ and the integer $n\geq 0$ are fixed.
A {\em realm} is a
saturated full sub-$n+1$-stack $\underline{R} \subset n\underline{STACK}$.
\end{definition}

A realm is completely determined by the subsheaf of sets
$$
\tau _{\leq 0} (\underline{R})\subset
\tau _{\leq 0} (n\underline{STACK})
$$
and conversely any such subsheaf of sets determines a realm.

This notion is quite close to Giraud's notion of ``lien'' \cite{Giraud}: a {\em
lien} is just a realm such that  $\tau _{\leq 0}
(\underline{R})\cong \ast$, thusly corresponding to a section of the sheaf of
sets  $\tau _{\leq 0} (n\underline{STACK})$. (The notion of ``lien'' in the
literature refers to the case of $1^{\rm gr}\underline{STACK}$ but we can
obviously make the same definition in general).

The specification of a realm is basically the same thing as the specification of
a property which is to be held by the elements of the realm. In order to obtain
a realm, the property has to be
\newline
---invariant under restriction of stacks (i.e. if
an $n$-stack $A$ on $\Gg /X$ has the property and if $Y\rightarrow
X$ is a morphism in $\Gg$ then $A|_{\Gg /Y}$ should also have the property);
and
\newline
---local (i.e. if $A$ is an $n$-stack on $\Gg /X$ and if $\{ U_{\alpha}
\rightarrow X\}$ is a covering family such that each $A|_{U_{\alpha}}$ has the
property, then $A$ should have the property).

Given a ``property'' satisfying the above two axioms, we obtain a realm by
setting
$$
\underline{R}(X) \subset n\underline{STACK}(X) = nSTACK (\Gg /X)
$$
equal to the saturated full sub-$n+1$-category consisting of the $n$-stacks $A$
on $\Gg /X$ which have the property in question. Conversely, given a realm
$\underline{R}$, then the property of being an element of $\underline{R}(X)$ is
a property which satisfies the above axioms.

\begin{parag}
\label{intersection}
It is clear that an arbitrary intersection of realms is again a realm. This
corresponds to concatenation of properties (i.e. the logical AND). In
particular, the intersection of any subset of the set of realms which we define
in this paper is again a realm!
\end{parag}

\subnumero{A few realms}

A first example of a realm is the {\em realm of $n$-groupoids} which we are
denoting
$$
\ngr \underline{STACK} \subset n\underline{STACK}.
$$
More generally, recall (\cite{HirschowitzSimpson} \cite{Simpson}(xi))
that we say that an $n$-stack $A$ is {\em $k$-groupic} if the values $A(X)$ are
$n$-categories in which the $i$-morphisms are invertible (up to equivalence)
for $i>k$. A $0$-groupic $n$-stack is the same thing as an $n$-stack of
groupoids.  It is clear that the property of being $k$-groupic is preserved by
pullback, and one can check (slightly less trivially) that it is local. We
obtain
the {\em realm of $k$-groupic $n$-stacks} denoted
$$
n^{k\, {\rm gr}}\underline{STACK} \subset
n\underline{STACK}.
$$
For $k=0$ this is the same as the realm $\ngr \underline{STACK}$.

Another example of a property which we can require comes from Grothendieck
\cite{Grothendieck}(iv), Breen \cite{Breen} and  Baez-Dolan
\cite{BaezDolan}, see also \cite{Simpson}(xiii). We say that an $n$-stack $A$ is
{\em $k$-connected} if $\tau _{\leq k}(A) = \ast$. Note that this does {\em not}
mean that the values $A(X)$ are $k$-connected $n$-categories, because  $\tau
_{\leq k}(A)$ is the $k$-stack associated to the $k$-prestack
$$
\tau ^{\rm pre}_{\leq k}(A)(X):= \tau _{\leq k}(A(X)).
$$
It does mean that the stack associated to this prestack, is trivial. Basically
this means that $A$ is ``locally'' $k$-connected. It is clear that this property
is preserved by pullback and that it is local.  We denote the {\em realm of
$k$-connected $n$-stacks} by
$$
n\underline{STACK}^{k\, {\rm conn}} \subset n\underline{STACK}.
$$

Generalizing the previous paragraph, we obtain the following construction: if
$\underline{R}\subset k\underline{STACK}$ is a realm of $k$-stacks, then the
inverse image of $\underline{R}$ by the truncation operation $\tau _{\leq k}$
is a realm of $n$-stacks which we denote
$$
(\tau _{\leq k})^{-1}(\underline{R}) \subset n\underline{STACK}.
$$

\subnumero{Index}

In this subsection we will list a number of realms which come into play further
on in the paper. Some of the concepts involved in the definitions we give here,
are only defined in future sections (notably \S 6, \S 7 and \ref{linear}). We
refer the reader to those sections for the precise definitions;  the
present list
is only intended to collect in one place the relevant notations.

The site in question is always $\Gg = Sch /k $ (which means the noetherian
schemes over $Spec (k)$), with the etale topology. The field $k$ is
assumed to be of characteristic zero. Fix a value of $n$.

\begin{parag}
\label{index1}
We now list our main examples.

\noindent
$\underline{PE}$ is the {\em realm of presentable $n$-stacks} (cf \S 6, p.
\pageref{pepage}); $T\in \underline{VP}(X)$ if $T$ is a presentable $n$-stack on
$\Gg /X$.

\noindent
$\underline{GE}$ is the {\em realm of geometric $n$-stacks of finite type}
(cf \S 7.1);
$T\in \underline{VG}(X)$ if $T$ is a geometric $n$-stack of finite type on $\Gg
/X$.

\noindent
$\underline{VP}$ is the {\em realm of very presentable $n$-stacks} (cf \S 6, p.
\pageref{vppage}); $T\in \underline{VP}(X)$ if $T$ is a very presentable
$n$-stack
on $\Gg /X$.

\noindent
$\underline{VG}$ is the {\em realm of geometric very presentable $n$-stacks}
(cf \S 7.1  and \S 6, p. \pageref{vppage});
$T\in \underline{VG}(X)$ if $T$ is a geometric very presentable $n$-stack
on $\Gg
/X$.

\noindent
$\underline{CV}$ is the {\em realm of connected very presentable
$n$-stacks};
$T\in \underline{CV}(X)$ if $T$ is a very presentable $n$-stack on $\Gg
/X$, with $\pi _0(T)= \ast _{\Gg /X}$.

\noindent
$\underline{CG}$ is the {\em realm of connected geometric very presentable
$n$-stacks};
$T\in \underline{CG}(X)$ if $T$ is a geometric very presentable $n$-stack
on $\Gg
/X$, with $\pi _0(T)= \ast _{\Gg /X}$.

\noindent
$\underline{FL}$ is the {\em realm of connected
$n$-stacks with $\pi _1$ a flat linear group scheme, and $\pi
_i$ locally free of finite rank}. Thus  $T\in \underline{FL}(X)$ if $T$ is an
$n$-stack on $\Gg /X$, with $\pi _0(T)= \ast _{\Gg /X}$, with $\pi _1(T, t)$
represented by a flat linear group scheme over $X$ (cf \ref{linear}), and with
the $\pi _i(T,t)$ being represented by finite-rank vector bundles over $X$
(here the basepoint $t$ exists
locally over $X$ in the etale topology so these conditions make sense).

\noindent
$\underline{FV}$ is the {\em realm of connected very presentable
$n$-stacks with $\pi _1$ a flat linear group scheme}. Thus  $T\in
\underline{FV}(X)$ if $T$ is a connected very presentable $n$-stack on $\Gg /X$
with $\pi _1(T,t)$ represented by a flat linear group scheme over
$X$ (cf \ref{linear}).

\noindent
$\underline{FG}$ is the {\em realm of connected geometric very presentable
$n$-stacks with $\pi _1$ a flat linear group scheme}. Thus  $T\in
\underline{FG}(X)$ if $T$ is a connected very presentable $n$-stack on $\Gg /X$
with $\pi _1(T,t)$ represented by a flat linear group scheme over
$X$ (cf \ref{linear}).
Note that $\underline{FG}= \underline{FV} \cap \underline{VG}$.

\noindent
$\underline{AL}$ is the {\em realm of $1$-connected
$n$-stacks with $\pi _i$ locally free of finite rank}. Thus  $T\in
\underline{AL}(X)$ if $T$ is a simply connected $n$-stack on $\Gg
/X$ with $\pi _i(T,t)$ represented by finite rank vector bundles
over $X$. These $1$-connected $n$-stacks were called {\em locally free} in the
previous chapter.

\noindent
$\underline{AV}$ is the {\em realm of $1$-connected very presentable
$n$-stacks}. Thus  $T\in
\underline{AV}(X)$ if $T$ is a simply connected very presentable $n$-stack
on $\Gg
/X$. Equivalently, if $T$ is a simply connected $n$-stack on $\Gg /X$ such that
the $\pi _i(T,t)$ are vector sheaves on $X$ (cf \S 6).

\noindent
$\underline{AG}$ is the {\em realm of $1$-connected geometric very presentable
$n$-stacks}. Thus  $T\in
\underline{AG}(X)$ if $T$ is a simply connected geometric very presentable
$n$-stack on $\Gg /X$. In Theorem \ref{criterion} we show that $T\in
\underline{AV}(X)$ is in $\underline{AG}(X)$ if and only if $H^{\ast}(T/X,
\Oo )$
is a residually perfect complex (\ref{residuallyperfect}) on $X$.
\end{parag}

\begin{parag}
\label{index2}
A somewhat different notation is the following: if $\Pp$ is a {\em directory of
Serre classes} (cf 4.8) then $\underline{M}^{\Pp}$ is the associated realm.
The realms $\underline{PE}$, $\underline{VP}$, $\underline{CV}$,
$\underline{FL}$, $\underline{FV}$, $\underline{AL}$ and $\underline{AV}$ are of
this form; whereas the realms of geometric $n$-stacks (those with a `G' in the
notation) don't come from directories of Serre classes.
\end{parag}

\begin{parag}
\label{index3}
We  have the following inclusions:
$$
\begin{array}{ccccc}
&& \underline{GE} &\subset &\underline{PE}\\
&& \cup && \cup \\
&& \underline{VG} &\subset &\underline{VP}\\
&& \cup && \cup \\
&& \underline{CG} &\subset &\underline{CV}\\
&& \cup && \cup \\
\underline{FL} &\subset& \underline{FG} &\subset &\underline{FV}\\
\cup && \cup && \cup \\
\underline{AL} &\subset& \underline{AG} &\subset &\underline{AV}.
\end{array}
$$
\end{parag}

\begin{parag}
\label{index4}
The values over $Spec(k)$ of many of the above realms are
all equal:
$$
\underline{CV}(Spec(k))=
\underline{CG}(Spec(k))=
\underline{FV}(Spec(k))=
\underline{FG}(Spec(k))=
\underline{FL}(Spec(k)).
$$
They are all equal to the $n+1$-category defined in \ref{degree3}, that of
connected $\ngr$-stacks $T$ on $Sch /k$ such that $\pi _1(T)$ is represented by
an affine algebraic group over $k$, and  $\pi _i(T)$ ($i\geq 2$) are represented
by finite dimensional vector spaces.  (We didn't put the basepoint into the
previous phrase and it may exist only over a finite extension $k'$ in which case
the conditions on $\pi _1$ and $\pi _i$ are meant to be taken over $k'$.) This
$n+1$-category is the category of coefficients which comes into play the most
often in our examples and applications.
\end{parag}

\subnumero{Realms of stacks not of finite type}

If $P$ is a property, then one habitually says that an object is ``locally $P$''
if it is covered by open subobjects having property $P$. The standard example is
the notion of ``locally of finite type''. We integrate this into our terminology
of realms. In this section we work over a site $\Gg$ of schemes. If $T$ is an
$n$-stack on some $\Gg /Y$, say that a morphism (of $n$-stacks on $\Gg /Y$)
$$
U\rightarrow T
$$
is an {\em open substack} if it is fully faithful, and if for
any scheme $X$ and any morphism $X\rightarrow T$ the fiber product $X\times _TU$
(\`a priori an $n$-stack) is a sheaf of sets represented by a Zariski open
subset
of $X$.

We say that a family $\{ U_{\alpha} \rightarrow T\}$ of open substacks {\em
covers $T$} if for any scheme $X$ and any morphism $X\rightarrow T$, the
family of Zariski open subsets $U_{\alpha} \times _TX$ covers $X$.

Suppose $\underline{R}$ is a realm. For $Y\in \Gg$ we say that an $n$-stack $T$
on $\Gg /Y$ is {\em locally of type $\underline{R}$} if there is a covering of
$T$ by open substacks $\{ U_{\alpha} \rightarrow T\}$ such that each
$U_{\alpha}$ is in $\underline{R}(Y)$. We define a realm $\underline{R}^{\rm
loc}$ by putting
$$
\underline{R}^{\rm loc}(Y)
$$
equal to the saturated full sub-$n+1$-category of $n\underline{STACK}(Y)$
consisting of the $n$-stacks $T$ on $\Gg /Y$ which are locally of type
$\underline{R}$ in the aforementionned sense.

In particular, refering to the index in the previous subsection, we obtain
realms
$$
\underline{PE}^{\rm loc},\;\;\;
\underline{GE}^{\rm loc},\;\;\;
$$
$$
\underline{VP}^{\rm loc},\;\;\;
\underline{VG}^{\rm loc}.
$$
The remaining realms of \ref{index1} are realms of connected $\ngr$-stacks,
so in
those cases $\underline{R}^{\rm loc} = \underline{R}$.

\subnumero{Closure properties}

Here we introduce some notations for certain basic properties that a
realm can have (but, let's not jump to ``meta-realms''\ldots ).
Fix a realm $\underline{R} \subset n\underline{STACK}$.

\begin{parag}
\label{closurelimits}
We say that it is {\em
closed under limits} if for every $X\in \Gg$, the full subcategory
$\underline{R}(X)\subset nSTACK (\Gg /X)$ is closed under taking finite limits
as defined in \cite{Simpson}(xi). In other words, we require that for any
$n+1$-category $I$ with functor $A: I \rightarrow \underline{R}(X)$,
the $n$-stack $\lim _{\leftarrow , I} A \in nSTACK (\Gg /X)$ is an element of
$\underline{R}(X)$.

We say that $\underline{R}$ is {\em closed under finite limits} if it satisfies
the above closure property for limits taken over finite $1$-categories $I$. (One
could certainly speak of finite $n+1$-categories here, but it isn't clear
whether that is useful.)

Closure under finite limits is equivalent to saying that $\underline{R}(X)$
contains the final object $\ast _X$, and is closed under fiber products.

\newparag{closurecolims}
We can similarly define the property of being {\em closed under colimits}
(resp. {\em finite colimits}).  Existence of limits and colimits in $nSTACK (\Gg
/X)$ was basically proven in \cite{Simpson}(xi) but in the stack case to be
rigorous some extra work is needed; we don't treat that question here.

\newparag{closuretruncation}
We say that a realm $\underline{R}$ is {\em closed under truncation} if
for every $X\in \Gg$, the full subcategory
$\underline{R}(X)\subset nSTACK (\Gg /X)$ is closed under the operations
$$
A\mapsto Ind _k^n(\tau _{\leq k}(A))
$$
where $Ind_k^n$ is the operation consisting of considering a $k$-stack as an
$n$-stack, for $0\leq k\leq n$.

\newparag{closureextension}
We say that a realm $\underline{R}$ is {\em closed under
extension} if the following property holds. Suppose that $B\in
\underline{R}(X)$ and suppose $A\in n\underline{STACK}(X)$ with a morphism of
$n$-stacks $f:A\rightarrow B$. Suppose that for every $Y\rightarrow X$ and for
every point $b\in B_0(Y)$, the fiber product
$$
Y\times _{B|_{\Gg /Y}} (A|_{\Gg /Y})
$$
is in $\underline{R}(Y)$. Then we want $A$ to be in $\underline{R}(X)$.

\newparag{yetanother}
Suppose $\underline{R}$ is any realm.
In view of the remark \ref{intersection}, we could look at the smallest realm
$\underline{R}'$ containing $\underline{R}$ and satisfying any given subset of
the above closure properties. In the definition of geometric
$n$-stack in \S 7 below, we'll do this with another closure property that
we haven't yet described.

\newparag{determined}
Suppose $\Gg$ has a final object.
A realm $\underline{R}$ which is closed under extensions and which contains the
representable objects $X\in \Gg$ is completely determined by its global
sections i.e. by the saturated full subcategory $\underline{R} (\ast )\subset
nSTACK (\Gg )$.
\end{parag}

\subnumero{Nonabelian cohomology of a locally constant $n$-stack}

To illustrate the usefulness of at least one of the above closure
properties, we
look at a concrete situation, basically the first that one
encounters.

\begin{theorem}
\label{cohconstant}
Suppose $K$ is a finite CW complex, and let $W$ be the constant
$n$-prestack whose
values are $\Pi _n(K)$. Then for any $n$-stack $T$, the cohomology stack
$\underline{Hom}(W, T)$ may be expressed as a finite limit of copies of $T$
(i.e. it is obtained by iterating constructions using fiber products of $T$ and
$\ast$). In particular, if $\underline{R}$ is a realm which is
closed under finite limits, and if $T\in \underline{R}(Z)$ then
$\underline{Hom}(W|_{\Gg /Z},T)\in \underline{R}(Z)$.
\end{theorem}
{\em Proof:}
We can represent $K$ as a homotopy colimit of copies of $\ast$ indexed by a
finite category $I$ (take a category whose nerve is homotopy equivalent to
$K$).
Then $\underline{Hom}(W,T)$ is the homotopy limit of the constant
functor $I\rightarrow nSTACK$ whose values are $T$.
\eop

We illustrate the above argument more
concretely with a few examples. If $W=\emptyset$ then
$\underline{Hom}(\underline{W}, T)=\ast$.  If $W=\ast$ then
$\underline{Hom}(\underline{W}, T)=T$. If $W=2\ast$ consists of two points then
writing
$$
W=\ast \cup ^{\emptyset} \ast
$$
yields
$$
\underline{Hom}(\underline{W}, T)=T\times T.
$$
If $W=S^1$ is the circle then writing it as the union of two contractible
intervals joined along a space homotopic to $2\ast$ i.e.
$$
W=\ast \cup ^{2\ast} \ast
$$
yields
$$
\underline{Hom}(\underline{W}, T)=T\times _{T\times T}T.
$$
The morphisms $T\rightarrow T\times T$ are the diagonal.
If $W=S^2$ is the $2$-sphere then writing it as the union of two contractible
hemispheres joined along the circle i.e.
$$
W=\ast \cup ^{S^1}\ast
$$
yields
$$
\underline{Hom}(\underline{W}, T)=T\times _{
(T\times _{T\times T}T)}T.
$$
From these examples one sees that
any $\underline{Hom}(\underline{W}, T)$ can be broken down into a succession of
fiber products of $T$.
The proof of Theorem  \ref{cohconstant} was just a fancy way of saying
this.

\subnumero{Directories of Serre classes}

We continue our study of the ``zoology'' of stacks, concentrating on
the case of realms of $n$-groupoids (i.e. realms contained in $\ngr
\underline{STACK}$). As one can see from the examples discussed above and below,
one good way of constructing a realm of $n$-groupoids is to impose conditions on
the homotopy group sheaves. We will first treat this idea in full generality,
introducing the notion of {\em directory of Serre classes} (it is a
generalization to the relative or ``stack'' case, of the classical notion of
``Serre class of groups''). In \S 6  we will specify the directory
of Serre classes which is most useful for nonabelian Hodge theory, leading to
the {\em realm of very presentable $n$-stacks of groupoids}. The present general
treatment is just  intended to motivate the idea of defining a realm by imposing
conditions on the homotopy group sheaves, and in particular to motivate why we
choose certain types of conditions.

Fix $n$. A {\em directory of Serre classes} $\Pp$ on a site $\Gg$ is a
collection $\Pp _i(X)$ indexed by $i=0,1,2,\ldots n$ and by $X\in \Gg$, where:
\newline
--$\Pp _0(X)$ is a full subcategory of the category of sheaves of sets on
$\Gg /X$;
\newline
--$\Pp _1(X)$ is a full subcategory of the category of sheaves of groups
on $\Gg /X$; and
\newline
--for $i\geq 2$,  $\Pp _I(X)$ is a full subcategory of the category of
sheaves of
abelian groups on $\Gg /X$.

\begin{parag}
We require:
\newline
--(compatibility with restrictions)\,  that if $G\in \Pp _i (X)$ and if
$Y\rightarrow X$ is a morphism in $\Gg$, then $G|_{\Gg /Y} \in \Pp _i(Y)$; and
\newline
--(locality) \, that if $\{ U_{\alpha}\rightarrow X\}$ is a covering of $X$, and
if $G$ is a sheaf on $\Gg /X$ (of sets for $i=0$, of groups for $i=1$, of
abelian
groups for $i\geq 2$) such that for each $\alpha$, $G|_{\Gg /U_{\alpha}} \in \Pp
_i(U_{\alpha})$ then $G\in \Pp _i(X)$.

\newparag{obtain}
We obtain the {\em realm defined by $\Pp$} denoted
$\underline{M}^{\Pp}$, by setting
$$
\underline{M} ^{\Pp}(X)\subset n\underline{STACK} (X)
$$
equal to the subcategory of $n$-stacks of groupoids $T$ on $\Gg /X$ such
that $\pi
_0(T)$ lies in $\Pp _0(X)$, and for $Y\in \Gg /X$ and each basepoint $t\in
T(Y)$, the homotopy group sheaves $\pi _i(T|_{\Gg /Y}, t)$ lie in $\Pp _i(Y)$.
Using the conditions of  compatibility with restriction and locality, it is
immediate that this defines a full $n+1$-substack $\underline{M}^{\Pp}$.
\end{parag}

Recall from above certain of the closure properties that one can ask of a realm
$\underline{M}$:
\newline
--(closure under truncation \ref{closuretruncation}) \, if $T\in
\underline{M}(Z)$ then the truncation $\tau _{\leq k}(T)$ is in
$\underline{M}(Z)$;
\newline
--(closure under finite limits \ref{closurelimits}) \, we have $\ast _Z\in
\underline{M}(Z)$, and if $R\rightarrow S \leftarrow T$ is a diagram in
$\underline{M}(Z)$ then $R\times _ST$ is in $\underline{M}(Z)$;
\newline
--(closure under extension \ref{closureextension}) \, if $E\rightarrow T$ is a
morphism of $n$-stacks of groupoids over $Z$ such  that $T$ is in
$\underline{M}(Z)$ and for every $Z'\in \Gg /Z$ and every point $t:\ast
_{Z'}\rightarrow T|_{\Gg /Z'}$, the fiber product
$$
F_t:= \ast _{Z'} \times _{T|_{\Gg /Z'}}E|_{\Gg /Z'}
$$
is in $\underline{M}(Z')$, then $E$ is in $\underline{M}(Z)$.

(Note that the closure under finite limits as we have stated it is equivalent to
the condition of \ref{closurelimits} because any limit over a finite
$1$-category can be broken down into a series of fiber products.)

\begin{parag}
\label{rmktrunc}
If $\Pp$ is a directory of Serre classes, then the associated realm
$\underline{M}^{\Pp}$ is obviously closed under truncation.
\end{parag}

\begin{lemma}
Suppose $\Pp$ is a directory of Serre classes which satisfies the following
properties (for every $X\in \Gg$):
\newline
---each $\Pp _i(X)$ is closed under extensions for $i\geq 0$;
\newline
---for $i\geq 1$, $\Pp _i(X)$ is closed under cokernels of maps from objects in
$\Pp _{i+1}(X)$;
\newline
---for $i\geq 2$, $\Pp _i(X)$ is closed under kernels of maps toward objects in
$\Pp _{i-1}(X)$;
\newline
---$\Pp _1(X)$ is closed under stabilizers of actions on objects in $\Pp _0(X)$;
\newline
---$\Pp _0(X)$ is closed under quotients by actions of groups in $\Pp _1(X)$.
\newline
Then the associated
realm $\underline{M}^{\Pp}$ is closed under extensions.
\end{lemma}
{\em Proof:}
Use the long exact sequence of homotopy group objects for a fibration sequence.
\eop

\begin{lemma}
Suppose $\Pp$ is a directory of Serre classes which satisfies the following
properties (for every $X\in \Gg$):
\newline
---$\Pp _i(X)$ is closed under kernel (or equalizer) for $i\geq 0$;
\newline
---the cokernel of a map in $\Pp _i(X)$ is in $\Pp _{i-1}(X)$, for $i\geq 1$;
\newline
---$\Pp _i(X)$ is closed under extension for $i\geq 0$.
\newline
Then the associated realm $\underline{M}^{\Pp}$ is closed under finite
limits. \end{lemma}
{\em Proof:}
We have a long exact sequence for the homotopy group sheaves of a fiber product
$R\times _ST$ as follows:
$$
\ldots \rightarrow \pi _i(R\times _ST, (r,t))\rightarrow
\pi _i(R,r)\times \pi _i(T,t) \rightarrow
$$
$$
\pi _i(S,s)\rightarrow \pi _{i-1}(R\times _ST, (r,t))\rightarrow \ldots
$$
which ends with the exact sequence
$$
\pi _2(R,r)\times \pi _2(T,t) \rightarrow \pi _2(S,s) \rightarrow \pi _1(R\times
_ST, (r,t))\rightarrow
$$
$$
\pi _1(R,r)\times \pi _1(T,t) \tworightarrows \pi _1(T,t)
$$
where the last part means taking the equalizers of the two morphisms. For the
next part of the exact sequence, we have an exact sequence
$$
\pi _0(R\times _ST) \rightarrow \pi _0(R)\times \pi _0(T) \tworightarrows
\pi _0(S),
$$
in other words the image of the first morphism is the equalizer of the second
pair of morphisms. Finally, given a point $(r,t)$ of $R\times _ST$, the
fiber of
$$
\pi _0(R\times _ST) \rightarrow \pi _0(R)\times \pi _0(T)
$$
over the image point $([r], [t])$ is an orbit of the group $\pi _1(S,s)$.

In the
above discussion, $s\in S$ denotes the image of $r$ and $t$ (they are the same).

From these long exact sequences and our hypotheses, one obtains the fact
that for
$R,S,T$ in $\underline{M}^{\Pp}$, the fiber product $R\times _ST$ is in
$\underline{M}^{\Pp}$.
\eop

The following theorem gives a converse to the construction $\Pp \mapsto
\underline{M}^{\Pp}$.

\begin{theorem}
Suppose $\underline{M} \subset n^{\rm gr}\underline{STACK}$ is a realm
which is closed under truncation, extension and finite limits.
Then there exists a directory of Serre classes $\Pp$ such that
$\underline{M}=\underline{M}^{\Pp}$.
\end{theorem}
{\em Proof:}
For $i=1$ (resp. $i\geq 2$) define $\Pp _i(Z)$ to be the set of sheaves of
groups
(resp. abelian groups) $G$ on $\Gg /Z$ such that $K(G,i)\in \underline{M}(Z)$.
Let $\Pp _0(Z)$ be the set of sheaves of sets $F$ on $\Gg /Z$ which (when
considered as $n$-stacks) are in $\underline{M}(Z)$. The fact that
$\underline{M}$
is a stack implies the locality and restriction conditions, so $\Pp$ is a
directory of Serre classes.

Suppose $T\in \underline{M}(Z)$, and suppose $t\in T(Z)$ is a basepoint. Then,
considering $t$ as a map $\ast _Z\rightarrow T$, we have
$$
\ast _Z \times _{\tau _{\leq i-1}T} \tau _{\leq i}T \cong K(\pi _i(T, t), i).
$$
Thus, the axioms of compatibility with truncations and fiber products together
with the definition of $\Pp _i(Z)$ imply that $\pi _i(T, t)\in \Pp _i(Z)$.
We obtain (by doing the argument over $Z'$) that for any $Z'\rightarrow Z$ and
basepoint $t\in T(Z')$, the $\pi _i(T|_{\Gg /Z'}, t)$ is in $\Pp _i(Z')$.  Note
also that (even without a basepoint) $\tau _{\leq 0}(T)= \pi _0(T)$ is in
$\underline{M}$ so again by definition $\pi _0(T)\in \Pp _0(Z)$. These
all together show that $T\in \underline{M}^{\Pp}(Z)$.

Suppose now on the other hand that $T\in \underline{M}^{\Pp}(Z)$. We show by
induction on $k$ that $\tau _{\leq k}T \in \underline{M}(Z)$. Note first that by
definition of $\Pp _0(Z)$,  $\tau _{\leq 0} T = \pi _0(T)$ is in
$\underline{M}(Z)$. Next, suppose that we know that $\tau _{\leq k-1}T\in
\underline{M}(Z)$. Apply the condition of closure under extensions, to
the morphism
$$
\tau _{\leq k}T \rightarrow \tau _{\leq k-1}T.
$$
For any basepoint $t:\ast _{Z'} \rightarrow \tau _{\leq k-1}T|_{\Gg /Z'}$,
the fiber $F_t$ is locally over $Z'$ equivalent to something of the form $K(\pi
_k(T|_{\Gg /Z'}, t), k)$. By definition of $\Pp _k$ and by the hypothesis that
$\pi _k(T|_{\Gg /Z'}, t)\in \Pp _k(Z')$, we get that $F_t\in \underline{M}(Z')$.
We can now imply the condition of closure under extensions to conclude that
$\tau _{\leq k}T\in \underline{M}(Z)$. The inductive statement for $k=n$ gives
that $T\in \underline{M}(Z)$.

We have now shown that $\underline{M}(Z)= \underline{M}^{\Pp}(Z)$ for all $Z\in
\Gg$.
\eop

\oldsubnumero{A few examples}

Define the following directories of Serre classes on the site $Sch /k$ over a
field $k$ of characteristic zero. They will be denoted $FlatLoc$ and $1ConLoc$
respectively. For $i\geq 2$ we put
$$
FlatLoc_i(X)= 1ConLoc_i(X):= \{ \mbox{loc. free sheaves} \}
$$
equal to the set of finite rank locally free sheaves of $\Oo _X$-modules
(considered as sheaves of abelian groups on $Sch /X$). For $i=0$ put
$$
FlatLoc_0(X)= 1ConLoc_0(X):= \{ \ast _X\}
$$
equal to the class consisting only of the trivial sheaf of sets $\ast _X$.
For $i=1$ put
$$
1ConLoc _1(X) := \{ \ast _X\}
$$
equal to the class consisting only of the trivial sheaf of groups $\ast _X$,
whereas
$$
FlatLoc _1(X) := \{ \mbox{flat affine group schemes} \}
$$
is the class of sheaves of groups represented by flat linear group schemes over
$X$; recall (cf \ref{linear}) that ``linear'' means a group scheme embedded as a
closed subgroup scheme of some $GL(E)$ with $E$ locally free of finite rank over
$X$.

These directories of Serre classes give rise to realms
$$
\underline{FL}:= \underline{M}^{FlatLoc},
$$
$$
\underline{AL}:= \underline{M}^{1ConLoc}.
$$

Similarly we can define (jumping ahead a  bit) the directories of Serre classes
$FlatVect$ and $1ConVect$, by the same classes of sheaves of sets or sheaves of
groups as above in degrees $i=0,1$, but allowing all of the vector sheaves (cf
\S 6) in degrees $i\geq 2$. These give rise to the realms
$$
\underline{FV}:= \underline{M}^{FlatVect},
$$
$$
\underline{AV}:= \underline{M}^{1ConVect}.
$$

\subnumero{Formalization of Postnikov-type arguments}

In this section we will fix a realm $\underline{M}$. If $X$ and $T$ are
$\ngr$-stacks (on $\Gg /Z$ for example), we would like to know when
$\underline{Hom}(X,T)$ lies in $\underline{M}(Z)$. We will approach this
question in the standard way  by using the  Postnikov tower for $T$.

For $T$ running through another realm $\underline{R}$, this will (sometimes)
answer the question of when $([X], \underline{R}, \underline{M})$ is
well-chosen.

It turns out that even in looking at the above case, a little bit more
generality comes into play. In order to get to the general form of the statement
we introduce the following notion. In this section, for simplicity we work only
with $\ngr$-stacks (some things we say would work in greater generality but the
coefficient stacks of cohomology must be $\ngr$-stacks in order to have a
Postnikov tower).

Suppose $\underline{M} \subset \ngr \underline{STACK}$ is a realm. We say that a
morphism of $\ngr$-stacks $U\rightarrow V$ is {\em of type $\underline{M}$} if
it satisfies the following condition: for every $Z\in \Gg$ and every point $v\in
V_0(Z)$, the fiber product
$$
Z\times _{V|_{\Gg /Z}} U|_{\Gg /Z}
$$
(which is considered as an $\ngr$-stack on $\Gg /Z$)
is in $\underline{M}(Z)$.

We say that an $\ngr$-stack $T$ is {\em of type $\underline{M}$} if the
structural morphism $T\rightarrow \ast$ is of type $\underline{M}$. If $\Gg$ has
a final object $\ast$ then this condition is equivalent to saying that $T\in
\underline{M}(\ast )$.

One can rewrite the condition of ``closure under extensions''
\ref{closureextension} as saying that if $U\rightarrow V$ is a morphism
of $\ngr$-stacks on $\Gg /X$ of type $\underline{M}|_{\Gg /X}$, and if $V\in
\underline{M}(X)$ then $U\in \underline{M}(X)$.

Throughout this section we will generally assume that $\underline{M}$ is closed
under extensions; and often we shall also assume closure under finite limits. In
any case these hypotheses will be explicitly mentionned in the statements.

\begin{lemma}
\label{compose}
Suppose $\underline{M}$ is closed under extensions.
Then the composition of two morphisms of type
$\underline{M}$ is again of type $\underline{M}$.
In particular if $T\rightarrow R$ is a morphism of type $\underline{M}$ and
if $R$ is of
type $\underline{M}$ then $T$ is of type $\underline{M}$.
\end{lemma}
{\em Proof:}
Left to the reader.
\eop

For the remainder of this section we will fix a realm $\underline{M}$ and
consider the following general situation: we have a diagram
$$
T\stackrel{f}{\rightarrow} R\rightarrow A\stackrel{p}{\rightarrow} B
$$
of $\ngr$-stacks. The {\em stack of sections of $R$ over $A$
relative to $B$} is defined as
$$
\underline{\Gamma}(A/B, R) := \underline{Hom}(A/B, R) \times
_{\underline{Hom}(A/B, A)}\{ 1_A\} .
$$
The morphism $f$ then induces a morphism on stacks of sections
$$
\underline{\Gamma} (A/B, f):
\underline{\Gamma} (A/B, T)\rightarrow \underline{\Gamma}(	A/B, R).
$$

We would like to know when this induced morphism
$\underline{\Gamma}(A/B, f)$ is of type $\underline{M}$. Using the previous
lemma,
if we know for some other reason that $\underline{\Gamma}(	A/B, R)$ is of type
$\underline{M}$, and if we can show that $\underline{\Gamma}(A/B, f)$ is of type
$\underline{M}$, then it will follow that
$\underline{\Gamma}(	A/B, T)$ is of type
$\underline{M}$.

Similarly, note that if $B=\ast$ and $T= T'\times A$ then
$$
\underline{\Gamma}(
A/B, T)= \underline{Hom}(A, T').
$$
Thus the statement about stacks of sections includes as a special case the
statement about $\underline{Hom}(-, -)$ being of type $\underline{M}$.

On the other hand, since we are going to be discussing Postnikov truncation
morphisms of the form $\tau _{\leq k} T\rightarrow \tau _{\leq k-1}T$,
it seems natural to look in general at a morphism $T\rightarrow R$.
The consideration of a relative $A/B$ in the domain is motivated by the fact
that many of our constructions related to Hodge theory will be relative (for
example, the Hodge filtration is an object relative to ${\bf A}^1$).

To close out our preliminary statements, we will avoid making too much of a
general treatment of degree $1$ nonabelian cohomology. On the one hand
particular cases of this were discussed at length in \cite{Simpson}, on the
other
hand doing this for general presentable group sheaves as coefficients presents
technical difficulties which, while treatable in the case of the de Rham
cohomology $H^1(X_{DR}, G)$, are not yet treated for a general formal category.
Thus, we will assume that we have already understood the degree $1$ case in the
sections  $\underline{\Gamma}(	A/B, R)$, so we assume that the
morphism $f: T\rightarrow R$ is relatively $1$-truncated. (Concretely if one
wants to look at $\underline{\Gamma}(	A/B, T)$, set $R:= \tau _{\leq 1}(T/A)$
(cf the definition of relative truncation below)
and suppose that we already understand $\underline{\Gamma}(	A/B, R)$, then apply
the subsequent discussion to climb up to $\underline{\Gamma}(	A/B, R)$.)

\begin{hypotheses}
The realm $\underline{M}\subset \ngr \underline{STACK}$ is closed under
extensions. We have a diagram of morphisms of $\ngr$-stacks on $\Gg$
$$
T\stackrel{f}{\rightarrow} R\rightarrow A\stackrel{p}{\rightarrow} B
$$
and we look at the morphism
$$
\underline{\Gamma}(A/B, f):
\underline{\Gamma} (A/B, T)\rightarrow \underline{\Gamma}(	A/B, R).
$$
We assume that the morphism $f: T\rightarrow R$ is relatively $1$-connected i.e.
that $\tau _{\leq 1}(T/R)=R$ (cf the definition of relative truncation below).
\end{hypotheses}

\oldsubnumero{Reduction to the Eilenberg-MacLane case}

In order to test whether $\underline{\Gamma}(A/B, f)$ is of type
$\underline{M}$,
we have to test over all objects $X$ mapping into $\underline{\Gamma} (A/B,
R)$.
A map from $X$ into this space of sections consists of a pair $(p,\eta )$ where
$p:X\rightarrow B$ is a morphism and $\eta$ is a section
$$
\eta : A_X:= A\times _BX \rightarrow R_X := R\times _BX .
$$

Given $(p,\eta )$ define
$$
F_{p,\eta }:= A_X \times _{R_X} T_X
$$
where $T_X := T\times _BX$ also.

\begin{lemma}
\label{identifyfiber}
With the above notations,
the fiber of $\underline{\Gamma}(A/B, f)$ over $(p,\eta )$ can
be identified by the formula
$$
\underline{\Gamma} (A/B, T)\times _{ \underline{\Gamma}(	A/B, R)}X
= \underline{\Gamma} (A_X/X, F_{p,\eta }).
$$
\end{lemma}
\eop

\begin{parag}
\label{postnik1}
We will apply the above discussion inductively using the Postnikov truncation.
First, if $U\rightarrow V$ is a morphism of topological spaces, define the {\em
relative Postnikov truncation}
$$
U\rightarrow \tau _{\leq k}(U/V) \rightarrow V
$$
by attaching cells to $U$ covering cells in $V$ so as to remove the homotopy
groups of the homotopy-fibers in degrees $>k$. Applying the correspondence
between $n$-groupoids and $n$-truncated spaces, we get a relative Postnikov
truncation for morphisms of $n$-groupoids.
If $T\rightarrow R$ is a morphism of $\ngr$-stacks,
define the {\em relative Postnikov truncation}
$$
T\rightarrow \tau _{\leq k}(T/R)\rightarrow R
$$
by setting
$$
\tau ^{\rm pre}_{\leq k}(T/R)(Z):= \tau _{\leq k}(T(Z)/R(Z))
$$
for each $Z\in \Gg$,
and by putting $\tau _{\leq k}(T/R)$ equal to the stack associated to
the prestack $\tau ^{\rm pre}_{\leq k}(T/R)$.

\newparag{postnik2}
The above form a tower, in that we have
$$
T =\tau _{\leq n}(T/R) \rightarrow \tau _{\leq n-1}(T/R) \rightarrow
\ldots
$$
$$
\ldots \rightarrow \tau _{\leq 1}(T/R) \rightarrow \tau _{\leq 0}(T/R)
\rightarrow R.
$$
Sometimes it will be convenient to establish the convention
$\tau _{\leq -1}(T/R):=R$.

\newparag{postnik3}
The relative Postnikov truncation is compatible with base-change over $R$:
if $R'\rightarrow R$ is a morphism and if we set $T':= T\times _RR'$ then
$$
\tau _{\leq k}(T'/R') = \tau _{\leq k} (T/R) \times _RR'.
$$

\newparag{postnik4}
Given a morphism $T\rightarrow R$, we say that $T/R$ is {\em relatively
$k$-truncated} (resp. {\em relatively $k$-connected})
if $T\cong \tau _{\leq k}(T/R)$ (resp. $\tau _{\leq k}(T/R)\cong R$).
An {\em Eilenberg-MacLane morphism of degree $k$} is a morphism $T\rightarrow R$
such that $T/R$ is relatively $k$-truncated and relatively $k-1$-connected.
\end{parag}

We now return to the previous discussion, where $A\rightarrow B$ is a morphism
and $T\rightarrow R\rightarrow A$ is a morphism of objects over $A$ (all objects
here are $\ngr$-stacks  on $\Gg$). Fix a morphism $p: X\rightarrow B$
with $X\in \Gg$, and denote by subscript $X$ the fiber-products over $B$ with
$X$. Note that
$$
\tau _{\leq k}(T_X /R_X) = \tau _{\leq k}(T/R)_X
$$
is again the Postnikov truncation (by \ref{postnik3}. Fix  $0\leq k\leq n$ and
look at a section
$$
\eta ^k: A_X \rightarrow \tau _{\leq k-1}(T_X /R_X)
$$
(recall that the latter maps to $R_X$ which maps to $A_X$ and ``section''
means a
section of the map to $A_X$). Set
$$
F^k_{p,\eta ^k} := \tau _{\leq k}(T_X /R_X)\times _{\tau _{\leq k-1}(T_X /R_X)}
A_X.
$$
Note that
$$
F^k_{p,\eta ^k} \rightarrow A_X
$$
is an Eilenberg-MacLane morphism of degree $k$. We call it the {\em
Eilenberg-MacLane morphism of degree $k$ coming from $\eta ^k$}.

The following proposition gives the reduction to consideration of these
Eilenberg-MacLane morphisms.

\begin{proposition}
\label{test}
Suppose that $\underline{M}$ is a realm closed under extension. Suppose
$A\rightarrow B$ is a morphism  of $\ngr$-stacks and $T\rightarrow
R\rightarrow A$
is a morphism of $\ngr$-stacks over $A$.
Suppose that for all $p:X\rightarrow B$, all $0\leq k \leq n$, and
all Eilenberg-MacLane morphisms
$$
F^k_{p,\eta ^k}\rightarrow A_X
$$
of degree $k$ arising from sections
$$
\eta ^k: A_X \rightarrow \tau _{\leq k-1}(T_X /R_X),
$$
we have that
$$
\underline{\Gamma}(A_X /X, F^k_{p,\eta ^k})\rightarrow X
$$
is of type $\underline{M}$. Then the morphism
$$
\underline{\Gamma}(A/B, T)\rightarrow \underline{\Gamma}(A/B, R)
$$
is of type $\underline{M}$.
\end{proposition}
{\em Proof:}
Using Lemma \ref{identifyfiber} and the hypotheses of the proposition, we obtain
that
$$
\underline{\Gamma}(A/B, \tau _{\leq k} (T/R) )\rightarrow
\underline{\Gamma}(A/B, \tau _{\leq k-1} (T/R) )
$$
is of type $\underline{M}$ for all $0\leq k \leq n$
(cf the convention \ref{postnik2} for
$k=0$). Now, applying Lemma \ref{compose} we obtain the claimed result.
\eop

In applying this proposition, one has sometimes to be careful to distinguish
which types of Eilenberg-MacLane morphisms can arise as $F^k_{p,\eta ^k}$ as
fibers over sections $\eta ^k$. Then one has to treat the sections of these
fibrations. The case $k=1$ is the most difficult; this is the subject of
Giraud's book \cite{Giraud} and as stated previously we try to avoid that as
much as possible by assuming that $f$ is relatively $1$-connected.

\oldsubnumero{Local systems}

Before treating the Eilenberg-MacLane case, we introduce the notion of ``local
system''. This terminology is employed in the same way as in \cite{Simpson}(iv),
(v). It should {\em not} be confused with the notion of local system over the
underlying topological space of a scheme (which was used uniquely in \S 2);
in our
present terminology if $X$ is a scheme then a ``local system'' over $X$ just
means a sheaf of sets or groups on $X$.

\begin{definition}
\label{localsystem}
Suppose $R$ is an $\ngr$-stack. A {\em local system (of sets)} over $R$ is a
relatively $0$-connected morphism of $\ngr$-stacks $L\rightarrow R$. This is
equivalent to the data, for each $X\in \Gg$, of a local system in the usual
topological sense over the space corresponding to the $n$-groupoid $R(X)$.

A {\em local system of groups} is a local system $L\rightarrow R$ together
with an
$R$-morphism $L\times _RL\rightarrow L$ satisfying the usual axioms expressed in
terms of commutative diagrams. Note that the identity is a section $R\rightarrow
L$. A {\em local system of abelian groups} is a local system of groups
satisfying
the commutativity axiom too.
\end{definition}

In interpreting the above definition, note that if $L$ and $L'$ are two local
systems of sets over $R$ and if $f,g:L\rightarrow L'$ are two $R$-morphisms,
then if $f$ and $g$ are homotopic, the homotopy between them is unique (as are
all higher homotopies). Thus there is no need for ``homotopy coherence'' in the
definition of a group structure for example. Another way of saying this is that
the $n+1$-category of local systems over $R$ is actually a $1$-category.

Pullback along the morphism $R\rightarrow \tau _{\leq 1}R$ induces an
equivalence between the category of local systems over the $1$-stack $\tau
_{\leq 1}R$, and the category of local systems over $R$. The former is known to
form a topos.

Suppose $T\rightarrow R$ is a morphism of $\ngr$-stacks, and suppose that
$\rho : R\rightarrow T$ is a section. Then we obtain the {\em local systems of
homotopy groups} denoted
$$
\pi _i(T/R, \rho )\rightarrow R.
$$
For $i=0$ we don't need to fix a section.

If the morphism $T\rightarrow R$ is relatively $1$-connected, then we don't need
to choose a section (and in fact a globally defined section might not exist),
and we obtain local systems of homotopy groups
$$
\pi _i(T/R)\rightarrow R.
$$

If $n\geq 1$ and if $L\rightarrow R$ is a local system of groups (abelian for
$n\geq 2$) then we obtain the corresponding {\em Eilenberg-MacLane stack}
$$
K(L/R, n)\rightarrow R.
$$
It comes with a section $0: R\rightarrow K(L/R, n)$, and is characterized by the
conditions
$$
\pi _i(K(L/R, n)/R, 0) = \{ 0\} , \;\;\;\; i\neq  n,
$$
$$
\pi _n(K(L/R, n)/R, 0) = L.
$$

\begin{parag}
\label{cohomology}
Suppose $T\rightarrow R$ is a morphism of $\ngr$-stacks, and suppose
that $L\rightarrow T$ is a local system. Then we define
$$
H^i(T/R, L):= \pi _0(\underline{\Gamma} (T/R, K(L/T, i))/R).
$$
This exists for any $i\geq 0$ if $L$ is a local system of abelian groups, and
for $i=0,1$  if $L$ is a local system of gropus. It is a local system (of
abelian groups, groups or eventually sets) over $R$.
\end{parag}

In general the cohomology is different from the higher direct image coming from
the morphism of topoi of local systems corresponding to $T\rightarrow R$. The
morphism of topoi will give rise to the cohomology for $\tau _{\leq 1} T
\rightarrow \tau _{\leq1}R$; thus one could use this type of formulation in case
$R$ and $T$ are $1^{\rm gr}$-stacks.

\oldsubnumero{Treating the Eilenberg-MacLane case}

From now on we make the additional hypothesis that $\underline{M}$ is closed
under finite limits.

Recall the assumption that $f:T\rightarrow R$ is relatively
$1$-connected, implying that the relative homotopy group objects
$$
L^i:=\pi _i(T/R)\;\;\;\; (i\geq 2)
$$
are well defined as local systems of abelian
groups on $R$.
Furthermore we have the classical ``classification'' of Eilenberg-MacLane
fibrations: the morphism
$$
\tau _{\leq k}(T/R)\rightarrow \tau _{\leq k-1}(T/R)
$$
is classified by a morphism (section over $R$)
$$
\tau _{\leq k-1}(T/R)\rightarrow K(L^k/R, k+1)
$$
in the sense that we have the following cartesian diagram:
$$
\begin{array}{ccc}
\tau _{\leq k}(T/R)&\rightarrow & R \\
\downarrow && \downarrow \\
\tau _{\leq k-1}(T/R)&\rightarrow &K(L^k/R, k+1)\\
\end{array}
$$
where the bottom arrow is the classifying map and the right vertical arrow is
the basepoint section.

Because of the occurrence of $k+1$ in the above diagram, if we want a nice
statement staying within the world of $\ngr$-stacks we make the seemingly
somewhat artificial hypothesis that $f$ is $n-1$-truncated in the following
statement. This is related to our standard observation that the truncation level
of the base should be one more than that of the map. In practice this doesn't
pose any problem because we can just increase $n$ (our realms are really going
to be defined for any $n$).

\begin{proposition}
\label{postniprop}
Suppose that the realm $\underline{M}$ is closed under extensions and under
finite limits. Suppose that
$$
T\stackrel{f}{\rightarrow} R\rightarrow A\stackrel{p}{\rightarrow} B
$$
is a diagram with $f$ relatively $1$-connected and relatively
$n-1$-truncated. Put
$L^k := \pi _k(T/R)$ as a local system on $R$. Suppose that for every $X\in \Gg$
and every section $\eta : A_X \rightarrow R$, and for every $2\leq k\leq
n-1$, the
$\ngr$-stack
$$
\underline{\Gamma}(A_X/X, K(\eta ^{\ast}L^k /A_X, k+1))
$$
on $\Gg /X$ is in $\underline{M}(X)$.
Then $\underline{\Gamma}(A/B, f)$ is of type $\underline{M}$.
\end{proposition}
{\em Proof:}
Use the same notations as in Proposition \ref{test}. Note that any section $\eta
^k$ projects to a section $\eta$ as in the present statement.  We have a
cartesian
diagram (deduced from the previous one above)
$$
\begin{array}{ccc}
F^k_{p,\eta ^k}&\rightarrow & A_X \\
\downarrow && \downarrow \\
A_X&\rightarrow &K(\eta ^{\ast} L^k/A_X, k+1)\\
\end{array}
$$
where the right vertical arrow is the zero-section and the bottom arrow is the
classifying section. This gives a cartesian diagram on stacks of sections
$$
\begin{array}{ccc}
\underline{\Gamma}(A_X/X, F^k_{p,\eta ^k})&\rightarrow &
\underline{\Gamma}(A_X/X,A_X) \\
\downarrow && \downarrow \\
\underline{\Gamma}(A_X/X,A_X)&
\rightarrow &\underline{\Gamma}(A_X/X,K(\eta ^{\ast} L^k/A_X, k+1))\\
\end{array}
$$
which we rewrite (in view of the fact that the sections of $A_X$ itself are
trivial) as
$$
\begin{array}{ccc}
\underline{\Gamma}(A_X/X, F^k_{p,\eta ^k})&\rightarrow &
X \\
\downarrow && \downarrow \\
X&\rightarrow &\underline{\Gamma}(A_X/X,K(\eta ^{\ast} L^k/A_X, k+1))\\
\end{array} .
$$
This is a cartesian diagram of $\ngr$-stacks on $\Gg /X$. The fact that
$\underline{M}(X)$ is closed under finite limits implies in particular that the
final object of $\ngr\underline{STACK}(X)$, which is $X$, is in
$\underline{M}(X)$. The hypothesis of the present proposition is that the stack
of sections on the bottom right is in $\underline{M}(X)$, so again by closure of
$\underline{M}(X)$ under finite limits, we get that
$$
\underline{\Gamma}(A_X/X, F^k_{p,\eta ^k})\in \underline{M}(X).
$$
This verifies the necessary hypotheses to put into Proposition \ref{test} and by
that proposition we obtain the conclusion
that $\underline{\Gamma}(A/B, f)$ is of type $\underline{M}$.
\eop

In order to get a nice statement, we rewrite the condition of the above
proposition entirely in terms of usual abelian cohomology (using the derived
category of complexes). In order to shorten a whole bunch of notation involving
shifting a complex, taking the canonical truncation and then applying
Dold-Puppe, we introduce the following notation (as was done in Illusie
\cite{Illusie}). If $C^{\cdot}$ is a complex of sheaves of abelian groups on
$X\in \Gg$  supported in positive degrees, then we set
$$
{\bf K}(C^{\cdot}/X, m)
$$
equal to the $m$-stack on $\Gg /X$ obtained by shifting the complex so that it
is supported in degrees $\geq -m$; taking the canonical truncation so that it is
supported in the interval $[-m, 0]$; and then applying the Dold-Puppe
construction (cf \cite{Illusie}) to obtain a presheaf of spaces. Finally, take
the associated stack (noting that if the complex consisted of injectives then
the resulting object is already a stack).

Denote by ${\bf R} p_{X,\ast} (-)$ the ``higher derived direct image functor''
for $p_X: A_X \rightarrow X$. It is defined by the formula
$$
\underline{\Gamma}(A_X/X, K(\eta ^{\ast}L^k /A_X, k+1))
= {\bf K}({\bf R} p_{X,\ast} (\eta ^{\ast} L^k)/X, k+1).
$$
As was stated above, this is not in general a higher derived direct image for a
morphism of topoi; that will however be the case if $A_X$ is $1$-truncated.

We can rewrite the above proposition as follows.

\begin{corollary}
\label{postnicor}
Suppose that the realm $\underline{M}$ is closed under extensions and under
finite limits. Suppose that
$$
T\stackrel{f}{\rightarrow} R\rightarrow A\stackrel{p}{\rightarrow} B
$$
is a diagram with $f$ relatively $1$-connected and relatively $n-1$-truncated.
Put
$L^k := \pi _k(T/R)$ as a local system on $R$. Suppose that
for every $X\in \Gg$
and every section $\eta : A_X \rightarrow R$, and for every $2\leq k\leq
n-1$, the
$\ngr$-stack
$$
{\bf K}({\bf R} p_{X,\ast} (\eta ^{\ast} L^k)/X, k+1)
$$
on $\Gg /X$ is in $\underline{M}(X)$.
Then $\underline{\Gamma}(A/B, f)$ is of type $\underline{M}$.
\end{corollary}
\eop

\begin{exercise}
Rewrite the above discussion in terms of a ``cartesian family'' $A\rightarrow
\underline{B}$ as defined in \S 5 below.
\end{exercise}

\subnumero{The non-groupoid case: extended directories of Serre classes}

It can be useful to consider $n$-stacks not necessarily of groupoids.
We can make a definition like that of ``directory of Serre class'' but
appropriate for this situation. This subsection is very optional.

Recall that the {\em interior} of an
$n$-category $A$ is the maximal sub-$n$-category which is an $n$-groupoid; in
other words, it is the $n$-groupoid consisting of all $i$-morphisms which
are invertible up to equivalence in $A$. It is denoted $A^{\rm int}$.

Recall that $A_{1^k/}$ denotes the $n-k$-category of
$k$-arrows in $A$. For example with $k=1$ $A_{1/}$ is the $n-1$-category of
arrows in $A$, and inductively $A_{1^k/}= (A_{1^{k-1}/})_{1/}$.
Recall that if $x,y$ are objects of $A$ then we denote by $A_{1/}(x,y)$ the
$n-1$-ctaegory of arrows from $x$ to $y$. We can extend this notation in the
following way. Suppose $u,v$ are $k-1$-arrows, sharing a common source and
target. Then let
$$
A_{1^k/}(u,v)\subset A_{1^k/}
$$
be the full sub-$n-k$-category of $k$-arrows whose source and target are
respectively $u$ and $v$.

The above notations extend to the case of $n$-stacks.

An {\em extended directory of Serre classes} is a collection
$$
\Pp = \{ \Pp _{i,j}(X),\;\;\; 0\leq i \leq j \leq n, \;\;\; X\in \Gg\}
$$
where for each $i$, $\Pp _{i,\cdot -i} (-)$ is a directory of Serre classes for
$n-i$-stacks on $\Gg$.  Let $\Pp _{i/}$ denote this directory of Serre classes.
We define the full substack
$$
\underline{M}^{\Pp} \subset n\underline{STACK}
$$
as follows. An $n$-stack $A$ on $\Gg /X$ is in $\underline{M}^{\Pp}(X)$
if and only if the $n$-stack of groupoids $A^{\rm int}$ is in
$\underline{M}^{\Pp _{0/}}(X)$, and for any $X'\rightarrow X$, any $1\leq i \leq
n$ and any two $i-1$-morphisms $u,v$ sharing the same sources and targets in
$A(X')$, we have that  the $n-i$-stack of groupoids  on $X'$
$$
(A_{1^i/}(u,v))^{\rm int} \;\;\;\; \mbox{is in}\;\;\;\;
\underline{M}^{\Pp _{i/}}(X').
$$

We leave as a problem for further study, to find natural conditions which
should be imposed on the $\Pp _{i,j}(X)$.

\numero{Cartesian families and base change}
\label{cartesianpage}

We define a notion of ``cartesian family of $n$-stacks'' parametrized by a base
$n+1$-stack. This gives an improved approach to nonabelian cohomology in a
relative situation, and allows us to formulate a base-change statement. The
notion
of cartesian family fits into a natural picture for defining a canonical fibrant
replacement for $nCAT$: thus, in some sense, a weak morphism $B\rightarrow nCAT$
is the same thing as a cartesian family over $B$. In other words, this says what
universal property is satisfied by $nCAT$ (or similarly
$n\underline{STACK}$). It
also allows us to define the ``arrow family'' for a general $n+1$-category.

We describe the idea of the notion of cartesian family first, then we give the
applications to base-change and the arrow family, and only at the end do we give
the complete technical discussion of the construction.

\subnumero{A canonical fibrant replacement for $nCAT$}

One of the main problems with the construction of the $n+1$-category
$nCAT$ given in \cite{Simpson}(x) is that the result is not fibrant. Thus,
when speaking of ``weak morphisms'' to $nCAT$ one must choose a fibrant
replacement $nCAT \rightarrow nCAT'$ and look at morphisms to $nCAT'$. The
process of choosing a fibrant replacement is not very concrete and it is
consequently hard to get a hold of what a morphism to $nCAT'$ really means.
This led to difficulties in \cite{Simpson}(xi) and
\cite{HirschowitzSimpson}.

Here, we will sketch a method for constructing a canonical fibrant replacement
which we denote
$$
{\bf u} : nCAT \rightarrow nFAM.
$$
The $n+1$-category $nFAM$ parametrizes ``cartesian families'' of
$n$-categories in a sense which will be explained precisely below. The morphism
${\bf u}$ comes from a universal cartesian family over $nCAT$.
This picture gives a reasonably good understanding (which could of course be
much improved with time) of the ``universal property'' of $nCAT$.

One can realize that the universal property of $nCAT$ might not be extremely
simple to describe, by making the preliminary observation that $nCAT$ is an
$n+1$-category, and indeed, families of $n$-categories are most naturally seen
as having parameter space which is an $n+1$-category. For example, when one
looks at a fibration of topological spaces where the fiber is $n$-truncated, the
fibration is pulled back from a fibration over the $n+1$-truncation of the base
space---in other words, the classifying space for fibrations with an
$n$-truncated fiber, is $n+1$-truncated. Keeping to the topological case, the
previous phrase can be made more precise by observing that for a given
$n$-truncated fiber $F$, the ``group'' (i.e. space with $1$-fold delooping
structure) $Aut(F)$ is again $n$-truncated, and its classifying space $B\,
Aut(F)$ is $n+1$-truncated. This latter space classifies fibrations with fiber
$F$. Specializing a bit more (to $n=1$), we get the following concrete and
well-known example: if $F$ is a $K(G,1)$ for a group $G$ then $B\, Aut(F)$ is a
$2$-truncated space with $\pi _2=Z(G)$ and $\pi _1= Out(G)$; here we recover
Giraud's the nonabelian $2$-cohomology picture \cite{Giraud} where
fibrations with
fiber a ``gerb'' for group $G$, are classified by a nonabelian $2$-cohomology
whose Postnikov tower includes $H^2(-, Z(G))$ and $H^1(-, Out(G))$.

The first guess for  a universal property of $nCAT$ would be to say that a
morphism $E\rightarrow nCAT'$ is the same thing as a morphism of $n$-categories
$F\rightarrow E$. However, in this picture, $E$ itself is an $n$-category and
from the previous paragraph, this doesn't look like the right thing to do.
Later on, we will be able to obtain a morphism $E\rightarrow nCAT'$ for certain
families $F\rightarrow E$, but in fact this will not work under all possible
circumstances. The required condition is probably related to the
``compatible with
change of base'' (ccb) condition of \cite{HirschowitzSimpson} also see above.
We restrict our attention to the case where the base is an $\ngr$-stack, in
which case the (ccb) condition is automatic  \cite{HirschowitzSimpson}.

Here is a better approach. We will look at an $n+1$-category (or even
$n+1$-precat) $E$ and try to define a reasonable type of ``family of
$n$-categories parametrized by $E$''. The basic idea is to say that such a
family should consist of an $n$-category $F_x$ for each object $x\in E_0$, plus
for two objects $x,y\in E_0$, a morphism of $n$-categories
$$
F_x \times E_{1/}(x,y) \rightarrow F_y.
$$
We of course need compatibility with composition together with higher-order
homotopy coherence, as shall be encoded in the actual definition below.
Before getting there, we  note that this idea is reasonably well in accord with
the discussion of topological classifying spaces a couple of paragraphs ago: if
$E$ has only one object $x$, a map from $E$ into $B\, Aut (F)$  (sending the
object $x$ to the unique object ``$F$'') should consist essentially of a map
$E_{1/}(x,x)\rightarrow  Aut (F)$ compatible with the composition law. In this
point of view, the level of truncation of all of the spaces in question is the
right one.  An added advantage of the resulting approach will be to avoid the
need for cumbersome additional conditions such as (ccb) when applying the
universal property. (Note however that we have to refer to that condition in
some of our proofs.)

Before going on, we describe briefly the replacement $nFAM$. It is defined by  a
universal property: we actually describe the {\em morphisms} from an
$n+1$-precat $E$ into $nFAM$. These morphisms are by definition the {\em
cartesian families} over $E$. We will now briefly describe the notion of {\em
precartesian family}; the notion of cartesian family is obtained by adding an
additional fibrancy condition which will be explained in detail later. In
practice, this fibrancy condition is imposed by choosing a fibrant replacement,
so morally speaking a morphism $E\rightarrow nFAM$ is obtained by having a
precartesian family over $E$.

As suggested by the prior introductory remarks, we proceed as follows. An
$n+1$-precat $E$ can be considered as a simplicial object in the category $nPC$
of $n$-precats; we denote this ``reduction in dimension'' by ${\bf r}E$. Thus
$({\bf r}E)_p:= E_{p/}$. Note in particular that the term of degree zero is the
discrete set $E_0$ of objects of $E$. Now a {\em precartesian family}
over $E$ is a morphism of simplicial objects in $nPC$,
$$
\Ff \rightarrow {\bf r} E,
$$
which satisfies a certain Segal-type condition (see \ref{precartesiandef}
below).
The $n$-category in the family corresponding to an object $x\in E_0$ is the
fiber of $\Ff _0$ over $x$. Denote this by $\Ff _0(x)$.
The precartesian condition of \ref{precartesiandef} puts into practice the idea
described a few paragraphs ago: for example on the first level, we require
(in obvious notations) that
the morphism
$$
\Ff _1(x,y) \rightarrow {\bf r}E_1(x,y) \times \Ff _0(x)
$$
be an equivalence; and then the morphism
$$
\Ff _1(x,y) \rightarrow \Ff _0(y)
$$
is what amounts to the morphism
$$
{\bf r}E_1(x,y) \times \Ff _0(x)\rightarrow \Ff _0(y)
$$
in the heuristic description above.

In \ref{correspondence} we will explain how to interpret a morphism of
$n$-categories $F\rightarrow E$ as a cartesian family over $E$, in the case
where
$E$ is an $n$-groupoid.

\subnumero{Nonabelian cohomology in a relative situation and base-change}

Before getting to the technical details, we give some Hodge-theoretic
motivation.
In \cite{Navarro-Aznar}, V. Navarro Aznar gave a definition of the ``family'' of
de Rham homotopy types, corresponding to a family of simply connected varieties
$X\rightarrow S$. He used an explicit way of
``integrating'' the sheafified de Rham complex to obtain a global d.g.a. defined
over the base $S$. He showed that the family of de Rham homotopy types admits a
``Gauss-Manin connection'' up to coherent homotopy. He calculated the homotopies
of coherence up to a preliminary level of coherence which was sufficient to
obtain a flat connection on the homotopy group sheaves; however, as this got
computationally complicated, he didn't treat the full set of higher-order
coherencies. His theory also works  for the nilpotent completion of the
fundamental group, in the case of a family of non-simply connected varieties.

We will consider this type of question in a context where the ``d.g.a.
approach''
to defining rational homotopy type is replaced by the shape-theoretic approach.
One of the results of our study (see \S 9 and \S 10 below) will be a
definition of
the Gauss-Manin connection on the homotopy type representing the shape, in the
simply connected case, with all higher-order homotopy coherencies. Thus, even in
this case we extend somewhat the result of \cite{Navarro-Aznar}. However,
we don't
show the conjectured compatibility of our definition with the definition of
\cite{Navarro-Aznar}, and until that is done one cannot strictly speaking say
that we generalize \cite{Navarro-Aznar}.

One of the interesting aspects about the construction of \cite{Navarro-Aznar} is
that it is compatible with ``base change'' i.e. with pulling back the family
along a morphism $S'\rightarrow S$. What we will explain in the present chapter,
is how to obtain a similar base-change property for the ``de Rham shape'' or
more generally, for any type of nonabelian cohomology.

\oldsubnumero{Base change for relative nonabelian cohomology}

We start with a preliminary version that doesn't require the use of ``cartesian
families''. Suppose $X\rightarrow S$ is a morphism of $n$-stacks, and suppose
that $S$ is an $n$-stack of groupoids; thus it is ``compatible with change of
base'' (ccb) cf \cite{HirschowitzSimpson}, also \S 3.5 above. Suppose $T$ is an
$n$-stack. Then we obtain the relative nonabelian cohomology
$$
\underline{Hom}(X/S, T)\rightarrow S.
$$

\begin{theorem}
\label{basechangeA}
If $S'$ is another $\ngr$-stack with morphism $S'\rightarrow S$ and if we set
$X':= S'\times _SX$, then
$$
\underline{Hom}(X'/S', T)= \underline{Hom}(X/S, T)\times _S S'.
$$
\end{theorem}
{\em Proof:}
This follows tautologically from the definition of $\underline{Hom}(X/S, T)$.
\eop

In subsubsection 5.7 below we will see how this works when these formulas are
considered as functors in the variable $T$. For this, it is more convenient to
change point of view and start with a cartesian family $X\rightarrow
\underline{B}$; thus we first need to know more precisely what a cartesian
family is.

\subnumero{The notion of cartesian family and construction of $nFAM$}

We now start in on the technical details of the notion of cartesian family.

\oldsubnumero{Piecewise fibrations}

An {\em indecomposable object} is a $W\in nPC$  with the
property that for any pushout diagram $B\leftarrow A \rightarrow C$, we have
$$
Hom (W, B\cup ^AC) = Hom (W, B\cup ^AC).
$$
(Quillen calls these the ``small projective generators'' \cite{Quillen}.)
The basic indecomposables are the representable $n+1$-precats of the form $h(M)$
for $M\in \Theta ^{n}$. Note that any $n$-precat is a colimit of
indecomposable ones.

To introduce the subsequent ideas, we say that a morphism $U\rightarrow E$ in
$nPC$ is  {\em piecewise fibrant} if, for every
indecomposable $W$ mapping to $E$, the morphism $U\times _EW\rightarrow W$ is
fibrant.

\begin{parag}
\label{glueing}
This has the following glueing property: if
$$
\begin{array}{ccccc}
U' & \leftarrow & U & \rightarrow & U'' \\
\downarrow && \downarrow && \downarrow \\
E' & \leftarrow & E & \rightarrow & E''
\end{array}
$$
is a diagram such that the vertical morphisms are piecewise fibrant, and both
squares are fiber-products (i.e. cartesian), then the morphism $$
U'\cup ^UU'' \rightarrow E' \cup ^E E''
$$
is also piecewise fibrant.  (The usual global notion of ``fibrant  morphism''
doesn't have this glueing property.)
\end{parag}

There is a companion property which fits into a lifting property. Given a
sequence
$$
V' \stackrel{i}{\rightarrow} V \rightarrow E,
$$
we say that $i$ is a {\em trivial cofibration universally over $E$} if
for every morphism $W\rightarrow E$, we have that
$$
V'\times _EW \rightarrow V \times _EW
$$
is a trivial cofibration. It suffices to test this on indecomposable $W$.

\begin{lemma}
\label{unilifting}
Suppose
$$
\begin{array}{ccc}
V ' & \rightarrow & F \\
\downarrow && \downarrow \\
V &\rightarrow & E
\end{array}
$$
is a diagram such that the vertical morphism on the left is a trivial
cofibration universally over $E$, and the morphism on the right is
piecewise fibrant. Then a lifting $V\rightarrow F$
exists. \end{lemma}
{\em Proof:}
Express $E$ as a colimit of indecomposables (viewed as ``cells'') and construct
the lifting over one cell at a time.
\eop

\begin{example}
In the case of simplicial sets, the notion of ``piecewise fibration'' is
the same
as that of ``fibration'' since one can test for being a fibration using only
horn-extensions over indecomposable objects (simplices) mapping into the
base. In
the categorical situation (say, with $n=1$) this is no longer true, as is
shown by
the following example. Consider the diagram
$$
\begin{array}{ccc}
\ast & \rightarrow & \ast \\
\downarrow && \downarrow \\
\overline{I} &\rightarrow & \overline{I}
\end{array} .
$$
The horizontal morphisms are the identity and the vertical morphisms are the
inclusion of the vertex $0$ in the ``interval'' $\overline{I}$. The vertical
morphism is a piecewise fibration. It is also a trivial cofibration, but lifting
doesn't occur; thus it is not fibrant. However, it is not universally a trivial
cofibration, so the lifting property of the previous lemma is not contradicted.
\end{example}

We don't actually use the above definitions {\em per se} in what follows, but
variants as the reader shall see.

\oldsubnumero{Precartesian and cartesian families}

If $E$ is an $n+1$-precat, we may
think of it as a simplicial object in the category of $n$-precats, denoted
${\bf r}(E)$ (the letter ``r'' is for ``reduction'' as in reduction of the
dimension). This is defined by
$$
{\bf r}(E)_p:= E_{p/}.
$$
Note that ${\bf r}(E)_0$ is a discrete set.
This gives a functor from the model category $(n+1)PC$ of $n+1$-precats, to
the model category of $n$-prestacks over $\Delta$ (with the
coarse topology on $\Delta$), which we denote by $nPS(\Delta )$.

If $p\in \Delta$ then, thinking of $p$ as being the ordered set
$$
p=\{ 0', 1',\ldots , p'\} ,
$$
let the {\em initial} morphism $\iota : 0\rightarrow p$ be the one which
sends the unique element $0'$ of $0$, to the first element $0'$ of $p$.

({\em Note:} it is in the  definition of ``initial'' that we specify the
``directions'' of our arrows: choosing instead the last element would result in
putting opposites in various places below, or equivalently changing our
conventions about directions of arrows, or again equivalently, looking at
contravariant rather than covariant functors into $nCAT$. As it is, the
convention is that $A_{1/}(x,y)$ is the space of arrows {\em from} $x$ and {\em
to} $y$.)

\begin{parag}
\label{precartesiandef}
Suppose $E$ is an $n+1$-precat. A {\em precartesian family} over
$E$ is a morphism of $n$-prestacks over $\Delta$
$$
F \rightarrow {\bf r}(E)
$$
with the following property (which we call the ``cartesian'' property):
for any $p\in \Delta$ and for the initial morphism $\iota : 0\rightarrow p\in
\Delta$, the diagram
$$
\begin{array}{ccc}
F_p & \rightarrow & F_0 \\
\downarrow && \downarrow \\
E_{p/} & \rightarrow & E_{0}
\end{array}
$$
is homotopy-cartesian, in that it induces a weak equivalence
$$
F_p \cong E_{p/} \times _{E_0} F_0.
$$
Note here that $E$ satisfies the ``constancy condition'' that $E_0$ is a set;
however we {\em do not} require this condition of $F_0$, and indeed the whole
point is that $F_0\rightarrow E_0$ is a family of $n$-precats indexed by the set
$E_0$. For $x\in E_0$ denote by $F(x)$ the fiber of $F_0$ over $x$. It is an
$n$-precat, and is to be considered as the ``element'' of the family $F$,
corresponding to the object $x\in E_0$.
\end{parag}

The fact that $E_0$ is a set means that the fiber product appearing above has a
correct homotopical meaning.

We say that a precartesian family $F\rightarrow {\bf r}(E)$ is {\em piecewise
fibrant} if for every indecomposable $n+1$-precat $W$ mapping to $E$, the
pullback
$$
F\times _{{\bf r}(E)}{\bf r}(W) \rightarrow {\bf r}(W)
$$
is a fibrant morphism of $n$-prestacks on $\Delta$.

This might be slightly
different from the corresponding definition just in the world of $n$-prestacks
on $\Delta$, because we ask the condition only over things of the form ${\bf
r}(W)$ for $W$ an indecomposable $n+1$-precat (although it might also be
equivalent, I don't know).  Note, in any case, that if $W$ is an indecomposable
$n+1$-precat then ${\bf r}(W)$ is an indecomposable $n$-prestack. The properties
\ref{glueing} and \ref{unilifting} still work here.

A {\em cartesian family} $F \rightarrow {\bf r}(E)$ is a precartesian family
that is piecewise fibrant.

Notice as a first consequence that the fibers $F(x)$ for $x\in E_0$, are fibrant
$n$-categories. Also, the morphisms of $n$-precats
$$
F_p \rightarrow {\bf r}(E)_p = E_{p/}
$$
are piecewise fibrant. Thus the homotopy-cartesian diagrams appearing in the
definition of ``precartesian'' involve piecewise fibrant vertical morphisms.

The basic part of the structure of a cartesian family are the diagrams
$$
F(x) \times E_{1/}(x,y)\stackrel{\cong}{\leftarrow} F_1(x,y) \rightarrow F(y).
$$
Thus, this notion corresponds to what was said previously. The rest of the
simplicial structure corresponds to the compatibility with composition and
higher homotopy-coherence for this structure.

\oldsubnumero{The construction of $nFAM$}

Set
$\Phi (E)$
equal to the set of cartesian families $F\rightarrow {\bf r}(E)$.
We ignore set-theoretic difficulties, just saying that to be correct one should
look at the set of cartesian families such that the underlying set is contained
in a given (large) fixed  set.

\begin{parag}
\label{contravariant}
Note that $\Phi$ is a contravariant functor from $(n+1)PC$ to $Sets$: if
$f:E'\rightarrow E$ is a morphism of $n$-precats then the pullback morphism
$\Phi (f)$ takes a cartesian family
$$
F\rightarrow {\bf r}(E)
$$
to the fiber-product family
$$
F\times _{{\bf r}(E)} {\bf r}(E').
$$
This is again a cartesian family. To prove that, note first of all that the
piecewise fibrant condition is preserved under pullbacks; after that we just
have to verify $(\ast )$ that the initial map $\iota : 0\rightarrow p$
induces an
equivalence  $$
F_p \times _{E_{p/}} E'_{p/} \cong E'_{p/}\times _{E_0} F_0.
$$
To do this, note that we can consider $E_{p/}$ as a colimit of indecomposable
$W$. Then $F_p$ is a colimit of the corresponding
$$
F_p(W):=W\times _{E_{p/}} F_p
$$
and similarly $E'_{p/}$ is a colimit of the
$$
E'_{p/}(W):=W\times
_{E_{p/}}E'_{p/}.
$$
The fiber product of these two colimits over $E_{p/}$ is the colimit of the
fiber products of the elements over $W$. Thus it suffices to prove that for any
indecomposable $W\hookrightarrow E_{p/}$,
$$
F_p(W)\times _W E'_{p/}(W)\cong
E'_{p/}(W)\times _{E_0}F_0
$$
is a weak equivalence (this suffices because the model category is left proper,
i.e. colimits preserve levelwise weak equivalences).
Now, as input we know that
$$
F_p(W) \cong W \times _{E_0}F_0,
$$
and also that $F_p(W)\rightarrow W$ is a fibration.

One has to be careful that the closed model category $nPC$ is not right proper,
i.e. fiber products don't necessarily preserve weak equivalences. However, in
our case the
morphism $W\times _{E_0} F_0\rightarrow W$ is a disjoint union of trivial
fibrations, so (by \cite{Simpson}(x) Theorem 5.1) it is ``compatible with
change of base'' in the sense of (\cite{HirschowitzSimpson} 11.12); thus
$F_p(W)\rightarrow W$ is compatible with change of base, and we get the desired
property that
$$
F_p(W)\times _W E'_{p/}(W)\cong
E'_{p/}(W)\times _{E_0}F_0.
$$
This completes the proof that
$F\times _{{\bf r}(E)} {\bf r}(E')$ is a cartesian family, so it completes the
proof that $\Phi$ is a contravariant functor.
\end{parag}

We would like to define $nFAM$ by setting
$$
Hom _{(n+1)PC}(E, nFAM):= \Phi (E).
$$
This suffices to define an $n+1$-precat $nFAM$, in view of the
following lemma.

\begin{lemma}
\label{nFAMexists}
The functor $E\mapsto \Phi (E)$ takes colimits to limits, so there is an
essentially unique $n+1$-precat $nFAM$ with functorial isomorphism
$Hom _{(n+1)PC}(E, nFAM)\cong \Phi (E)$.
\end{lemma}
{\em Proof:}
We write the proof in terms of a single coproduct but the same argument works
for any colimit. If
$$
E' \leftarrow E \rightarrow E''
$$
is a diagram of $n+1$-precats, and if
$$
F'\rightarrow {\bf r}(E'), \;\;\;\; F'' \rightarrow {\bf r}(E'')
$$
are cartesian families such that the pullbacks to $E$ are the same cartesian
family
$$
F\rightarrow {\bf r}(E),
$$
then we claim that
$$
F' \cup ^F F'' \rightarrow {\bf r}(E' \cup ^E E'')
$$
is again a cartesian family. First of all, it is piecewise fibrant because of
the glueing property \ref{glueing} for  piecewise fibrations (which works in the
same way in our present situation as for the original situation discussed in
\ref{glueing}). On the other hand, it is a precartesian family, by
essentially the
same argument as in \ref{contravariant} above. Thus, it is a cartesian family.
We obtain that
$$
\Phi (E' \cup ^E E'') = \Phi (E')\times _{\Phi (E)} \Phi (E'').
$$
By a Yoneda-type argument there exists an $n+1$-precat $nFAM$ such that a map
$E\rightarrow nFAM$ is the same thing as an element of $\Phi (E)$ i.e. as a
cartesian family over $E$.
\eop

The description of maps $E\rightarrow nFAM$ as corresponding to cartesian
families over $E$ is particularly natural and should allow a reasonable facility
of construction of such maps.  This will be seen in applications such as the
construction of the ``arrow family'' below.

Our next task is to prove that $nFAM$ is in fact a fibrant $n+1$-precat.

\begin{lemma}
\label{fibrant}
If $E'\rightarrow E$ is a trivial cofibration of $n+1$-precats and if $F' \in
\Phi (E')$ is a cartesian family over $E'$ then $F'$ extends to a cartesian
family $F\in \Phi (E)$ over $E$. In other words, the $n+1$-precat $nFAM$ is
fibrant.
\end{lemma}
{\em Proof:}
The equation $Hom (E, nFAM)=\Phi (E)$ provides the equivalence between the
statement that $nFAM$ is fibrant, and the extension condition of the first
sentence of the lemma. Thus, we prove the extension condition.

{\bf Step 1.}
We first consider the case when $E'\rightarrow E$ induces an
isomorphism on the set of objects, and induces trivial cofibrations
$$
E'_{p/} \hookrightarrow E_{p/}.
$$
In other words, ${\bf r}(E')\rightarrow {\bf r}(E)$ is a trivial cofibration of
$n$-prestacks over $\Delta$. We are given a cartesian family $F'\rightarrow {\bf
r}(E')$. Choose a factorization
$$
F' \hookrightarrow A \stackrel{a}{\rightarrow} {\bf r}(E)
$$
with $a$ a fibration (of $n$-prestacks over $\Delta $). We can't directly use
$A$ as our extension $F$ because it doesn't restrict to $F'$ over ${\bf r}(E')$.
Thus, we first note that
$$
i:F' \hookrightarrow A\times _{{\bf r}(E)}{\bf r}(E')
$$
is a trivial cofibration universally relative to
${\bf r}(E')$. To see this, note that the morphism $a$ is homotopically
equivalent, over each component ${\bf r}(E)_p=E_{p/}$, to a disjoint union of
trivial fibrations (pulled back from $F'_0\rightarrow E'_0=E_0$) and in
particular
it is compatible with change of base.  Again using the expressions of
equivalence with trivial cofibrations, we find that the morphism $i$ is a weak
equivalence.

Next, recall that by definition the morphism
$$
F'\rightarrow {\bf r}(E')
$$
is piecewise a fibration (at least over the indecomposable pieces which make up
${\bf r}(E')$, coming from indecomposable pieces of $E$). Therefore the lifting
property \ref{unilifting} (twisted slightly to our situation in the present
subsection, but again the proof is the same) implies that there is a retraction
$$
d:A\times _{{\bf r}(E)} {\bf r}(E')\rightarrow F'
$$
over ${\bf r}(E')$. Use this retraction to define
$$
F:= A \cup ^{A\times _{{\bf r}(E)} {\bf r}(E')} F'.
$$
Now we have a morphism of $n$-prestacks $F\rightarrow {\bf r}(E)$ which
restricts over ${\bf r}(E')$ to $F'$. One can easily check that it is cartesian
(using the glueing property \ref{glueing} for the piecewise fibrant condition).
This completes the proof of Step 1.

{\bf Step 2.}
We look at what a cartesian family over $E=\Upsilon (U)$ means (refer to
\cite{Simpson}(xi) for the notation $\Upsilon$). Here $U$ is an $n$-precat.
Note that
$$
E_{p/} = \ast \sqcup \coprod _{i= 1}^p U \sqcup \ast
$$
where the factor $U$ corresponding to a given index $i$ is equal to
$E_{p/}(x_0,\ldots , x_p)$ with $x_0=\ldots = x_{i-1} = 0$ and
$x_i=\ldots = x_p=1$.  (Recall that the objects of $\Upsilon (U)$ are denoted
$0$ and $1$). The first $\ast$ corresponds to the sequence $x_0= \ldots = x_p =
0$ and the second $\ast$ corresponds to $x_0=\ldots = x_1$.

If $F$ is a cartesian family over $E$, let $F(0)$ (resp. $F(1)$) denote the
pieces of $F_0$ lying over $0$ (resp. $1$) in $E_0$, and let  $F(0,1)$ be the
piece of $F_1$ lying over $U\subset E_{1/}$.
We have an equivalence $F(0,1)\stackrel{\cong}{\rightarrow} F(0)\times U$.
Define the ``strictification'' of
$F$ in the following way. First, define a morphism of precartesian families
$$
Str^1(F) \stackrel{a}{\rightarrow} F
$$
by putting
$$
Str^1(F)_p := F(0) \sqcup \coprod _{i= 1}^p (F(0,1)) \sqcup F(1).
$$
The simplicial structure comes from structural maps which are determined by the
condition of compatibility with the morphism $a$ which we now define. The maps
$$
Str^1(F)_p\stackrel{a_p}{\rightarrow} F_p
$$
come from the degeneracy maps
$$
F(0)\rightarrow F(0,\ldots , 0),\;\;\;\; F(1)\rightarrow F(1,\ldots , 1),
$$
and
$$
F(0,1)\rightarrow F(0,\ldots , 0,1, \ldots , 1).
$$
The precartesian condition for $F$ garanties that $a$ is an equivalence.

Next, use the equivalence
$$
F(0,1)\stackrel{\cong}{\rightarrow} F(0)\times U
$$
and the condition that $F(1)$ is fibrant, to choose a map
$$
F(0)\times U \rightarrow F(1)
$$
such that the composition
$$
i': F(0,1)\rightarrow F(0)\times U \rightarrow F(1)
$$
is homotopic to the original structural map denoted
$$
i: F(0,1)\rightarrow F(1).
$$
The homotopy may be
realized as a map
$$
F(0,1) \times \overline{I}\rightarrow F(1).
$$
Using this map (and the structural map $F(0,1)\rightarrow F(0)$ composed with
the first projection), define a precartesian family $Str^2(F)$ by the formula
$$
Str^2(F)_p := F(0) \sqcup \coprod _{i= 1}^p (F(0,1)\times \overline{I}) \sqcup
F(1).
$$
The inclusion at the start of the homotopy (i.e. the end which corresponds to
the original structural map $i$) gives an equivalence of precartesian families
$$
b: Str^1(F) \rightarrow Str ^2(F).
$$
On the other hand, define the precartesian family $Str^3(F)$ by the same formula
as for $Str^1(F)$, but using the morphism $i'$ in place of $i$ to define the
simplicial structure. Inclusion at the other end of our homotopy gives an
equivalence of precartesian families
$$
c: Str^3(F) \rightarrow Str ^2(F).
$$
Finally, define a precartesian family $Str^4(F)$ by the formula
$$
Str^4(F)_p:= F(0) \sqcup \coprod _{i= 1}^p (F(0)\times U) \sqcup F(1).
$$
The factorization of the structural map $i'$ into
$$
F(0,1)\rightarrow F(0)\times U \rightarrow F(1)
$$
yields an equivalence of precartesian families
$$
d:Str^3(F) \rightarrow Str^4(F).
$$
Set $Str(F):= Str^4(F)$; note that it is a strict precartesian family
coming from a morphism $U\times F(0)\rightarrow F(1)$; we call it the
{\em strictification of $F$}. Together, $a$, $b$, $c$ and $d$ provide a chain of
equivalences of precartesian families linking $F$ to the
strictification $Str(F)$.

{\bf Step 3:}
We now consider the trivial cofibrations of the following type
(again see \cite{Simpson}(xi) for the notations):
$$
E'= \Upsilon (U_1)\cup ^{\{ 1\} } \cdots \cup ^{\{ k-1\}} \Upsilon (U_k)
\hookrightarrow
\Upsilon ^k(U_1,\ldots , U_k) = E.
$$
If $F'$ is a cartesian family over $E'$ then we can view $F'$ as a coproduct of
cartesian families $F^{(i)}$ over the $\Upsilon (U_i)$. Let $Str(F')$ denote the
corresponding coproduct of the $Str(F^{(i)})$. Applying Step 2, we get a
chain of
equivalences (levelwise over $\Delta$ and relatively with respect to ${\bf
r}(E')$) $$
F' \stackrel{\cong}{\leftarrow} \ldots \stackrel{\cong}{\rightarrow} Str(F').
$$
Using Step 1 and standard model-category-theoretic arguments (including, where
appropriate, replacement of precartesian families by their equivalent cartesian
replacements using a fibrant replacement of $n$-prestacks \ldots ), it
suffices in
order to extend $F'$ to $F$, to be able to extend $Str(F')$. But now
$Str(F')$ is
determined by the data $F(0), \ldots , F(k)$ and $$ U_i \times F(i-1)\rightarrow
F(i) $$
for $i=1,\ldots , k$. Using these data, we can directly construct a strict
family $Str(F)$ over $E$ extending $Str(F')$. This shows that any cartesian
family over $E'$ extends over $E$.

{\bf Step 4:}
Consider the cofibration $E'=\ast \rightarrow E = \overline{I}$.
Since $E$ retracts back onto $E$, it is immediate that any cartesian family over
$E'$ extends to a cartesian family over $E$.

{\bf Step 5:} Now, combining Steps 1, 3 and 4 and using a standard
model-category-theoretic argument (also using Lemma \ref{nFAMexists})
we obtain the required extension along any trivial cofibration
$E'\hookrightarrow E$. Indeed, any such can be embedded into a larger one
$$
E'\hookrightarrow E \rightarrow E''
$$
where $E'\hookrightarrow E''$ is obtained as a sequential colimit of things of
the form treated in Steps 1, 3, and 4. Thus, the family $F'$ extends to $E''$
and the restriction to $E$ is the extension we are looking for.
This completes the proof of the lemma.
\eop

\oldsubnumero{The map $nCAT\rightarrow nFAM$}

We now define the map ${\bf u}: nCAT\rightarrow nFAM$ which will be a weak
equivalence. As $nFAM$ is fibrant by Lemma \ref{fibrant}, this will constitute a
concrete, canonical fibrant replacement for $nCAT$. Unfortunately, the map
${\bf u}$ itself is not canonical but depends on a choice of fibrant
replacement
(this is because of the extra fibrant condition in the definition of
``cartesian'').

To define the map, it suffices to construct a ``universal'' cartesian family
over $nCAT$. Recall how $nCAT$ was defined: the objects are the fibrant
$n$-categories $A$, and for a sequence of objects
$A_0,\ldots , A_p$ we put
$$
nCAT_{p/}(A_0,\ldots , A_p) := \underline{Hom}(A_0,A_1)\times
\ldots \times \underline{Hom}(A_{p-1}, A_p).
$$
Note that we can write
$$
nCAT_{p/}= \coprod _{A_0,\ldots , A_p} nCAT_{p/}(A_0,\ldots , A_p).
$$
The morphisms of the simplicial structure are obtained using the (strictly
associative) composition of internal $\underline{Hom}$'s. Define a cartesian
family
$$
\Uu \rightarrow {\bf r}(nCAT)
$$
by setting
$$
\Uu _p:= \coprod _{A_0,\ldots , A_p}A_0 \times nCAT_{p/}(A_0,\ldots , A_p)
$$
$$
= \coprod _{A_0,\ldots , A_p}A_0 \times \underline{Hom}(A_0,A_1)\times
\ldots \times \underline{Hom}(A_{p-1}, A_p).
$$
Again, the morphisms of the simplicial structure are obtained using the
compositions and, eventually, the evaluation morphisms
$$
A_0\times \underline{Hom}(A_0,A_i)\rightarrow A_i .
$$
(these evaluation maps share strict associativity with respect to the
compositions
so this works).
The map to ${\bf r}(nCAT)_p= nCAT_{p/}$ is the obvious projection.

The family $\Uu \rightarrow {\bf r}(nCAT)$ is obviously precartesian.

The morphisms $\Uu _p \rightarrow nCAT_{p/}$ are fibrant, but we don't know if
the full morphism $\Uu \rightarrow {\bf r}(nCAT)$ is piecewise fibrant; thus,
let
$$
\Uu \hookrightarrow \Uu '\rightarrow {\bf r}(nCAT)
$$
be a fibrant replacement in the world of $n$-prestacks over $\Delta$.
This is again precartesian, and it is fibrant thus piecewise fibrant, so $\Uu '$
is a cartesian family over $nCAT$. This corresponds to a morphism
$$
{\bf u}:nCAT\rightarrow nFAM.
$$

\begin{exercise}
We leave as an exercise to the reader to state and prove a
unicity-up-to-coherent-homotopy for ${\bf u}$ (we have made a choice of $\Uu '$
but this choice is unique up to coherent homotopy).
\end{exercise}

\oldsubnumero{Functors of $n$-precats}

Before getting to the next theorem which states that ${\bf u}$ is an
equivalence, we need to develop one more technique. The reason is that we have
made an arbitrary choice in defining ${\bf u}$, due to the fact that the
standard universal family over $nCAT$ is precartesian but not cartesian.
Thus it would be more convenient to work with precartesian families; but these
don't have the glueing property \ref{glueing}, so the functor in question
doesn't
come from an $n+1$-precat. Thus, instead, we shall study this type of functor.
As one can see by looking in the proof of Theorem \ref{famequiv} below, we
actually  come across this problem for $n$-precats rather than $n+1$-precats.

Let $Funct(nPC)$ denote the category of functors
$$
F: nPC^o \rightarrow Sets .
$$
The map taking an object to its associated representable functor gives a full
embedding $nPC \subset Funct(nPC)$. Thus, if $U\in nPC$ and
$F\in Funct(nPC)$ we can speak of the set of morphisms from $U$ to $F$.
Using the same definitions as in Quillen \cite{Quillen}, we can define the
notion of {\em homotopy} between two morphisms $U\tworightarrows F$.
Furthermore, we can define a notion of {\em derived morphism from $U$ to $F$}
as being a diagram
$$
U\leftarrow U' \rightarrow F
$$
where $U'\rightarrow U$ is a weak equivalence in $nPC$. Finally, we can define
the set of {\em derived morphisms from $U$ to $F$ up to homotopy} (it is left
to the reader to write down the necessary diagrams). If $F\in nPC$ too then the
set of derived morphisms  up to homotopy is just the set of morphisms from $U$
to $F$ in the homotopy category $Ho(nPC)$.

Suppose $f:A\rightarrow B$ is a morphism in $Funct(nPC)$. We say that
$f$ has the {\em weak homotopy lifting property} if the following holds: for any
square
$$
\begin{array}{ccc}
U& \rightarrow & A \\
\downarrow && \downarrow \\
U' & \rightarrow & B
\end{array}
$$
where the vertical morphisms are a morphism $U\rightarrow U'$ in $nPC$, and $f$
respectively, where the horizontal morphisms are derived morphisms, and where
the diagram commutes up to homotopy; then there exists a derived morphism
$U'\rightarrow A$ such that the two triangles commute up to homotopy. Note that
we don't require compatibility between the homotopies involved.

In the case of a morphism $f: A\rightarrow B$ in $nPC \subset Funct(nPC)$,
the weak homotopy lifting property is equivalent to saying that $f$ projects to
an isomorphism in the homotopy category $Ho(nPC)$. By \cite{Quillen}
this is equivalent to $f$ being a weak equivalence. Recall that weak
equivalences satisfy the ``three for the price of two'' property. The following
lemma states that the same is true of morphisms in $Funct(nPC)$ which have the
weak homotopy lifting property.

\begin{lemma}
\label{threefortwo}
Suppose $A\stackrel{f}{\rightarrow} B \stackrel{g}{\rightarrow}C$
is a sequence of morphisms in $Funct(nPC)$. Then if any two of $f$, $g$  and
$gf$ satisfy the weak homotopy lifting property, the third one does too.
\end{lemma}
{\em Proof:} Left to the reader, using the techniques of \cite{Quillen}.
\eop

We can now state the theorem which completes the construction of a fibrant
replacement for $nCAT$. The reader will note that the proof we give here is
somewhat incomplete at the end.

\begin{theorem}
\label{famequiv}
The morphism ${\bf u}:nCAT\rightarrow nFAM$ is a weak equivalence.
\end{theorem}
{\em Proof:}
It is clearly surjective on objects (the objects of $nFAM$ are themselves
exactly the fibrant $n$-categories). Since we already know that both sides are
$n+1$-categories, it suffices to prove that for two fibrant $n$-categories
$A,B$, the morphism
$$
\underline{Hom}(A,B) = nCAT_{1/}(A,B) \rightarrow nFAM _{1/}(A,B)
$$
is a weak equivalence of $n$-categories.

Recall that if $U$ is an $n$-precat and $C$ an $n+1$-category with objects
$x,y\in C_0$, a morphism  $E\rightarrow C_{1/}(x,y)$ is the same thing as a
morphism $\Upsilon (E)\rightarrow C$ sending $0$ to $x$ and $1$ to $y$.
Applying this to $C=nFAM$, we get that a morphism
$f:U\rightarrow nFAM_{1/}(A,B)$ is the same thing as a cartesian family $F$ over
$\Upsilon (U)$, with $F(0)=A$ and $F(1)=B$.

Define the following functors in $Funct(nPC)$:
$$
{\bf cartesian}(\Upsilon (\cdot ); A,B); \;\;\;\; {\bf precartesian}(\Upsilon
(\cdot ); A,B).
$$
These are defined by setting, for $E\in nPC$,
${\bf cartesian}(\Upsilon (\cdot ); A,B)(E)$ equal to the set of cartesian
families over $\Upsilon (E)$ restricting to $A$, $B$ on the endpoints;
respectively ${\bf precartesian}(\Upsilon (\cdot ); A, B)(E)$ equal to the
set of
precartesian families over $\Upsilon (E)$ restricting to $A$, $B$ on the
endpoints. We have a morphism of functors
$$
{\bf cartesian}(\Upsilon (\cdot ); A,B)\rightarrow {\bf precartesian}(\Upsilon
(\cdot ); A,B).
$$
On the other hand, by definition ${\bf cartesian}(\Upsilon (\cdot ); A,B)$ is
the functor represented by $nFAM_{1/}(A,B)$.

Set $H:= \underline{Hom}(A,B)= nCAT_{1/}(A,B)$. The map ${\bf u}$ induces a map
$$
H\rightarrow {\bf cartesian}(\Upsilon (\cdot ); A,B),
$$
whereas the standard
precartesian family over $nCAT$ induces a map
$$
H\rightarrow {\bf precartesian}(\Upsilon (\cdot ); A,B).
$$
We claim three things:
\newline
(1) \, the triangle composed of the above three maps commutes up to homotopy;
\newline
(2) \, the morphism
$$
{\bf cartesian}(\Upsilon (\cdot ); A,B)\rightarrow {\bf precartesian}(\Upsilon
(\cdot ); A,B)
$$
satisfies the weak homotopy lifting property; and
\newline
(3) \, the morphism
$$
H\rightarrow {\bf precartesian}(\Upsilon (\cdot ); A,B)
$$
satisfies the weak homotopy lifting property.

With these three, Lemma \ref{threefortwo} implies that the first morphism
$$
H\rightarrow {\bf cartesian}(\Upsilon (\cdot ); A,B) = nFAM _{1/}(A,B)
$$
satisfies weak homotopy lifting, but this is a morphism in $nPC$ so it is a weak
equivalence. This is what we needed to show for the theorem.

We will not give the detailed proofs of (1)-(3) here. We just note that
(1) and (2) are general types of statements related to the use of the fibrant
replacement in the construction of ${\bf u}$, and that for (3), one can use the
``strictification'' procedure which was explained in Step 2 of the proof of
Lemma
\ref{fibrant}. \eop

\subnumero{Cartesian families of $n$-stacks}

The above discussion also works in the situation of $n$-stacks over a site
$\Gg$.
We briefly indicate how this goes; the proofs are left to the reader.

Suppose $E$ is an $n+1$-prestack on $\Gg$ (i.e. it is a presheaf of
$n+1$-precats). Then applying the construction ${\bf r}$ object-by-object, we
obtain an $n$-prestack ${\bf r}(E)$ over $\Delta \times \Gg$. A {\em
precartesian (resp. cartesian) family of $n$-prestacks over $E$} is a morphism
$$
F\rightarrow {\bf r}(E)
$$
of $n$-prestacks over $\Delta \times \Gg$, such that for each object
$X\in \Gg$, the value
$$
F(X) \rightarrow {\bf r}(E) (X)
$$
is a precartesian (resp. cartesian) family of $n$-precats in the sense of the
previous subsection.

{\em Caution:} If $E$ is an $n+1$-stack, the reduction ${\bf r}(E)$ will not in
general be an $n$-stack on $\Delta \times \Gg$; for example the presheaf of
objects $E_0$ (which becomes ${\bf r}(E)_0$) is not in general a sheaf.

Suppose $F\rightarrow {\bf r}(E)$ is a precartesian (resp. cartesian) family of
$n$-prestacks over $E$. We say that it is a {\em precartesian (resp. cartesian)
family of $n$-stacks over $E$} if, for every $X\in \Gg$ and every $e\in E_0(X)$,
the restriction of $F_0$ over the point $$
e: X\rightarrow {\bf r}(E)_0 = E_0
$$
is an $n$-stack on $\Gg /X$.

If $E$ is an $n+1$-prestack over $\Gg$, let $\underline{\Phi}(E)$ be the set of
cartesian families of $n$-stacks over $E$.

\begin{lemma}
\label{cartesianstacks1}
The functor $\Phi$ is represented by an $n+1$-prestack $n\underline{FAM}$
over $\Gg$.
\end{lemma}
{\em Proof:}
Left to the reader, following the proof of \ref{nFAMexists}.
\eop

There is a canonical precartesian family of $n$-stacks over
$n\underline{STACK}$ constructed as in the previous subsection for $nCAT$.
Choosing a fibrant replacement (for the coarse topology on $\Delta \times
\Gg$) gives a cartesian family of $n$-stacks over  $n\underline{STACK}$, hence a
morphism of $n+1$-prestacks $$
\underline{{\bf u}}:  n\underline{STACK}\rightarrow n\underline{FAM}.
$$

\begin{theorem}
\label{cartesianstacks2}
The $n+1$-prestack $n\underline{FAM}$ is fibrant for the Grothendieck topology of
$\Gg$ (in particular it is an $n+1$-stack).
The morphism $\underline{{\bf u}}$ is an equivalence of $n+1$-stacks.
Thus it constitutes a fibrant replacement for $n\underline{STACK}$.
\end{theorem}
{\em  Proof:}
First, prove the theorem
for the case where $\Gg$ has the coarse topology. In this case, ``weak
equivalence'' means equivalence object-by-object over $\Gg$. Thus the extension
problem refered to in the proof of Theorem \ref{cartesianstacks1} is
essentially the same as that of Lemma \ref{fibrant}, and can safely be left
to the
reader.

Now treat the case where $\Gg$ has a different Grothendieck topology. The
essential problem is to show that if $E\rightarrow E'$ is a trivial cofibration
of $n+1$-prestacks on $\Gg$, then any cartesian family over $E$ extends to one
over $E'$. This is somewhat delicate since, for example, the corresponding
morphism of object presheaves $E_0\rightarrow E'_0$ need not be a weak
equivalence.

Instead of treating this problem directly, we use the glueing result of
\cite{HirschowitzSimpson} (and in fact, the extension problem refered to above
is essentially a problem of glueing $n$-stacks).

Use the
superscript $\Gg$ to denote the topology of $\Gg$, and the superscript
``coarse'' to denote the coarse topology on $\Gg$. By \cite{HirschowitzSimpson},
$$
n\underline{STACK}^{\Gg} \subset n\underline{STACK}^{\rm coarse}
$$
is the full sub-prestack consisting of the coarse-topology stacks which are
$\Gg$-stacks. On the other hand,
$$
n\underline{FAM}^{\Gg} \subset n\underline{FAM}^{\rm coarse}
$$
is by definition the full sub-prestack consisting of the coarse-topology
cartesian families which are cartesian families of $\Gg$-stacks. The morphism of
Theorem \ref{cartesianstacks2} for the coarse topology is an equivalence
$$
\underline{{\bf u}}^{\rm coarse}
n\underline{STACK}^{\rm coarse}\stackrel{\cong}{\rightarrow}
n\underline{FAM}^{\rm coarse}.
$$
This morphism sends $n\underline{STACK}^{\Gg}$ to
$n\underline{FAM}^{\Gg}$, and since both of these are full sub-prestacks,
$\underline{{\bf u}}^{\rm coarse}$ induces an equivalence of $n+1$-prestacks
$$
n\underline{STACK}^{\Gg}
\stackrel{\cong}{\rightarrow}
n\underline{FAM}^{\Gg}.
$$
Theorem 20.1 of \cite{HirschowitzSimpson} says that the $n+1$-prestack
$n\underline{STACK}^{\Gg}$ is an $n+1$-stack. This implies that
$n\underline{FAM}^{\Gg}$ is an $n+1$-stack. On the other hand, the result of
the present theorem for the case of the coarse topology says that
$n\underline{FAM}^{\Gg}$ is fibrant in the coarse topology. Finally, Lemma 9.2
of \cite{HirschowitzSimpson} says that an $n+1$-prestack which is fibrant
for the
coarse topology and which is an $n+1$-stack, is fibrant for the topology $\Gg$.
Thus $n\underline{FAM}^{\Gg}$ is $\Gg$-fibrant.

For the second statement of the theorem, one has to verify that the restriction
of $\underline{{\bf u}}^{\rm coarse}$ to the full sub-prestack
$n\underline{STACK}^{\Gg}$, is equal to
$\underline{{\bf u}}^{\rm \Gg}$. This is of course true on the level of the
canonical precartesian families; but note that we used a fibrant replacement for
the coarse topology to define $\underline{{\bf u}}^{\rm \Gg}$, so it is
the same as $\underline{{\bf u}}^{\rm coarse}$. With what was said previously,
we now obtain the second statement of the theorem.
\eop

In view of this theorem, we may define as our fibrant replacement
$$
n\underline{STACK}'= n\underline{FAM}.
$$
Then, in view of our convention that a ``morphism to $n\underline{STACK}$''
really means a morphism to the fibrant replacement $n\underline{STACK}'$,
we can say that a cartesian family of $n$-stacks over an $n+1$-stack $E$,
corresponds to a morphism $E\rightarrow n\underline{STACK}$ (and vice-versa).

Similarly, we will make a confusion between realms $\underline{R}\subset
n\underline{STACK}$ and their corresponding full sub-stacks
$\underline{R}'\subset n\underline{FAM}$ (noting that
$\underline{R}'$ will be a fibrant replacement for $\underline{R}$).
Thus, for example, when we speak of
$$
\underline{Hom}(\underline{R},
\underline{A})
$$
this will really mean
$$
\underline{Hom}(\underline{R}',
\underline{A}')
$$
or even $\underline{Hom}(\underline{R}',
\underline{A}') $.

In the remainder of the paper for clarity of notation we will keep the original
$n\underline{STACK}$.

\subnumero{Construction of cartesian families}

A crucial question which arises from the above point of view is: under what
circumstances can a morphism of $n$-stacks $X\rightarrow B$ be considered as
a cartesian family, interpreting $B$ as an $n+1$-stack?

Our basic answer is that this works if $B$ is an $n$-stack of groupoids. In
general, some condition such as ``ccb'' will almost certainly be necessary but
we don't attack the general situation.

Before getting to the statement of that construction, we do a ``warm-up
exercise'' with the following lemma. For both this lemma and the subsequent
proposition, we use Tamsamani's construction of the Poincar\'e $n$-groupoid of a
space $\Pi _n(X)$, and his theorem that any $n$-groupoid is equivalent to one
of the form  $\Pi _n(X)$ \cite{Tamsamani}. It would be nice to have a more
algebraic alternative to this type of argument, but I haven't yet found one.

\begin{lemma}
\label{opposite}
Suppose $B$ is an $n$-groupoid (resp. $n$-stack of groupoids). Then
``inverting $1$-morphisms''  is an equivalence (which exists as a morphism if we
assume that $B$ is fibrant)
$$
B^o\stackrel{\cong}{\rightarrow} B.
$$
\end{lemma}
{\em Proof:}
Treat first the punctual case where $B$ is an $n$-groupoid, and assume that it
is fibrant as an $n$-precat. Tamsamani constructs a space
$$
Y=\Re (B)
$$
and an equivalence of $n$-groupoids
$$
B\stackrel{\cong}{\rightarrow} \Pi _n(Y)
$$
where $\Pi _n$ is the ``Poincar\'e $n$-groupoid'' construction predicted in
\cite{Grothendieck} and defined in \cite{Tamsamani}.

Now, reversing the time direction for the $1$-morphisms (but not the
$i$-morphisms, $i\geq 2$) gives an isomorphism
$$
\Pi _n(Y) ^o \cong \Pi _n(Y).
$$
Recall from \cite{Tamsamani} that
the objects of $\Pi _n(Y)$ are the points of $Y$, and
$$
\Pi _n(Y)_{p/} (y_0, \ldots , y_p)
= \Pi _{n-1}(Hom _{\rm top}^{y_0,\ldots , y_p}(R^{p}, Y))
$$
where
$$
R^p := \{ (t_0,\ldots , t_p) \in \rr ^{p+1},\;\;\; t_i \geq 0, \;\;\; \sum t_i
= 1\}
$$
is the standard $p$-simplex and
$$
Hom _{\rm top}^{y_0,\ldots , y_p}(R^{p}, Y)
$$
denotes the space of maps from $R^p$ to $Y$ which send the vertices to $y_i$.
There is a reflection of $R^p$ which inverts the order of the vertices, giving
an isomorphism
$$
Hom _{\rm top}^{y_0,\ldots , y_p}(R^{p}, Y)
\cong
Hom _{\rm top}^{y_p,\ldots , y_0}(R^{p}, Y).
$$
This gives an isomorphism
$$
\Pi _n(Y)_{p/} (y_0, \ldots , y_p) \cong
\Pi _n(Y)_{p/} (y_p, \ldots , y_0)=
\Pi _n(Y)^o_{p/} (y_0, \ldots , y_p).
$$
The set
of these reflections is compatible with the face and degeneracy maps in the
appropriate way, so we obtain an isomorphism of $n$-categories
$$
\Pi _n(Y) \cong \Pi _n(Y)^o.
$$
Now if $B$ was fibrant, we could invert the original equivalence
and compose to obtain an equivalence
$$
B^o\rightarrow  \Pi _n(Y)^o\cong \Pi _n(Y)
\rightarrow B
$$
as desired. The map $B^o\rightarrow B$ depends on a number of choices but one
can verify that, up to homotopy, it is independent of the choices.

The above construction is
functorial in $B$, except for the choice of inverse at the last step.
Thus we obtain, for an $n$-stack of groupoids (i.e. presheaf of $n$-groupoids
over $\Gg$) a diagram of the form
$$
B^o\rightarrow A^o \cong A \leftarrow B
$$
where $A$ is the $n$-prestack $X\mapsto \Pi _n (\Re (B(X)))$, and where the
arrows are object-by-object equivalences. If $B$ is fibrant over $\Gg$ then the
last arrow can be inverted and we obtain the desired equivalence $B^o \cong B$.
\eop

\begin{proposition}
\label{correspondence}
Suppose that $B$ is an $n$-stack of groupoids, and let $\underline{B}$ denote
$B$ considered as an $n+1$-stack of groupoids. Then any morphism of $n$-stacks
$\Ff \rightarrow B$ gives rise to a natural cartesian family
$$
\Ff ^{{\rm fam}/B}\rightarrow {\bf r}(\underline{B}).
$$
Conversely, any such cartesian family comes from a morphism of $n$-stacks
$\Ff \rightarrow B$.
\end{proposition}
{\em Construction-Proof:}
The main problem is to get a hold of ${\bf r}(\underline{B})$. We do this in the
following way. Assume that $B$ is fibrant. Let $\overline{I}^{(k)}$ denote the
$1$-groupoid with $k+1$ isomorphic objects denoted $0,\ldots , k$ (and no
nontrivial automorphisms). These are organized into a cosimplicial object in the
category of $1$-groupoids. We can also think of them as constant prestacks of
$1$-groupoids over $\Gg$. Now put
$$
\rho '(B)_k:= \underline{Hom}(\overline{I}^{(k)}, B),
$$
where the internal $\underline{Hom}$ is that of $n$-prestacks on $\Gg$. Thus
$\rho '(B)$ is a simplicial object in the category of $n$-prestacks on $\Gg$.
Change this slightly, noting that $\rho '(B)_0 =B$. Put
$$
\rho (B)_k:= \rho '(B) _k \times _{B^{k+1}} B_0^{k+1},
$$
where the first structural map for the fiber product is restriction
of a map $\overline{I}^{(k)}\rightarrow B$, to the ``vertices'' i.e. the
objects of $\overline{I}^{(k)}$. The second structural map comes from
$B_0\rightarrow B$.

If $x_0,\ldots , x_k$ are objects of $B_0(Y)$, set (for any $Z/Y$)
$$
\rho (B)_k(x_0,\ldots , x_k)(Z)
$$
equal to the preimage of $(x_0,\ldots ,x_k)|_Z\in B_O^{k+1}(Z)$, in $\rho
(B)_k(Z)$.
Thus $\rho (B)_k(x_0,\ldots , x_k)$ is an $n$-prestack on $\Gg /Y$.

{\bf Claim:} We claim that
$\rho (B)_k(x_0,\ldots , x_k)$ is $n-1$-truncated so it can be considered as an
$n-1$-prestack on $\Gg /Y$, and that as such there is a natural equivalence
(or chain of equivalences)
$$
\rho (B)_k(x_0,\ldots , x_k) \cong B_{k/}(x_0,\ldots , x_k);
$$
and that as the sequence of objects varies, this gives
an equivalence of simplicial $n$-prestacks
$$
\rho (B) \cong {\bf r}(\underline{B}).
$$

To prove the claim, it suffices to treat the punctual case, i.e. we can work
over a site with only one object, thus considering everything as $n$-categories
rather than $n$-stacks. As in the proof of the previous lemma, the only way I
currently see to prove this claim is to  pass through the topological situation
via Tamsamani's $\Pi _n$ construction. Again, it would be nice to have a more
``algebraic'' proof of this claim.

Setting $X:= \Re (B)$, we have seen above that there is a natural equivalence
$$
B \cong \Pi _n(X).
$$
Note that $B_0$ is a subset of the set of objects of $\Pi _n(X)$, namely
consisting of only the points of $X$ corresponding to objects of $B$. Note also
that $\Pi _n(X)$ is fibrant, so it may be used as a target in an internal
$\underline{Hom}$.

By the
adjunction between $\Re$ and $\Pi _n$, and the fact that the realization $\Re
(\overline{I}^{(k)})$ is the standard $k$-simplex which we denote $R^k$, we have
an isomorphism  $$
\underline{Hom}(\overline{I}^{(k)}, \Pi _n(X))
\cong \Pi _n\underline{Hom}_{\rm top}(R^k, X).
$$
These are equivalent to $\Pi _n(X)$. Thus we can write
$$
\rho (B) _k (x_0,\ldots , x_k) = \Pi _n\left( \underline{Hom}_{\rm top}(R^k, X)
\times _{X^{k+1}} \{ (x_0,\ldots , x_k\} \right) .
$$
For topological reasons it is easy to see that this is $n-1$-truncated, in fact
it is equivalent to
$$
Path ^{x_0,x_1} (X) \times \ldots \times Path ^{x_{k-1},x_k}(X).
$$
This proves that $\rho (B) _k$ is $n-1$-truncated. Finally we note that
by construction,
$$
\Pi _n(X)_{k/}(x_0,\ldots , x_k)
=
\Pi _{n-1} \left( \underline{Hom}_{\rm top}(R^k, X)
\times _{X^{k+1}} \{ (x_0,\ldots , x_k\} \right) .
$$
The fact that the space in question is $n-1$-truncated means that the
truncation
morphism
$$
\Pi _n\left( \underline{Hom}_{\rm top}(R^k, X)
\times _{X^{k+1}} \{ (x_0,\ldots , x_k\} \right)
\rightarrow
\Pi _{n-1}\left( \underline{Hom}_{\rm top}(R^k, X)
\times _{X^{k+1}} \{ (x_0,\ldots , x_k\} \right)
$$
is an equivalence. Therefore we obtain a functorial morphism which is an
equivalence,
$$
\rho (B) _k (x_0,\ldots , x_k)\rightarrow
\Pi _{n}(X)_{k/}(x_0,\ldots , x_k) .
$$
On the other hand we have a morphism which is again an equivalence,
$$
B_{k/} (x_0,\ldots , x_k)
\rightarrow
\Pi _{n}(X)_{k/}(x_0,\ldots , x_k).
$$
This gives a chain of equivalences to prove the claim. Note that as $k$ and
$x_0,\ldots , x_k$ vary, the $\Pi _{n}(X)_{k/}(x_0,\ldots , x_k)$ fit together
to form a simplicial $n-1$-category
intermediate in the chain of equivalences. As everything is functorial, it also
works in the situation of stacks over $\Gg$.

Now continue with the proof of the proposition, 	assuming that $\Ff \rightarrow
B$ is a fibrant morphism of $n$-stacks.  In view of the claim it suffices to
construct a precartesian family
$$
\rho (\Ff  /B) \rightarrow \rho (B).
$$
To do this, set
$$
\rho '(\Ff  /B)_k:= \underline{Hom}(\overline{I}^{(k)}, \Ff ),
$$
thus $\rho '(\Ff  /B)_k \cong \Ff$.
Then set
$$
\rho (\Ff  /B)_k:= \rho '(\Ff  /B)_k \times _{\rho '(B)_k}\rho (B)_k.
$$
This is a simplicial $n$-stack, mapping to $\rho (B)$. We have to show that the
map
$$
\rho (\Ff  /B)\rightarrow \rho (B)
$$
is a cartesian family.

To see this, fix an initial morphism $0\rightarrow k$ in $\Delta$, and look
at the
cube whose top face is
$$
\begin{array}{ccc}
\rho (\Ff /B)_k & \rightarrow & \rho '(\Ff /B)_k\\
\downarrow && \downarrow \\
\rho (\Ff /B)_0 & \rightarrow & \rho '(\Ff /B)_0,
\end{array}
$$
whose bottom face is
$$
\begin{array}{ccc}
\rho (B)_k & \rightarrow & \rho '(B)_k\\
\downarrow && \downarrow \\
\rho (B)_0 & \rightarrow & \rho '(B)_0,
\end{array}
$$
and whose side faces contain (A1 $+$ A2):
$$
\begin{array}{ccccc}
\rho (\Ff /B)_k & \rightarrow & \rho '(\Ff /B)_k& \rightarrow & \rho '(\Ff
/B)_0\\
\downarrow && \downarrow && \downarrow \\
\rho (B)_k & \rightarrow & \rho '(B)_k
& \rightarrow & \rho '(B)_0,
\end{array}
$$
and (B1 $+$ B2):
$$
\begin{array}{ccccc}
\rho (\Ff /B)_k & \rightarrow & \rho (\Ff /B)_0& \rightarrow & \rho '(\Ff
/B)_0\\
\downarrow && \downarrow && \downarrow \\
\rho (B)_k & \rightarrow & \rho (B)_0
& \rightarrow & \rho '(B)_0.
\end{array}
$$
(The reader is urged to draw the cube.)

Now the sides (B2)
$$
\begin{array}{ccc}
\rho (\Ff /B)_0& \rightarrow & \rho '(\Ff
/B)_0\\
 \downarrow && \downarrow \\
 \rho (B)_0
& \rightarrow & \rho '(B)_0
\end{array}
$$
and (A1)
$$
\begin{array}{ccc}
\rho (\Ff /B)_k & \rightarrow & \rho '(\Ff /B)_k\\
\downarrow && \downarrow \\
\rho (B)_k & \rightarrow & \rho '(B)_k
\end{array}
$$
are cartesian, by construction. On the other hand, the side (A2)
$$
\begin{array}{ccccc}
\rho '(\Ff /B)_k& \rightarrow & \rho '(\Ff
/B)_0\\
\downarrow && \downarrow \\
\rho '(B)_k
& \rightarrow & \rho '(B)_0
\end{array}
$$
is cartesian because the top and bottom arrows are equivalences (the top
elements are equivalent to $\Ff$, the bottom elements are equivalent to $B$)
and since $B$ is an $n$-stack of groupoids, the morphism $\Ff \rightarrow B$ is
automatically ccb \cite{HirschowitzSimpson}, so the pullback of an equivalence
is again an equivalence.

It follows that the composition of A1 and A2 is a cartesian square; but this is
also the same as the composition of B1 and B2; and since we know that B2 is
cartesian, it follows that the square (B1)
$$
\begin{array}{ccc}
\rho (\Ff /B)_k & \rightarrow & \rho (\Ff /B)_0 \\
\downarrow && \downarrow \\
\rho (B)_k & \rightarrow & \rho (B)_0
\end{array}
$$
is homotopy-cartesian. This property is the same as the property of $\rho (\Ff
/B) \rightarrow \rho (B)$ being a precartesian family.

If one wants to get a cartesian family, make a fibrant replacement (in the world
of $n$-stacks on $\Delta \times \Gg$). The fact that $n\underline{FAM}$ is
fibrant means that we obtain an equivalent cartesian family over the base ${\bf
r}(\underline{B})$ since this latter is equivalent to $\rho (B)$. This completes
the construction.

For the converse statement in the proposition, we note that
$$
B \cong {\rm hocolim}_{\Delta} \rho (B).
$$
If $\Ee \rightarrow \rho (B)$ is a cartesian family, set
$$
\Ff := {\rm hocolim}_{\Delta} \Ee ;
$$
with its map to ${\rm hocolim}_{\Delta} \rho (B)\cong B$.
We leave to the reader to verify that this is (up to equivalence) an inverse of
the preceding construction ({\em Caution:} the author has not made this
verification).
\eop

\subnumero{The arrow family}

We can apply the above notion of cartesian family to construct the ``arrow
family'' for any $n+1$-category $A$. This answers a question posed in
\cite{Simpson}(xi) and allows us to define the notion of ``representable
functor'' (asked for by A. Hirschowitz).

The ``arrow family'' has always played an important role in enriched category
theory and I am not competent to give a complete list of references. One place
to start is Cordier and Porter \cite{CordierPorter}. The recent paper of
Carboni, Kelly, Verity, Wood \cite{CKVW} also seems to be relevant.

Fix an $n+1$-category $A$. We will define a precartesian family
$$
ArrFam(A) \rightarrow {\bf r}(A^o \times A).
$$
A sequence of objects of $A^o\times A$ may be denoted as
$$
(x_0,y_0),\ldots , (x_p, y_p).
$$
Then
$$
(A^o\times A)_{p/}((x_0,y_0),\ldots , (x_p, y_p))=
A_{p/}(x_p,\ldots , x_0) \times A_{p/}(y_0,\ldots , y_p).
$$
We put
$$
ArrFam(A)_{p}((x_0,y_0),\ldots , (x_p, y_p)):=
A_{2p+1 /}(x_p,\ldots , x_0,y_0,\ldots , y_p).
$$
The morphism from here to $(A^o\times A)_{p/}((x_0,y_0),\ldots , (x_p, y_p))$
is composed of the two structural morphisms of $A$ for the two morphisms
$$
p\tworightarrows (2p+1)
$$
in $\Delta$, these two morphisms being
$$
\{ 0,\ldots , p\} \mapsto \{ 0,\ldots , p\}\subset \{ 0,\ldots , 2p+1\}
$$
and
$$
\{ 0,\ldots , p\} \mapsto \{ p+1,\ldots , 2p+1\}\subset \{ 0,\ldots , 2p+1\}
$$
respectively. We put
$$
ArrFam(A)_{p}:= \coprod _{(x,y)} ArrFam(A)_{p}((x_0,y_0),\ldots , (x_p, y_p))
$$
where the coproduct is taken over all sequences of objects
$(x_0,y_0),\ldots , (x_p, y_p)$. This maps to
$$
{\bf r}(A^o \times A)_p =
\coprod _{(x,y)}(A^o\times A)_{p/}((x_0,y_0),\ldots , (x_p, y_p))
$$
by the above map. We leave it to the reader to specify the structural morphisms
for the simplicial structure of $ArrFam(A)$ mapping to the structural
morphisms of ${\bf r}(A^o \times A)$.

The object $ArrFam(A)$ is a precartesian family. First note that
$$
ArrFam(A)_1 = \coprod _{(x,y)} A_{1/}(x,y)
$$
so the $n$-categories indexed by this family are indeed the ``$Hom$''
$n$-categories $A_{1/}(x,y)$ of $A$. The precartesian property comes from the
fact that the morphisms
$$
A_{2p+1 /}(x_p,\ldots , x_0,y_0,\ldots , y_p) \rightarrow
$$
$$
A_{p/}(x_p,\ldots , x_0) \times A_{1/}(x_0, y_0) \times A_{p/}(y_0,\ldots , y_p)
$$
are equivalences of $n$-categories (which in turn follows easily from the Segal
condition for the $n+1$-category $A$).

Now replace $ArrFam(A)$ by an equivalent fibrant morphism of $n$-stacks over
$\Delta$,
$$
ArrFam(A)'\rightarrow {\bf r}(A^o \times A).
$$
This is a cartesian family so it corresponds to a morphism
$$
{\bf Arr}_A: A^o \times A \rightarrow nFAM = nCAT'.
$$

We don't treat here the problem of showing that this is uniquely defined
up to coherent higher homotopies by the process of choosing the fibrant
replacement $ArrFam(A)'$. (But in order to be completely rigorous one needs to
do this.)

\begin{exercise}
Show that for $A=nCAT$, the precartesian family $ArrFam(A)$ is exactly the
pullback of the universal precartesian family over $nCAT$ for the explicit
(strict) arrow family
$$
{\bf Arr}_{nCAT}: nCAT^o \times nCAT \rightarrow nCAT
$$
which was constructed previously (\ref{arrownstack2}). In particular the present
general construction of ${\bf Arr}_A$ coincides (up to coherent homotopy)
with the
previous specific ${\bf Arr}_{nCAT}$ in this case.
\end{exercise}

\begin{exercise}
Verify that the above construction of the ``arrow
family'' works {\em mutatis mutandis} to construct  ``arrow families'' for
$n+1$-stacks, namely if $\underline{A}$ is an $n+1$-stack then we obtain
$$
{\bf Arr} _{\underline{A}}: \underline{A}^o\times \underline{A} \rightarrow
n\underline{FAM} = n\underline{STACK}'.
$$
\end{exercise}

\subnumero{Base change for relative shapes}

We now get back to the question of relative nonabelian cohomology and base
change.
Fix the site $\Gg$. Suppose we have fixed
realms
$$
\underline{R}, \underline{A} \subset n\underline{STACK},
$$
(recall that a ``realm'' is just a saturated full substack of
$n\underline{STACK}$).

Suppose $\underline{B}$ is an $n+1$-stack and suppose that $\Ff \rightarrow {\bf
r}\underline{B}$ is a cartesian family of $n$-stacks over $\underline{B}$.
We obtain a classifying morphism
$$
[\Ff ]: \underline{B} \rightarrow n\underline{STACK}'.
$$
This morphism is tautological if we take $n\underline{FAM}$ as the fibrant
replacement $n\underline{STACK}'$ (but is essentially well-defined for any
fibrant
replacement, by comparing with the fibrant replacement $n\underline{FAM}$).

We have an arrow family
$$
{\bf Arr}_{n\underline{STACK}'} :
(n\underline{STACK}')^o\times n\underline{STACK}'
\rightarrow n\underline{STACK}';
$$
this can be obtained either by homotopy invariance from the
arrow family ${\bf Arr}_{n\underline{STACK}}$ which was constructed in
\ref{arrownstack2}, or from the general construction of \S 5.6 above.

Composing the classifying map $[\Ff ]$ with the arrow family we obtain
$$
{\bf Arr}([\Ff ], -): \underline{B}^o \times n\underline{STACK} \rightarrow
n\underline{STACK}'.
$$
As before, in view of the definition of internal $\underline{Hom}$, this map
corresponds to a map
$$
{\bf Shape}(\Ff /B) : \underline{B}^o\rightarrow
\underline{Hom}(n\underline{STACK}, n\underline{STACK}').
$$
We call this map the {\em relative shape} for the cartesian family $\Ff /B$. It
may also be seen as the composition of the absolute shape map
$$
{\bf Shape} : (n\underline{STACK}')^o\rightarrow
\underline{Hom}(n\underline{STACK}', n\underline{STACK}'),
$$
with the cartesian family $[\Ff ]:B^o \rightarrow (n\underline{STACK}')^o$.

If the base is the final object $B=\ast$ then the relative shape map coincides
with the absolute shape map defined previously.

The relative shape map ${\bf Shape}(\Ff /B) $ is supposed to be related to
relative nonabelian cohomology in the same way that the absolute shape ${\bf
Shape}(\Ff )$ is related to nonabelian cohomology. Unfortunately, I haven't had
the time (nor does the present paper permit the space) to verify this
compatibility, so we formulate this as a problem rather than a theorem.

\begin{problem}
\label{compatibility}
Suppose $\Ff \rightarrow B$ is a morphism of $n$-stacks, where $B$ is an
$\ngr$-stack. Let $\underline{B}$ be $B$ considered as an $n+1$-stack of
groupoids, and let $\Ff ^{fam /B} \rightarrow {\bf r}\underline{B}$
be the cartesian family associated to $\Ff /B$ by Proposition
\ref{correspondence}.   By composing the relative shape map for $\Ff ^{fam
/B}/\underline{B}$ with the morphism of evaluation on $T$ in the first variable
$$
{\bf eval}(T):
\underline{Hom}(n\underline{STACK}',n\underline{STACK}')\rightarrow
n\underline{STACK}' ,
$$
we get a morphism
$$
{\bf Shape}(\Ff ^{fam /B} /\underline{B})\circ {\bf eval}(T):
\underline{B}^o \rightarrow
n\underline{STACK} '.
$$
The problem is to show that this morphism is the same as the classifying
morphism for the cartesian family corresponding (again by \ref{correspondence})
to the structural morphism for the relative nonabelian cohomology of $\Ff
/B$ with
coefficients in $T$,
$$
\underline{Hom}(\Ff /B, T)\rightarrow B.
$$
\end{problem}

One remark is that for an object $b\in B_0(Z)$, the fiber of
$\underline{Hom}(\Ff /B, T)$ over $b$ is the same as $\underline{Hom}(\Ff _b,
T|_{\Gg /Z})$. On the other hand, the restriction of the morphism $[\Ff ^{fam
/B}]$ to $b$ is just the classifying morphism for $\Ff _b$, so (in view of the
definition of ${\bf Arr}_{n\underline{STACK}}$), the restriction
$$
{\bf Shape}(\Ff ^{fam /B} /\underline{B})\circ {\bf eval}(T)|_{\{ b\}}
$$
is point of $n\underline{STACK} '(Z)$ corresponding to
$\underline{Hom}(\Ff _b, T|_{\Gg /Z})$. Thus, on the level of fibers over
objects of $B$, the compatibility of the above problem is easy. The problem  is
to verify that the ``action'' (plus higher homotopy coherencies) of the
morphisms of $B$ on these fibers, are the same in the two cases.

In what follows, we will tacitly assume that the above compatibility is known.

\begin{parag}
\label{relativezoo1}
Suppose now that $\underline{R} \subset n\underline{STACK}$ is a realm.
Suppose as above that $\Ff \rightarrow {\bf r} B$ is a cartesian family of
$n$-stacks indexed by a base $n+1$-stack $B$.  We obtain a restricted shape map,
which we call the {\em $\underline{R}$-shape of $\Ff /B$},
$$
{\bf Shape}_{\underline{R}}(\Ff /B): B^o \rightarrow
\underline{Hom}(\underline{R}, n\underline{STACK}).
$$
We can again pose a general type of problem analogous to \ref{zoo2} in this
situation: given another realm $\underline{A}$, when does the
$\underline{R}$-shape of $\Ff /B$ have answers in $\underline{A}$, which means,
when does the relative shape map factor through a map
$$
{\bf Shape}_{\underline{R}}^{\underline{A}}(\Ff /B): B^o \rightarrow
\underline{Hom}(\underline{R}, \underline{A})\,  ?
$$
As before, we shall say that
{\em ${\bf Shape}_{\underline{R}}^{\underline{A}}(\Ff /B)$ exists} for this
condition; and we can also write this statement as
$$
\underline{R} \stackrel{{\bf Shape}(\Ff /B)}{\longrightarrow}
\underline{A}.
$$

\newparag{relativezoo2}
If $\underline{P}$ is the largest realm such that $(\underline{P},
\underline{R}, \underline{A})$ is well-chosen, then the above condition is
equivalent to the condition that the classifying map for the cartesian family
takes values in $\underline{P}$,
$$
[\Ff ] : B\rightarrow \underline{P}' \subset n\underline{STACK}'.
$$
\newparag{relativezoo3}
This condition can be rewritten more concretely as follows: for any $Z\in \Gg$
and any $b\in B_0(Z)$ we obtain the fiber of our family $\Ff _b$ which is an
$n$-stack over $\Gg /Z$. We need this $n$-stack to be contained in
$\underline{P}(Z)$; thus we need to know that for any $Y\rightarrow Z$ and any
$n$-stack $T\in \underline{R}(Y)$, the nonabelian cohomology
$\underline{Hom}(\Ff _b|_{\Gg /Y}, T)$ should lie in $\underline{A}(Y)$.

Note that the restriction of $\Ff _b$ to $\Gg /Y$ is just $\Ff _{b_Y}$ for $b_Y$
the restriction of $b$ to $b_Y\in B_0(Y)$. Thus we can rewrite the above
condition combining together the choice of $Z$ and $Y$. We state this as a
lemma.
\end{parag}

\begin{lemma}
\label{relativezoo4}
Suppose $\Ff \rightarrow {\bf r}B$ is a cartesian family of $n$-stacks over an
$n+1$-stack $B$. Suppose that $\underline{R}$ and $\underline{A}$ are realms.
Then the relative $\underline{R}$-shape of $\Ff /B$ takes values in
$\underline{A}$, i.e.
${\bf Shape}_{\underline{R}}^{\underline{A}}(\Ff /B)$ exists or which we
write as
$$
\underline{R} \stackrel{{\bf Shape}(\Ff /B)}{\longrightarrow}
\underline{A},
$$
if and only if the following condition holds: for any $Z\in \Gg$, any object
$b\in B_0(Z)$, and any $n$-stack $T\in \underline{R}(Z)$, the nonabelian
cohomology $n$-stack of the fiber
$$
\underline{Hom}(\Ff _b, T)
$$
lies in $\underline{A}(Z)$.
\end{lemma}
\eop

\begin{parag}
\label{relativezoo5}
We will often want to apply the lemma to the case where the cartesian family
$\Ff\rightarrow {\bf r}B$ comes from a morphism of $n$-stacks $F\rightarrow B$
via the correspondence of Proposition \ref{correspondence}. In this case, an
object $b\in B_0(Z)$ may be viewed as a morphism $Z\rightarrow B$ and the fiber
is then $\Ff _b= F\times _BZ$.
\end{parag}

Finally we come to the statement of our base-change theorem for the relative
shape maps.

\begin{theorem}
\label{basechange2}
The relative shape map satisfies base change: suppose $\underline{R}$ is a
realm, and suppose $\Ff \rightarrow {\bf r} \underline{B}$ is a
cartesian family. If $\varphi : \underline{B}'\rightarrow \underline{B}$ is a
morphism of $n+1$-stacks and if  $\Ff '\rightarrow {\bf r} \underline{B}'$
denotes
the base change of the cartesian family $\Ff$, then
$$
{\bf Shape}_{\underline{R}}(\Ff '/\underline{B}')=
{\bf Shape}_{\underline{R}}(\Ff /\underline{B})\circ \varphi ^o.
$$
\end{theorem}
{\em Proof:}
This again is tautological from the definition (all of the work has already
been done in the construction of the relative shape map itself).
\eop

Modulo the compatibility of Problem \ref{compatibility}, this theorem contains
the statement of Theorem \ref{basechangeA}.

\numero{Very presentable $n$-stacks}
\label{presentablepage}

\subnumero{Presentability}

Work on the site $Sch /k$ of noetherian schemes over a field $k$ of
characteristic zero, with the etale topology. Fix $n$ where necessary.

A couple of particular examples of directories of Serre classes were defined in
\cite{Simpson}(vii). These give rise to the notions of {\em
presentable $n$-stack} and {\em very presentable $n$-stack}. The main reference
for these notions is \cite{Simpson}(vii), and we give only a brief
discussion here. There are three classes of sheaves involved. In degrees $i\geq
2$, we put the {\em vector sheaves}, 	a notion which we discuss in more detail
in subsection 6.2 below. In degree $1$ we put the {\em affine presentable group
sheaves} which we discuss briefly in the  paragraph after next. For the
degree $0$
part, we refer directly to \cite{Simpson}(vii).

Very briefly, the category of {\em vector sheaves} on an affine noetherian
scheme $X$ is the smallest subcategory of the category of sheaves of abelian
groups on $Sch /X$ which contains the finite rank vector bundles on $X$ and
which
is closed under kernel and cokernel. (It is also closed under extensions but
this is not part of the definition.) This notion will be recalled in greater
detail in 6.2 below.

In \cite{Simpson}(vii) is defined a notion of {\em presentable group sheaf}
over a scheme $X$ (in characteristic $0$). This is a certain type of sheaf of
groups on the big site $Sch /X$.  We don't give the definition here, but note
some of the properties. First, the category of presentable group sheaves is
closed under kernels, cokernels (by normal subgroups) and extensions. A
presentable group sheaf $G$ has a Lie algebra sheaf $Lie(G)$ which is a vector
sheaf. It also has a ``connected component'' $G^{0}$, and the behavior of
morphisms for example, on the connected component, is determined by their
restriction to $Lie(G)$. If $X=Spec (k)$ is a point, then a presentable group
sheaf is just a group scheme of finite type. If $G$ is a presentable group sheaf
over $X=Spec (A)$ with $A$ artinian of finite type (hence, of finite
length) over
a field $k$ then  the projection $p:X\rightarrow Spec (k)$ induces a
push-forward
presentable group sheaf $p_{\ast}(G)$ on $Spec (k)$, in other words
$p_{\ast}(G)$
is a group scheme of finite type. In this sense, one can at least calculate what
is going on over an artinian base.

We say that a presentable group sheaf $G$ over $X$ is {\em affine presentable}
if, for every finte type artinian scheme $Z\rightarrow X$ denoting the
projection
by $p:Z\rightarrow Spec (k)$, the group scheme $p_{\ast}(G|_Z)$ is affine.

Note that a vector sheaf is itself an abelian affine presentable group
sheaf, and one could conjecture that any abelian affine presentable group sheaf
such that the $p_{\ast}(G|_Z)$ are unipotent, is a vector sheaf. It might
be reasonable to replace the notion of ``vector sheaf'' by a notion of
``unipotent affine presentable group sheaf'' but instead of proceeding like
this,
we just use vector sheaves.

Recall also from \cite{Simpson}(vii) that we have defined a condition
P$3\frac{1}{2}$ for sheaves of sets over $Sch/X$. We don't go into the
definition
here at all.

Define a directory of Serre classes $Pres$ by setting $Pres_i(X)$ equal to
the abelian presentable group sheaves for $i\geq 2$; the presentable group
sheaves for $i=1$; and the sheaves of sets satisfying condition
P$3\frac{1}{2}$ for $i=0$.
We obtain the realm
\label{pepage}
$$
\underline{PE} := \underline{M}^{Pres}
$$
of {\em
presentable $n$-stacks of groupoids}.

Define a subdirectory of Serre classes $VPres \subset Pres$ as follows:
$VPres_0=Pres_0$. For $i=1$, $VPres_1(X)$ is the class of affine presentable
group sheaves over $X$; and for $i\geq 2$, $VPres_i(X)$ is the class of vector
sheaves over $X$. We obtain the realm
\label{vppage}
$$
\underline{VP}:= \underline{M}^{VPres}
$$
of {\em
very presentable $n$-stacks of groupoids}.

Define a further subdirectory of Serre classes $Vect \subset VPres$ as follows:
For $i=0,1$, $Vect_i(X)$ is the class containing only the trivial punctual sheaf
$\ast _X$. For $i\geq 2$, $Vect_i(X)= VPres_i(X)$ is the class of vector sheaves
over $X$.  We obtain the realm
$$
\underline{AV}:= \underline{M}^{Vect}
$$
of {\em
$1$-connected very presentable $n$-stacks of groupoids}. This realm will be
important in our discussion of representability of shape below.

(The realm $\underline{CV}$ of connected very presentable $n$-stacks of
groupoids
also comes from a directory of Serre classes $CVPres$ with $CVPres_0(X) = \{
\ast _X\}$ and $CVPres_i(X) = VPres_i(X)$ for $i\geq 1$.)

Note also that we obtain the necessary ingredients to define the realm
$\underline{FV}$, see \ref{index1}.

All of these realms give rise to shape theories and theories of nonabelian
cohomology. Furthermore, the realms of
coefficients  $\underline{PE}$ and $\underline{VP}$
are also appropriate as realms of answers, as we now state.

\begin{lemma}
\label{vpcompatible}
Both realms $\underline{PE}$ and $\underline{VP}$ are compatible with
truncations, finite limits, and extensions. In particular, if $K$ is a finite
CW complex and $W$ is the constant prestack with values $\Pi _n(K)$, then
${\bf Shape}(W):= \underline{Hom}(W, -)$ takes  $\underline{PE}$ to
$\underline{PE}$ (resp. $\underline{VP}$ to $\underline{VP}$).
\end{lemma}
{\em Proof:}
We refer to \cite{Simpson}(vii) for these results. Note first of all that
the notions of presentability (resp. very presentability) which we use here are
the same as those of \cite{Simpson}(vii) (whereas these notions differ very
slightly---in the condition on $\pi _0$---from the preliminary version defined
in \cite{Simpson}(v)). This compatibility can be seen from Theorem 10.1
of \cite{Simpson}(vii).

Compatibility with truncations is tautological cf \ref{rmktrunc}.

Compatibility with finite limits \ref{closurelimits} is the statement of
Corollary
10.6 of \cite{Simpson}(vii).

Finally, compatibility with extensions \ref{closureextension}
comes from the
statement of Corollary 10.9 of \cite{Simpson}(vii). To see this note that
if one takes the hypothesis of \ref{closureextension} (we use the same
notations)
and applies it to any morphism $X\rightarrow B$ inducing a surjection on
$\pi _0$
(noting that such a morphism exists by the definition of presentability) then
Corollary 10.9 implies that $A$ is presentable (resp. very presentable).
\eop

\begin{parag}
As was stated in \ref{index4}, the connected very presentable
$n$-stacks over a point are those which were considered in \ref{degree3},
i.e. the $n$-stacks $T$ with $\pi _0(T)=\ast$, with $\pi _1(T)$ an affine
algebraic group of finite type, and $\pi _i(T)$ finite dimensional vector spaces
for $i\geq 2$.
\end{parag}

\subnumero{More on vector sheaves}
We discuss in a bit more detail the category of vector sheaves
over a base scheme $X$. The references for this notion are as follows.
Auslander \cite{Auslander} defined a closely related notion of ``coherent
functor'' which meant functor of $A$-modules (rather than of $A$-algebras).
Apparently M. Artin wrote a letter to A. Grothendieck in the 1960's
containing the version where functors of $A$-algebras are considered, but this
never seems to have been made public. Hirschowitz gave the first publicly
available version of the definition of vector sheaf which he called
``$U$-coherent sheaf'' in \cite{Hirschowitz}. Hirschowitz obtained most of the
important properties such as involutivity of duality. Jaffe later considered
this notion in \cite{Jaffe}. Some further properties of vector sheaves were
proven in \cite{Simpson}(v), (vii).

For
affine $X$, the category of {\em vector sheaves} over $X$ is the smallest
subcategory of sheaves on the big site $Sch /X$ which contains the vector
bundles
over $X$ and which is closed under kernel and cokernel. For general $X$ the
category of {\em vector sheaves over $X$} is defined by the condition of
locality
plus compatibility with the affine case.

The category of vector sheaves over $X$ is an abelian subcategory of the
category of sheaves of groups on $Sch/X$.  The coherent sheaves on $X$,
which are
defined as cokernels
$$
\Oo ^a \rightarrow \Oo ^b \rightarrow \Ff \rightarrow 0,
$$
are vector sheaves. Coherent sheaves are injective objects. Their duals, which
we call {\em vector schemes}, are the group-schemes with vector space structure
(but not necessarily flat) over $X$. These admit dual presentations as
kernels of maps $\Oo ^b\rightarrow \Oo ^a$. They are projective objects (at
least
if the base $X$ is affine). If $X$ is affine, then any vector sheaf $U$ admits
resolutions
$$
0\rightarrow V\rightarrow V' \rightarrow V'' \rightarrow U
\rightarrow 0 $$
with $V$, $V'$ and $V''$ vector schemes; and
$$
0\rightarrow U \rightarrow \Ff \rightarrow \Ff '\rightarrow \Ff '' \rightarrow 0
$$
with $\Ff$, $\Ff '$, $\Ff ''$ coherent sheaves.

The most surprising property \cite{Hirschowitz} is
that the duality functor $U^{\ast} := Hom (U,\Oo )$ is exact, and is an
involution. This is due to the fact that we take the big site $Sch /S$ rather
than the small Zariski or etale sites.

Another interesting point is that there are two different types of tensor
products of vector sheaves: the {\em tensor product}
$$
U\otimes _{\Oo} V := Hom (U, V^{\ast} )^{\ast},
$$
and the {\em cotensor product}
$$
U\otimes ^{\Oo} V := Hom (U^{\ast}, V).
$$
These are not the same (although they coincide for coherent sheaves) and in
particular they don't have the same exactness properties. Neither of them is
equal to the tensor product of sheaves of $\Oo$-modules. See \cite{Simpson}(v)
and (vii) for further discussion.

\subnumero{The Breen calculations for vector sheaves in characteristic zero}

The above facts work in any characteristic and depend on the
$\Oo$-module structure. However, in characteristic zero, vector sheaves have the
additional property that the morphisms $U\rightarrow V$ of sheaves of abelian
groups over $Sch /\cc $ are automatically morphisms of $\Oo$-modules, see
\cite{Simpson}(v), (vii). Similarly, extensions of sheaves of abelian
groups, between two vector sheaves, are again vector sheaves. These properties
persist for the higher $Ext^i$, see Corollary \ref{ext} below. These properties
do not remain true in characteristic $p$: the basic problem is that Frobenius
provides a morphism $\Oo \rightarrow \Oo$ of sheaves of abelian groups, which is
not a morphism of sheaves of $\Oo$-modules.

We need to be able to calculate the Postnikov invariants of the
$n$-stacks in a realm such as $\underline{VP}$.  For example, we would like
to calculate $H^i(K(V , m), W )$, when $V$ and $W$ are vector sheaves. In the
case $V=W=\Oo$, this calculation is the algebraic analogue of the classical
Eilenberg-MacLane calculations. The algebraic version is the subject of
Breen's work \cite{Breen}. Breen concentrated on the case of
characteristic $p$ in \cite{Breen}(i) (where the answer is very
complicated because of Frobenius-type elements). The characteristic $0$
version is implicit in \cite{Breen} because it is strictly easier than the
characteristic $p$ case.

Here we work only in characteristic $0$.
For the proofs of statements in this section, the reader may refer to the two
appendices of \cite{Simpson}(xii) (as well as refering to \cite{Breen}).

\begin{theorem}
\label{bc}
Suppose $S$ is  a scheme over $Spec (\qq  )$.
Suppose $V$ is a vector scheme over $X$ and suppose $\Ff$ is a coherent sheaf
over $S$. Then for $n$ odd we have
$$
H^i(K(V/S, n)/S, \Ff ) = \Ff \otimes _{\Oo} \bigwedge _{\Oo}^{i/n} (V^{\ast}).
$$
For $n$ even we have
$$
H^i(K(V/S, n)/S, \Ff ) = \Ff \otimes _{\Oo} Sym _{\Oo}^{i/n} (V^{\ast}).
$$
In both cases the answer is $0$ of $i/n$ is not an integer. The multiplicative
structures on the left sides, in the case $\Ff = \Oo$, coincide with the
obvious ones on the right sides. In the case of arbitrary $\Ff$, the
natural structures
of modules over the cohomology with coefficients in $\Oo$, on both sides,
coincide. \end{theorem}

The main ingredient in the proof is the following result.

\begin{proposition}
\label{somecomplexes}
Suppose $S$ is a scheme and $V$ is a vector sheaf on $S$. Then the complexes
$$
\ldots \Ff \otimes _{\Oo}\bigwedge _{\Oo} ^j \otimes _{\Oo} Sym _{\Oo}^k V
\rightarrow \Ff \otimes _{\Oo}\bigwedge _{\Oo} ^{j-1}
\otimes _{\Oo} Sym _{\Oo}^{k+1} V
\ldots $$
and
$$
\ldots \Ff \otimes _{\Oo}\bigwedge _{\Oo} ^j
\otimes _{\Oo} Sym _{\Oo}^k V \rightarrow
\Ff \otimes _{\Oo}\bigwedge _{\Oo} ^{j+1}
\otimes _{\Oo} Sym _{\Oo}^{k-1} V \ldots
$$
are exact as sequences of vector sheaves (i.e. as sequences of sheaves on the
site $Sch /\cc $).
\end{proposition}
{\em Proof:}
See \cite{Simpson}(xii).
\eop

With this proposition, the proof of Theorem \ref{bc} is essentially just a
rewriting of the standard spectral-sequence proof of the Eilenberg-MacLane
calculations. See \cite{Simpson}(xii).
\eop

\begin{remark}
\label{specseqcalx}
Using the result of Theorem \ref{bc} one can envision the calculation of
$H^i(K(U/S, n)/S, W)$ for any vector sheaves $U$ and $W$ over $S$. In particular
this would lead to a calculation for the case where $U$ and $W$ are coherent
sheaves (this question was mentionned in \S 3). To treat $W$, resolve by
coherent sheaves
$$
0\rightarrow W \rightarrow \Ff \rightarrow \Ff ' \rightarrow \Ff '' \rightarrow
0
$$
which gives a pair of long exact sequences relating
$H^i(K(U/S, n)/S, W)$ to
$$
H^i(K(U/S, n)/S, \Ff), \;\; H^i(K(U/S, n)/S, \Ff '), \;\;
H^i(K(U/S, n)/S, \Ff '').
$$
Thus we reduce to the case of coefficients in a coherent sheaf $\Ff$.
To treat $U$ (for cohomology with coefficients in a coherent sheaf $\Ff$),
resolve
by vector schemes
$$
0\rightarrow V'' \rightarrow V' \rightarrow V \rightarrow U\rightarrow 0
$$
which breaks into two short exact sequences
$$
0\rightarrow V'' \rightarrow V' \rightarrow U' \rightarrow 0
$$
and
$$
0\rightarrow U' \rightarrow V \rightarrow U\rightarrow 0.
$$
These give fibration sequences (relative to $S$)
$$
K(V'/S, n+1) \rightarrow K(U'/S, n+1) \rightarrow K(V''/S , n+2),
$$
and
$$
K(V/S, n) \rightarrow K(U/S,n) \rightarrow K(U'/S, n+1).
$$
We get Leray spectral sequences, the first going from something for $V''$
and $V'$
to the the answer for $U'$, and the second going from this plus the answer for
$V$, to the answer for $U$. For example, if $n$ is odd the first spectral
sequence is $$
E_2^{(n+2)p , (n+1)q} = \Ff \otimes  \bigwedge ^p(V'')^{\ast}
\otimes Sym ^q(V')^{\ast} \Rightarrow H^{(n+2)p+(n+1)q}(K(U'/S, n+1), \Ff ).
$$
The first nonzero differential comes from the map $V''\rightarrow V'$.
I don't know if there are higher differentials.

To get the $E_2$ term of the second spectral sequence (for cohomology of
$K(U/S,n)$ with coefficients in a coherent sheaf $W$) one has to apply the first
spectral sequence's calculation of the cohomology of $K(U'/S, n+1)$ with
coefficients in the coherent sheaf $H^j(K(V/S, n), \Ff )$.

It would be interesting to carry out the calculation but I haven't done that.
In any case, it is clear from the fact that the category of vector sheaves is
closed under kernel, cokernel and extension, that the answer $H^i(K(U/S, n), W)$
will be a vector sheaf on $S$.
\end{remark}

\begin{corollary}
\label{cohVPisVS}
Suppose $S$ is a scheme and $T\rightarrow S$ is a relatively $1$-connected very
presentable $n$-stack. Then for any vector sheaf $V$ on $S$,
$$
H^i(T/S, V)
$$
is a vector sheaf.
\end{corollary}
{\em  Proof:}
Use the Leray spectral sequence for the Postnikov tower of $T$, and Theorem
\ref{bc}. See \cite{Simpson}(xii).
\eop

\begin{corollary}
\label{homVP}
If $S$ is a scheme and $T\rightarrow S$ and $T'\rightarrow S$ are relatively
$1$-connected very presentable $n$-stacks then $\underline{Hom} (T/S, T'/S)$ is
very presentable.
\end{corollary}
Note that $Aut(V)$ is a very presentable group
sheaf when $V$ is a vector sheaf, see \cite{Simpson}(vii).
\eop

The following application  was the original motivation for
Breen's calculations of the cohomology of the Eilenberg-MacLane
sheaves \cite{Breen}. From our version Theorem \ref{bc}, we obtain
the corresponding result in the relative case in characteristic zero. Similar
corollaries were stated for example for cohomology with coefficients in the
multiplicative group $\Gm$, in \cite{Breen}.

\begin{corollary}
\label{ext}
{\rm (\cite{Simpson}(v) Corollary 3.11)}
Suppose $U,V$ are vector sheaves over a scheme $S$. Let $Ext^i_{\rm gp} (U,V)$
denote the $Ext$ sheaves between $U$ and $V$ considered as sheaves of abelian
groups on $Sch /S$, let $Ext^i_{\rm vs}(U,V)$ denote the $Ext$ sheaves between
$U$ and $V$ considered as vector sheaves on $S$. Then the natural morphisms are
isomorphisms
$$
Ext^i_{\rm vs} (U,V)
\stackrel{\cong}{\rightarrow}Ext^i_{\rm gp} (U,V).
$$
The $Ext^i$ vanish for $i>2$.
\end{corollary}
{\em Proof:}
See \cite{Simpson}(xii).
\eop

\subnumero{Representability of very presentable shape}

Under certain circumstances, we can hope that the
very presentable shape
$$
{\bf Shape} _{\underline{VP}}(\Ff )
$$
be {\em
representable}. By this we mean that there should be a morphism $\Ff \rightarrow
R$ from $\Ff$ to a very presentable $n$-stack $R$ such that for any
other very presentable $n$-stack $T$ we have
$$
Hom (R , T)\stackrel{\cong}{\rightarrow} Hom (\Ff , T).
$$
The main case where the representing stack $R$ will exist, is when
the very presentable shape of $\Ff$ is simply connected; in this case the
representing stack $R$ will be simply connected. Also in this case, the
representing $R$ can be viewed as representing the shape of $\Ff$ on the
smaller realm $\underline{AV}$ of $1$-connected very presentable $\ngr$-stacks.
However, it turns out that the condition needed in order to be able to obtain
representability of ${\bf Shape} _{\underline{AV}}(\Ff )$ is much weaker,
basically just vanishing of $H^1(\Ff , V)$ for vector sheaves $V$. Thus we can
sometimes obtain representability of the $\underline{AV}$-shape
even when the $\underline{VP}$-shape is not representable. In this situation,
the representing object for the $\underline{AV}$-shape
is somewhat akin to Quillen's ``plus'' construction.

Because of what was said in the previous paragraph, consideration of
representability of the the $\underline{VP}$-shape has no advantages over
consideration of representability of the $\underline{AV}$-shape, so we restrict
our attention to the latter.

\begin{parag}
Recall from above that $\underline{AV}$ denotes the realm of $1$-connected
$\ngr$-stacks whose higher homotopy group sheaves are vector sheaves. If $\Ff
\rightarrow \Ee$ is a morphism of $\ngr$-stacks then we obtain the {\em
$\underline{AV}$-shape of $\Ff $ relative to $\Ee$} which is a morphism
$$
{\bf Shape}_{\underline{AV}}(\Ff /\Ee ): \Ee \rightarrow
\underline{Hom}(\underline{AV}, \ngr\underline{STACK}).
$$
We say that this shape is {\em representable} if there is a morphism
$$
{\rm rep}(\Ff /\Ee ): \Ee \rightarrow \underline{AV}\subset
\ngr\underline{STACK}
$$
and a morphism of cartesian families over $\Ee$
$$
\Ff \rightarrow {\rm rep}(\Ff /\Ee ),
$$
such that this morphism induces an equivalence of shape functors
$$
{\bf Shape}_{\underline{AV}}(\Ff /\Ee ) \stackrel{\cong}{\rightarrow}
{\bf Shape}_{\underline{AV}}({\rm rep}(\Ff /\Ee ) ).
$$
The cartesian family corresponding to the morphism ${\rm rep}(\Ff /\Ee )$
also corresponds by Proposition \ref{correspondence} to a  morphism of type
$\underline{AV}$ of $\ngr$-stacks
$$
{\bf rep}(\Ff /\Ee )\rightarrow \Ee ,
$$
and the morphism from $\Ff$ corresponds to a diagram plus homotopy of
commutativity
$$
\begin{array}{ccc}
\Ff & \rightarrow &{\bf rep}(\Ff /\Ee )\\
\downarrow && \downarrow \\
\Ee & = & \Ee .
\end{array}
$$
The universal property is again that this morphism induces an equivalence of
shapes. This can be rephrased as follows: that for any morphism from a scheme
$$
Z\rightarrow \Ee
$$
and any $T\in \underline{AV}(Z)$ (i.e. $\ngr$-stack $T$ over $Z$ which is
$1$-connected and whose higher homotopy group sheaves are vector sheaves over
$Z$), we have an equivalence of nonabelian cohomology
$$
\underline{Hom}({\bf rep}(\Ff /\Ee )\times _{\Ee}Z/Z, T)
\stackrel{\cong}{\rightarrow}
\underline{Hom}(\Ff\times _{\Ee}Z/Z, T).
$$

\newparag{sflksf}
We make a first reduction.
Suppose for the moment that the base is already a scheme, i.e. we have a
morphism of $\ngr$-stacks $\Ff \rightarrow Z$ with $Z$ a scheme. Suppose
$$
R\rightarrow Z
$$
is an $\ngr$-stack over $Z$, in $\underline{AV}(Z)$ (which means that it is
simply connected and the higher homotopy group sheaves are vector sheaves), and
with a morphism of $\ngr$-stacks over $Z$
$$
\Ff \rightarrow R.
$$
We say that this situation is a {\em weak representing object for the
$\underline{AV}(Z)$-shape of $\Ff /Z$} if for any $T\in \underline{AV}(Z)$,
it induces an equivalence of nonabelian cohomology stacks over $Z$,
$$
\underline{Hom}(R/Z,T) \stackrel{\cong}{\rightarrow}
\underline{Hom}(\Ff /Z, T).
$$
In this case we write $R={\bf rep}^{\rm wk}(\Ff /Z)$. If it exists, the
representing object is unique up to coherent homotopy of all orders.
\end{parag}

\begin{lemma}
\label{reduction1}
Suppose that $\Ff \rightarrow \Ee$ is a morphism of $\ngr$-stacks
such that for every morphism from a scheme
$Z\rightarrow \Ee$, the weak representing object ${\bf rep}^{\rm wk}(\Ff \times
_{\Ee}Z/Z)$ exists.
Then the representing object ${\bf rep}(\Ff /\Ee )$ exists
(and its pullback to any $Z$ is the same as the weak representing object).
\end{lemma}
{\em Proof (new in {\tt v2}):}
Note first of all that if there exists an object $R\rightarrow \Ee$ with
morphism $\Ff \rightarrow R$ over $\Ee$, such that for any morphism from
a scheme $Z\rightarrow \Ee$, the pullback $R\times _{\Ee}Z$ together with its
morphism from $\Ff \times _{\Ee}Z$ is equivalent to
${\bf rep}^{\rm wk}(\Ff \times
_{\Ee}Z/Z)$, then $R$ is a representing object $R={\bf rep}(\Ff /\Ee )$.
This is clear from the definitions.
The problem is to construct $R$.

{\em Claim:} if
$Z'\rightarrow Z \rightarrow \Ee$ then the natural morphism
$$
{\bf rep}^{\rm wk}(\Ff \times _{\Ee}Z'/Z')
\rightarrow
{\bf rep}^{\rm wk}(\Ff \times _{\Ee}Z/Z) \times _ZZ'
$$
is an equivalence.

Here is the proof of the claim. For brevity,
write
$$
P':=
{\bf rep}^{\rm wk}(\Ff \times _{\Ee}Z'/Z')
$$
and
$$
P:={\bf rep}^{\rm wk}(\Ff \times _{\Ee}Z/Z).
$$
We have
$$
H^i(P\times _ZZ'/Z', \Oo )
= H^i(P/Z, \Oo )|_{\Gg /Z'};
$$
this is in fact tautological because $H^i(P/Z, \Oo )$ is the sheafification
of the functor $Y\mapsto H^i(P\times _ZY, \Oo )$ on $Y\in \Gg /Z$.
On the other hand, recall that
$$
H^i(P/Z, \Oo )= H^i(\Ff \times _{\Ee}Z/Z, \Oo )
$$
and
$$
H^i(P'/Z', \Oo )= H^i(\Ff \times _{\Ee}Z'/Z', \Oo )
$$
since $P$ is a weak representing object for $\Ff \times _{\Ee}Z/Z$
and $P'$ is a weak representing object for $\Ff \times _{\Ee}Z'/Z'$.
By the same tautological argument,
$$
H^i(\Ff \times _{\Ee}Z'/Z', \Oo ) = H^i(\Ff \times _{\Ee}Z/Z, \Oo )|_{\Gg /Z'}.
$$
The diagram
$$
\Ff \times _{\Ee}Z' \rightarrow P' \rightarrow P\times _ZZ'
$$
induces a diagram
$$
H^i(P\times _ZZ'/Z', \Oo )\rightarrow
H^i(P'/Z', \Oo ) \rightarrow
H^i(\Ff \times _{\Ee}Z'/Z', \Oo ).
$$
In this diagram, the second morphism is an isomorphism, and combining what
was said above the composition is an isomorphism, so the first
morphism is an isomorphism.
As the
cohomology functors (of $P'/Z'$ and $P\times _ZZ'/Z'$) are anchored, it follows
that our morphism
$$
P'\rightarrow P\times _ZZ'
$$
induces an isomorphism on cohomology with any vector sheaf
coefficients. This implies that it  is an equivalence
(note that both sides are relatively $1$-connected
and very presentable). This proves the claim.

The claim implies that ${\bf rep}^{\rm wk}(\Ff \times
_{\Ee}Z/Z)$ is a representing object for
$\Ff \times
_{\Ee}Z/Z$ (not just a weak one). In particular, we obtain the existence of $R$
in the case where the base is a scheme.

For the case of a more general base object $\Ee$, we will give two
constructions.
The first is quicker and requires less verification of details, whereas the
second is more conceptual and uses the notion of cartesian family we have
introduced above. However, in order to obtain a fully rigorous version of the
second construction, one would have to prove a certain number of
compatibility and
uniqueness statements (such as the statement which was left open in Proposition
\ref{correspondence} above).

{\em First Construction:}
As noted above, if $\Ee$ is represented by a scheme $Z$ then the weak
representing
object ${\bf rep}^{\rm wk}(\Ff /Z)$ is in fact a representing object which we
take as $R$.

If $\Ee$ is a disjoint union of schemes, then over each component we obtain
the representing object, and we can take for $R$ the disjoint union of these
component representing objects. Thus $R$ exists if the base is a disjoint union
of schemes.

The existence of $R$ is invariant under equivalences of $\Ee$. Furthermore, $R$
is essentially unique if it exists. More precisely, suppose that $R$ and $R'$
are two $\ngr$-stacks over $\Ee$ with morphisms (relative to $\Ee$)
$$
u:\Ff \rightarrow R, \;\;\;\;\; u':\Ff \rightarrow R'
$$
both satisfying the condition in the first paragraph of the proof. Assume that
$R$ and $R'$ are fibrant over $\Ee$. Let
$$
H:= \underline{Hom}(R/\Ee , R' /\Ee )\times _{
\underline{Hom}(\Ff /\Ee , R' /\Ee )} \Ee
$$
where the first structural morphism is composition with $u$,
and the second structural morphism is the section corresponding to $u'$. Now
$$
H\rightarrow \Ee
$$
is a morphism of $\ngr$-stacks with the property that for any
map from a scheme
$Z\rightarrow \Ee$,  the projection
$$
H\times _{\Ee}Z\rightarrow Z
$$
is an equivalence (this property comes from the universal property of the
representing object $u_Z: \Ff\times _{\Ee}Z \rightarrow R\times _{\Ee}Z$).
It follows that the morphism $H\rightarrow \Ee$ is an equivalence. The
projection
$$
H\rightarrow \underline{Hom}(R/\Ee , R' /\Ee )
$$
therefore corresponds (in an essentially unique way) to a section of
$$
\underline{Hom}(R/\Ee , R' /\Ee )\rightarrow \Ee ,
$$
thus to a morphism $f:R\rightarrow R'$. In view of the construction, we have
$f\circ u \cong u'$. Similarly going back in the other direction we obtain an
essentially unique morphism $R'\rightarrow R$ and again applying the same
argument for $(R,R)$ and for $(R', R')$ we obtain that these morphisms are
inverses. Thus $R$ is essentially unique.

Now suppose that $\Ee$ can be written as a coproduct
$$
\Ee = A\cup ^BC
$$
of two cofibrations $B\hookrightarrow A$ and $B\hookrightarrow C$,
and suppose that we know existence of the representing objects
$$
{\bf rep}(\Ff \times _{\Ee} A/A), \;\;\;
{\bf rep}(\Ff \times _{\Ee} B/B), \;\;\;
{\bf rep}(\Ff \times _{\Ee} C/C).
$$
From the uniqueness of the previous paragraph, we may choose equivalences
$$
{\bf rep}(\Ff \times _{\Ee} A/A)\times _AB \cong
{\bf rep}(\Ff \times _{\Ee} B/B)
$$
and
$$
{\bf rep}(\Ff \times _{\Ee} C/C)\times _CB \cong
{\bf rep}(\Ff \times _{\Ee} B/B).
$$
Replace the coproduct $\Ee$ by a coproduct
$$
\Ee ' := A \cup ^{B\times \{ 0\}} (B\times \overline{I}) \cup ^{B\times \{
1\}}C.
$$
Using the above equivalences, we
may replace ${\bf rep}(\Ff \times _{\Ee} B/B)$ by
an object
$$
R_{B\times \overline{I}}\rightarrow B\times \overline{I}
$$
whose restriction to $B\times \{ 0\}$ is
$$
{\bf rep}(\Ff \times _{\Ee} A/A)\times _AB
$$
and whose restriction to $B\times \{ 1\}$ is
$$
{\bf rep}(\Ff \times _{\Ee} C/C)\times _CB.
$$
We can furthermore assume that this object recieves a map from $\Ff \times
_{\Ee} (B\times \overline{I})$.  Now we can set
$$
R' := {\bf rep}(\Ff \times _{\Ee} A/A)
\cup ^{{\bf rep}(\Ff \times _{\Ee} A/A)\times _AB}
R_{B\times \overline{I}}
\cup ^{{\bf rep}(\Ff \times _{\Ee} C/C)\times _CB}
{\bf rep}(\Ff \times _{\Ee} C/C).
$$
Setting $\Ff ':= \Ff \times _{\Ee} \Ee '$,
we have a map
$$
\Ff ' \rightarrow R'
$$
relative to $\Ee '$. This map has the property required in the first paragraph
of the proof (since any map from a scheme $Z$ into $\Ee '$ has to factor
through one of the components of the coproduct). Note here that $R'$ is not
necessarily fibrant over $\Ee '$ but it is a ``quasifibration'' in the
sense that
the fibrant replacement has the same fibers over the objects of $\Ee '(Z)$.

Thus we have constructed the representing object for $\Ff '\rightarrow \Ee '$,
and since $\Ee '$ is equivalent to $\Ee $ and $\Ff '$ equivalent to $\Ff$,
we obtain (by invariance under equivalence of the base) existence of a
representing object for $\Ff /\Ee$.

The construction of $R$ in general follows from the above cases, because any
$\ngr$-stack $\Ee$ is equivalent to one obtained from disjoint unions of
schemes, by a finite number of coproducts. For example, $\Ee$ is equivalent to
the realization of a simplicial object whose components are disjoint unions of
schemes (this technique was pointed out to me by C. Teleman). The
$n$-truncation  of the realization of this simplicial object may be obtained
from various disjoint unions of the component objects, by a finite number of
coproducts.

{\em Second Construction:}
We now explain a more conceptual way to approach the problem of existence
of $R$,
using the notion of cartesian family. The drawback of this approach is that it
requires several compatibility results for which we don't give very many
details.

Refer to \cite{Simpson}(xi) for the notation $\Upsilon$ and surrounding
techniques.

We have a subcategory (realm)
$$
\underline{RBLE}(Z)\subset \ngr\underline{STACK}(Z)
$$
whose objects are those $\ngr$-stacks over $Z$ for which a representing
object exists. These fit together to form a realm
$\underline{RBLE}$. Fix a fibrant replacement
$\ngr\underline{STACK}'$ and let $\underline{RBLE}'$ denote the
full substack of $\ngr\underline{STACK}'$ with the same objects as
$\underline{RBLE}$.

Next, look at
$$
W(Z)\subset \underline{Hom}(I, \ngr\underline{STACK}'(Z)),
$$
the full sub-$n+1$-category of morphisms which are equivalent to representing
objects relative to $Z$. (Note that a map $I\rightarrow
\ngr\underline{STACK}'(Z)$
is the  same thing as a morphism in $\ngr\underline{STACK}'(Z)$).

We claim that the evaluation on $0\in I$ is a morphism
$$
{\bf ev}_0:W(Z)\rightarrow \ngr\underline{STACK}'(Z)
$$
which is fully faithful, and whose essential image is
$\underline{RBLE}'(Z)$. It is clear that the objects in the image
are exactly those for which the representing object exists;
thus we just have to prove that it is fully faithful.

Suppose $v,w:I\rightarrow \ngr\underline{STACK}'(Z)$ are two objects of $W(Z)$.
Now
$$
W(Z)_{1/}(v,w) = \left( E \mapsto Hom ^{v,w}(\Upsilon (E) \times I,
\ngr\underline{STACK}'(Z)) \right) .
$$
We will show that the evaluation (which is clearly a fibration) is a trivial
fibration, by showing that it satisfies lifting for any cofibration.
This means that we suppose we have a cofibration $E\hookrightarrow E'$,
that we have a map
$$
\alpha ': \Upsilon (E') \rightarrow \ngr\underline{STACK}'(Z),
$$
and that we have a map
$$
A:\Upsilon (E) \times I\rightarrow \ngr\underline{STACK}'(Z)
$$
restricting to $u$ on $0\times I$, restricting to $v$ on $1\times I$,
and restricting to $\alpha '|_{\Upsilon (E)}$ on $\Upsilon (E) \times 0$.
We would like to extend $A$ to a map
$$
A' : \Upsilon (E') \times I\rightarrow \ngr\underline{STACK}'(Z)
$$
agreeing with $\alpha '$ on $\Upsilon (E')\times 0$.

Dividing up the square $\Upsilon (E')\times I$ into two triangles
$\Upsilon ^2(E ', \ast )$ and $\Upsilon ^2( \ast , E')$,
the extension to one of the triangles $\Upsilon ^2(E', \ast )$
is automatic. This gives an extension along the diagonal. Thus we are
reduced to the following problem: we have a morphism
$$
B: \Upsilon ^2(\ast , E)\rightarrow \ngr\underline{STACK}'(Z)
$$
plus an extension of the $02$ edge $B_{02}$ to a morphism
$$
B'_{02}: \Upsilon (E')\rightarrow \ngr\underline{STACK}'(Z);
$$
also the first edge $B_{01}$ is the  morphism $v$ to a representing object
(the target object $B_1$ being in $\underline{AV}(Z)$).
The third object $B_2$ is also in $\underline{AV}(Z)$ (it is the target of the
other morphism $w$). We would like to obtain an extension to a morphism
$$
B': \Upsilon ^2(\ast , E')\rightarrow \ngr\underline{STACK}'(Z).
$$
The functor
$$
E\mapsto Hom ^{v,B_2} (\Upsilon ^2(\ast , E), \ngr\underline{STACK}'(Z))
$$
is representable by an $n$-precat we denote $H$. This functor (and hence $H$)
maps by restriction to the edge $12$ to the functor
$$
E\mapsto Hom ^{B_1,B_2} (\Upsilon ( E), \ngr\underline{STACK}'(Z));
$$
this latter functor is just
$$
\ngr\underline{STACK}'(Z)_{1/}(B_1,B_2).
$$
The resulting map
$$
H\rightarrow \ngr\underline{STACK}'(Z)_{1/}(B_1,B_2)
$$
is a weak equivalence, because the inclusion
$$
\Upsilon (\ast ) \cup ^{\{ 1\} }\Upsilon (E) \hookrightarrow
\Upsilon ^2(\ast , E)
$$
is a trivial cofibration.

On the other hand, restriction to the third edge $02$ gives a map
$$
H\rightarrow \ngr\underline{STACK}'(Z)_{1/}(B_0,B_2).
$$
The diagram
$$
\ngr\underline{STACK}'(Z)_{1/}(B_1,B_2)
\stackrel{\cong}{\leftarrow}H\rightarrow \ngr\underline{STACK}'(Z)_{1/}(B_0,B_2)
$$
is homotopy equivalent to the morphism
$$
\ngr\underline{STACK}(Z)_{1/}(B_1,B_2)
\rightarrow
\ngr\underline{STACK}(Z)_{1/}(B_0,B_2)
$$
of composition by $v$.  (Note here that when we take off the fibrant
replacements,
$$
\ngr\underline{STACK}(Z)_{1/}(B_1,B_2)=\underline{Hom}(B_1,B_2)(Z)
$$
and
$$
\ngr\underline{STACK}(Z)_{1/}(B_1,B_2)=\underline{Hom}(B_0,B_2)(Z),
$$
so composition with $v$ is well-defined).
Now the facts that $v:B_0\rightarrow B_1$ is a morphism to a representing
object in $\underline{AV}(Z)$, and that $B_2$ is also in $\underline{AV}(Z)$,
imply that the morphism of composition with $v$
$$
\underline{Hom}(B_1,B_2)(Z)\rightarrow \underline{Hom}(B_0,B_2)
$$
is  an equivalence. Thus, going back to the situation of our diagrams into the
fibrant replacement $\ngr\underline{STACK}'(Z)$, we get that
$$
H\rightarrow \ngr\underline{STACK}'(Z)_{1/}(B_0,B_2)
$$
is an equivalence. It is a fibration, so it satisfies lifting for
any cofibration. This lifting property exactly gives the existence of the
extension of $B$ to
$$
B': \Upsilon ^2(\ast , E')\rightarrow \ngr\underline{STACK}'(Z)
$$
that we were looking for. This completes the proof that
$$
{\bf ev}_0:W(Z)\rightarrow \ngr\underline{STACK}'(Z)
$$
is fully faithful.

On the other hand, evaluation at $1\in I$ gives a map
$$
W (Z)\rightarrow \underline{AV}'(Z).
$$
Here $\underline{AV}\subset \ngr \underline{STACK}$ is the realm that we are
using, and $\underline{AV}'$ denotes the corresponding full substack of the
fibrant replacement $\ngr \underline{STACK}'$.

Now we let $Z$ vary. The diagram
$$
\underline{RBLE}'(Z)\stackrel{\cong}{\leftarrow}W(Z)
\rightarrow
\underline{AV}'(Z)
$$
is functorial in $Z$. This is because the pullback along a morphism
$Z'\rightarrow
Z$ of any representing object over $Z$, is a representing object over
$Z'$; i.e. pullback maps $W(Z)$ into $W(Z')$.

Thus we get a diagram of $n+1$-stacks
$$
\underline{RBLE}'\stackrel{\cong}{\leftarrow}W
\rightarrow
\underline{AV}'.
$$
Since $\underline{AV}'$ is fibrant we can choose a morphism
$$
\rho : \underline{RBLE}'
\rightarrow
\underline{AV}'
$$
which represents, up to homotopy, the composition of the second arrow with the
inverse of the first arrow in the previous diagram. The morphism $\rho$ is
unique up to coherent homotopy (see \cite{Simpson}(iv) and (xi) Theorem 2.5.1
for uniqueness of the inverse of an equivalence).

Furthermore, note that $\underline{AV}'\subset \underline{RBLE}'$.
Standard closed model category arguments show that we can choose
the morphism $\rho$ to be the identity on $\underline{AV}'$. Also, if we let
$i: \underline{AV}'\hookrightarrow \underline{RBLE}'$ denote the inclusion,
then the technique of \cite{HirschowitzSimpson} \S 13 allows construction of
a natural transformation (i.e. a $2$-morphism between $1$-morphisms in
$(n+1)STACK$)
$$
\eta : 1 \rightarrow i\circ \rho
$$
which, on each object, is the morphism from a representable $n$-stack
to its representing object. This situation can be interpreted from a
shape-theoretical point of view, see the references for shape theory given
above.

We  now turn to the problem at hand, that of constructing the representing
object for $\Ff /\Ee $ under the hypotheses of the present lemma. Via
Proposition \ref{correspondence} we may consider $\Ff /\Ee$ as coming from a
cartesian family $\Ff ^{{\rm fam} /\Ee}\rightarrow {\bf r}(\underline{\Ee })$.
This family in turn corresponds to a morphism
$$
[\Ff /\Ee ]: \Ee \rightarrow \ngr \underline{STACK}'.
$$
The hypothesis of our lemma says that this morphism has image in the
realm $\underline{RBLE}'$ (the claim at the start of the present proof
implies that the weak representing objects in the hypothesis of the
lemma, are in fact representing objects). Composing with the morphism $\rho$
constructed above we get $$ \rho \circ [\Ff /\Ee ] : \Ee \rightarrow
\underline{AV}', $$
and going back in the other direction in Proposition \ref{correspondence},
this gives a morphism
$$
R\rightarrow \Ee ,
$$
which is of type $\underline{AV}$. The natural transformation
$\eta$ gives a morphism of $\ngr$-stacks over $\Ee$, $\Ff \rightarrow R$.

The restriction of this over any object $a\in \Ee (Z)$ (i.e. the pullback by $a:
Z\rightarrow \Ee$) yields the weak representing object ${\bf rep}^{\rm wk}(\Ff
\times _{\Ee}Z/Z)$ (to show this, note that the weak representing family
is one possible choice, but the construction is canonical). This
implies that $\Ff \rightarrow R\rightarrow \Ee$ is a representing object.
\eop

We make a second reduction. This reduction will not be used {\em per se} in
what we say afterward, but it is implicit in \cite{Simpson}(xii) and seems
important enough to mention here.

\begin{lemma}
\label{reduction2}
In the situation where $\Ff \rightarrow
Z$ is over a base scheme, a morphism of $\ngr$-stacks
$\Ff \rightarrow R$ over $Z$, is a weak representing object over $Z$ if and only
if for any vector sheaf $V$ on $Z$, it induces an isomorphism of cohomology
sheaves over $Z$,
$$
H^i(R/Z, V)\stackrel{\cong}{\rightarrow}
H^i(\Ff /Z, V).
$$
\end{lemma}
{\em Proof:} This is immediate from the Postnikov tower for a general $T\in
\underline{AV}(Z)$.
\eop

\oldsubnumero{A criterion for representability}

We will now recall the criterion for representability which was proved in
\cite{Simpson}(xii). It is of course just an adaptation to our
present situation of a technique well-known from long ago in algebraic topology.

The notion of representability used in \cite{Simpson}(xii) was what we are
calling ``weak representability'' here, so we will transfer the criterion to our
present situation using Lemma \ref{reduction1}.

\begin{parag}
\label{defanchored}
Introduce the following terminology. Let $Z$ be a scheme. We say that a
covariant
endofunctor $F$ from the category of vector sheaves on $Z$ to itself, is {\em
anchored} if the natural map
$$
F(U)\rightarrow Hom (Hom (F(\Oo ),\Oo ), U)
$$
is an isomorphism for any coherent sheaf $U$ (recall that the coherent sheaves
are the injective objects in the category of vector sheaves). The above natural
map comes from the trilinear map
$$
F(U)\times Hom (F(\Oo ), \Oo ) \times Hom (U,\Oo ) \rightarrow \Oo
$$
defined by $(f,g,h)\mapsto g( F(h)(f))$.
\end{parag}

Recall the following lemmas concerning this property.

\begin{lemma}
\label{anchored1}
(A) If
$$
0\rightarrow F' \rightarrow F \rightarrow F'' \rightarrow 0
$$
is a short exact sequence of natural transformations between covariant
endofuncturs on the category of vector sheaves over $Y$, then if any two of the
three endofunctors is anchored, so is the third.
\newline
(B) If $F$ is an anchored endofunctor which is also left exact, then $F$ is
representable $F(V)= Hom (W, V)$ for a vector sheaf $W=Hom (F(\Oo ), \Oo )$.
\end{lemma}
{\em Proof:} See \cite{Simpson}(xii).
\eop

\begin{lemma}
\label{anchored2}
Suppose $V$  is a vector sheaf. Then the endofunctor on the category of vector
sheaves defined by
$$
U\mapsto H^i(K(V,m), U)
$$
is anchored.
\end{lemma}
{\em Proof:} See \cite{Simpson}(xii).
\eop

\begin{corollary}
\label{anchored3}
Suppose $T$ is a relatively $1$-connected very presentable $n$-stack over
a scheme $Z$. Then the endofunctor
$$
U\mapsto H^i(T/Z, U)
$$
is anchored.
\end{corollary}
{\em Proof:}
Decompose $T$ into a Postnikov tower where the pieces are of the form $K(V/Z,m)$
for vector sheaves $V$
\eop

Here is the weak-representability criterion.

\begin{theorem}
\label{representable1}
Suppose $Z$ is a scheme and $\Ff \rightarrow Z$ is a morphism of $\ngr$-stacks.
Suppose that  for any vector sheaf $V$ over $Z$, the cohomology  $H^i(\Ff
/Z, V)$
is again a vector sheaf over $Z$.  Suppose furthermore that
$H^0(\Ff /Z, V) = V$ and $H^1(\Ff /Z, V)= 0$ for any vector sheaf $V$.
Finally suppose that the functors $V\mapsto H^i(\Ff /Z,  V)$ are anchored.
Then
there is a weak representing object
$$
\Ff \rightarrow {\bf rep}^{\rm wk}(\Ff /Z)\rightarrow Z
$$
for the $\underline{AV}(Z)$-shape of $\Ff$.
\end{theorem}
{\em Proof:}
See \cite{Simpson}(xii).
\eop

We obtain the following corollary about representability as we have defined it
above.

\begin{corollary}
\label{representable2}
Suppose $\Ff \rightarrow \Ee$ is a morphism of $\ngr$-stacks satisfying the
following properties: that for any morphism from a scheme $Z\rightarrow \Ee$ and
any vector sheaf $V$ on $Z$, the cohomology
$$
H^i(\Ff \times _{\Ee}Z/Z, V)
$$
is again a vector sheaf, with
$$
H^0(\Ff \times _{\Ee}Z/Z, V) = V,\;\;\mbox{and}\;\; H^1(\Ff \times _{\Ee}Z/Z,
V) =0;
$$
and that the endofunctor $V\mapsto H^i(\Ff \times _{\Ee}Z/Z, V)$ on the category
of vector sheaves over $Z$ is anchored.
Then there is a representing object
$$
\Ff \rightarrow {\bf rep}(\Ff /\Ee )\rightarrow \Ee
$$
for ${\bf Shape}_{\underline{AV}}(\Ff /\Ee )$.
\end{corollary}
{\em Proof:}
Apply Lemma \ref{reduction1} and Theorem \ref{representable1}.
\eop

If the base object is the final object $\ast$ then we omit it from the notation;
thus
$$
{\bf rep}(\Ff ):= {\bf rep}(\Ff /\ast ).
$$
The above theorem and corollary give criteria for existence of ${\bf rep}(\Ff )$
(put $Z = \ast$ in the theorem, for example).

\begin{remark}
\label{basechangerep}
Suppose $\Ff \rightarrow \Ee$ is a morphism of $\ngr$-stacks.
If a representing object ${\bf rep}(\Ff /\Ee )$ for the $\underline{AV}$-shape
of $\Ff /\Ee$ exists, then it is compatible with base change in the following
sense: for any morphism of $\ngr$-stacks $\Ee '\rightarrow \Ee$, if we put
$$
\Ff ':= \Ff \times _{\Ee} \Ee '
$$
then
$$
\Ff  ' \rightarrow {\bf rep}(\Ff /\Ee )\times _{\Ee} \Ee '
$$
is a representing object for $\Ff '/\Ee '$, in other words
$$
{\bf rep}(\Ff '/\Ee ')={\bf rep}(\Ff /\Ee )\times _{\Ee} \Ee '
$$
with this formula implying existence of the object on the left.
\end{remark}

\subnumero{The very presentable suspension}

In this section we fix a base scheme $Z$ and work with $\ngr $-stacks over $Z$.
Where convenient, we consider these as $\ngr$-stacks on the site $\Gg /Z$ of
schemes over $Z$, but for clarity we keep notations such as $\pi _i(\Ff /Z)$,
${\bf rep}(\Ff /Z)$ etc.

We can apply the above notions of representability to define a very presentable
``suspension'' of a very presentable $\ngr$-stack. Suppose $\Ff$ is a relatively
$1$-connected very presentable $\ngr$-stack over a base scheme $Z$

Fix $N$. Define the {\em $N$-stack suspension} $Susp _N(\Ff /Z)$ to be the
$N$-fold truncation of the objectwise suspension of the $\ngr$-stack $\Ff$;
it is
an $N^{\rm gr}$-stack, which can be defined (in the world of $N$-stacks on
$\Gg/Z$) by the homotopy-pushout which we denote
$$
Susp_N(\Ff /Z):= \tau _{\leq N}(\ast _Z\cup ^{\Ff } \ast  _Z).
$$
Technically speaking, at least one of the copies of $\ast$ has to be replaced by
something contractible which recieves a cofibration from $\Ff$. One must first
take the pushout (as a presheaf of $N$-precats) and then take the associated
stack. We could write, for example:
$$
Susp_N(\Ff /Z):= \tau _{\leq N}\left( \Ff \times \overline{I}_Z \cup ^{\Ff
\times
\{ 0,1\} _Z }) \{ 0,1\} _Z /Z \right) .
$$

If we think of our $\ngr$-stacks as being simplicial presheaves on $\Gg
/Z$, then
$Susp _N(\Ff /Z)$ is obtained by taking the objectwise suspension and then
$N$-truncating.

\begin{lemma}
If $L$ is any sheaf of abelian groups on $Sch /Z$ then
$$
H^{i+1} (Susp _N(\Ff /Z)/Z, L) = H^i(\Ff /Z, L)
$$
for $1\leq i \leq N-1$. We have
$$
H^1(Susp _N(\Ff /Z)/Z, L) = H^0(\Ff /Z, L) /L.
$$
\end{lemma}
{\em Proof:}
This follows from the Mayer-Vietoris exact sequence.
\eop

\begin{corollary}
If $\Ff$ is a $1$-connected and very presentable $\ngr$-stack on $\Gg /Z$, then
the suspension $Susp _N (\Ff /Z)$ satisfies the conditions of
\ref{representable1}. In particular, there exists a representing $1$-connected
very presentable $N$-stack ${\bf rep} Susp _N(\Ff /Z)$ over $Z$, with
$$
H^{i+1} ({\bf rep}Susp _N(\Ff /Z)/Z, V) = H^i(\Ff /Z, V)
$$
for any vector sheaf $V$, and any $1\leq i \leq N-1$.
\end{corollary}
{\em Proof:}
The functor $V\mapsto H^{i}(Susp _N(\Ff /Z)/Z, V)$ is anchored, and
$$
H^0(Susp
_N(\Ff /Z)/Z, V)=V, \;\;\; H^1(Susp _N(\Ff /Z)/Z,V)=0.
$$
These follow from the
previous lemma. Applying Theorem \ref{representable1} with base $Z=\ast$ gives
existence of  ${\bf rep}Susp _N(\Ff /Z)$. Note that the cohomology of
${\bf rep}Susp _N(\Ff/Z )$ with vector sheaf coefficients is the same as that of
$Susp _N(\Ff /Z)$ so the formula in the present statement is also a
consequence of
the previous lemma.
\eop

We think of ${\bf rep} Susp _N(\Ff /Z)$ as the suspension of $\Ff$ in the
world of
very presentable $N^{\rm gr}$-stacks on $\Gg /Z$.

We can think of the suspension as being pointed (relative to $Z$) by one of the
endpoints so ${\bf rep} Susp _N(\Ff /Z)$ has a natural basepoint section we
denote $0$.

Let $\Omega$ denote the loop-space functor (relative to $Z$). Its input is a
pointed $N$-stack, the output being a pointed $N-1$-stack. In particular we can
apply it to the very presentable suspension above. We obtain the ``loops on the
suspension'' of a very presentable $N$-stack, which is a very presentable
$N-1$-stack denoted  $$
\Omega {\bf rep} Susp _N (\Ff /Z).
$$
We will often want to get back to an $n$-stack here, in which case one should
take $N= n+1$.

We also want to do the iterated loops on the iterated suspension. One can
iterate the suspension operation several times without having to take ${\bf
rep}$ each time (taking it only once at the end). Thus we obtain the {\em
iterated suspension} of a very presentable $n$-stack, which is
$$
{\bf rep} Susp _N^k (\Ff /Z).
$$
(Here we omit one copy of ``$/Z$'' which should go into the notation.)
Again, we would like to then apply the loop space functor, to obtain the
$N-k$-stack
$$
\Omega ^k{\bf rep} Susp _N^k (\Ff /Z).
$$
We usually want to get back to an $n$-stack here, so one should take
$N=n+k$. We make the {\em convention} that if $N$ is supressed from the
notation, it means to take $N=n+k$:
$$
\Omega ^k{\bf rep} Susp ^k (\Ff /Z):=
\Omega ^k{\bf rep} Susp _{n+k}^k (\Ff /Z).
$$
This is a very presentable $n$-stack.

We have a natural morphism of very presentable $\ngr$-stacks
$$
\Ff \rightarrow \Omega ^k{\bf rep} Susp ^k (\Ff /Z).
$$

\begin{lemma}
This construction stabilizes for $k\geq n+2$. In other words, the natural
map
$$
\Omega ^k{\bf rep} Susp ^k (\Ff /Z)\rightarrow
\Omega ^{k'}{\bf rep} Susp ^{k'} (\Ff /Z)
$$
is an equivalence for $k'\geq k \geq n+2$.
\end{lemma}
{\em Proof:}
Apply the usual stabilization theorems from topology.
\eop

We denote by $\Omega ^{\infty}{\bf rep} Susp ^{\infty} (\Ff /Z)$
the stabilized version (recall that here we stay within the world of
$n$-stacks).
To be precise, fix some $k\geq n+2$ and define
$$
\Omega ^{\infty}{\bf rep} Susp ^{\infty} (\Ff /Z):=
\Omega ^{k}{\bf rep} Susp ^{k} (\Ff /Z).
$$
The usefulness of the morphism
$$
\Ff \rightarrow \Omega ^{\infty}{\bf rep} Susp ^{\infty} (\Ff /Z)
$$
comes from the fact that it provides a canonical
alternative to the Postnikov truncation. This alternative truncation will
preserve the geometricity property which we shall consider below. For now we
just note:

\begin{proposition}
\label{nilpotence}
Suppose $\Ff$ is a $1$-connected very presentable $\ngr$-stack over $Z$.
The fiber of the morphism of $\ngr$-stacks
$$
\Ff \rightarrow \Omega ^{\infty}{\bf rep} Susp ^{\infty} (\Ff /Z)
$$
over the basepoint section, is an $\ngr$-stack which we denote $\Phi (\Ff /Z)$.
If $\Ff$ is $i$-connected relative to $Z$ then $\Phi  (\Ff /Z)$ is
$i+1$-connected relative to $Z$.
\end{proposition}
{\em Proof:}
In what follows, we use the {\em homology sheaves} of an $\ngr$-stack $\Ff$.
These are defined as the sheaves $H_j(\Ff /Z)$ associated to the presheaves
$Y\mapsto H_j(\Ff (Y))$. Note that the $\Ff (Y)$ are rational spaces (i.e. the
homotopy groups are rational vector spaces) so the homology sheaves are sheaves
of rational vector spaces.

Suppose $\Ff$ is $i$-connected. Choose a large enough value of $k$ so that
$$
\Omega ^{\infty}{\bf rep} Susp ^{\infty} (\Ff /Z)=
\Omega ^{k}{\bf rep} Susp ^{k} (\Ff /Z).
$$
Note that $H_i(\Ff /Z) = \pi _i(\Ff /Z)$ is a vector sheaf, and
$$
H^{i+k} ({\bf rep}Susp ^k \Ff /Z, V) = H^i (\Ff /Z, V)
$$
for any vector sheaf. Note that $\pi _i(\Ff /Z)$ is identified by the functor
$V\mapsto H^i(\Ff /Z, V)$. Using the fact that $\pi _{i+k}({\bf rep}Susp ^k \Ff
/Z)$ is a vector sheaf, it is identified by the same functor $V\mapsto
H^{i+k}({\bf rep}Susp ^k \Ff  /Z, V)$, so we get
$$
\pi _{i+k} ({\bf rep}Susp ^k \Ff /Z) = \pi _i(\Ff /Z).
$$
Thus
$$
\pi _i(\Omega ^k{\bf rep}Susp ^k \Ff /Z) = \pi _i(\Ff /Z).
$$
To finish the proof we need to show that the morphism
$$
\pi _{i+1} (\Ff /Z) \rightarrow \pi _{i+1}(\Omega ^k{\bf rep}Susp ^k \Ff /Z)
$$
is a surjection.

If $i\geq 3$ then for this statement we are in the stable range and the morphism
in question is an isomorphism. Thus we may suppose $i=2$. Similarly we may take
$k=1$.

The sheaf $\pi _2(\Ff /Z)$ contributes to rational
homotopy in degrees $2$ and $\geq 4$. We obtain that $H_2(\Ff /Z)$ and $H_3(\Ff
/Z)$ are vector sheaves. Thus $H_3(Susp ^1\Ff /Z)$ and $H_4(Susp ^1\Ff /Z)$ are
vector sheaves, and in view of the stable range these are the same as
$\pi _3(Susp ^1\Ff /Z)$ and $\pi _4(Susp ^1\Ff /Z)$ respectively. Thus
$$
\tau _{\leq 4} ({\bf rep}Susp ^1( \Ff /Z)/Z) =
\tau _{\leq 4} (Susp ^1 (\Ff /Z) /Z),
$$
as the latter is already very presentable. Finally we get back to
$$
\pi _3(\Omega ^1{\bf rep}Susp ^1 \Ff /Z) = H_4(Susp ^1\Ff /Z)=
H_3(\Ff /Z).
$$
Thus it suffices to show that the morphism
$$
\pi _3(\Ff /Z)\rightarrow H_3(\Ff /Z)
$$
is surjective. But this can be verified on the level of each object $Y \in \Gg
/Z$: the morphism
$$
\pi _3(\Ff (Y))\rightarrow H_3(\Ff (Y))
$$
is surjective, since $\Ff (Y)$ is a $1$-connected rational space.

This completes the proof (using the long exact sequence for the homotopy groups
of the fiber $\Phi (\Ff )$).
\eop

\numero{Geometric $n$-stacks}
\label{geometricpage}

The realms defined in the previous section,
of presentable or very presentable $n$-stacks, take us fairly far from the realm
of ``geometry'' in that the stacks in question can't be pictured as geometric
objects in a very nice way. For example, a vector sheaf $V$ over a base $X$ with
the property that the dimensions of the vector space fibers $V(x)$ for every
closed point $x\in X$ are all the same, need not be locally free. Thus it is
natural to look for another condition. It should be noted that we were forced to
go to some types of objects like presentable group sheaves, if we wanted to have
compatibility with truncation and with finite limits. As we have seen, the
compatibility with finite limits is quite important, so we could throw out
compatibility with truncation instead. This means that we will be looking at
objects characterized by a condition which mixes up the various homotopy group
sheaves $\pi _i$ for different $i$.

Throughout this section we work only with $n$-stacks of groupoids but for
clarity of notations we write ``$n$-stack'' rather than ``$\ngr$-stack''.

One prototypical example of this type of behavior is the notion of {\em perfect
complex}: this is a complex of sheaves of $\Oo$-modules on $X$ which is locally
quasiisomorphic to a complex of vector bundles. Note that the cohomology sheaves
of a perfect complex (which correspond to the homotopy group sheaves of a stack)
are not themselves vector bundles, and similarly if we truncate (in the
canonical way preserving cohomology) a perfect complex, the result is no longer
necessarily perfect.

For $1$-stacks, we have Artin's notion of {\em algebraic
stack} \cite{Artin}: this is a $1$-stack (of groupoids) $F$ such that there
exists
a smooth surjective morphism $X\rightarrow F$ from a scheme $X$, and such that the
fiber product $R:= X\times _FX$ is itself an algebraic space. Here the pair
$(X,R)$ forms a category-object in the category of algebraic spaces, and $F$ is
the $1$-stack which it represents. Thus one can think of algebraic stacks as
being the stacks represented by groupoids in the category of algebraic spaces
$(X,R)$ with the property that the two projections $R\rightarrow X$ are smooth.

One obtains an analogous notion by replacing ``smooth'' in the above definition,
by any other  condition. Amazingly enough, this definition for ``flat''
groupoids gives rise to the same class of objects as for ``smooth'' groupoids
(\cite{Artin}). Unfortunately, this nice result doesn't seem to persist for
$n$-stacks.

We can formalize this situation as follows, making it applicable to $n$-stacks.
Fix a full saturated
$n+1$-subcategory $M\subset
nSTACK$ (recall that $nSTACK$ is the $n+1$-category of global sections of
$n\underline{STACK}$, i.e. the $n+1$-category of $n$-stacks on our site $\Gg$).
Fix also an $n+1$-subcategory $L\subset M$ with the property that every
object of
$M$ is in $L$, that for any $x,y\in M_0$, $L_{1/}(x,y)\subset M_{1/}(x,y)$ is a
full saturated sub-$n$-category, and that all equivalences between $x$ and
$y$ are
in $L_{1/}(x,y)$. A {\em groupoid in $M$} is a morphism of $n+1$-categories $$
G:\Delta \rightarrow M
$$
such that $G(0)$ is a $0$-stack; such that the {\em Segal maps}
$$
G(m)\rightarrow G(1)\times _{G(0)} \ldots \times _{G(0)}G(1)
$$
are equivalences of $n$-stacks; and such that the resulting $n+1$-stack
${\bf i} G$ obtained by considering each simplicial $n$-category as an
$n+1$-category, is an $n+1$-stack of groupoids.

We say that a groupoid in $M$ is {\em of $L$-type} if the morphisms
$G(1)\rightarrow G(0)$ are in $L$.

We say that $G$ is a {\em groupoid of $n-1$-stacks in $M$} if each $G(m)$ is
$n-1$-truncated. Then ${\bf i} G$ is an $n$-stack.

Finally, we say that {\em $M$ is closed under integration of groupoids of
$n-1$-stacks of $L$-type} if for any groupoid of $n-1$-stacks $G$ of $L$-type,
the integral ${\bf i} G$ is in $M$.

We say that the pair $(M, L)$ is {\em closed under
integration} if $M$ is closed  under integration of groupoids of
$n-1$-stacks of $L$-type, and if a morphism
$$
{\bf i} G \rightarrow {\bf i} G'
$$
is in $L$ whenever $G\rightarrow G'$ is a morphism of
groupoids of $n-1$-stacks of $L$-type, such that the morphism
$$
G(0)\times _{{\bf i} G'} G'(0) \rightarrow G'(0)
$$
is in $L$.

Suppose now that $\underline{M}\subset n\underline{STACK}$ is a realm (i.e.
saturated full substack), and $\underline{L} \subset \underline{M}$ is a
substack
which over each object has the properties  required in the first paragraph
concerning $L\subset M$. If, for each $Z\in \Gg$, the subcategory
$\underline{M}(Z) \subset n\underline{STACK}(Z)$ is closed under integration of
groupoids of $n-1$-stacks of $\underline{L}(Z)$-type, then we say that {\em
$\underline{M}$ is closed  under integration of groupoids of $n-1$-stacks of
$\underline{L}$-type}. If for each $Z\in \Gg$ the pair $(\underline{M}(Z),
\underline{L}(Z))$ is closed under integration, then we say that {\em the pair
$(\underline{M}, \underline{L})$ is closed under integration}.

It is clear that an arbitrary intersection of pairs which are closed under
integration, is again closed under integration.
Therefore, if  we fix a preliminary pair
$\underline{K}\subset \underline{N} \subset n\underline{STACK}$ then we can take
the minimal pair $(\underline{M}, \underline{L})$  containing $(\underline{N},
\underline{K})$ and being closed under integration.

\subnumero{Geometric $n$-stacks---the concrete definition}

Start with $\underline{N}$ equal to the stack of schemes
and $\underline{K}$ equal to the substack of schemes with only smooth
morphisms.
One obtains for $n=1$ the notion of Artin algebraic stack (i.e. $\underline{M}$
is the $2$-stack of algebraic stacks) and for general $n$ one obtains the notion
of geometric $n$-stack of \cite{Simpson}(ix).

We recall here the more concrete version of this notion of geometric $n$-stack;
work on the site of noetherian schemes with the etale topology.

The definition is inductive on $n$. For $n=0$, we say that a {\em geometric
$0$-stack} is a sheaf of sets which is represented by an Artin algebraic space.
Recall that an algebraic space $A$ has a chart which is a surjective smooth
morphism from  a scheme $X\rightarrow A$. Recall that a morphism
$A\rightarrow B$
of algebraic spaces is smooth if, for smooth charts $X\rightarrow A$ and
$Y\rightarrow B$, the morphism of schemes $X\times _BY\rightarrow Y$ is smooth.

Now suppose $n\geq 1$ and suppose that we have defined the notion of geometric
$n-1$-stack as well as the notion of smooth morphism between geometric
$n-1$-stacks. An $n$-stack $A$ is {\em geometric} if there is a surjective
morphism from a scheme $X\rightarrow A$, which has the following property:
\newline
---if $Y$ is any scheme mapping to $A$ then $X\times _AY$ is a geometric
$n-1$-stack, and the map $X\times _AY\rightarrow Y$ is smooth.

We say that $X\rightarrow A$ is a ``smooth surjection'' or a {\em chart}.
If $A\rightarrow B$ is a morphism of geometric $n$-stacks then we say that $f$
is {\em smooth} if for any (or all) charts $X\rightarrow A$ and $Y\rightarrow
B$, the morphism of geometric $n-1$-stacks $X\times _BY\rightarrow Y$ is smooth.
Note that for verification of this latter condition at level $n-1$, one can take
$Y$ as a chart for itself and take some chart $W\rightarrow X\times _BY$, so the
condition becomes just that $W\rightarrow Y$ be a smooth map of schemes.

This completes the inductive definition. See \cite{Simpson}(ix)
for a number of further remarks, like the independence of these properties
under different choices of charts, etc.

\begin{parag}
\label{induction}
It is clear that if $A$ is a geometric $n$-stack and if $N\geq n$, then $A$
considered as an $N$-stack (which we sometimes write $Ind _n^N(A)$) is again
geometric.
\end{parag}

Generally speaking we shall work over a base field $k$ and only consider schemes
of finite type over $k$; and in this case we will look at geometric $n$-stacks
which are of finite type (i.e. where all the charts in question are schemes of
finite type) without further mention. When the first chart $X\rightarrow A$
is not
necessarily of finite type, we say that $X$ is {\em locally geometric}, but even
in this case we assume that the further charts which enter inductively into
the definition are of finite type.

\oldsubnumero{Compatibility with finite limits and extensions}

Let $\underline{GE}$ (resp. $\underline{VG}$) denote the realm defined by
setting $\underline{GE}(X)$ (resp. $\underline{VG}(X)$) equal to the class of
$\ngr$-stacks $T$ on $\Gg /X$ which are geometric (resp. geometric and very
presentable) when considered as $\ngr$-stacks on $\Gg$.

\begin{lemma}
\label{vgcompatible}
The realms
$\underline{GE}$ and $\underline{VG}$ are compatible with finite limits and with
extensions.
\end{lemma}
{\em Proof:}
We refer to \cite{Simpson}(ix) for these statements. In view of the
definition of geometric morphism in \cite{Simpson}(ix), compatibility with
extensions is the statement of Corollary 2.6 of \cite{Simpson}(ix).
Compatibility with finite limits is the statement of Proposition 2.1 of
\cite{Simpson}(ix). The previous remarks are for $\underline{GE}$. We obtain
the same for
$$
\underline{VG} = \underline{GE} \cap \underline{VP}
$$
in view of Lemma \ref{vpcompatible}.
\eop

\begin{remark}
Geometric $n$-stacks are presentable; see \cite{Simpson}(ix) Proposition 5.1.
Thus we can write
$$
\underline{GE} \subset \underline{PE}.
$$
So, in the definition of $\underline{VG}$, the new conditions being put on a
geometric $n$-stack  are the conditions of ``very-presentability'', i.e.
affineness for $\pi _1$ and the condition that the $\pi _i$ be vector sheaves
for $i\geq 2$.
\end{remark}

\subnumero{Doing geometry with geometric $n$-stacks}

One can do a lot of standard things from algebraic geometry, for geometric
$n$-stacks. This was pointed out in Laumon-Moret Bailly \cite{LaumonMoretBailly}
for algebraic (i.e. geometric) $1$-stacks, and our remarks here are exactly the
same.

One can define the properties of a morphism
$f:A\rightarrow B$ in the following way. Say $P$ is a property (for example,
from the list of properties in \cite{LaumonMoretBailly}).
Then we say that $f$ has property $P$ if, for the smooth surjections
$$
Y\rightarrow B\;\;\; \mbox{and} \;\;\; X\rightarrow A\times _BY,
$$
the morphism $X\rightarrow Y$ has property $P$.

Note that our previous definition of smooth morphism fits into this framework.

We obtain in particular a notion of {\em flat} morphism between two geometric
$n$-stacks.

Similarly we obtain a notion of {\em birational} morphism between two geometric
$n$-stacks.

An {\em open set} $U\subset A$ is a full substack such that for any map from a
scheme $Y\rightarrow A$, the pullback $U\times _AY$ is represented by an open
subset of $Y$.

If $f:B\rightarrow A$ is a smooth morphism between geometric $n$-stacks,
then there is an open set $U\subset A$ which is the ``image'' of $f$, in the
sense that for any morphism $Y\rightarrow A$ from a scheme, and for the smooth
surjection  $$
X\rightarrow Y\times _AB,
$$
we have that the image of $X\rightarrow Y$ in $Y$ (which is open because this is
a smooth map of schemes) is equal to $Y\times _AU$. We denote $U$ by $im(f)$.

An open set $U\subset A$ is {\em dense} if for any smooth map from a scheme
$Y\rightarrow A$, the pullback $U\times _AY$ is dense in $Y$.

A geometric $n$-stack $A$ is {\em irreducible} if there is a smooth morphism
from an irreducible scheme $Y\rightarrow A$, with dense image.

A {\em closed substack} is a morphism $Z\rightarrow A$ such that for any map
from a scheme $Y\rightarrow A$, the pullback $Z\times _AY$ is represented by a
closed subscheme of $Y$.

A closed subscheme is a {\em Weil divisor} if for any smooth map from a scheme
$Y\rightarrow A$, the pullback $Z\times _AY$ is a Weil divisor on $Y$.

A morphism $f: B\rightarrow A$ is {\em representable} if for any morphism from a
scheme $X\rightarrow A$, the fiber product $B\times _AX$ is an algebraic space;
and {\em proper representable} if in the same circumstances, the morphism
$B\times _A\rightarrow X$ is a proper morphism of algebraic spaces.

\oldsubnumero{Resolution of singularities}

Canonical resolution of singularities implies resolution of singularities for
geometric $n$-stacks. As usual we are working in characteristic zero, so
canonical resolution is available by Bierstone-Milman \cite{BierstoneMilman}. In
characteristic $p>0$ this should be possible in the future by Spivakovsky's
program.

Let $\Ll _n\subset nSTACK (Sch /k)$ be the sub-$n+1$-category whose objects are
the reduced geometric $n$-stacks, and whose morphisms are those generated by the
smooth morphisms and the morphisms which are locally sections of smooth
morphisms. We include all $i$-morphisms for $i\geq 2$.

\begin{theorem}
There is an $n+1$-endofunctor $Res$ on $n+1$-category $\Ll _n$.
There is a  morphism of $n$-stacks (usually not in $\Ll _n$)
$Res(A)\rightarrow A$
which is natural transformation with respect to morphisms of $\Ll_n$ in the
variable $A$. These satisfy:
\newline
---$Res(A)$ is smooth;
\newline
---the map $Res(A)\rightarrow A$ is a representable proper birational morphism,
isomorphism on the open set where $A$ is smooth;
\newline
---if $A$ is smooth then $Res(A)\rightarrow A$ is an isomorphism;
\newline
---if $Y$ is a smooth scheme then $Res(A\times Y)=Res(A)\times Y$.
\end{theorem}
{\em Sketch of Proof:}
For schemes this result is due to Bierstone and Milman \cite{BierstoneMilman}.
For
algebraic spaces it follows by the same argument as we are about to do. Thus we
may assume that it is known for $n=0$. We prove it by induction on $n$; so
assume
that it is known for $n-1$.

Let $A$ be a reduced geometric $n$-stack, and choose a smooth surjective
morphism from a scheme $X\rightarrow A$. Form the simplicial object in the
$n$-category of geometric $n-1$-stacks $\Ll _{n-1}$,
$$
G_{\cdot} := X\times _	A \ldots \times _A X, \;\;\;\; {\bf i} G_{\cdot} \cong A.
$$
This may be chosen as a strict functor from $\Delta ^o$ to the model category of
$n-1$-prestacks. Note that the face
morphisms are smooth, and the degeneracy morphisms are sections of smooth
morphisms. Let $F_k:= Res (G_k)$ (which exists by the inductive hypothesis).
This again is a simplicial object in $\Ll _{n-1}$ and we have a morphism of
simplicial $n-1$-stacks $F_{\cdot} \rightarrow G_{\cdot}$. Set
$$
Res(A):= {\bf i} F_{\cdot}
$$
which has a map to $A\cong {\bf i} G_{\cdot}$. The $F_k$ are smooth, so ${\bf i}
F_{\cdot}$ is a smooth geometric $n$-stack. Note that $G_0$ and hence $F_0$ are
algebraic spaces; they are equipped with smooth maps
$$
G_0\rightarrow A, \;\;\;\; F_0 \rightarrow Res(A),
$$
and the map $F_0\rightarrow G_0$
is birational; thus $Res(A)\rightarrow A$ is birational.
Note also that $F_0 = G_0 \times _A Res(A)$; the fact that $F_0\rightarrow G_0$
is proper and representable implies the same for $Res(A)\rightarrow A$.

If $A$ is smooth then $G_k$ are smooth so $F_k=G_k$ and $Res(A)=A$.

For the geometric $n$-stack $A\times Y$ we can choose a smooth
surjection of the form $X\times Y\rightarrow A\times Y$, and then
the resulting simplicial object becomes $G_{\cdot} \times Y$, its resolution is
$F_{\cdot}\times Y$ and we obtain
$$
Res(A\times Y)= {\bf i} (F_{\cdot} \times Y)= ({\bf i} F_{\cdot} )\times Y
=Res(A)\times Y.
$$

We should show that $Res(A)$ is independent of the choice of smooth surjection
$X\rightarrow A$ in a homotopy-coherent way (i.e. with all higher homotopies),
which then gives the functoriality of $Res$. We don't do this here, which is why
it is only a ``sketch of proof''.
\eop

One can also resolve the singularities of a Weil divisor, turning it into a
normal crossings divisor; or of the complement of an open set $U\subset A$. In
the latter case, if $U$ is smooth then there is a surjective proper birational
morphism $p:B\rightarrow A$ which is an isomorphism over $U$, and such that
$B$ is smooth and the complement of $p^{-1}(U)$ in $B$ is a divisor with normal
crossings.

\oldsubnumero{Deformation theory}

One can of course envision a deformation theory for geometric $n$-stacks
analogous to that of schemes. I don't have any concrete results in this
direction but we can give the obvious definitions and pose some questions.

Suppose $A_0$ is a geometric $n$-stack, and suppose $D$ is the spectrum of an
artinian local ring (everything of finite type over a base field $k$, let's
say). A {\em deformation of $A$ parametrized by $D$} is a triple $(A_D, f,
\xi )$
consisting of a flat morphism of geometric $n$-stacks
$$
f: A_D \rightarrow D,
$$
together with an equivalence $\xi : A_0\cong A_D \times _D Spec (k)$.

We obtain an $n+1$-category of deformations of $A$ indexed by $D$, which fit
together into a functor
$$
Def _{A_0}: Art /k \rightarrow (n+1)CAT.
$$
To be precise, this is defined as
$$
Def_{A_0}(D):= nSTACK ^{\rm geom, flat}(Sch /D) \times
_{nSTACK ^{\rm geom, flat}(Sch /Spec (k))} (\ast )
$$
where the second morphism in the fiber product corresponds to $A_0$.

Taking $D= Spec (k[\epsilon ]/(\epsilon ^2 ))$ we obtain an $n+1$-category
which we denote by $TDef _{A_0}$.

\begin{lemma}
The $n+1$-categories $Def_{A_0}(D)$ are $n+1$-groupoids.
\end{lemma}
\eop

We formulate as a conjecture some of the ``standard-type'' statements about
deformation theory that one would like to prove. Recall that we have defined the
{\em tangent spectrum} $TA$ of a geometric $n$-stack, see \cite{Simpson}(ix).

Let $\Omega ^{-1}TA$ denote the $1$-fold delooping of $TA$.

\begin{conjecture}

\noindent
(1)\,
The $n+1$-groupoid $TDef _{A}$ has a structure of ``spectrum'' i.e. it is the
$N$-fold looping of  an $N+n+1$-groupoid.
\newline
(2) \, If $A$ is a smooth geometric $n$-stack then there is a natural
equivalence
(of $n+1$-categories)
$$
TDef _A \cong \Gamma (A, \Omega ^{-1} TA),
$$
which in the case where $A$ is a smooth scheme reduces to the well-known
isomorphism $TDef_A \cong H^1(A, TA)$.
\newline
(3) \, In general there is a ``cotangent complex'' $Cot(A)$ which is a
spectrum over $A$, and we have
$$
TDef _A \cong \Gamma (A, \Omega ^{-1} Cot(A)).
$$
Again, this should generalize the classical cotangent complex in the case where
$A$ is a scheme, and should coincide with the cotangent complex defined by
Laumon-Moret Bailly in the case $n=1$.
\end{conjecture}

One would also like to have an ``obstruction theory'', construct versal
deformation stacks, and so on.

\subnumero{A criterion}

Suppose $Z$ is a scheme of finite type over a field $k$ of characteristic zero.
We would like to study geometric $n$-stacks $T\rightarrow Z$ which are
relatively $1$-connected, and very presentable (i.e. the $\pi _i(T/Z)$ are
vector sheaves on $Z$).
In this subsection we will prove a result which makes a connection between the
notion of geometric $n$-stack and the notion of perfect complex.
Before getting there, we need some preliminary definitions.

\oldsubnumero{Formal smoothness}

We say that a morphism of $n$-stacks $f:A\rightarrow B$ is {\em formally smooth}
if it satisfies the following strong form of the infinitesimal lifting property:
for any Artinian local ring $R$ with ideal $I$ such that ${\bf m}_RI=0$, and any
diagram
$$
\begin{array}{ccc}
Spec (R/I)& \rightarrow & 	A\\
\downarrow && \downarrow \\
Spec (R) & \rightarrow & B
\end{array}
$$
by which we mean a square of morphisms plus homotopy of commutativity
(i.e. a morphism $I\times I \rightarrow nSTACK (\Gg )'$ to the fibrant
replacement of $nSTACK (\Gg )$), there exists a lifting
$Spec (R)\rightarrow A$ together with homotopies of commutativity for the
resulting two triangles plus a $2$-homotopy between the composition of the two
triangles, and the original homotopy of the square (again all of these can be
represented by saying that there is a certain extension of the above morphism
into $nSTACK (\Gg )'$).

One could say that $f$ is {\em weakly formally smooth} if the same condition as
above holds, but without requiring the existence of the $2$-homotopy. It isn't
clear whether this notion is useful, though.

Any base-change of a formally smooth morphism is again formally smooth, and
conversely if $Y\rightarrow B$ is a  formally smooth surjective morphism
such that $A\times _BY\rightarrow Y$ is formally smooth, then $A\rightarrow B$
is formally smooth.

It follows immediately that
if $f:A\rightarrow B$ is a morphism between geometric $n$-stacks (of finite type
over a field) then $f$ is smooth if and only if it is formally smooth.

Furthermore, we can rewrite the definition of ``geometric $n$-stack'' as
follows: an $n$-stack $A$ is geometric if and only if there is a formally
smooth surjective morphism from a scheme $X\rightarrow A$ such that
$X\times _AX$
is a geometric $n-1$-stack.

\oldsubnumero{Almost-geometric $n$-stacks}

For the purposes of one of the proofs below, the following weakening of the
notion of geometricity is useful. We work only with schemes of finite type over
a base (say, a base field, although another type of reasonable base scheme
would probably also work). We say that a $0$-stack $A$ is {\em almost geometric}
if there is a formally smooth surjective morphism from a scheme
$X\rightarrow A$.
(Note that we don't require $X\times _AX$ to be a scheme, so this is more
general
than an algebraic space.) For $n\geq 1$ we say that an $\ngr$-stack $A$ is  {\em
$n$-almost-geometric} if there is a formally smooth surjective
morphism  from a scheme $X\rightarrow A$ such that $X\times _AX$ is an
$n-1$-almost-geometric $n-1$-stack.

\begin{lemma}
\label{aggeo}
Suppose $A$ is an $n$-stack which, when considered as an $n+1$-stack, is
$n+1$-almost geometric. Then $A$ is a geometric $n$-stack.
\end{lemma}
{\em Proof:}
We prove this by induction on $n$. Suppose $A$ is a $0$-stack which is
$1$-geometric. Then there is a formally smooth surjection $X\rightarrow A$
from a scheme. Furthermore, $R:= X\times _AX$ is a $0$-geometric $0$-stack,
so there is a formally smooth surjection $Y\rightarrow R$. Note that
$$
Y\times _RY = Y\times _{(X\times X)}Y
$$
is a scheme, and the projection morphisms $Y\times _RY\rightarrow Y$
are formally smooth, hence smooth. Thus $R$ is an algebraic space (quotient of a
scheme by a smooth equivalence relation). Now the projection morphisms
$R\rightarrow X$ are formally smooth hence smooth. Furthermore, $R$ is an
equivalence relation because $A$ is a $0$-stack; thus $A$ is an algebraic space
i.e. a geometric $0$-stack.

Now we treat the inductive step. Suppose the statement is known for $n-1$.
Suppose $A$ is an $n$-stack which is $n+1$-almost-geometric. Choose a formally
smooth surjection $X\rightarrow A$. Then $X\times _AX$ is an $n-1$-stack which
is $n$-almost-goemetric, so by induction $X\times _AX$ is a geometric
$n-1$-stack. The projection morphisms from here to $X$ are formally smooth which
implies smooth, and this shows that $A$ is a geometric $n$-stack.
\eop

On the other hand, the notion of almost-geometric is compatible with truncation:

\begin{lemma}
\label{agtruncation}
Suppose $A$ is an $n$-almost geometric $\ngr$-stack, and suppose $0\leq k
\leq n$.
Then $\tau _{\leq k}(A)$ is a $k$-almost geometric $k^{\rm gr}$-stack.
\end{lemma}
{\em Proof:}
Prove this by induction on $n$.
The formally smooth surjection $X\rightarrow A$ provides a formally smooth
surjection $X\rightarrow \tau _{\leq k}(A)$, and we have
$$
X\times _{\tau _{\leq k}(A)}X = \tau _{\leq k-1} ( X\times _AX).
$$
thus if $X\times _AX$ is $n-1$-almost geometric, then by the inductive
statement, we get that $\tau _{\leq k-1} ( X\times _AX)$ is
$k-1$-almost-geometric. Thus, directly from the definition,
$\tau _{\leq k}(A)$ is $k$-almost-geometric.
\eop

{\em Caution:} the $k$-almost-geometric $k$-stack $\tau _{\leq k}(A)$ will in
general {\em not} be an $n$-almost-geometric $n$-stack, if $k<n$. This can
provide
an example showing that
\ref{induction} no longer holds for the notion of
``almost-geometric''.

The property of being almost-geometric also is compatible with extension.

\begin{lemma}
\label{agextension}
The notion of ``almost-geometric $n$-stack'' satisfies the following  property
of closure under extension. Suppose that $f:A\rightarrow B$ is a morphism of
$n+1$-stacks such that $B$ is an almost geometric $n+1$-stack. Suppose that for
every scheme $Y$ with map $Y\rightarrow B$, the fiber product $A\times _BY$ is
an almost-geometric $n$-stack. Then $\tau _{\leq n}(A)$ is an
almost-geometric $n$-stack.
\end{lemma}
{\em Proof:}
Use the same type of arguments as in \cite{Simpson}(ix) and (vii); this
is somewhat long so we don't give the details here.
\eop

\oldsubnumero{A criterion for geometricity}

We work in characteristic zero.

\begin{definition}
\label{residuallyperfect}
A complex of sheaves $C^{\cdot}$ on the big etale site of a scheme $X$ is {\em
residually perfect} if for any $N$ and locally on $X$, there is a perfect
complex
$L^{\cdot}$ (i.e. a complex of finite rank vector bundles) and a
quasiisomorphism in the derived category
$$
\tau ^{[-N,N]}(C^{\cdot}) \cong
\tau ^{[-N,N]}(L^{\cdot})
$$
where $\tau ^{[-N,N]}$ means the intelligent truncation (preserving cohomology)
to a complex supported in $[-N,N]$.
\end{definition}

Note of course that a perfect complex is residually perfect.

\begin{theorem}
\label{criterion}
Suppose $p:T\rightarrow Z$ is a relatively $1$-connected $n$-stack of groupoids,
which is very presentable (i.e. the $\pi _i(T/Z)$ are vector sheaves on $Z$).
Then $T$ is geometric if and only if the relative cohomology complex ${\bf R}
p_{\ast} \Oo$ is a residually perfect complex on $Z$.
\end{theorem}

The proof of this theorem requires several steps which take up the rest of this
subsection.

We suppose in what follows that the base $X$ is affine, writing $X=Spec (A)$.
If $M$ is an $A$-module, then the {\em quasicoherent sheaf $\tilde{M}$} is
defined by
$$
\tilde{M}(Spec (B)):= M\otimes _AB.
$$
This extends to any scheme $Z\rightarrow X=Spec (A)$ by the sheaf condition.
Say that $\tilde{M}$ is {\em flat} if $M$ is a flat $A$-module (note that
$M=\tilde{M}(X)$ so this condition is well defined.

Suppose $V$ is a vector scheme over $X$. We can write $V=Spec (E)$ with
$E=\bigoplus _i E_i$ decomposed into homogeneous components. Thus, $E_1$ is the
dual coherent sheaf to $V$. The set of morphisms of sheaves of groups (or
equivalently, sheaves of $\Oo$-modules) from $V$ to $\tilde{M}$ is equal to
$M\otimes _A E_1$.  This is the subset of sections of homogeneity $1$, in
the space of
sections $H^0(V, \tilde{M})=M\otimes _AE$.
(For the proof of this, use the arguments of \cite{Simpson}(vii).)

Suppose
$$
\tilde{M} \rightarrow \tilde{N} \rightarrow \Ff \rightarrow 0
$$
is an exact sequence with the first two terms being quasicoherent. If $V$ is a
vector scheme, then any morphism $V\rightarrow \Ff$ lifts to a morphism
$V\rightarrow \tilde{N}$. To see this, we can work in the Zariski topology over
$V$ itself. The fact that $V\rightarrow X$ is
an affine map implies that the given section in $H^0(V, \Ff )$ lifts to a
section in $H^0(V, \tilde{N})$. Decomposing this latter section according to its
homogeneous components and taking the component of degree $1$, gives the desired
lifting.

We remark that if $p:Y\rightarrow X$ is a morphism of schemes and $\Ff$ is a
coherent sheaf on $Y$, flat over $X$, then ${\bf R}p_{\ast}(\Ff )$ is
quasiisomorphic to a complex, concentrated in positive degrees, whose components
are quasicoherent flat sheaves on $X$. For this it suffices to calculate the
cohomology by a \v{C}ech complex with respect to an affine open covering.

An argument analogous to that given in the next lemma, implies that the same is
true for a morphism $p:Y\rightarrow X$ with $Y$ an algebraic space.

\begin{lemma}
\label{mumford1}
Suppose $X$ is a scheme and $p:T\rightarrow X$ is a morphism of geometric
$n$-stacks. Suppose that $p$ is flat. Then ${\bf R} p_{\ast}(\Oo )$ (which is a
complex of sheaves on the big etale site over $X$) is quasi-isomorphic to a
complex made up of flat quasicoherent sheaves of $\Oo_X$-modules,
concentrated in degrees $\geq 0$.
\end{lemma}
{\em Proof:}
We prove this by induction on $n$. It is well-known for $n=0$ (here $p$ is a
flat morphism of algebraic spaces). Suppose it is known for $n-1$. Choose a
smooth surjective morphism $Y\rightarrow T$ from a scheme $Y$. The condition
that $p$ is flat means that $Y$ is flat over $X$. Define the simplicial
$n-1$-stack
$$
G_k:= Y\times _T \ldots \times _T Y \;\;\; (k+1 \,\, \mbox{times}),
$$
which resolves $T$. The elements $G_k$ are $n-1$-stacks. This simplicial object,
{\em a priori} a weak functor $\Delta ^o \rightarrow (n-1)STACK(Sch /X)$,
may be strictified and considered as a simplicial object in the category of
$n-1$-prestacks on $Sch /X$. As such we can integrate it to an $n$-stack
${\bf i} G_{\cdot}$ which is equivalent to $T$. Now we can use this object to
calculate cohomology: let $q^k: G_k \rightarrow X$ denote the projections; then
$$
{\bf R} p_{\ast} (\Oo ) = \int _k{\bf R}q^k_{\ast} (\Oo )
$$
where here $\int _k$ means taking the associated total complex of the
cosimplicial complex of sheaves ${\bf R}q^{\cdot}_{\ast} (\Oo )$ on $X$.
By the inductive hypothesis, each of the terms ${\bf R}q^k_{\ast} (\Oo )$ is
quasiisomorphic to a complex of flat quasicoherent sheaves of $\Oo _X$-modules.
These choices can be made in a strictly coherent way (to prove this, use a
Reedy-type argument as in \cite{HirschowitzSimpson}), so that the total complex
itself is quasiisomorphic to a complex of flat quasicoherent sheaves of
$\Oo_X$-modules. \eop

The following lemma is a slight variation on Mumford's standard argument
\cite{Mumford}, and is also related to Illusie's use of Deligne-Lazard
\cite{Illusie}.

\begin{lemma}
\label{mumford2}
Suppose $C^{\cdot}$ is a complex of sheaves (supported in positive degrees)
on the
big etale site over a scheme $X$, such that each $C^k$ is a flat quasicoherent
sheaf on a noetherian scheme $X$. Suppose also that the $H^i(C^{\cdot})$ are
vector sheaves. Then $C^{\cdot}$ is a residually perfect complex.
\end{lemma}
{\em Proof:}
This proof is based on totally standard techniques. We give it in detail, just
so as to verify correctness of the statement. The proof is composed of several
steps.

\noindent
{\bf Step 1.} \, We show that the first nonzero cohomology sheaf (say $V=H^i$)
is a vector scheme. More precisely, we show the following: if
$$
0\rightarrow V \rightarrow E \rightarrow F
$$
is a complex of sheaves on the big site over a scheme $X$, such that $E$ and
$F$ are flat quasicoherent sheaves and $V$ is a vector sheaf, then $V$ is a
vector scheme.

We may assume that the base $X=Spec (\Oo (X))$ is affine; thus $E$ and $F$ are
respectively of the form $\tilde{M}$ and $\tilde{N}$ for flat $\Oo (X)$-modules
$M$ and $N$.

Choose a surjection $V'\rightarrow V$ from a vector scheme $V'$ to $V$.
Then we can consider the map $V'\rightarrow  E$ as being a section of $E|_{V'}$,
in other words as being given by an element of $\Oo (V')\otimes _{\Oo (X)}M$.
The element in question is a tensor product involving a finite number of
generators of $M$. Thus there is a submodule $M_1\subset M$ of finite type such
that our map factors through
$$
V'\rightarrow \tilde{M_1} \rightarrow \tilde{M}.
$$
Let $V''$ be a vector scheme surjecting to
the kernel of $V'\rightarrow V$. Thus the map $V'\rightarrow E$ restricts to
zero on $V''$. Again, we can express this fact using a finite number of
generators of the module $M$ (since going to zero in a tensor product is
expressed using a finite number of the relations in the definition of tensor
product). Thus, by possibly increasing the size of $M_1$ but keeping the finite
type condition, we may assume that $V''$ goes to zero in $\tilde{M_1}$. Thus we
have a factorization
$$
V\rightarrow \tilde{M_1} \rightarrow \tilde{M}.
$$
By Deligne-Lazard, we may assume that $M_1$ is flat, thus $\tilde{M_1}$
is a finite rank vector bundle which we denote $L$. We obtain a factorization
$$
0\rightarrow V \rightarrow L \rightarrow E
$$
with $L$ a finite rank vector bundle.
Note that the map $V\rightarrow L$ is
injective because the composition into $E$ is injective.

Now the fact that $V$ maps to zero in $F$ means that we have a map
$$
L/V \rightarrow F.
$$
The same argument as above shows that there is a finite rank vector bundle
$K$ with a factorization
$$
L/V\rightarrow K\rightarrow F.
$$
Let $W$ be the kernel of $L\rightarrow K$. Being the kernel of a map of finite
rank vector bundles, $W$ is a vector scheme. But $W$ certainly maps to zero in
$F$, so the map $W\rightarrow E$ takes $W$ to $V$. This is a retraction of $W$
onto $V$. Therefore we can write
$$
W=V\oplus V^{\perp}.
$$
In particular, $V$ is the kernel of the projection onto the second factor
$W\rightarrow W$, so $V$ is a vector scheme. This completes the proof of Step 1.

\noindent
{\bf Step 2.} \, If $V$ is a vector scheme with a morphism
$V\rightarrow \tilde{M}$ to a flat quasicoherent sheaf, then there is a vector
bundle $K$ (i.e. locally free sheaf of finite rank) and factorization
$V\hookrightarrow K \rightarrow \tilde{M}$.
This was shown in the course of the proof of Step 1.

\noindent
{\bf Step 3.} \, Suppose $V$ is a vector scheme with a morphism to a vector
bundle
$V\rightarrow L$ followed by $L\rightarrow \tilde{M}$ where $\tilde{M}$ is a
flat quasicoherent sheaf. Suppose that the composed morphism $V\rightarrow
\tilde{M}$ is zero. Then there is another vector bundle $L'$ and factorization
$L\rightarrow L' \rightarrow \tilde{M}$ such that the composition $V\rightarrow
L'$ is zero.

\noindent
{\bf Step 4.} \, We now construct a sequence of complexes $C^{\cdot}_k$ of flat
quasicoherent sheaves, all quasiisomorphic to $C^{\cdot}$ (in the derived
category), together with complexes of vector
bundles
$$
U^0_k \rightarrow U^1_k \rightarrow \ldots
$$
with morphisms $\varphi _k: U_k^{\cdot} \rightarrow C^{\cdot}_k$ such that
$\varphi _k$ induces an isomorphism on $H^i$ for $i<k$ and an injection on
$H^k$.
For $k=0$ we set $C^{\cdot}_0=C^{\cdot}$ and it suffices to take the complex
$U_k:= 0$. Thus we may now assume we have constructed
$C^{\cdot}_{k-1}$ and $U^{\cdot}_{k-1}$,  and we try to construct
$C^{\cdot}_{k}$ and $U^{\cdot}_{k}$.  Set
$$
W^{\cdot} := Cone
(U^{\cdot}_{k-1}\stackrel{\varphi _{k-1}}{\rightarrow} C^{\cdot}_{k-1}),
$$
in other words
$$
W^i :=U^i_{k-1}\oplus  C^{i-1}_{k-1}
$$
with differential $d_{U_{k-1}} + d_{C_{k-1}} + \varphi _{k-1}$,
with morphisms
$$
C^{\cdot}_{k-1}[-1] \rightarrow W^{\cdot} \rightarrow U^{\cdot}_{k-1}
$$
(these form a ``triangle'' in the derived category). Note that
$W^{\cdot}$ is a complex of flat quasicoherent sheaves, and
$$
H^i(W^{\cdot}) = 0 \;\; \mbox{for} \;\; i\leq k-1.
$$
The other cohomology sheaves are vector sheaves (recall that the category of
vector sheaves is closed under kernel, cokernel and extension
\cite{Hirschowitz} \cite{Simpson}(vii)).

We will construct (in Step 5 below) a complex of finite-rank vector bundles
$L^{\cdot}$ (supported in degrees $\geq k$)
with a morphism  $L^{\cdot}\rightarrow W^{\cdot}$ which induces an isomorphism
on $H^{k}$ and an injection on $H^{k+1}$. This will suffice to complete step 4,
since $C^{\cdot}_{k-1}$ is quasiisomorphic to $Cone(W^{\cdot} \rightarrow
U^{\cdot} _{k-1})$, so we
can set
$$
C^{\cdot} _{k}:=Cone(W^{\cdot} \rightarrow U_{k-1}^{\cdot})
$$
and set
$$
U^{\cdot} _{k}:= Cone(L^{\cdot} \rightarrow U_{k-1}^{\cdot}) \rightarrow
C^{\cdot}_k.
$$
The cone of this map has cohomology vanishing in degrees $\leq k$ so
$U_k^{\cdot}$ is the complex we are looking for at stage $k$, and this will
complete the proof of Step 4.

{\bf Step 5:}
Thus we are now reduced to constructing $L^{\cdot}\rightarrow W^{\cdot}$
for $W^{\cdot}$ as in the previous step, at stage $k$
(basically, we have reduced to the problem of constructing the first few
terms of
our complex of vector bundles).
By step 1 we get that the first nonzero cohomology object
$V:=H^k(W^{\cdot})$ is a
vector scheme. Pulling back the complex $W^{\cdot}$ to the scheme $V$ itself
and using the fact that it is a complex of quasicoherent sheaves and $V$ is
affine, we can lift the isomorphism $V=H^k(W^{\cdot})$ to a section
$$
V\rightarrow W^k.
$$
By looking at degrees of homogeneity, we can interpret this map as a morphism of
sheaves of groups on $X$.
By step 2 there is a factorization
$$
V\rightarrow N^k \rightarrow W^k
$$
with $N^k$ a vector bundle. Applying Step 3 to the morphism $N^k\rightarrow
W^{k+1}$ we obtain a diagram
$$
\begin{array}{ccc}
W^k & \rightarrow & W^{k+1} \\
\uparrow && \uparrow \\
N^k & \rightarrow & N^{k+1}
\end{array}
$$
and continuing in this way we obtain a complex of vector bundles
$N^{\cdot}$ with morphism to $W^{\cdot}$ and with a factorization
$$
V \rightarrow H^k(N^{\cdot})\rightarrow H^k(W^{\cdot}).
$$
In particular, the morphism on $k$-th cohomology (on the right in the display)
is surjective.
Also we get from here that $V\rightarrow N^k$ is injective.

Since $V$ is a
vector scheme, there is a vector bundle $B$ and morphism $N^k\rightarrow B$
having $V$ as kernel. We now set $M^{k+1}:= N^{k+1} \oplus B$, with the
morphism
$N^k\rightarrow M^{k+1}$ having the two existing maps as components. Note that
$$
V= \ker (N^k\rightarrow M^{k+1} ).
$$
We have a morphism
$$
N^k/V \rightarrow N^{k+1},
$$
and on the other hand an injection
$$
N^k/V \hookrightarrow B.
$$
The fact that $N^{k+1}$ is a coherent sheaf implies that it is an injective
object in the category of vector sheaves \cite{Simpson}(vii), so there exists
a morphism
$$
B\rightarrow N^{k+1}
$$
compatible with the map from $N^k$. This gives a retraction $M^{k+1} \rightarrow
N^{k+1}$ compatible with the map from $N^k$. We can now set $M^k:= N^k$, and
continue the choice of $M^{k+2}\rightarrow W^{k+2}$ etc. using Step 3, to obtain
a complex of finite-rank vector bundles $M^{\cdot}$ mapping to $W^{\cdot}$ with
$H^k(M^{\cdot}) = V$. In particular the map
$$
H^k(M^{\cdot} )\rightarrow H^k(W^{\cdot})
$$
is an isomorphism.

Finally we need to arrange for injectivity on $H^{k+1}$. Let
$$
U:= \ker \left( H^{k+1}(M^{\cdot})\rightarrow H^{k+1}(W^{\cdot}) \right) .
$$
It is the first nonzero cohomology of $Cone (M^{\cdot} \rightarrow W^{\cdot})$
and this cone is a complex of flat quasicoherent sheaves, on the other hand $U$
is a vector sheaf so by Step 1, $U$ is a vector scheme. Choose a locally free
sheaf $B'$ with injection $U\hookrightarrow B'$. The morphism
$$
U\rightarrow M^{k+1} /im(M^k)
$$
is an injective morphism of vector sheaves (because $U$ is by definition
contained in $H^{k+1}(M^{\cdot})$). Therefore, since $B'$ is an injective
object in the category of vector sheaves, we obtain a morphism
$$
g:M^{k+1} /im(M^k)\rightarrow B'
$$
extending the given morphism on $U$ (in particular, $g$ induces an injection on
$U$). Now put
$$
B'' := B' \oplus M^{k+2}
$$
with morphism
$$
h:M^{k+1}\rightarrow B''
$$
obtained by adding $g$ plus the
differential of $M^{\cdot}$. Note that
$$
\ker (h) \subset \ker (M^{k+1}\rightarrow M^{k+2})
$$
because one of the factors of $h$ is the morphism to $M^{k+2}$. Therefore
there is a morphism
$$
M^{k+1} /\ker (h) \rightarrow M^{k+2},
$$
and the fact that $M^{k+2}$ is an injective
object in the category of vector sheaves means that we can extend this to a
morphism on $B''$ giving a factorization of the differential
$$
M^{k+1} \rightarrow B'' \rightarrow M^{k+2}.
$$
Set
$$
L^k:= M^k,\;\; L^{k+1}:= M^{k+1}\;\;
$$
(with the differential of $M^{\cdot}$ as differential $L^k\rightarrow L^{k+1}$),
and
$$
L^{k+2} := M^{k+2} \oplus B'',
$$
with differential $L^{k+1}\rightarrow L^{k+2}$ equal to the sum of the
differential of $M^{\cdot}$ plus the morphism $h$. Above, we have constructed a
retraction $L^{k+2}\rightarrow M^{k+2}$ compatible with the differentials.
Finally, put $L^i:= M^i$ for $i\geq k+3$ and use the above retraction followed
by the differential of $M^{\cdot}$ to define the differential
$L^{k+2}\rightarrow L^{k+3}$. We now have a complex $L^{\cdot}$ with map
$$
L^{\cdot} \rightarrow M^{\cdot}.
$$
Note that $H^k(L^{\cdot})= H^k(M^{\cdot})$. The morphism of complexes is in fact
surjective and its kernel (which we shall note $K^{k+2}$ is concentrated in
degree $k+2$ with $K^{k+2} = B''$. We obtain a long exact sequence of
cohomology,
one part of which is
$$
0\rightarrow H^{k+1}(L^{\cdot})\rightarrow H^{k+1}(M^{\cdot} ) \rightarrow B'' .
$$
In particular, the morphism $U\rightarrow B''$ is injective so the image of
$$
H^{k+1}(L^{\cdot})\rightarrow H^{k+1}(M^{\cdot} )
$$
doesn't meet $U$. This implies that the composed morphism
$$
H^{k+1}(L^{\cdot})\rightarrow H^{k+1}(W^{\cdot} )
$$
is injective. We have now constructed the desired complex $L^{\cdot}$, which
completes the proof of Step 5.

{\bf Step 6:} The approximations $U_k^{\cdot}$ fit into Definition
\ref{residuallyperfect} to show that $C^{\cdot}$ is residually perfect.
\eop

We now obtain one half of Theorem \ref{criterion}:

\begin{corollary}
If $p:T\rightarrow Z$ is a relatively $1$-connected geometric $n$-stack which is
very presentable, then
${\bf R} p_{\ast}
\Oo$ is a residually perfect complex on $Z$.
\end{corollary}
{\em Proof:}
By Corollary \ref{cohVPisVS} the cohomology
sheaves of  ${\bf R} p_{\ast} \Oo$ are vector sheaves. On the other hand, by
Lemma \ref{mumford1}, ${\bf R} p_{\ast} \Oo$ is quasiisomorphic to a complex of
flat quasicoherent sheaves. Lemma \ref{mumford2} now implies that ${\bf R}
p_{\ast} \Oo$ is a residually perfect complex.
\eop

\oldsubnumero{Proof of the other half of Theorem \ref{criterion}}

For the other direction in the theorem, we use the ``very presentable
suspension'' ${\bf rep} Susp_N$ defined in \S 6. We suppose now
that we are in the situation of Theorem \ref{criterion} and that we know that
${\bf R}p_{\ast}\Oo$ is a residually perfect complex.
Put
$$
C^{\cdot}:= [{\bf R}
p_{\ast} \Oo ]^{\ast}
$$
which is a residually perfect complex supported in negative degrees.
The
$N$-truncated Dold-Puppe of $C^{\cdot}$
is equal to
the very presentable suspension:
$$
\Omega ^{\infty} {\bf rep} Susp _N^{\infty} (T/Z)=
\tau _{\leq N}DP (C^{\cdot}) .
$$
Note that $\tau _{\leq N}DP (C^{\cdot})$ is quasiisomorphic to the
Dold-Puppe of the truncation of a perfect complex \ref{residuallyperfect}. The
Dold-Puppe of a perfect complex is geometric \cite{Simpson}(ix),
so the truncation is an $N$-almost-geometric $N^{\rm gr}$-stack.

We obtain a morphism of $N$-stacks
$$
T\rightarrow \tau _{\leq N}DP (C^{\cdot}).
$$
Let $T'$ be the fiber of this morphism (relative to $Z$), and let $p':
T'\rightarrow Z$ denote the projection. Note that $T'$ is very presentable. We
claim that the truncated cohomology complex
$$
\tau _{\leq N-1}({\bf R} p'_{\ast} (\Oo ))
$$
is quasiisomorphic to the $N-1$-truncation of a perfect
complex.

To prove the claim, note that we have a fibration sequence
$$
\Omega \tau _{\leq N}DP (C^{\cdot}) \rightarrow T' \rightarrow T.
$$
But
$$
\Omega \tau _{\leq N}DP (C^{\cdot}) = \tau _{\leq N-1}DP (C^{\cdot}[1])
$$
is the $N-1$-truncation of a perfect complex. The cohomology of this fiber is
given by the graded-symmetric algebra on the complex
$$
\tau _{\leq N-1}C^{\cdot}[1],
$$
so the $N-1$-truncation of the higher direct image complex for $T'\rightarrow T$
is an $N-1$-truncation of a perfect complex locally over $T$. Using the Leray
spectral sequence we get that $\tau _{\leq N-1}({\bf R} p'_{\ast} (\Oo ))$
is the
$N-1$-truncation of a perfect
complex.

Define a sequence of  integers $N_i$ by $N_0=N$ and
$N_{i+1}=N_i-1$. Define a sequence of $N_i$-stacks
$$
T_i\stackrel{p_i}{\rightarrow}Z,
$$
and complexes $L_i$, starting with $T_0=T$; setting
$$
L_i:=  \tau _{\leq N_i}[{\bf R}
p_{i,\ast} \Oo ]^{\ast};
$$
and setting $T_{i+1}$ equal to the $N_i-1$-truncation of the fiber of
$$
T_i \rightarrow DP(L_i).
$$
From the above, we obtain that the $L_i$ are
$N_i$-truncations of perfect complexes; thus the $DP(L_i)$ are $N_i$-almost
geometric $N_i$-stacks. On the other hand, by Proposition \ref{nilpotence},
there is a $j$ such that $T_j=\ast _Z$.

Now we go backwards from $j$ to $0$.  Set $M_j=N_j$ and define $M_i$ for
$0\leq i
\leq j$ by $M_{i-1} = M_i -1$.
We know that $T_j$ is an
$M_j$-almost-geometric $M_j^{\rm gr}$-stack. Applying Lemma
\ref{agextension} we obtain
by descending induction on $i$ (starting with $i=j$), that $\tau _{\leq
M_i}(T_i)$
is an $M_i$-almost-geometric $M_i^{\rm gr}$-stack.
Finally we get back to the statement that
$$
\tau _{\leq M_0} (T)
$$
is an $M_0$-almost-geometric $M_0$-stack.

By choosing $N$ large enough at the start, we may assume that $M_0 \geq n+1$.
Now, applying Lemma \ref{aggeo}, we obtain that $T$ is a geometric $n$-stack.
This completes the proof of the second half of Theorem \ref{criterion}.
\eop

\subnumero{A semicontinuity-type result}

Related to the above criterion is the following result. It generalizes a similar
statement for perfect complexes.

\begin{theorem}
\label{semicon}
Suppose $f:T\rightarrow Z$ is a morphism from a very presentable geometric
$n$-stack to an integral (i.e. reduced and irreducible) scheme $Z$, such
that $f$
is relatively $1$-connected. If $z\in Z$ is a closed point, let $\pi _i(T_z)$
denote the ${\bf k}(z)$-vector space which represents the sheaf
$$
\pi _i(T\times _Z Spec ({\bf k}(z) )).
$$
(In order to get a base point to define this, it may be necessary to extend the
ground field but this doesn't affect the notion of the dimension of the vector
space in question.) If the functions
$$
z\mapsto dim _{{\bf k}(z)} \pi _i(T_z)
$$
are locally constant in the Zariski topology of $Z$, then the $\pi _i(T/Z)$ are
locally free coherent sheaves over $Z$, i.e. $T$ is ``locally free'' in the
terminology of \S 3.
\end{theorem}
{\em Proof:}
We prove this when $T$ is $k$-connected, by descending induction on $k$. It is
clearly true if $k=n$, so we may assume that $1\leq k <n$ and we may assume
that it is known for $k+1$-connected $n$-stacks. By the Hurewicz theorem,
(version for very presentable $n$-stacks) we have an isomorphism of vector
sheaves over $Z$,
$$
\pi _k(T/Z) \cong H^k(T/Z, \Oo )^{\ast}.
$$
One can see from the arguments used in the previous section that
$H^k(T/Z, \Oo )$ is a vector scheme, thus its dual is a coherent sheaf. In other
words, under the hypotheses of the theorem the first nonvanishing homotopy group
sheaf is a coherent sheaf. Semicontinuity for coherent sheaves over an
integral scheme implies that if the dimensions of the fibers are constant
then it
is locally free. Thus we obtain $\pi _k(T/Z)$ locally free. Now let $T'$ be the
fiber of $$
T\rightarrow K(\pi _k(T/Z)/Z, k).
$$
It shares the same homotopy group sheaves in degrees $i>k$, and is
$k+1$-connected. Furthermore, $T'$ is geometric. Thus the inductive statement of
the theorem for $k+1$-connected $n$-stacks implies that the $\pi _i(T'/Z)$ are
locally free. We now have that $\pi _i(T/Z)$ are locally free for all
$i\geq k$.
\eop

{\em Remark:} One cannot say, however, that the dimensions $\pi _i(T_z)$ are
semicontinuous as functions of $z$. From the above proof one gets only that the
first (in order of increasing $i$) function which is not locally constant, is
semicontinuous.

\subnumero{Criteria for geometricity of representing objects}

We can apply our criterion to the object ${\bf rep}(\Ff /\Ee
)$ representing the $\underline{AV}$-shape
of $\Ff /\Ee$, if it exists.

\begin{corollary}
Suppose $p:\Ff\rightarrow Z$ is a
morphism from an $\ngr$-stack to an integral scheme, and suppose that the weak
representing object $q:{\bf rep}^{\rm wk}(\Ff /Z )\rightarrow Z$ for the
$\underline{AV}(Z)$-shape of $\Ff /Z$ exists. Suppose that the higher direct
image complex ${\bf R}p_{\ast} (\Oo )$ is a perfect complex on $Z$. Then
${\bf rep}^{\rm wk}(\Ff /Z )$ is a geometric $\ngr$-stack over $Z$.
\end{corollary}
{\em Proof:}
The morphism
$$
\Ff \rightarrow {\bf rep}^{\rm wk}(\Ff /Z )
$$
induces  a quasiisomorphism of complexes over $Z$
$$
{\bf R}q_{\ast} (\Oo ) \rightarrow {\bf R}p_{\ast} (\Oo )
$$
(cf \ref{reduction2} for
example). The hypothesis therefore implies that
${\bf R}q_{\ast} (\Oo )$ is a perfect complex on $Z$. Then applying Theorem
\ref{criterion} we get that  ${\bf rep}^{\rm wk}(\Ff /Z )$ is a very presentable
geometric $\ngr$-stack over $Z$.
\eop

\begin{corollary}
\label{geometricity}
Suppose $p:\Ff\rightarrow \Ee$ is a morphism of $\ngr$-stacks, such that the
representing object
${\bf rep}(\Ff /\Ee )$ exists. Suppose that for every map from a scheme
$Z\rightarrow \Ee$, the higher direct image complex
${\bf R}p_{2,\ast} (\Oo )$ is a perfect complex on $Z$, where
$$
p_2: \Ff \times _{\Ee} Z\rightarrow Z
$$
is the second projection.  Then the
morphism
$$
{\bf rep}(\Ff /\Ee )\rightarrow \Ee
$$
is of type $\underline{AG}$, in other words the corresponding cartesian family
(cf \ref{correspondence}) corresponds to a morphism
$$
\Ee \rightarrow \underline{AG}.
$$
\end{corollary}
{\em Proof:}   Immediate. \eop

If $p:\Ff\rightarrow \Ee$ is a morphism of $\ngr$-stacks, then for any ``point''
i.e. morphism $z: Spec (k)\rightarrow \Ee$, we obtain an $\ngr$-stack
$\Ff _z$ over $Spec(k)$. If the representing object
${\bf rep}(\Ff /\Ee )$ exists, then by Remark \ref{basechangerep} we have
$$
{\bf rep}(\Ff _z/Spec(k))= {\bf rep}(\Ff /\Ee )\times _{\Ee}Spec(k).
$$

\begin{corollary}
\label{locfreeness}
In the situation of Corollary \ref{geometricity},
suppose that $\Ee$ receives a surjection from
a disjoint union of integral schemes, and suppose that the dimensions of the
homotopy group vector-spaces
$$
\pi _i({\bf rep}(\Ff _z/Spec(k)))
$$
are constant as functions of the point $z$. Then
the homotopy group local systems over $\Ee$ are local systems of
finite-rank vector bundles over $\Ee$, and the morphism
$$
{\bf rep}(\Ff /\Ee )\rightarrow \Ee
$$
is of type $\underline{AL}$, in other words the corresponding cartesian family
(cf \ref{correspondence}) corresponds to a morphism
$$
\Ee \rightarrow \underline{AL}.
$$
\end{corollary}
{\em Proof:}
Apply \ref{semicon}.
The hypothesis that there is a surjection $Y\rightarrow \Ee$ with
$Y$ a disjoint union of integral schemes, implies that for any map from a scheme
$Z\rightarrow \Ee$ there is (locally on $Z$ in the etale topology) a
factorization $Z\rightarrow Y\rightarrow \Ee$. Thus we may apply
\ref{semicon} for
the integral schemes (components of $Y$) to get the local-freeness result over
$Z$. Then in view of the definition of $\underline{AL}$ this gives the desired
statement. \eop

\numero{Formal groupoids of smooth type}
\label{formalpage}

We now turn to the problem of finding ``domains'', i.e. objects for
which the nonabelian cohomology is interesting. This discussion is by no means
exhausitive! We concentrate on the type of construction which will be most
useful
for defining things in Hodge theory: formal categories.  Work over a field
$k$ of
characteristic zero.

This section is mostly taken from Berthelot \cite{Berthelot} and Illusie
\cite{Illusie}.

\subnumero{Definitions}

Recall  that a {\em formal scheme} $N$ is a locally ringed space which is
locally
isomorphic to the formal completion of a scheme. There is a unique ideal
$I$ cutting out a reduced subscheme with the same underlying topological space,
and we put $N^{(k)}$ equal to the subscheme cut out by $I^k$ (note that this is
a scheme rather than a formal scheme). We define the {\em functor represented by
$N$} on the category $Sch /k$ of all $k$-schemes, by setting
$$
N(Y):= \lim _{\rightarrow , k} N^{(k)}(Y).
$$
This is most useful on the subcategory of noetherian schemes or even schemes of
finite type over $k$.

A {\em formal category} is a sextuple $(X, N, e,s,t,m)$ where $X$ is a
scheme and
$N$ is a formal scheme, and the rest are morphisms
$$
e: X\rightarrow N,\;\;\; s,t: N\rightarrow X, \;\;\; m: N\times _XN\rightarrow N
$$
together satisfying the usual axioms for a category in the category of formal
schemes, and subject to the additional condition that $e(X)$ is the topological
space underlying the formal scheme $N$. To restate this, a formal category is a
category in the category of formal schemes, such that the object object is a
scheme and the morphism object is supported along the diagonal $e(X)$.

Let $I$ denote the ideal defining $e(X)$  (which is a closed subscheme because
$e$ is a section of either of the projections $s,t$). Note that if $X$ is not
reduced, this will be different from the ideal previously denoted by $I$ (but it
can play the same role).

A formal category $(X,N,e,s,t,m)$ represents a functor
$$
\Ff ^{\rm pre}= \Ff ^{\rm pre}_{(X,N)}: Sch /k \rightarrow Cat
$$
since both $X$ and $N$ represent functors to sets (for $N$ use the construction
described above). We call $\Ff ^{\rm pre}$ the {\em $1$-prestack associated to
the formal category} and call the associated $1$-stack $\Ff$ the {\em $1$-stack
associated to the formal category $(X,N)$}.

The data of a formal category $(X,N)$ are equivalent to the data of a $1$-stack
$\Ff$ and a surjective morphism from a scheme $X\rightarrow \Ff$. Given $\Ff$
and the morphism, we obtain $N:= X\times _{\Ff} X$ (which automatically has the
required morphisms $e,s,t,m$). One must require of $(X\rightarrow \Ff )$
that the
resulting $N$ be represented by a formal scheme concentrated along the diagonal.
We often use interchangeably the notations $(X,N)$ and $(X,\Ff )$.

\subnumero{Formal categories of smooth type}

We say that a formal category $(X,N,e,s,t,m)$ is {\em of smooth type} if $N$
is locally the formal completion of a scheme of finite type and if the
projections $s$ and $t$ are formally smooth. This condition
 implies several things (from
\cite{Illusie} (ii), see also \cite{Berthelot}):

\begin{parag}
\label{smooth1}
Locally over $X$, $N$ is isomorphic to the formal completion of
$X\times {\bf A}^n$ along the zero section, for some $n$.

\newparag{smooth2}
That $I/I^2$ is a locally free $\Oo _X$-module and
$$
Gr _I^{\cdot}(\Oo _N) \cong \bigoplus _k Sym ^k(I/I^2).
$$
We denote $I/I^2$ by $\Omega ^1_{X/\Ff}$ and the dual vector bundle
$(I/I^2)^{\ast}$ on $X$ by $T(X/\Ff )$.

\newparag{smooth3}
In particular, $\Oo _N / I^k$ is locally free as an $\Oo _X$-module for either
of the two structures induced by the projections $s,t$.

\newparag{smooth4}
We obtain (see \cite{Illusie} (ii) for the argument) that the presheaf
of categories $\Ff ^{\rm pre}_{(X,N)}$ is in fact a presheaf of groupoids, thus
the stack associated to $(X,N)$ is a $1$-stack of groupoids. In view of this, we
can use interchangeably the terminology ``formal category of smooth type'' or
``formal groupoid of smooth type''.

\newparag{smooth5}
Set
$$
\Lambda ^r:= (\Oo _N /I^{r+1})^{\ast} :=  \underline{Hom}_{\Oo _X} (\Oo _N
/I^{r+1}, \Oo _X)
$$
where, for taking the dual we fix the structure of $\Oo_X$-module on
$\Oo _N$ determined by the projection $s$. The two projections $s,t$ induce
structures of left and right $\Oo _X$-module on $\Lambda ^r$ (which are not
in general the same). The multiplication $m$ induces morphisms
$$
\Lambda ^r \otimes _{\Oo _X} \Lambda ^s \rightarrow \Lambda ^{r+s}
$$
and if we set $\Lambda := \bigcup \Lambda ^r$ then $\Lambda$ is a sheaf of rings
of differential operators on $X$.

\newparag{smooth6}
The projections $s,t$ induce splittings
$$
\Lambda ^1 \cong \Oo _X \oplus (I/I^2)^{\ast},
$$
so $\Lambda$ is a {\em split almost-polynomial sheaf of rings of differential
operators on $X$} in the sense of \cite{Simpson}(iii).

\newparag{smooth7}
For a split almost-polynomial sheaf of rings of differential
operators on $X$ to come from a formal category of smooth type $(X,N)$,
one needs to have a compatible cocommutative coalgebra structure on $\Lambda$
(this corresponds to multiplication of functions in $\Oo _N$).
\end{parag}

\subnumero{Calculation of cohomology via the de Rham complex for $\Ff$}

We now fix a formal category of smooth type $(X,N)$ with associated
$1^{\rm gr}$-stack $\Ff = \Ff_{(X,N)}$ and associated sheaf of rings
of differential operators $\Lambda$. The cohomology of $\Ff $ can be calculated
by a de Rham complex for $\Lambda$. Again this is from \cite{Illusie} (ii), see
also \cite{Berthelot}. In this section we will derive this result without
refering to rings of differential operators such as $\Lambda$.

Recall the terminology that a {\em local system} over $\Ff $ is a relatively
$0$-connected morphism of stacks $L\rightarrow \Ff $. We get local systems of
groups, rings, modules, and so on. There is of course the structure sheaf
considered as local system of rings, which we denote
$$
\Oo \times \Ff = :\Oo _{\Ff}
\rightarrow \Ff .
$$
A {\em local system of $\Oo$-modules} is a local system
$L\rightarrow \Ff$ with relative structure of $\Oo _{\Ff }$-module.
We say that a local system of $\Oo$-modules $L$ is {\em $X$-locally free} if the
pullback $L_X:= L\times _{\Ff} X$ is a locally free sheaf of $\Oo
_X$-modules over
$X$. We say that a local system (of abelian groups or of $\Oo$-modules) $L$
is a {\em local system of vector sheaves over $\Ff$} if $L\times _{\Ff} X$ is a
vector sheaf over $X$.

\begin{parag}
\label{lambda}
If $L$ is a local system of vector sheaves on $\Ff$ then $L_X$ has a
structure of $\Lambda$-module. Conversely, a vector sheaf $L_X$
together with structure of $\Lambda$-module comes from a unique $X$-locally free
local system on $\Ff$.
\end{parag}

In the remainder of this subsection we fix a local system of vector sheaves $L$
over $\Ff$. We will see how to calculate the cohomology $H^i(\Ff , L)$.
This cohomology may be defined as
$$
H^i(\Ff , L):= \pi _0\underline{\Gamma}(\Ff , K(L/\Ff , i)),
$$
or it may be viewed as the cohomology in the topos of local systems of sets over
$\Ff$.

Let $X^{\rm zar}$ denote the Zariski topological space underlying $X$,
thought of as a site (i.e. it is a category). We have a functor ${\bf zar}$ from
$X^{\rm zar}$ to the category of open substacks of $\Ff$. (This should be
thought of as a ``morphism of objects'' from $\Ff$ to $X^{\rm zar}$). In
particular, if  $L\rightarrow \Ff$ is a local system, then we obtain a pullback
sheaf, which should be thought of as a direct image and which we
consequently denote ${\bf zar} _{\ast}(L)$, on $X^{\rm zar}$. This situation has
all of the same aspects as for the direct image via a usual morphism: we have
the higher derived direct image ${\bf R}{\bf zar} _{\ast}(L)$ and
$$
H^i(\Ff , L) = {\bf H}^i(X^{\rm zar}, {\bf R}{\bf zar} _{\ast}(L)).
$$

\oldsubnumero{The \v{C}ech-Alexander and de Rham complexes for $X/\Ff $}

The cohomology of $\Ff$ with coefficients in local system of vector sheaves
$L$ is calculated by a {\em de Rham complex} of sheaves on $X^{\rm zar}$,
denoted
$$
L_X\otimes \Oo _X \stackrel{d_{L/\Ff}}{\rightarrow}
L_X\otimes \Omega ^1_{X/\Ff} \stackrel{d_{L/\Ff}}{\rightarrow}  \ldots .
$$
The components are in fact sheaves of the form indicated, in other words
we have a sheaf of relative differential forms $\Omega ^k_{X/\Ff}$ which is
locally free as a sheaf of $\Oo _X$-modules on $X^{\rm zar}$, and the components
of the de Rham complex are of the form
$$
L_X \otimes _{\Oo _X} \Omega ^k_{X/\Ff}.
$$
However, the differentials $d_{L/\Ff}$ are not morphisms of $\Oo _X$-modules
(the standard example of this being the case $\Ff = X_{DR}$ where the
differentials are the usual de Rham differentials). The differentials are
morphisms of sheaves of abelian groups on $X^{\rm zar}$. Because of this, it is
not altogether easy to define this complex (or to see why it calculates the
cohomology of $L$). Our approach is based on the observation that the
Hodge-to de
Rham spectral sequence comes from the Hodge filtration on the de Rham complex;
but the Hodge filtration can actually be defined before knowing what the de Rham
complex is, for example it can be defined as a filtration on the
\v{C}ech-Alexander complex. Then the spectral sequence for this filtration, with
respect to the direct image functor ${\bf zar}_{\ast}$, gives the above de Rham
complex.

Our use of the \v{C}ech-Alexander complex is inspired by
Berthelot \cite{Berthelot} and Illusie \cite{Illusie}, and also by Constantin
Teleman \cite{CTeleman}. Note that our use of the ``Hodge-to-de Rham spectral
sequence'' is a mundane one and shouldn't be confused with the interesting
recent results of C. Teleman and I. Grojnowski.

The stack $\Ff$ has a resolution by the simplicial presheaf of sets
(whose components are formal schemes)
$$
\ldots X\times _{\Ff } X \tworightarrows X \rightarrow \Ff .
$$
The general term is $X\times _{\Ff} \ldots \times _{\Ff} X$.
Suppose $L$ is a local system of vector sheaves on $\Ff$. Then we obtain the
cosimplicial sheaf on $X^{\rm zar}$ which we denote
$$
{\bf zar}_{\ast} (X; L) \tworightarrows {\bf
zar}_{\ast} (X\times _{\Ff } X; L) \ldots .
$$
The general term is
$$
\check{C}A^{k}(X/\Ff /X^{\rm zar}; L):=
{\bf zar}_{\ast} (X\times _{\Ff }  \ldots \times _{\Ff} X;
L|_{X\times _{\Ff }  \ldots \times _{\Ff} X})
$$
where ${\bf zar}_{\ast}$ denotes the ``direct image'' construction analogous to
that described previously, here for
$$
X\times _{\Ff }  \ldots \times _{\Ff} X\rightarrow X^{\rm zar}.
$$
The {\em \v{C}ech-Alexander complex} is the complex of sheaves on $X^{\rm zar}$
associated to this cosimplicial sheaf. It is denoted
$\check{C}A^{\cdot}(X/\Ff /X^{\rm zar}; L)$. It calculates the cohomology of
$\Ff$ with coefficients in $L$.

\begin{parag}
\label{hf}
Let $\underline{e}(X)$ denote the constant simplicial subobject
 of the simplicial presheaf of sets considered above (all of its
terms are $X$, the inclusion being given by the iterated degeneracy map).
In each
component it is a closed sub-formal scheme. Let $I$ denote its defining ideal
(this means the collection consisting of the ideal defining the closed subscheme
in each component). Then we obtain a filtration of $\check{C}A^{k}(X/\Ff /X^{\rm
zar}; L)$ by setting
$$
F_k\check{C}A^{\cdot }(X/\Ff /X^{\rm zar}; L) :=
I^k \check{C}A^{\cdot }(X/\Ff /X^{\rm zar}; L).
$$
This filtration is the {\em Hodge filtration} of the \v{C}ech-Alexander complex
(corresponding to a trivial filtration on $L$).
\end{parag}

The spectral sequence for a filtered complex, for the cohomology sheaves of
the complex of sheaves, is
$$
\underline{H}^i(Gr ^{k}_F\check{C}A^{\cdot }(X/\Ff /X^{\rm zar}; L))
\Rightarrow
\underline{H}^i\check{C}A^{\cdot }(X/\Ff /X^{\rm zar}; L).
$$
The first differential is
$$
d: \underline{H}^i(Gr ^{k}_F\check{C}A^{\cdot }(X/\Ff /X^{\rm zar}; L))
\rightarrow \underline{H}^{i+1}(Gr ^{k-1}_F\check{C}A^{\cdot }(X/\Ff
/X^{\rm zar}; L)).
$$

Let $\Ff _0$ be the ``normal cone'' of $\Ff$. It comes from a formal category
$(X, N_0)$ where
$$
\Oo _{N_0} = Gr ^{\cdot}_I(\Oo _N)
$$
(cf \ref{smooth2}). In terms of our later discussion of the Hodge
filtration, the
formal category $\Ff _0$ is the fiber of ${\bf Hodge}(\Ff )\rightarrow \af$ over
$0\in \af$. In concrete terms, $N_0= \widehat{T(X/\Ff )}$ is the completion of
the ``tangent bundle'' $T(X/\Ff )$ (cf \ref{smooth2}) along the zero-section,
and both projections $s$ and $t$ are the same projection of $N_0$ to $X$.
We have
$$
\Ff
_0= K(\widehat{T(X/\Ff )}/X, 1).
$$
Note in particular that there is a projection $\Ff _0\rightarrow X$
(which doesn't exist for a more general formal category such as $X_{DR}$).

There is a local system $L_0$ on $\Ff _0$ corresponding to the trivial
filtration on $L$, which is automatically compatible with the Hodge filtration
on $\Ff$. Since it is the trivial filtration, $L_0:= Gr (L) \cong L$ as a vector
sheaf on $X$, and the action of $N_0$ is trivial. Thus the local system $L_0$ on
$\Ff _0$ is the pullback of a vector sheaf on $X$ by the
projection $\Ff _0\rightarrow X$.

We
have
$$
Gr ^{\cdot} _F\check{C}A^{k }(X/\Ff /X^{\rm zar}; L) =
\check{C}A^{k}(X/\Ff _0/X^{\rm zar}; L_0)
$$
where the degree in the grading of the associated-graded corresponds to the
degree
of the decomposition of the $\Gm$-action on $\check{C}A^{k}(X/\Ff _0/X^{\rm
zar};
L_0)$ (action by homotheties on $N_0=\widehat{T(X/\Ff )}$ with trivial action on
$L_0$).

Thus, in order to calculate the terms in the previous spectral sequence,
it suffices to calculate the cohomology of the \v{C}ech-Alexander complex in the
case of $(X, \Ff _0, L_0)$. Recall that the coefficient local system $L_0$ comes
from a vector sheaf $L_X$ on $X$ with trivial action of the gerb $\Ff _0$.

For this case, we suppose in general that $V\rightarrow X$ is a vector bundle.
Let $\widehat{V}$ be the completion of $V$ along the zero-section, and let
$V_{DR/X}:= V_{DR} \times _{X_{DR}}X$.
Then we have a fibration diagram
$$
K(\widehat{V}/X, 1)\rightarrow K(V/X, 1) \rightarrow K(V_{DR/X}/X, 1).
$$
Shifting gives
$$
V_{DR/X}\rightarrow K(\widehat{V}/X, 1)\stackrel{p}{\rightarrow} K(V/X, 1).
$$
The morphism $V_{DR/X} \rightarrow X$ is acyclic for vector sheaves pulled back
from $X$. This is the ``Poincar\'e lemma'' in our treatment. For the proof one
can reduce to the case $rk(V)=1$ and then do an explicit calculation with the
\v{C}ech-Alexander complex to prove the acyclicity.

From the acyclicity in the previous paragraph we get that if $L'_0$ denotes the
pullback of $L_X$ to $K(V/X, 1)$ then
$$
H^i( K(\widehat{V}/X, 1)/X, L_0) =
H^i( K(V/X, 1)/X, L'_0).
$$
The {\em Eilenberg-MacLane-Breen calculations}
\cite{Breen} (cf Theorem \ref{bc})
say that
$$
H^i( K(V/X, 1)/X, L'_0) = L_X\otimes _{\Oo _X}\bigwedge ^i_{\Oo _X}(V^{\ast}).
$$
Thus we obtain
$$
H^i( K(\widehat{V}/X, 1)/X, L_0)=
L_X\otimes _{\Oo _X}\bigwedge ^i_{\Oo _X}(V^{\ast}).
$$
Note that since $V$ is locally free, there is no other contribution from the
``universal coefficients theorem'' for passing from a coherent sheaf as
coefficients as in \ref{bc}, to a general vector sheaf $L_X$.

Put this back into the previous situation with $V= T(X/\Ff )$. Note that
$$
V^{\ast} = \Omega ^1_{X /\Ff _0}=\Omega ^1_{X /\Ff }
$$
so we get
$$
H^i(\Ff _0 /X, L_0) = L_X \otimes _{\Oo _X} \Omega ^i_{X/\Ff }.
$$
Finally, note that the homothety action of $\Gm$ gives an action of pure degree
$-i$ on  $H^i(\Ff _0 /X, L_0)$.
Putting this all together we get that
$$
\underline{H}^i(Gr ^{k}_F\check{C}A^{\cdot}(X/\Ff /X^{\rm zar}; L))
= {\bf zar} _{\ast} (L_X \otimes _{\Oo _X} \Omega ^i_{X/\Ff })
$$
for $k=-i$ and it is zero otherwise. The first differential of the spectral
sequence thus defines a differential
$$
d_{L/\Ff}: {\bf zar} _{\ast} (L_X \otimes _{\Oo _X} \Omega ^i_{X/\Ff })
\rightarrow {\bf zar} _{\ast} (L_X \otimes _{\Oo _X} \Omega ^{i+1}_{X/\Ff }).
$$
This complex is the {\em de Rham complex} for $X/\Ff$ with coefficients in $L$.
The spectral sequence degenerates and we immediately get that the cohomology
sheaves of the \v{C}ech-Alexander complex are isomorphic to those of the de Rham
complex. However, in fact we can use the degenerate form of the spectral
sequence to directly define a quasiisomorphism between the \v{C}ech-Alexander
complex and the de Rham complex. Let
$$
B^{\cdot}\subset \check{C}A^{\cdot}(X/\Ff /X^{\rm zar};
L)
$$
denote the subcomplex defined by
$$
B^i := \delta ^{-1}(F_{-i-1}) \subset F_{-i}\check{C}A^{\cdot}(X/\Ff
/X^{\rm zar};
L).
$$
Here $\delta$ denotes the differential of $\check{C}A^{\cdot}$.
We have a  morphism
$$
B^{\cdot} \rightarrow \bigoplus _{i} \underline{H}^ i(
Gr ^{-i}_F   \check{C}A^{\cdot}(X/\Ff /X^{\rm zar};
L))
$$
which sends the differential of $B^{\cdot}$ to the differential coming from the
spectral sequence. The spectral sequence differential was by definition our
differential for the de Rham complex, so this morphism can be rewritten as
$$
B^{\cdot} \rightarrow ({\bf zar}_{\ast}L_X \otimes _{\Oo _X} \Omega ^i_{X/\Ff },
d_{L/\Ff }). $$
The degenerate form of the first term in the spectral sequence implies that this
morphism as well as the inclusion
$B^{\cdot}\hookrightarrow \check{C}A^{\cdot}$ are quasi-isomorphisms. Thus
we have
established a chain of two quasiisomorphisms linking
the de Rham complex and the \v{C}ech-Alexander complex:
$$
({\bf zar}_{\ast}L_X \otimes _{\Oo _X} \Omega ^i_{X/\Ff }, d_{L/\Ff })
$$
$$
\leftarrow
B^{\cdot} \hookrightarrow
$$
$$
\check{C}A^{\cdot}(X/\Ff /X^{\rm zar}; L).
$$

Since the \v{C}ech-Alexander complex calculates the cohomology of $\Ff$, we
obtain---by globalizing over $X^{\rm zar}$---an isomorphism
$$
H^i(\Ff , L) \cong {\bf H}^i(X^{\rm zar},
({\bf zar}_{\ast}L_X \otimes _{\Oo _X} \Omega ^i_{X/\Ff }, d_{L/\Ff }) ).
$$

Now the global ``Hodge to de Rham'' spectral sequence coming from the ``stupid
filtration'' on the de Rham complex (but which could also be seen directly as
coming from the Hodge filtration on the \v{C}ech-Alexander complex) is
$$
H^q(X^{\rm zar}, {\bf zar}_{\ast}L_X \otimes _{\Oo _X} \Omega ^p_{X/\Ff })
\Rightarrow H^{p+q}(\Ff , L).
$$
The beginning term here is of course the same thing as a cohomology group of
$X$, so we can write the spectral sequence as
$$
H^q(X, L_X \otimes _{\Oo _X} \Omega ^p_{X/\Ff })
\Rightarrow H^{p+q}(\Ff , L).
$$
This proves finite-dimensionality of the $H^{i}(\Ff , L)$, for example.

It is this latter utilization which interests us for the relative case.
Suppose $X\rightarrow \Ff \rightarrow S$ is a morphism from a formal category to
a scheme $S$. We assume that $X$ is flat over $S$, and that $X/\Ff$ is of smooth
type. Implicitly we are assuming that the formal category $(X, N)$ corresponding
to $\Ff$, maps to the formal category corresponding trivially to $S$ (this is
what is meant by the existence of the morphism $\Ff \rightarrow S$). Define in
this setting a {\em higher derived \v{C}ech-Alexander complex}
$$
{\bf R} \check{C}A^{\cdot}(X/\Ff /S; L),
$$
which is a complex of sheaves on $S$ obtained by integrating the
cosimplicial complex of sheaves, higher derived direct image of the
components of
the usual $\check{C}A^{\cdot}$. This avoids trying to define a site $X^{\rm zar}
/S$ etc. From here we can employ exactly the same argument as above:
filter the \v{C}ech-Alexander complex by the ``Hodge filtration'' and take the
spectral sequence for this filtered complex as above. We obtain a
quasiisomorphism
between ${\bf R} \check{C}A^{\cdot}(X/\Ff /S; L)$ and the higher direct image of
the de Rham complex.
This gives a relative Hodge-to-de Rham spectral sequence of the form
$$
H^q(X/S, L_X \otimes _{\Oo _X} \Omega ^p_{X/\Ff })
\Rightarrow H^{p+q} (\Ff /S, L).
$$
On the other hand, if $L$ is $X$-locally free, then we also get (from the usual
Mumford-type argument using flatness of $X/S$) that
$$
{\bf R} \check{C}A^{\cdot}(X/\Ff /S; L),
$$
is a {\em perfect complex} on $S$. In particular, the cohomology
sheaves $H^i(\Ff /S, L)$ are vector sheaves on $S$. Perfectness implies that the
$n$-stack
$$
\underline{\Gamma} (\Ff /S, K(L/\Ff , n)) \rightarrow S
$$
is {\em geometric}.

It follows (still in the case that $L$ is $X$-locally free) that if the
dimensions
of the cohomology are constant on  $S$, then the relative cohomology sheaves
$H^i(\Ff /S, L)$ are locally
free over $S$ (compatible with base-change).

In \S 10 we extend this discussion to the case where $L_X$ is locally the
pullback
of a vector sheaf from $S$, showing then that the $H^i(\Ff /S, L)$ are vector
sheaves on $S$.

\subnumero{Blowing up formal categories}

The following construction will be of primary importance in many of our
examples.
It will be used in the next section below. Suppose $(X,N)$ is a formal category
of smooth type, with corresponding $1$-stack $\Ff$. Suppose $D\subset X$ is a
Cartier divisor. Recall that, according to the definition of Cartier divisor,
the local defining equations of $D$ are not zero-divisors on $X$.
We define
as follows the {\em formal category obtained by blowing up along $D$}
denoted $(X, {\bf BL}_D(N))$, with associated $1$-stack which we shall denote
(abusing notation) also by  ${\bf BL}_D(\Ff )$. The underlying scheme is the
same $X$. Let $\tilde{N}$ denote the blow-up of $N$ along $e(D)$ (recall that
$e:X\rightarrow N$ is the  morphism ``identity''). Let $\tilde{e}: X\rightarrow
\tilde{N}$ denote the strict transform of $X$ (see below for why it is
well-defined) and let  ${\bf BL}_D(N)$ denote the formal completion of
$\tilde{N}$ along $\tilde{e}(X)$. The morphisms $s$ and $t$ give, by
composition,
morphisms  $\tilde{s}, \tilde{t}: \tilde{N}\rightarrow X$ and by restriction
these give the morphisms $s',t'$ for $(X, {\bf BL}_D(N))$. Similarly, the
morphism $\tilde{e}$ gives the morphism $e'$ for $(X, {\bf BL}_D(N))$.

We need to show several things: first, why the ``strict transform of $X$'' is
well-defined; then why the multiplication map $m'$ is well-defined and
associative; and finally, why the resulting formal category is of smooth type.

\begin{parag}
\label{bld1}
The first task is to determine the effect of blowing up
as above. Use the projection $s: N\rightarrow X$ to pull back the
normal bundle of $X$ in $N$, to $N$; call this bundle $E$. Assume that $D$ is
defined by a regular function $z$ on $X$. Denote also by $z$ the pullback of
this function to a function on $N$.
 We can choose  a section $\eta$ of $E$ with $e(X)$ the zero-subscheme of
$\eta$,
and such that $d\eta$ is the identity isomorphism between the normal bundle of
$N$ in $X$ and $E|_X$. (The previous two sentences require that we restrict to a
Zariski open subset in $X$.)

\newparag{bld2}
Let $\tilde{N}$ be the blow-up of $N$ with center $D$, which we
can write down explicitly using the equations $(\eta , z)$ for $D$. Let
$\tilde{N}'$ denote the standard coordinate neighborhood of  $\tilde{N}$
containing the set-theoretic strict transform of $X$. If $Z$ is any scheme, then
a morphism $Z\rightarrow {\bf BL}_D(N)'$ is the same thing as a pair $(f, \alpha
)$ where $f: Z\rightarrow N$, and $\alpha $ is a section of $f^{\ast}(E)$ such
that
$$
f^{\ast}(\eta ) = z\alpha .
$$
Recalling that $e:X\rightarrow N$ is the morphism ``identity'' (whose image is
the set-theoretic support of the formal scheme $N$) we obtain a morphism
$\tilde{e}:X\rightarrow \tilde{N}'$ by setting
$$
\tilde{e}:= (e, 0).
$$
This formula for $\tilde{e}$ shows that the ``strict
transform'' of $X$ is well-defined (one can check that it is independent of the
choices we have made). Let ${\bf BL}_D(N)$ be the formal completion of
$\tilde{N}'$ along
$\tilde{e}(X)$.

\newparag{bld3}
A morphism $Z\rightarrow {\bf BL}_D(N)$ is a pair
$(f, \alpha )$ where $f: Z\rightarrow N$ and $\alpha$ is a section of
$f^{\ast}(E)$  such that $f^{\ast}(\eta ) = z\alpha $ and
$\alpha |_{Z^{\rm red}}=0$ (this latter condition meaning that $\alpha$ is
infinitesimally near to the zero-section).

\newparag{bld4}
We would like to put this into a more invariant formulation. Suppose $Z$ is a
scheme with a morphism $f: Z\rightarrow N$. We would like to classify the
liftings to $Z\rightarrow {\bf BL}_D(N)$. For this, suppose that the
composition $p:=s\circ f: Z\rightarrow X$ is flat. Then the section $\alpha$ in
the above description will be a section denoted
$$
a\in H^0(X, p_{\ast} (\Oo _Z) \otimes _{\Oo _X} E)
$$
subject to the equation $za = b$ where $b$ is the section in the same space
corresponding to $f^{\ast}(\eta )$. In this flat case, the fact that $z$ is not
a zero-divisor in $\Oo _X$  implies that multiplication by  $z$ induces an
injection $$ z\cdot -: H^0(X, p_{\ast} (\Oo _Z) \otimes _{\Oo _X}
E)\hookrightarrow  H^0(X, p_{\ast} (\Oo _Z) \otimes _{\Oo _X} E).
$$
Therefore, if $\alpha$ exists, it is uniquely defined. Thus, if a lifting of
$Z\rightarrow N$ into $Z\rightarrow {\bf BL}_D(N)$ exists, it is uniquely
defined. In particular, when we define our multiplication $m'$, it will
automatically be associative.

Now observe that the lifting $\alpha$ exists if and only if the inverse image
$p^{-1}(D)$ (pullback of $D$ from $X$ to $Z$),
maps to $D\subset N$ under the map
$Z\rightarrow N$. Using this, we can see that the restriction of
the multiplication $m$ to a map
$$
{\bf BL}_D(N)\times
_X{\bf BL}_D(N)\rightarrow N,
$$
lifts to a multiplication $m'$. Indeed, the inverse image of $D\subset X$
in ${\bf BL}_D(N)$ under either of the two projections, is the exceptional
divisor (recall that ${\bf BL}_D(N)$ is the formal neighborhood of the strict
transform of $X$ in the blow-up $\tilde{N}$, and in this neighborhood the
equation defining the exceptional divisor is the same as the equation of $D$
pulled back from $X$). It follows that the inverse image of $D$ in ${\bf
BL}_D(N)\times _X{\bf BL}_D(N)$, by either of the two end projections, is the
product of two copies of the exceptional divisor over $D$; and this maps  by $m$
to $D\subset N$. Thus by the previously-mentionned criterion, the lifting $m'$
exists.

Finally, note that the projection $s: {\bf BL}_D(N)\rightarrow X$ is formally
smooth. To see this, note that locally the projection looks like what is
obtained by blowing up the divisor $D$ in the zero section of a vector bundle
over $X$, and then taking the completion around the strict transform of $X$; but
the answer is again the completion of the zero section in a new vector bundle
obtained from the original one by twisting by $\Oo _X(-D)$, in particular it is
formally smooth.

By making the same calculations as above but using the projection $t$ from the
beginning instead of $s$, we get that the projection $s$ is formally smooth.
Thus, the formal category $(X, {\bf BL}_D(N))$ is of smooth type.
This completes the construction.

\newparag{bld6}
We have the formula
$$
\Omega ^1_{X/{\bf BL}_D(\Ff )} = \Omega ^1_{X/\Ff }(D),
$$
where the latter means sections having one pole along $D$.
The first differential of the de Rham complex comes from the composition
$$
d: \Oo _X \rightarrow \Omega ^1_{X/\Ff }\hookrightarrow
\Omega ^1_{X/\Ff }(D).
$$
\end{parag}

\subnumero{The Hodge filtration, revisited}

We now look at the  Hodge filtration of \ref{hf}, from a slightly
different point of view. Recall from \cite{Simpson}(i) that a {\em filtered
scheme} is a family $X\rightarrow \af $ with action of $\Gm$ covering the
standard
action on $\af $.

We can make the same definition for a filtered $n$-stack, using the
homotopy-coherent notion of ``group action'' \cite{Simpson}(vi).
The homotopy-coherent notion of group action is
the same as the following: a group $G$ acting on an $n$-stack $T$ is a fibration
sequence of $n$-stacks
$$
T\rightarrow E \rightarrow BG.
$$
With these notations, a {\em filtered $n$-stack} is a diagram
$$
\begin{array}{ccccc}
W & \rightarrow & E & \rightarrow & B\Gm \\
\downarrow && \downarrow && \downarrow \\
{\bf A} ^1 &\rightarrow & {\bf A^1}/\Gm & \rightarrow & B\Gm
\end{array}
$$
where the rows are fibration sequences.
This data is obtained from just $E\rightarrow \af /\Gm$.

The {\em underlying $n$-stack} is the fiber of $W\rightarrow \af $
over the point $1\in \af $. This is the same as the fiber of $E$ over the
point $[1]: Spec (k)\rightarrow \af /\Gm$. The {\em associated-graded $n$-stack}
is the fiber of $W$ over $0\in \af $, which is the same as the fiber of $E$ over
$0\rightarrow  \af /\Gm$. Note that this latter morphism is in fact a
morphism $[0]: B\Gm \rightarrow \af /\Gm$, so the associated-graded
naturally comes with an action of $\Gm$.

In fact, one can abstract the definition of \cite{Simpson}(i) even more, by
saying that a {\em filtered point} in a $1$-stack $\underline{M}$, is a morphism
of stacks $\af / \Gm \rightarrow \underline{M}$. The definitions of filtered
objects given in \cite{Simpson}(i) are obtained by applying this
definition to the parametrizing $1$-stacks $\underline{M}$ for the type of
object considered. The same definition works for an $n+1$-stack $\underline{M}$:
a {\em filtered point of $\underline{M}$} is a morphism of $n+1$-stacks
$$
 f: \af / \Gm \rightarrow \underline{M},
$$
where the $1$-stack $\af / \Gm$ is
considered as an $n+1$-stack. The above notion of ``filtered $n$-stack'' is
obtained by applying this definition with $\underline{M} = n\underline{STACK}$.
(For this translation, use the notion of cartesian family and Proposition
\ref{correspondence}, noting that $\af / \Gm$ is a stack of groupoids.)
Again, the
``underlying point'' of $f$ is the restriction of $f$ to $[1]$, and the
``associated-graded'' (with its $\Gm$-action) is the restriction of $f$ to
$[0] =
B\Gm$.

Apply these definitions to the $1$-stacks associated to a formal groupoid.
Suppose $(X,N)$ is a formal groupoid of smooth type. Then we define the {\em
Hodge filtration} of this formal groupoid to be
$$
(X\times \af , {\bf BL}_{X\times 0}(N))
$$
together with its natural action of $\Gm$. This yields the {\em Hodge
filtration} on the associated stack, which we can denote by
$$
{\bf Hodge}(\Ff ) \rightarrow \af /\Gm .
$$
The fiber over $[1]$ is just $\Ff$. The fiber over $0$ is the formal groupoid
$(X, N_0)$ where $N_0$ is the formal completion of the zero-section in the
normal bundle of $N$ to $e(X)$.

The relative \v{C}ech-Alexander complex for $X\times \af \rightarrow {\bf
Hodge}(\Ff
)$ relative to $\af$, is actually that obtained by applying the construction
$\xi$ of \cite{Simpson}(i) to the Hodge filtration of the  \v{C}ech-Alexander
complex for $X\rightarrow \Ff$. Thus, the filtration induced by ${\bf Hodge}(\Ff
)$ on the cohomology of $\Ff$, is the same as that defined in
\ref{hf}.

It would be nice to be able to describe how to obtain the spectral sequence
directly from this version of the notion of filtration. This would give a
geometric interpretation of the spectral sequence using ${\bf Hodge}(\Ff )$.
This should permit one to  conclude, for example, that the differentials in the
spectral sequence respect the product structure, i.e. satisfy the Leibniz rule.
In turn, this would determine the de Rham complex once the first
differential $d:
\Oo _X \rightarrow \Omega ^1_{X/\Ff}$ is known.
We leave this problem open for further research.
(But of course it should be said that the arguments of Berthelot
\cite{Berthelot} and Illusie \cite{Illusie} serve to calculate the differentials
and to prove that they satisfy the Leibniz rule; the question raised here is how
to see this geometrically.)

\subnumero{Morphisms of smooth type}

Suppose
$X\rightarrow \Ff$ and
$S\rightarrow \Ee$ are formal categories. A {\em morphism of formal categories}
is a commutative diagram of the form
$$
\begin{array}{ccc}
X & \rightarrow & \Ff \\
\downarrow{\scriptstyle f_X} && \downarrow {\scriptstyle f}\\
S& \rightarrow & \Ee .
\end{array}
$$
Such a morphism induces a morphism of locally free sheaves on $X$,
$$
d(f):T(X/\Ff )\rightarrow f_X^{\ast}T(S/\Ee ).
$$
We say that $f$ is {\em of smooth type}
if $X$ is flat over $S$ and if the above morphism $d(f)$ is surjective.

Note that even if the base $\Ee =S$ is already a scheme, the condition that $f$
be of smooth type is non-vacuous; it says that $f:X\rightarrow S$
should be flat.

As a sidelight (useful later on but not related to the remainder of the present
subsection) we say that a morphism $f$ is {\em projective} if $f_X: X\rightarrow
S$ is projective.

\begin{proposition}
\label{usualsmooth}
Suppose $f:(X,\Ff )\rightarrow (S, \Ee )$ is  a morphism of smooth type
between formal categories of smooth type. If $(S',\Ee ')$ is any formal
category of smooth type with a morphism
$(S',\Ee ')\rightarrow (S, \Ee )$ then the fiber product
$$
(S' \times _S X, \Ee ' \times _{\Ee } \Ff )
$$
is again a formal category of smooth type and the projection to $(S', \Ee ')$
is a morphism of smooth type. In particular applying that with $S'=\Ee ' = S$
we get that
$$
(X, \Ff \times _{\Ee} S)
$$
is a formal category of smooth type.
\end{proposition}
{\em Proof:} Let $M,M', N$ denote the morphism objects of $\Ee , \Ee ', \Ff$
respectively. The morphism object $N'$ of the formal category  $\Ff ':= \Ee '
\times _{\Ee } \Ff$ is given by the formula
$$
N'= (S' \times _S X)\times _{\Ee ' \times _{\Ee } \Ff}(S' \times _S X)
$$
$$
=(S'\times _{\Ee '}S')\times _{S\times _{\Ee}S} (X\times _{\Ff}X)
$$
$$
= M'\times _{M} N = (M'\times _S'X') \times _{M\times _XX'} (N\times _XX')
$$
where $X':= S'\times _SX$. The three terms
$$
M'\times _S'X', \;\;\; M\times _XX', \;\;\; N\times _XX'
$$
are formal schemes mapping to $X'$ (say by the projection $s$), which are
formally smooth over $X'$ with relative tangent spaces along the unique section,
respectively
$$
T(S'/\Ee ')|_{X'},\;\;\; T(S/\Ee )|_{X'},\;\;\;
T(X/\Ff )|_{X'}.
$$
The surjectivity of the morphism
$$
d(f)|_{X'}: T(X/\Ff )|_{X'}\rightarrow T(S/\Ee )|_{X'}
$$
implies that the
pullback
$$
(M'\times _S'X') \times _{M\times _XX'} (N\times _XX')
$$
is formally smooth over $X'$, which gives that $X'\rightarrow \Ff '$ is a formal
category of smooth type. It is easy to see that the projection to $\Ee '$ is of
smooth type. The last statement is obtained by applying the main statement to
the formal category $S'=\Ee ' = S$ (which is of smooth type with
morphism object trivial equal to $e(S)$).
\eop

The notion of morphism of smooth type is compatible with taking the ``Hodge
filtration'' as viewed above.

\begin{lemma}
\label{hodgesmooth}
Suppose $f: (X,\Ff ) \rightarrow (S, \Ee )$ is a morphism of formal categories
which is of smooth type. Then the morphism
$$
{\bf Hodge}(f): (X\times \af , {\bf Hodge}(\Ff )) \rightarrow (S\times \af ,
{\bf Hodge}(\Ee
)) $$
is of smooth type.
\end{lemma}
{\em Proof:}
We have
$$
\Omega ^1_{X\times \af /{\bf Hodge}(\Ff )} = \Omega ^1_{X/\Ff} |_{X\times \af}
$$
functorially. The morphism induced by ${\bf Hodge} (f)$ is just the
pullback of the
morphism induced by $f$, so if the latter is surjective then so is the former.
\eop

{\em Remark:}
The formula given in the previous proof is functorial but is not compatible with
the $\Gm$-action on $\af$. If you want a formula compatible with this
action then
you have to tensor by  $\Oo _{\af}(1)$.

\numero{Formal categories related to Hodge theory}
\label{hodgepage}

Our nonabelian Hodge theory for nonabelian cohomology requires the construction
of various different formal groupoids of smooth type, and morphisms of smooth
type between them. By making these constructions, we obtain: the
nonabelian de Rham cohomology, the nonabelian Dolbeault cohomology, the Hodge
filtration (relating these two via a deformation), the Gauss-Manin connection on
the above, logarithmic homotopy type and regular singularities of the
Gauss-Manin
connection, and Griffiths transversality.  In order to obtain these properties
for the nonabelian cohomology and shape, one should apply the general results
of \S 10 below to the formal categories in the present section.

Most of these examples of formal categories are classical. In my own works  on
nonabelian Hodge theory, they appeared as split almost polynomial sheaves of
rings of differential operators in \cite{Simpson}(iii), and then as formal
groupoids used for degree $1$ nonabelian Hodge theory in \cite{Simpson}(vi). Of
course they were all well-known before that.

The deformation between de Rham and Dolbeault homotopy type in the simply
connected (compact complex) case was defined by Neisendorfer and Taylor
\cite{NeisendorferTaylor}. On the level of cohomology, the interpretation of the
Hodge filtration in terms of this deformation was basically the subject of
the article of Deninger \cite{Deninger}. For degree $1$ nonabelian cohomology,
this deformation is Deligne's space of {\em $\lambda$-connections} cf
\cite{Simpson}(i) and (vi).

Throughout this section the site in question is that of noetherian schemes over
a field $k$ of characteristic zero, with the etale topology. We fix a value of
$n$ and look at $n$-stacks.

\subnumero{De Rham theory}
Our first task is to construct the formal groupoid $X\rightarrow X_{DR}$.
Suppose $X$ is a smooth projective variety.
Define the $0$-stack $X_{DR}$ by the formula
$$
X_{DR}(Y) = X(Y^{\rm red}).
$$
It is easy to see that this is represented by the formal groupoid $(X,N_{DR})$
where $N_{DR}$ is the formal completion of the diagonal in $X\times X$.
Note that $N_{DR} \times _X N_{DR} = N_{DR}$ so the identity map serves as
multiplication. The fact that $X$ is smooth means that the ideal $I_{DR}$
defining the diagonal in $N_{DR}$ has $I/I^2 = \Omega ^1_X$ locally free
(note that the diagonal is $e(X)$). The de Rham complex as defined above is thus
of the form
$$
\Oo _X \stackrel{\delta _{DR}}{\rightarrow} \Omega ^1_X
\stackrel{\delta _{DR}}{\rightarrow} \Omega ^2_X
\stackrel{\delta _{DR}}{\rightarrow} \ldots .
$$
The differential $\delta _{DR}$ is equal to the usual de Rham
differential $d$ (\cite{Illusie} \cite{Berthelot}).

If $T$ is an $n$-stack then we define the {\em nonabelian de Rham cohomology of
$X$ with coefficients in $T$} to be
$$
\underline{Hom}(X_{DR}, T).
$$
Here of course, as remarked in \S 3 above, $T$ is assumed to be fibrant.

For the purposes of the paragraphs that follow we refer to the coefficient realm
$\underline{VP}$ of very presentable $n$-stacks. One could insert other realms
from \ref{index1} instead.

If $X^{\rm top}$ is simply connected then the de Rham shape
$$
{\bf Shape}_{\underline{VP}}(X_{DR}): T\mapsto \underline{Hom}(X_{DR}, T)
$$
is representable (cf the more general discussion in \S 10 below), and we
can think
of the representing $n$-stack  ${\bf rep}(X_{DR})$ as being the
``de Rham homotopy type'' of $X$. This object is a simply connected very
presentable $n$-stack, so the  higher homotopy sheaves
$$
\pi _i({\bf rep}(X_{DR}))
$$
are represented by finite dimensional vector spaces over the base field $k$.
We can call these vector spaces the ``de Rham homotopy groups of $X/k$''.

In the simply connected case,
Navarro Aznar \cite{Navarro-Aznar} and Wojtkowiak \cite{Wojtkowiak} had
independently and with methods completely different from the present, defined
``de Rham homotopy groups'' of $X/k$ which are
$k$-vector spaces. We conjecture that there is a natural isomorphism between
the vector spaces we have constructed above, and those constructed in
\cite{Navarro-Aznar} and \cite{Wojtkowiak}. (In fact I am not even sure if
anyone has verified the compatibility between these latter two constructions.)
Over $\cc$, the de Rham homotopy groups are identified with the complexified
rational homotopy groups, so in this case one obtains all of the compatibilities
in question.

If $X$ is not simply connected, the de Rham shape
${\bf Shape}_{\underline{VP}}(X_{DR})$ will not in general be representable, and
the shape itself seems to be the only object which could justifiably be called
the ``de Rham homotopy type'' of $X/k$.

\subnumero{Dolbeault theory}

Suppose again that $X$ is a smooth projective variety.  Put
$$
X_{Dol}:= K(\widehat{TX}/X, 1)\rightarrow X
$$
equal to the classifying $1$-stack for the sheaf of groups $\widehat{TX}$ over
$X$ (which means the formal completion of the tangent bundle $TX$ along the
zero-section). We have
$$
\Omega ^1_{X/X_{Dol}} = \Omega ^1_X
$$
so the terms in the de Rham complex for $X_{Dol}$ are again $\Omega ^i_X$,
but now the differential is equal to zero. Thus the de Rham complex for
$X_{Dol}$ is what is habitually called the {\em (algebraic) Dolbeault
complex} for
$X$.

The case of nonabelian Dolbeault cohomology was discussed at some length in
\cite{Simpson}(xii) and we refer the reader there.

If $X$ is simply connected, then the representing $n$-stack ${\bf
rep}(X_{Dol})$ exists, and it is a $1$-connected very
presentable $n$-stack. In particular, the homotopy groups
$\pi _i({\bf rep}(X_{Dol}))$ are $k$-vector spaces.

For $k=\cc$ the {\em
Dolbeault homotopy groups} of $X$ were defined
By Neisendorfer and Taylor \cite{NeisendorferTaylor}. Again, there is a
compatibility question to show that their definition is the same as ours.

\subnumero{Relative de Rham theory: the Gauss-Manin connection}

Now suppose that $X\rightarrow S$ is a smooth morphism of smooth varieties.
We obtain a commutative square
$$
\begin{array}{ccc}
X & \rightarrow & X_{DR} \\
\downarrow && \downarrow \\
S & \rightarrow & S_{DR} .
\end{array}
$$
We define the $0$-stack
$$
X_{DR/S} := X_{DR} \times _{S_{DR}} S.
$$
We claim that this is again of smooth type. In fact, it comes from the formal
category $N_{DR /S}$ which is the formal completion of the diagonal
in $X\times _SX$. The fact that the map is smooth means that this fiber product
is smooth. The ideal defining the diagonal has
$$
I_{DR/S} /I^2_{DR/S} = \Omega ^1_{X/S}.
$$
By looking at the morphism $N_{DR/S} \rightarrow N_{DR}$ we calculate that the
differential of the de Rham comples of $N_{DR/S}$ is again the standard de Rham
differential of the relative de Rham complex $\Omega ^{\cdot}_{X/S}$.

If $T$ is a very presentable $n$-stack on $Sch /\cc$, then we obtain the {\em
relative de Rham cohomology}
$$
\underline{Hom}(X_{DR}/S_{DR}, T)\rightarrow S_{DR}.
$$
Using the interpretation of $X_{DR}\rightarrow S_{DR}$ as a cartesian family
over $S_{DR}$ (\ref{correspondence}), the functor $T\mapsto
\underline{Hom}(X_{DR}/S_{DR}, T)$ may be viewed as a morphism of $n+1$-stacks
$$
{\bf Shape}_{\underline{VP}}(X_{DR}/S_{DR}): S_{DR} \rightarrow
\underline{Hom}(\underline{VP} , \underline{VP}).
$$
The pullback to $S$ is the family of nonabelian de Rham cohomology functors
of the
fibers, and the fact that this morphism descends to $S_{DR}$ is the {\em
Gauss-Manin connection}.

If the fibers are simply connected then the very presentable shape is
representable, so we obtain an $n$-stack
$$
{\bf rep}(X_{DR}/S_{DR})\rightarrow S_{DR}
$$
representing the shape. Again, the fact that this is defined over $S_{DR}$ is
the Gauss-Manin connection on the de Rham homotopy type. In particular, taking
homotopy groups we obtain vector bundles with integrable connection
$$
\pi _i ({\bf rep}(X_{DR}/S_{DR})/S_{DR})\rightarrow S_{DR}
$$
(the fact that these are vector bundles follows from \ref{last2} below, using
the comparison between de Rham homotopy groups and complexified rational
homotopy groups cf \cite{Simpson}(v) to verify that the dimensions are
constant).

These are the de Rham homotopy groups with their integrable connections. We
recover the construction of Navarro-Aznar \cite{Navarro-Aznar}, up to the  fact
that we have not verified that our construction is the same as his. (Over $\cc$,
this compatibility follows {\em a posteriori} from the fact that these
connections
have regular singularities, since a given flat bundle---in this case, the bundle
of complexified rational homotopy groups---has a unique algebraic structure with
regular singularities. However, this type of argument doesn't show that the
structures of definition over a smaller field $k\subset \cc $, which exist here
as well as in \cite{Navarro-Aznar}, are the same.)

\subnumero{Logarithmic de Rham theory}

Suppose first that $X$ is a smooth curve, and let $D\subset X$ be a divisor
consisting of a collection of isolated points. Then we define
$$
X_{DR}(\log D):= BL_D(X_{DR}).
$$
According to \ref{bld6}, the de Rham
complex in this case is given by
$$
\Omega ^1_{X/X_{DR}(\log D)}= \Omega ^1_X(D).
$$

Now we use the curve case, direct products and glueing to define
$X_{DR}(\log D)$
for any smooth variety $X$ and divisor $D$ with normal crossings.
We work in the etale topology. If $P$ is  a point of $X$,
we will  construct $X_{DR}(\log D)$ in an etale neighborhood of $P$; it will be
clear that these glue together to give $X_{DR}(\log D)$ as a formal algebraic
space. If $P\not\in D$ then $X_{DR}(\log D)=X_{DR}$ even in a Zariski
neighborhood of $P$. Suppose $P$ is a point contained in $a$ components of the
divisor $D$. For this we first treat the case $Y= {\bf A}^n$, and $D$ is the
union of the first $a$ coordinate hyperplanes. Let
$$
(Y_i,D_i)= ({\bf A}^1, 0)
$$
for $i=1,\ldots , a$ and let $Y_i={\bf A}^1$ with $D_i$ empty for $i=a+1,\ldots ,
n$. Set $Y=Y_1\times \ldots \times Y_n$, and define
$$
Y_{DR}(\log D):= Y_{1,DR}(\log D_1)\times \ldots \times Y_{a,DR}(\log D_a)
\times Y_{a+1,DR} \times \ldots \times Y_{n,DR}.
$$
For a general $P\in X$, there is an etale neighborhood $(X',D')$ of $P$
which is
also an etale neigborhood of $0\in Y$, and the above $Y_{DR}(\log D)$ pull back
to give a formal category on $X'$ which we can denote
by $X'_{DR}(\log D')$. These pieces
glue together (using etale glueing) to give $X_{DR}(\log D)$. The de
Rham complex for this formal category is the usual logarithmic de Rham complex.
The tangent bundle $T(X/X_{DR}(\log D))$ is the subsheaf of $T(X)$ of tangent
vector fields which are tangent to $D$.

\subnumero{Relative logarithmic de Rham theory}

Suppose $X$ and $Y$ are smooth, and $Y\subset X$ and $D\subset S$ are divisors
with normal crossings, and $f: X\rightarrow S$ is a morphism. We say that
{\em $f:
(X,Y)\rightarrow (S,D)$ is a relative normal crossings morphism} if
$f$ is flat, $f^{-1}(D)$ is contained in $Y$, all components of $Y$ which
don't surject onto $S$ are contained in $f^{-1}(D)$, and if there exist  local
holomorphic coordinates $x_i$ upstairs and $s_j$ downstairs such that the
divisors are defined by coordinate functions and such that the map $X\rightarrow
S$ is has coefficients $s_j$ which are  monomials in the $x_i$.

\begin{lemma}
\label{relnormcross}
If $f$ is a relative normal crossings morphism, then
the induced morphism of formal categories
$$
X_{DR}(\log Y) \rightarrow S_{DR}(\log D)
$$
is a morphism of smooth type.
\end{lemma}
{\em Proof:}
If we write
$$
s_j = \prod _i x_i^{a_{ij}}
$$
then the induced morphism on differentials is
$$
\Omega ^1_S(\log D)|_X\rightarrow \Omega ^1_X(\log Y)
$$
given by
$$
d\log s_j \mapsto \sum a_{ij} d\log x_i .
$$
This is a strict map of vector bundles because the matrix $a_{ij}$ is of
maximal rank (otherwise the morphism $X\rightarrow S$ wouldn't be onto and in
particular wouldn't be flat).
\eop

We call the relative nonabelian cohomology
$$
\underline{Hom}(X_{DR}(\log Y) / S_{DR}(\log D) , T)
$$
the {\em relative logarithmic de Rham cohomology}, and we call the shape
$$
{\bf Shape} (X_{DR}(\log Y) / S_{DR}(\log D)):
S_{DR}(\log D)\rightarrow \underline{Hom}(n\underline{STACK},
n\underline{STACK})
$$
the {\em relative logarithmic de Rham shape of $(X,Y) /(S,D)$}.

Note that all of our previous situations (except the Dolbeault example)
can be viewed as particular examples, when some of the divisors involved are
empty, when $S=\ast$, etc.

\subnumero{The Hodge filtration}

We can apply the construction $\Ff \mapsto {\bf Hodge}(\Ff )$ to all of the
formal
categories defined above. In fact, it also applies to all of the morphisms of
formal categories which we have defined and which are of  smooth type, by
Corollary \ref{hodgesmooth}.

We recapitulate what this means in a complete situation encompassing many of the
above examples. Suppose  $f: (X, Y)\rightarrow (S,D)$ is a relative normal
crossings morphism. Then we obtain the morphism
$$
{\bf Hodge}(X_{DR}(\log Y))\rightarrow {\bf Hodge}(S_{DR}(\log D))
$$
over $\af$. It is a morphism of smooth type, by Lemma \ref{hodgesmooth}.

The relative nonabelian cohomology
$$
\underline{Hom}({\bf Hodge}(X_{DR}(\log Y))/ {\bf Hodge}(S_{DR}(\log D)),
T)
$$
with coefficients in some $n$-stack $T$ (exactly which realm we take for the
coefficients $T$ will become more clear in \S 10 below), contains all at once
the Gauss-Manin connection with its regular singularities and with Griffiths
transversality for the Hodge filtration. Restricting over $0\in \af$ yields the
Kodaira-Spencer classes for deformation of the nonabelian Dolbeault cohomology.
Restricting over a singular point $s\in D$ yields the ``de Rham monodromy
action'' on the ``nearby cohomology''.

If we let the coefficient $n$-stack $T$ vary, then the above fits into the
relative shape  map. Again, all of the above information (Hodge filtration,
Gauss-Manin connection, regularity, Kodaira-Spencer classes etc.) for the de Rham
shape is contained in the relative shape map
$$
{\bf Shape}({\bf Hodge}(X_{DR}(\log Y))/ {\bf Hodge}(S_{DR}(\log D)):
$$
$$
{\bf
Hodge}(S_{DR}(\log D) \rightarrow  \underline{Hom}(n\underline{STACK},
n\underline{STACK}').
$$
If $\underline{R}$ is a realm, then we obtain similarly the relative
$\underline{R}$-shape map
$$
{\bf Shape}_{\underline{R}}({\bf Hodge}(X_{DR}(\log Y))/ {\bf
Hodge}(S_{DR}(\log D)),
$$
and if $\underline{A}$ is another realm such that the
relative $\underline{R}$-shape takes answers in $\underline{A}$ then we obtain
$$
{\bf Shape}_{\underline{R}}^{\underline{A}}({\bf Hodge}(X_{DR}(\log Y))/ {\bf
Hodge}(S_{DR}(\log D)) :
$$
$$
{\bf Hodge}(S_{DR}(\log D) \rightarrow
\underline{Hom}(\underline{R}, \underline{A}').
$$

In \S 10 we will give a few examples of pairs of realms $(\underline{R},
\underline{A})$ such that this shape map exists. For example (see \ref{index1}
for the notations of realms) we get
$$
{\bf Shape}_{\underline{FV}}^{\underline{VP}^{\rm loc}}({\bf Hodge}(X_{DR}(\log
Y))/ {\bf Hodge}(S_{DR}(\log D)) :
$$
$$
{\bf Hodge}(S_{DR}(\log D) \rightarrow
\underline{Hom}(\underline{FV}, \underline{VP}^{\rm loc})
$$
(Theorem \ref{fv}), and
$$
{\bf Shape}_{\underline{FL}}^{\underline{VG}^{\rm loc}}({\bf Hodge}(X_{DR}(\log
Y))/ {\bf Hodge}(S_{DR}(\log D)) :
$$
$$
{\bf Hodge}(S_{DR}(\log D) \rightarrow
\underline{Hom}(\underline{FL}, \underline{VG}^{\rm loc})
$$
(Theorem \ref{fl}).

In the case where the base is a point, the above is just the ``Hodge
filtration''.
(In the introduction, as well as in \cite{Simpson}(vi), ${\bf
Hodge}(X_{DR})$ was
denoted by $X_{Hod}$.)
The Hodge filtration is $\Gm$-equivariant so
it gives a relative shape map defined on $\af /\Gm$. Thus, for example, if $X$
is a smooth projective variety then we obtain the {\em Hodge filtration for the
de Rham shape of $X$}
$$
{\bf Shape}_{\underline{FV}}^{\underline{VP}^{\rm loc}}({\bf
Hodge}(X_{DR})/\af )
: \af /\Gm \rightarrow
\underline{Hom}(\underline{FV}, \underline{VP}^{\rm loc}),
$$
$$
{\bf Shape}_{\underline{FL}}^{\underline{VG}^{\rm loc}}({\bf
Hodge}(X_{DR})/\af )
: \af /\Gm \rightarrow
\underline{Hom}(\underline{FL}, \underline{VG}^{\rm loc}).
$$
The image of the point $[1]\in \af /\Gm$ is the shape of $X_{DR}$, whereas the
image of the point $[0]\in \af /\Gm$ is the shape of $X_{Dol}$.

Suppose that $X\rightarrow S$ is a smooth projective morphism with simply
connected fibers. Then the relative $\underline{AV}$-shape of ${\bf
Hodge}(X_{DR})\rightarrow {\bf Hodge}(S_{DR})$ is representable, cf \S 6.4.
We obtain a representing $n$-stack
$$
{\bf rep} ({\bf Hodge}(X_{DR})/{\bf Hodge}(S_{DR}))
\rightarrow {\bf Hodge}(S_{DR}).
$$
As we shall see in Corollary \ref{last1} using the criterion of geometricity
\ref{criterion}, this representing $n$-stack is a $1$-connected very presentable
geometric $n$-stack over ${\bf Hodge}(S_{DR})$.

Again when $S$ is a point, the representing $n$-stack (which is
$\Gm$-equivariant) is a $1$-connected very presentable geometric $n$-stack
$$
{\bf rep} ({\bf Hodge}(X_{DR})/\af )/\Gm \rightarrow \af /\Gm .
$$
This object is the {\em Hodge filtration on the complex homotopy type of
$X_{DR}$}. This object existed long ago, essentially as soon as the rational
homotopy type of $X$ was calculated using differential forms (because one has a
``Hodge filtration'' on the associated d.g.a.). However, it was not clear in the
d.g.a. point of view, exactly what type of object it was. The above description
gives a homotopically correct description of what type of object the Hodge
filtration is.

\numero{Presentability and geometricity results}
\label{resultspage}

We now give our main results about nonabelian cohomology of formal categories.
These results are phrased in a general way. They apply directly to all of the
formal categories constructed in the previous section, to give concrete results
concerning nonabelian Hodge theory, as was illustrated with a few statements
at the end of the previous section.

The site in question in this section is $\Gg = Sch /k$ the site of noetherian
schemes over a field $k$ of characteristic zero, with the etale topology.
Where necessary, we fix a positive integer $n$.

\subnumero{Notations}

\begin{hypothesis}
\label{situation}
Throughout this section, we will be considering the following situation: let
$$
\begin{array}{ccc}
X&\stackrel{p}{\rightarrow} &\Ff\\
\downarrow && \downarrow {\scriptstyle f} \\
S&\stackrel{q}{\rightarrow} &\Ee
\end{array}
$$
be a morphism of smooth type between two formal categories of smooth type.
We assume that the morphism $X\rightarrow S$ is projective.
This assumption is in effect throughout the present section unless specified
otherwise.
\end{hypothesis}

\begin{parag}
\label{meaning}
We will be looking at realms $\underline{R}$ and $\underline{A}$ and we would
like to show that
$$
{\bf Shape} _{\underline{R}}^{\underline{A}}(\Ff /\Ee )
$$
exists, i.e. that the $\underline{R}$-shape of $\Ff /\Ee$
takes values in $\underline{A}$, i.e. that
$$
{\bf Shape} _{\underline{R}}(\Ff /\Ee ) : \Ee \rightarrow \underline{Hom}
(\underline{R}, \underline{A}).
$$
For shorthand, we write this condition as
$$
\underline{R}
\stackrel{{\bf Shape} (\Ff /\Ee )}{\longrightarrow} \underline{A}
$$
Recall from \ref{relativezoo4}, \ref{relativezoo5} the meaning of this
statement. Consider $\Ff /\Ee$ as corresponding to a cartesian family,
coming from a morphism $[f ] : \Ee \rightarrow 1^{\rm gr}
\underline{STACK}$. The
arrow family composed with the map $[f]$ is
$$
{\bf Arr} \circ [f ]: \Ee ^o\times  \underline{R}\rightarrow \ngr
\underline{STACK}
$$
and we would like the image to lie in $\underline{A}$. If $Z\in \Gg$ then
it is a
map  $$
\Ee ^o(Z) \times \underline{R}(Z) \rightarrow \ngr \underline{STACK}(Z)
$$
defined as follows: for $T\in \underline{R}(Z)$ and $\varphi \in \Ee ^o (Z)$
we take the pullback family
$$
\varphi ^{\ast}(\Ff )= \Ff \times _{\Ee} Z
$$
(considered as a stack on $\Gg /Z$)
and look at
$$
\underline{Hom}(\varphi ^{\ast}(\Ff )/Z, T/Z) \in \ngr \underline{STACK}(Z).
$$
The problem is to see whether this $\ngr$-stack on $\Gg /Z$ is in
$\underline{A}(Z)$.
\end{parag}

We recapitulate the previous paragraph in the following proposition:

\begin{proposition}
\label{meaningprop}
Suppose $\underline{R}$ and $\underline{A}$ are realms, and keep our
notations of
\ref{situation}. Suppose that for every $Z\in \Gg$ with morphism $\varphi : Z
\rightarrow \Ee$, and for every $T\in \underline{R}(Z)$, the nonabelian
cohomology stack on $\Gg /Z$
$$
\underline{Hom}(\Ff \times _{\Ee}Z/Z, T/Z) \in \ngr \underline{STACK}(Z),
$$
lies in $\underline{A}(Z)$.
Then
$$
\underline{R}
\stackrel{{\bf Shape} (\Ff /\Ee )}{\longrightarrow} \underline{A}.
$$
\end{proposition}
\eop

\begin{parag}
\label{cases}
We will treat the cases where $\underline{R}$ is one of the following realms
(cf \ref{index1}):
$$
\underline{FL}\subset
\underline{FV},
$$
as well as the cases of the realms of simply connected stacks
$$
\underline{AL}\subset
\underline{AG}\subset
\underline{AV}.
$$
The answer realms for ${\bf Shape}(\Ff /S)$ will be
(cf Theorems \ref{fv} and \ref{fl} below)
$$
\underline{FV}\stackrel{{\bf Shape} (\Ff /\Ee )}{\longrightarrow}
\underline{VP}^{\rm loc},
$$
$$
\underline{FL}\stackrel{{\bf Shape} (\Ff /\Ee )}{\longrightarrow}
\underline{VG}^{\rm loc}.
$$
\end{parag}

The fact that the answers may not be of finite type
comes from the degree $1$ cohomology. Thus for the realms of simply connected
stacks we obtain slightly better statements:
$$
\underline{AV}\stackrel{{\bf Shape} (\Ff /\Ee )}{\longrightarrow}
\underline{VP},
$$
$$
\underline{AG}\stackrel{{\bf Shape} (\Ff /\Ee )}{\longrightarrow}
\underline{VG},
$$
$$
\underline{AL}\stackrel{{\bf Shape} (\Ff /\Ee )}{\longrightarrow}
\underline{VG}.
$$
We treat the proofs of these first, before recalling what is necessary about
degree $1$ cohomology.

Another remark to be made here is that in the case of simply connected
coefficients we treat $\underline{AG}$. This treatment is somewhat complicated
and we only give a sketch of proof of Theorem \ref{ag} below. It should be
possible to do the same thing for $\underline{FG}$ but we haven't included this
statement.

\begin{parag}
{\em Caution:} throughout this section, the term ``local system'' is used. See
\ref{localsystem}; this only refers to the possible action of $\pi _1$ of
the base
stack.   Note that $\pi _1$ here means e.g. the local automorphism groups
in a $1$-stack; it does {\em not} mean the topological $\pi _1(X^{\rm top})$.
Thus, for
example when the base is a $0$-stack (a scheme for example), a ``local system''
just means a sheaf.
This notion of ``local system'' should not be confused with
the differential-geometric notion of ``flat vector bundle'' over a space $X^{\rm
top}$ (which was used in \S 2).
\end{parag}

\subnumero{The abelian cohomology problem}

In the present paper we don't want to treat the problem of degree $1$ cohomology
with coefficients in a general affine presentable group sheaf (which is why we
restrict to e.g. $\underline{FV}$ rather than $\underline{CV}$), but
for future reference we would
like the treatment of higher-degree cohomology to be applicable even in the
presence of degree $1$ cohomology. Thus we will consider a morphism
$T\rightarrow
R$ (typically $R= \tau _{\leq 1}(T/S)$) and look at the morphism
$$
\underline{Hom}(\Ff /S, T/S)\rightarrow \underline{Hom}(\Ff /S, R/S).
$$
The question is when this morphism will be of type $\underline{A}$.

\begin{parag}
One important thing to point out is that it is impossible to have reasonable
results of the form we are looking for, for the stack of sections
$$
\underline{\Gamma} (\Ff /S, T)
$$
where $T\rightarrow \Ff$. Even in the case where $\Ff = X$ is a variety
projective over $S$, and $T= K(V/X, m)$ for some coherent sheaf $V$ on $X$,
the stack of sections (which is basically the higher derived direct image
complex ${\bf R} f_{\ast} (V)$ on $S$) doesn't have nice properties. For
example its cohomology sheaves on $S$ will not in general be vector sheaves
(Hirschowitz includes a counterexample which he ascribes to O. Gabber in
\cite{Hirschowitz}). Thus, it is crucial for the results of the present section
that the coefficient stack $T$ comes from a stack over $S$. In the somewhat
relative situation $T\rightarrow R \rightarrow S$ we will have to consider
sections over a morphism $\Ff \rightarrow R$, so for the same type of reason it
will be essential here that $R$ be $0$-connected relative to $S$.
\end{parag}

\begin{parag}
Both ``answer'' realms
$\underline{VG}$ and
$\underline{VP}$ are closed under extensions and finite limits.
Thus, when looking at a morphism $T\rightarrow R$ which is relatively
$1$-connected we can apply Proposition \ref{postniprop} (see also Corollary
\ref{postnicor}). When starting with realms of coefficients $\underline{FL}$ or
$\underline{FV}$, the Eilenberg-MacLane sheaves which enter in here are of the
form $K(L/R,m)$ for $L$ a local system of vector bundles (resp. local system of
vector sheaves) on $R$. The fact that $R$ is $0$-connected relative to $S$ means
that the pullback of $L$ via  $\eta :\Ff \rightarrow R$ will be locally in the
etale topology of $X$, a pullback of a vector sheaf from $S$.
\end{parag}

\begin{parag}
For the results starting with the realm $\underline{FG}$, which is not closed
under truncation, we cannot apply directly Proposition \ref{postniprop}. Instead
we apply the technique of \S 7.3 to decompose the coefficient
stack $T$ into a successive extension of $K(C^{\cdot}, m)$ for perfect
complexes $C^{\cdot}$; then these are themselves successive extensions of
$K(L,m)$ for vector bundles $L$. Thus in this case too, we are reduced to
considering Eilenberg-MacLane sheaves corresponding to local systems of
vector bundles $L$ over the base $R$. See the sketch of proof of Theorem
\ref{ag}
below.
\end{parag}

\begin{parag}
We are reduced to the following problem of abelian cohomology. Let $L= \pi
_k(T/R)$ as a local system on $R$. In all cases we want to treat, it is a local
system of vector sheaves on $R$. We suppose given an $S$-morphism $\eta : \Ff
\rightarrow R$ and then we have to look at the cohomology complex
$$
{\bf R} f_{\ast} (\eta ^{\ast} L)
$$
on $S$. The problem is to see what properties this has.

There are two cases: (1) where $L$ is locally free; then we would like the
higher direct image to be a perfect complex; and (2) where $L$ is a local system
of vector sheaves, in this case we would like the cohomology sheaves on $S$ to
be vector sheaves.
\end{parag}

So we get to the formulation of the following main statement.

\begin{theorem}
\label{mainab}
Suppose $f: \Ff \rightarrow S$ is a formal category over a base scheme $S$, with
$X\rightarrow S$ flat and projective as in \ref{situation}. Suppose
$R\rightarrow
S$ is a relatively $0$-connected $\ngr$-stack. Suppose $L$ is a local system of
vector sheaves on $R$. Suppose $\eta : \Ff \rightarrow R$ is an $S$-morphism.
Then the cohomology sheaves $$
H^i(\Ff /S, \eta ^{\ast} L)
$$
are vector sheaves on $S$. Furthermore, if $L$ is locally free then the higher
direct image complex
$$
{\bf R} f_{\ast} (\eta ^{\ast} L)
$$
is a perfect complex on $S$.
\end{theorem}
{\em Proof:}
Let $p: X\rightarrow \Ff$ be the projection from the underlying scheme to the
formal category $\Ff$. According to the discussion in the previous chapter, the
higher direct image complex ${\bf R} f_{\ast} (\eta ^{\ast} L)$ is
calculated by
the de Rham complex for $X/\Ff$. To be precise, let $Z\rightarrow S$ be a scheme
mapping to $S$, and let $X_Z^{\rm zar}$ denote the Zariski site of $X_Z:=
X\times _SZ$.

The de Rham complex with coefficients in $L$ pulls back to the de Rham complex
for $\Ff _Z:= \Ff \times _SZ$, which we denote by  $\Omega ^{\cdot} _{X_Z/\Ff
_Z}(\eta ^{\ast} _ZL)$ and which is a complex of sheaves on  $X_Z^{\rm zar}$.

Functorially in $Z$, the value of the higher direct image complex
is identified as the hypercohomology of $X_Z^{\rm zar}$
with coefficients in the de Rham complex:
$$
{\bf R} f_{\ast} (\eta ^{\ast} L)(Z) =
{\bf H} (X_Z^{\rm zar}, \Omega ^{\cdot} _{X_Z/\Ff
_Z}(\eta ^{\ast} _ZL)).
$$

The Hodge-to-de Rham spectral sequence, which is the same as the spectral
sequence for the de Rham complex with the stupid filtration, is
$$
H^i(X_Z^{\rm zar}, \Omega ^j _{X_Z/\Ff
_Z}(\eta ^{\ast} _ZL)) \Rightarrow {\bf H}^{i+j}(X_Z^{\rm zar}, \Omega
^{\cdot} _{X_Z/\Ff
_Z}(\eta ^{\ast} _ZL)).
$$
By functoriality in $Z$ this gives a spectral sequence $(\ast )$
$$
H^i(X/S, \Omega ^j _{X/\Ff }(\eta ^{\ast} L)) \Rightarrow
{\bf R} ^{i+j}f_{\ast} (\eta ^{\ast} L).
$$

An alternative method of obtaining the spectral sequence
$(\ast )$ is to look at the \v{C}ech-Alexander argument from the previous
chapter and redo the same for the higher direct image
${\bf R} f_{\ast}$ rather than for ${\bf zar} _{\ast}$.  The spectral sequence
resulting from the filtration given in the previous chapter for the
\v{C}ech-Alexander complex, is then exactly $(\ast )$.

A slightly more precise statement is to say that the ``Hodge filtration''
upstairs (e.g. the stupid filtration on the de Rham complex or alternatively the
filtration defined in the previous chapter on the \v{C}ech-Alexander complex)
yields a filtered complex $({\bf R} f_{\ast} (\eta ^{\ast} L), F^{\cdot})$
and the associated-graded is
$$
Gr ^j_F {\bf R} f_{\ast} (\eta ^{\ast} L) \cong
{\bf R} (pf)_{\ast}(\Omega ^j _{X/\Ff }(\eta ^{\ast} L)).
$$
Finally note that the component sheaves are identified by
$$
\Omega ^j _{X/\Ff }(\eta ^{\ast} L)=
(\eta ^{\ast} L|_X)\otimes _{\Oo _X}\Omega ^j _{X/\Ff }.
$$
To interpret this formula note that the $\Omega ^j _{X/\Ff }$
are locally free sheaves of $\Oo _X$-modules; and  $\eta ^{\ast} L|_X$ denotes
the pullback of $\eta ^{\ast} L$, which is a local system on $\Ff$, to $X$.

We can now prove the last statement of the theorem. If $L$ is locally free then
the component sheaves  $\Omega ^j _{X/\Ff }(\eta ^{\ast} L)$
are locally free sheaves on $X$. Included in the hypothesis that
$(X,\Ff )\rightarrow S$ is of smooth type, is the hypothesis that $X\rightarrow
S$ is flat. Therefore (by Mumford's argument \cite{Mumford}) the higher direct
images ${\bf R} (pf)_{\ast}(\Omega ^j _{X/\Ff }(\eta ^{\ast} L))$
are perfect complexes. Now
${\bf R} f_{\ast} (\eta ^{\ast} L)$ is a successive extension of these perfect
complexes, so it again is a perfect complex. This gives the last statement of
the theorem.

Suppose now that $R$ is relatively $0$-connected, that $L$ is a local system of
vector sheaves on $R$, and that we pull back by an $S$-map $\Ff \rightarrow R$.
Then
$$
\eta ^{\ast} L|_X = (\eta p)^{\ast} (L).
$$
The fact that $R\rightarrow S$ is relatively $0$-connected means that by going
to an etale covering of $S$ (and we can replace $S$ by this etale covering
without affecting the argument) we can assume that there is a section
$\rho : S\rightarrow R$. Furthermore we can assume that the morphism
$$
\eta p : X\rightarrow R
$$
is equal to $\rho fp$ over an etale covering $X'\rightarrow X$.
Thus we are reduced to proving the following statement:

\begin{proposition}
\label{statement}
If $X\rightarrow S$ is a
projective flat morphism, if $V$ is a vector sheaf on $S$  and if $L$ is a
local system over $X$  which is locally in the etale topology of $X$ equivalent
to the pullback of $V$, then the $H^i(X/S, L)$ are vector sheaves on $S$.
\end{proposition}

This
statement will be the subject of several lemmas which we now give; modulo the
proof of Proposition \ref{statement}, this completes the proof of the theorem.
\eop

\oldsubnumero{Lemmas for Proposition \ref{statement}}

\begin{hypotheses}
Fix now the following situation:  $f: X\rightarrow S$ is a projective flat
morphism of relative dimension $d$; and $V$ is a vector sheaf on $S$. Let
$$
\eta \in H^1(X^{\rm et}, Aut(V)|_X)
$$
be a cocycle on $X^{\rm et}$ with coefficients in the sheaf $Aut(V)$.
Let $L\rightarrow X$ be the associated local system (i.e. sheaf) which is a
vector
sheaf over $X$, locally (on $X^{\rm et}$) isomorphic to the pullback of $V$.
\end{hypotheses}

With these hypotheses, we would like to show that the $H^i(X/S,L)$ are vector
sheaves on $S$. The following lemma shows that it suffices to prove this for the
case of projective space over $S$.

\begin{lemma}
\label{reduction}
{\rm (reduction to the case of $\pp ^d$)}
Choose a finite flat projection
$$
p: X\rightarrow  Y:= \pp ^d_S.
$$
Let $r$ be the
degree of $p$. Then locally in the etale topology of
$Y$, $p_{\ast}(L)$ is isomorphic to the pullback of $V^{\oplus r}$. Furthermore,
$$
H^i(X/S, L)= H^i(Y/S, p_{\ast}(L)).
$$
\end{lemma}
{\em Proof:}
Choose a point $y\in Y$, and denote by $\tilde{Y}_y$ the henselization of $Y$ at
$y$. The inverse image  decomposes as a disjoint union
$$
p^{-1}(\tilde{Y}_y)=\bigcup _i U_i
$$
where $U_i$ are henselian in $X$, and where each projection
$$
p_i= U_i \rightarrow \tilde{Y}_y
$$
is finite and flat of degree $r_i$; and $\sum _i r_i =r$. Over each $U_i$,
$L$ is isomorphic to the pullback of $V$. For clarity, let $g: Y\rightarrow S$
denote the projection, and denote again by $g$ or $f_i$ the projections
restricted
to $\tilde{Y}_y$ or $U_i$ respectively (so $f_i=gp_i$). Thus
$$
L|_{U_i} \cong f_i^{\ast}(V).
$$
We get
$$
p_{i,\ast} (L|_{U_i})\cong
p_{i,\ast}(f_i^{\ast}(V)) =
p_{i,\ast}p_i^{\ast}(g^{\ast} V)
$$
$$
= [p_{i,\ast}p_i^{\ast}(\Oo _{Y})] \otimes _{\Oo _{Y}}
g^{\ast} V
$$
(for typographical reasons we write $\Oo _Y$ instead of
$\Oo _{\tilde{Y}_y}$ in this formula). Now
$$
p_{i,\ast}p_i^{\ast}(\Oo _{Y})\cong \Oo _{\tilde{Y}_y}^{\oplus r_i}
$$
so we get
$$
p_{i,\ast} (L|_{U_i})\cong g^{\ast} (V^{\oplus r_i}).
$$
Putting these together over all $U_i$ we get
$$
p_{\ast}(L)|_{\tilde{Y}_y} = p_{\ast}(L|_{p^{-1}(\tilde{Y}_y)}) \cong g^{\ast}
(V^{\oplus r}).
$$
This shows that locally in the etale topology of $Y$,
$p_{\ast}(L)$ is isomorphic to the pullback of $V^{\oplus r}$.

The last statement of the lemma comes from the fact that $p$ is finite: the
higher direct images of $p$ vanish (one can see that this holds even for vector
sheaf coefficients).
\eop

From now on we assume $X=\pp ^d_S$. Recall that $f$ is the projection to $S$.
The next step is to make use of the hypothesis that $L$ is locally isomorphic
to the pullback $f^{\ast}(V)$. The first step is to pass from triviality over an
etale covering of $X$, to triviality on a Zariski open covering of $X$.

Let
$$
\Ee := \underline{Hom}(L,f^{\ast} V)
$$
be the local system of vector
sheaves $\underline{Hom}$ on $X$. Let $\Ee ^{\rm in}$ be the coherent sheaf
mapping to $\Ee$ with the property that the sections of $\Ee$ over any scheme
which is etale over $X$, are the same as the sections of $\Ee ^{\rm in}$. One
can note that $\Ee$ and $\Ee ^{\rm in}$ also have the property of being, locally
in the etale topology of $X$, pullbacks (respectively of a vector sheaf, of a
coherent sheaf) from $S$.

 Fix an
etale surjection $X'\rightarrow X$ such that $V|_{X'}\cong L|_{X'}$, and let  $$
\varphi \in H^0(X', \Ee |_{X'})
$$
be the section giving our trivialization.  Fix   a point $P\in X$ lifting to
$P'\in X'$.

\begin{lemma}
There is an integer $k_0$ such that if $k\geq k_0$ and if $\psi
\in H^0(X', \Ee |_{X'})$ is a section such that
$$
\psi |_{Spec (\Oo _{X'} /{\bf m} _{P'}^k)} =
\varphi |_{Spec (\Oo _{X'} /{\bf m} _{P'}^k)}
$$
then $\psi$ gives a trivialization $V|_{X''}\cong L|_{X''}$ over a Zariski open
neighborhood $P\in X'' \subset X'$.
\end{lemma}
{\em Proof:}
Composition of endomorphisms is a bilinear map
$$
End(V)\times End(V)\rightarrow End(V).
$$
Let $\Ff$ denote the coherent sheaf associated to the module of global sections
of $End(V)$. The above composition map induces a morphism
$$
\Ff \otimes _{\Oo}\Ff \rightarrow \Ff .
$$
We will find $k$ such that for $\psi$ as in the lemma, then
$$
\varepsilon := 1-\psi\varphi^{-1} \in H^0(X', {\bf m} _{P'} \Ff |_{X'}).
$$
This suffices to prove the lemma, because $\varepsilon $
then restricts to a nilpotent element over any artinian scheme centered at $P'$,
so $\psi \varphi^{-1}= 1 -\varepsilon$ is invertible over any artinian scheme;
therefore it is invertible in the formal neighborhood of $P'$ and hence by Artin
approximation invertible in some etale neighborhood.

To find the required $k$, choose an injection of vector sheaves
$$
i:End(V)\hookrightarrow \Ff '.
$$
This induces a morphism of coherent sheaves $j: \Ff \rightarrow \Ff '$. Note
that, while $j$ is not injective as a morphism of vector sheaves, it does induce
an injection on modules of global sections, in other words it is injective as a
morphism of coherent sheaves in usual sense (on the Zariski topology).
The condition
$$
\psi
|_{Spec (\Oo _{X'} /{\bf m} _{P'}^k)} =  \varphi |_{Spec (\Oo _{X'} /{\bf m}
_{P'}^k)}
$$
means that $\varepsilon$ restricts to $0$ on $Spec (\Oo _{X'} /{\bf m}
_{P'}^k)$. Since $i$ is an injection of vector sheaves, this condition is
equivalent to saying that $i(\varepsilon )$ restricts to $0$ on
$Spec (\Oo _{X'} /{\bf m}
_{P'}^k)$, or equivalently that $i(\varepsilon )\in {\bf m}^k_{P'}\Ff '$.

Thus we have reduced the problem to the following one (we now consider
$\varepsilon$ as a section of $\Ff$). Given a morphism of coherent sheaves
$$
j:\Ff \rightarrow \Ff '
$$
which is injective in the usual (Zariski) sense, we would like to find an
integer $k$ such that
$$
j(\varepsilon ) \in {\bf m}^k_{P'}\Ff ' \Rightarrow \varepsilon
\in {\bf m}_{P'}\Ff .
$$
The existence of such a $k$ is an easy consequence of Krull's lemma:
apply Krull to the module $\Ff ' / {\bf m}_{P'} j(\Ff )$. We obtain that
$$
\bigcap _k ({\bf m}^k_{P'}\Ff ' + {\bf m}_{P'} j(\Ff )) =
{\bf m}_{P'} j(\Ff ).
$$
Therefore,
$$
\bigcap _k [({\bf m}^k_{P'}\Ff ' + {\bf m}_{P'} j(\Ff ))\cap j(\Ff )] =
{\bf m}_{P'} j(\Ff ).
$$
But this latter intersection is an intersection of vector spaces  lying between
${\bf m}_{P'} j(\Ff )$ and $j(\Ff )$; since the dimension of the quotient vector
space $j(\Ff )/{\bf m}_{P'} j(\Ff )$ is finite, it follows that the intersection
stops at a finite stage, in other words there is $k_0$ such that for $k\geq k_0$
we have
$$
[({\bf m}^k_{P'}\Ff ' + {\bf m}_{P'} j(\Ff ))\cap j(\Ff )] \subset
{\bf m}_{P'} j(\Ff ) .
$$
This gives the desired property.
\eop

\begin{corollary}
\label{zariskitriv}
There is a Zariski open neighborhood $P\in U\subset X$ and a trivialization
$V|_{U}\cong L|_{U}$.
\end{corollary}
{\em Proof:}
Note that $\varphi$ is actually a section of the coherent sheaf $\Ee ^{\rm in}$
over $X'$. Thus it can be approximated by sections of $\Ee ^{\rm in}$ on $X$. In
other words, we can choose a Zariski open neighborhood $P\in U'
\subset X$ and a section $\psi \in H^0(U, \Ee |_U)$ approximating $\varphi$ to
order $k$ at $P'$, i.e. with $$ \psi |_{Spec (\Oo _{X} /{\bf m} _{P}^k)} =
\varphi |_{Spec (\Oo _{X'} /{\bf m} _{P'}^k)}
$$
Note here that
$$
Spec (\Oo _{X} /{\bf m} _{P}^k)\cong Spec (\Oo _{X'} /{\bf m} _{P'}^k)
$$
because $X'\rightarrow X$ is etale at $P'$. By the above lemma, $\psi$ restricts
to a trivialization on some  neighborhood $P\subset X'' \subset X'\cap
p^{-1}(U')$. Let $U$ be the image of $X''$ in $U'$; then $\psi$ gives a
trivialization over $U$.
\eop

This corollary says that $L$ is locally in the Zariski topology isomorphic to
$f^{\ast}(V)$. Continuing with the notation of this corollary, we can
assume that
$U= X-D$ where $D$ is a Cartier divisor relative to $S$. Our section $\psi$
may be
considered as a meromorphic section of a coherent sheaf $\Ee ^{\rm in}$ on $X$,
and as such it has poles of a finite order along $D$.

Similarly, the inverse of $\psi$ also may be considered as a meromorphic section
having poles of finite order along $D$.  Therefore, $L$ is defined by a cocycle
in the Zariski topology, with coefficients in $Aut(V)$, and having poles of a
finite order.

\begin{lemma}
\label{boundedness}
There is an integer $r$ and a bounded family $\Bb$ of vector bundles of rank
$\leq r$ such that the following holds. If $A\rightarrow S$ is any morphism from
an artinian scheme of finite type over $k$ and if $a: A\rightarrow Spec (k)$
denotes the projection (here $k$ is our base field which we assume to be
algebraically closed of characteristic zero), then
$$
a_{\ast} (L_A) \rightarrow \pp ^d_k
$$
is a vector bundle on $\pp ^d_k$ with the property that
it has a filtration whose subquotients are vector bundles of rank $\leq r$
members of the bounded family $\Bb$.
\end{lemma}
{\em Proof:}
For any $A$, the  bundle $a_{\ast}(L_A)$ is a twisted version of the vector
space $V(A)$, defined by a cocycle with coefficients in $Aut(V)(A)$. This
cocycle is the image of the cocycle with coefficients in $Aut(V)$ which defines
$L$. From the previous discussion, we may assume that this cocycle has component
functions which have divisors of poles of degree $\leq m_0$. This implies, for
example, that the bundle $a_{\ast}(L_A)$ and its dual, are generated by global
sections after twisting by $\Oo _X(m)$ for $m\geq m_0$.

To get the filtration in question, we will construct a filtration of the module
$V_A=V(A)$ which is invariant under $Aut(V)(A)$, and whose subquotients have
rank $\leq r$. This filtration then  induces a filtration of $a_{\ast}(L_A)$,
since $L_A$ is defined by a cocycle in $Aut(V)(A)$. Furthermore, the
cocycle is defined in the Zariski topology by functions whose divisors of poles
have degree less than $m$, so the same will hold for the cocycles defining
the subquotients in the filtration.
Note that the open covering for
the cocycle description is fixed independent of $A$; and the vector bundles of
rank $r$ described by cocycles having poles of degree less than $m$, with
respect to a fixed covering, form our bounded family $\Bb$ which contains the
subquotients in the filtrations.

We now  construct the required filtration of $V(A)$. Choose a sequence of ideals
$I_j\subset \Oo (A)$ with $I_1={\bf m}_A$ and $I_{j+1} \subset I_j$
such that $I_j/I_{j+1} \cong \Oo (A) /{\bf m}_A$ as an $\Oo (A)$-module.
Choose an injection
$$
0\rightarrow V \rightarrow \Ff
$$
of vector sheaves, with $\Ff$ a coherent sheaf (this injection is chosen
independently of $A$). Recall from \cite{Simpson}(vii) and (xii) that the
cotensor product $\otimes ^{\Oo}$ is left exact; but that on coherent sheaves it
coincides with the usual tensor product. Thus we obtain an injection
(of vector sheaves over $A$)
$$
0\rightarrow V|_A\otimes ^{\Oo} (\Oo _A/I_j)
\rightarrow \Ff |_A\otimes _{\Oo} (\Oo _A/I_j).
$$
Let $V_{A,j}$ denote the kernel of
$$
V(A) \rightarrow [V|_A\otimes ^{\Oo} (\Oo _A/I_j)] (A).
$$
Similarly, let
$\Ff _{A,j}$ denote the kernel of
$$
\Ff (A) \rightarrow [\Ff |_A\otimes _{\Oo} (\Oo _A/I_j)] (A).
$$
Note however that
$$
\Ff _{A,j} = I_j\cdot \Ff (A) \subset \Ff (A).
$$
We have
$$
V_{A,j} = V(A) \cap \Ff _{A,j} \subset \Ff (A).
$$
The subquotients of the filtration $V_{A,j}$ are subobjects of the subquotients
of the filtration $\Ff _{A,j}$; but this latter filtration has subquotients
which are of length $\leq r_A$ where $r_A= dim (\Ff / {\bf m}_A\Ff )$.
Let $r$ be the maximum of the $r_A$; it is the maximal dimension of a
fiber of $\Ff$, which is bounded independently of $A$. The subquotients
of the filtration $V_{A,j}$ on the module $V(A)$, are of length $\leq r$.
Finally
we note that the filtration $V_{A,j}$ is intrinsically defined, so it is
invariant under $Aut(V)(A)$. (This is in spite of the fact that the injection
$V\rightarrow \Ff$ is not intrinsic). Thus we have constructed the required
filtration.
\eop

\begin{lemma}
In the above situation (and fixing $D\subset X$ equal to the hyperplane at
infinity), there is $m_0$ so that for  $m\geq m_0$ and all $i>0$
we have
$$
H^i(X/S, L(mD)) = 0.
$$
Furthermore, for any $m\geq m_0$ there is $m_1$ such that for $m'\geq m_1$,
we have an inclusion
$$
0\rightarrow H^0(X/S, L) \rightarrow H^0(mD/S, L|_{mD})
$$
whose image is the same as the image of
$$
H^0(m'D/S, L|_{m'D})\rightarrow H^0(mD/S, L|_{mD}).
$$
\end{lemma}
{\em Proof:}
By standard arguments with Artin approximation and Krull's lemma, it suffices to
obtain these properties (with $m_0$ and $m_1(m)$ uniform) for all vector bundles
of the form  $a_{\ast}(L_A)$ for artinian schemes of finite type $a:A\rightarrow
Spec(k)$ with morphisms $A\rightarrow S$. Furthermore, these properties are
invariant under taking extensions. Lemma \ref{boundedness} says
that there is a bounded family $\Bb$ of vector bundles such that any
$a_{\ast}(L_A)$ has a filtration with subquotients in this bounded family.  We
can find $m_0$ and $m_1(m)$ as necessary for the bounded family $\Bb$, and then
the required properties follow for the $a_{\ast}(L_A)$ which are successive
extensions of bundles in $\Bb$.
\eop

\oldsubnumero{Proof of Proposition \ref{statement}}
As stated above we may assume that $X=\pp ^d_S$. Proceed by induction on $d$
(the statement being automatic if $d=0$). In the situation of the preceding
lemma, the inductive hypothesis (plus the fact that vector sheaves are closed
under kernel, cokernel and extension) implies that the terms
$H^0(mD/S, L|_{mD})$ and $H^0(m'D/S, L|_{m'D})$ which appear are themselves
vector sheaves. The image of the map between them is therefore a vector sheaf,
so this shows that $H^0(X/S, L)$ is a vector sheaf. It follows that the same is
true of $H^0(X/S, L(mD))$ (because the same hypotheses apply to $L(mD)$). Now a
descending induction on $m$ starting with $m\geq m_0$, using the hypothesis that
the statement is known in dimension $d-1$ and
preservation of the category of vector sheaves
under kernel, cokernel and extension, shows that the $H^i(X/S, L(mD))$ are
vector sheaves for all $m$.
\eop

We reiterate that the proposition is not true if $L$ is an arbitrary
vector sheaf over $X$, see O. Gabber's counterexample in Hirschowitz
\cite{Hirschowitz}.

\subnumero{The main results for simply connected coefficient stacks}

Recall from \ref{index1} the definition of the realm $\underline{AV}$.
An $\ngr$-stack $T$ on $\Gg /X$ is in $\underline{AV}(X)$ if and only if
$\pi _0(T)=\pi _1(T)=\ast$, and the $\pi _i(T)$ (which are now well-defined as
sheaves on $X$) are vector sheaves on $X$ for $i\geq 2$. Note that this
means that
for any $Y\rightarrow X$ the $\pi _i(T|_{\Gg /Y})$ are also vector sheaves
on $Y$
because they are just the pullbacks of the $\pi _i(T)$ from $X$ (this type of
statement is not true in general but it is true here due to the fact that $T$ is
connected so there is no choice of basepoint).

\begin{theorem}
\label{av}
Suppose $(X,\Ff )\rightarrow (S, \Ee )$ is a morphism of smooth type
between formal categories of smooth type, such that $X\rightarrow S$ is
projective. Then the $\underline{AV}$-shape of $\Ff /\Ee$ takes values in
$\underline{VP}$, i.e.
$$
{\bf Shape}_{\underline{AV}}(\Ff /\Ee ): \Ee \rightarrow
\underline{Hom}(\underline{AV},\underline{VP}),
$$
which we have written above as
$$
\underline{AV}
\stackrel{{\bf Shape} (\Ff /\Ee )}{\longrightarrow}
\underline{VP}.
$$
\end{theorem}
{\em Proof:}
Apply Proposition \ref{meaningprop}.
Suppose $Z\rightarrow \Ee$ is a morphism from a scheme. Then $\Ff \times _{\Ee}
Z\rightarrow Z$ is a formal category of smooth type over the base scheme $Z$,
itself satisfying the hypotheses \ref{situation}.
Suppose $T\in \underline{AV}(Z)$. This means that $T$ is a relatively
$1$-connected $\ngr$-stack on $\Gg /Z$, with higher homotopy groups being vector
sheaves. Since $T$ is simply connected, the higher homotopy local systems are
actually vector sheaves pulled back from $Z$.  We apply Corollary
\ref{postnicor}
with
$$
R= A=\Ff \times _{\Ee}Z
$$
and $B = Z$.
What is called $T$ in Corollary \ref{postnicor} is, in terms of our notations,
$T\times _ZR$.  The local systems $L^k$ which
enter into the statement of Corollary \ref{postnicor} are pullbacks to $R$
of the $\pi
_i(T/Z)$, thus (by the definition of $\underline{AV}$) the $L^k$ are pullbacks
to $R$ of  vector sheaves on $Z$.

The realm $\underline{VP}$
is closed under extension and finite limits, by Lemma \ref{vpcompatible}.

Note also that pulling back to another $Y\in \Gg$ as needed in
Corollary \ref{postnicor} (what was denoted by $X$ in Corollary
\ref{postnicor}, we denote by $Y$ in
the present phrase so as not to confuse it with the $X$ of \ref{situation})
amounts just to another pullback of the same form as $Z$ so we can ignore it
(i.e. assume $Y=Z$). Furthermore the section $\eta : A:\rightarrow R=A$ is
unique, equal to the identity. By  Theorem \ref{mainab}, the cohomology
sheaves of the direct image complex
$$
{\bf R} f_{Z,\ast}(L^k)
$$
are vector sheaves on $Z$. This implies that the Dold-Puppe of this complex is a
very presentable $n$-stack on $\Gg /Z$, i.e. is in $\underline{VP}(Z)$. By
Corollary \ref{postnicor}, we get that the morphism
$$
\underline{\Gamma}(A/B, T\times _ZR)= \underline{Hom}(\Ff \times _{\Ee}
Z/Z, T/Z)
\rightarrow \underline{\Gamma}(A/B, R) = Z
$$
is of type $\underline{VP}$. In other words,
$\underline{Hom}(\Ff \times _{\Ee} Z/Z, T/Z)$ is a very presentable
$n$-stack over
$Z$.

Plugging this statement back into Proposition \ref{meaningprop}, we obtain the
desired statement that
$$
\underline{AV}
\stackrel{{\bf Shape} (\Ff /\Ee )}{\longrightarrow}
\underline{VP}.
$$
\eop

Recall from \ref{index1} the definition of the realm $\underline{AL}$ which is
similar to $\underline{AV}$ but with vector sheaves replaced by locally free
sheaves. An $\ngr$-stack $T$ on $\Gg /X$ is in $\underline{AL}(X)$ if and only
if $\pi _0(T)=\pi _1(T)=\ast$ and the $\pi _i(T)$ are locally free sheaves on
$X$ for $i\geq 2$. Again these latter are well-defined because $T$ is simply
connected, and being locally free sheaves on $X$ implies that the $\pi
_i(T|_{\Gg /Y})$ are locally free on $Y$ for any $Y\rightarrow X$.

\begin{theorem}
\label{al}
Suppose $(X,\Ff )\rightarrow (S, \Ee )$ is a morphism of smooth type
between formal categories of smooth type, such that $X\rightarrow S$ is
projective. Then the $\underline{AL}$-shape of $\Ff /\Ee$ takes values in
the realm of geometric very presentable $n$-stacks $\underline{VG}$, i.e.
$$
{\bf Shape}_{\underline{AL}}(\Ff /\Ee ): \Ee \rightarrow
\underline{Hom}(\underline{AL},\underline{VG}),
$$
which we have written above as
$$
\underline{AL}
\stackrel{{\bf Shape} (\Ff /\Ee )}{\longrightarrow}
\underline{VG}.
$$
\end{theorem}
{\em Proof:}
The same as the proof of Theorem \ref{av}.
Use the statement of Lemma \ref{vgcompatible} that
$\underline{VG}$ is closed under extensions and finite limits.
By the definition of $\underline{AL}$,
the local systems $L^k$ which enter into Corollary \ref{postnicor} are pullbacks of
vector bundles on $Z$; and Theorem \ref{mainab} then gives that the  higher
direct image complexes are perfect complexes, so their Dold-Puppe are geometric
(and also very presentable) $n$-stacks. Again, plugging into Proposition
\ref{meaningprop}, we obtain the statement that
$$
\underline{AL}
\stackrel{{\bf Shape} (\Ff /\Ee )}{\longrightarrow}
\underline{VG}.
$$
\eop

We now give the statement for coefficients in $\underline{AG}$. Recall that
this realm is the sub-realm of $\underline{AV}$ consisting of objects which
are geometric (it contains $\underline{AL}$). The proof the following theorem is
somewhat more complicated than the previous ones, and we only sketch the idea
here.

\begin{theorem}
\label{ag}
Suppose $(X,\Ff )\rightarrow (S, \Ee )$ is a morphism of smooth type
between formal categories of smooth type, such that $X\rightarrow S$ is
projective. Then the $\underline{AG}$-shape of $\Ff /\Ee$ takes values in
the realm of geometric very presentable $n$-stacks $\underline{VG}$, i.e.
$$
{\bf Shape}_{\underline{AG}}(\Ff /\Ee ): \Ee \rightarrow
\underline{Hom}(\underline{AG},\underline{VG}),
$$
which we have written above as
$$
\underline{AG}
\stackrel{{\bf Shape} (\Ff /\Ee )}{\longrightarrow}
\underline{VG}.
$$
\end{theorem}
{\em Sketch of proof:}
We can't directly apply a Postnikov argument such as Corollary
\ref{postnicor} because the
geometric $n$-stacks aren't closed under truncation. However, the technique of
\S 7.3 does permit us approximately to decompose a simply connected geometric
$n$-stack $T$ over $Z$, into components which are Eilenberg-MacLane stacks for
vector bundles. We have obtained a sequence of morphisms  giving fibration
sequences of the form
$$
T_{i+1}\rightarrow T_i\rightarrow \Sigma _i,
$$
starting with $T_0=T$. The $\Sigma _i$
are Dold-Puppe of perfect complexes over $Z$. Eventually the $T_i$ become highly
connected. On the other hand, $\Ff /\Ee$ has finite cohomological dimension, so
eventually the $T_i$ become irrelevant. (To make this argument precise one again
has to do something similar to what was done in \S 7.3 using truncation, the
notion of ``almost-geometric'', and Lemma \ref{aggeo}.

In sum, it suffices to
treat the cohomology of $\Ff /\Ee$ with coefficients in $\Sigma _i$. But now we
can decompose a perfect complex into a successive extension of vector bundles,
which gives a decomposition of $\Sigma_i$ into extensions of Eilenberg-MacLane
stacks of the form $K(L,m)$ for vector bundles $L$ on $Z$. Now, as shown by
Theorem \ref{mainab}, the cohomology of $\Ff /\Ee$ with coefficients in $K(L,m)$
is the Dold-Puppe of a perfect complex, so it is geometric. Again, apply closure
of $\underline{VG}$ under extensions to complete the proof.
\eop

\subnumero{Degree $1$ cohomology with flat linear coefficients}

We next want to extend the previous results to the case of coefficient
stacks $T$
in the realms $\underline{FV}$, $\underline{FL}$, and (again with only a
sketch of
the proof) $\underline{FG}$. The main additional ingredient necessary to do this
is an analysis of degree $1$ cohomology with coefficients in affine flat group
schemes. We discuss this briefly here.

\begin{parag}
\label{linear}
We start with the following definition. If $S$ is a scheme, a {\em linear group
scheme over $S$} is a group scheme $G\rightarrow S$ such that locally in the
etale topology of $S$, $G$ embeds as a closed subgroup scheme of some $GL(E)$
for a finite rank vector bundle $E$ over $S$.
\end{parag}

Note that a linear group scheme is automatically affine and of finite type over
$S$. Generally we will impose the extra condition that it be flat.

\begin{conjecture}
\label{maybeDL}
Any flat affine group scheme of finite type $G\rightarrow S$ is linear.
\end{conjecture}

This is of course true if $S=Spec (k)$.
In order to prove it in general one would seem to need an equivariant version of
Deligne-Lazard and I don't know whether such a theorem is available. Thus, in a
relative situation we will work only with linear group schemes.

\begin{parag}
\label{situationhere}
Suppose $X\rightarrow \Ff \rightarrow S$ is a projective morphism from a formal
category of smooth type, to a base scheme $S$. Suppose $G$ is a flat
linear group scheme over $S$. We will look at the {\em moduli $1$-stack}
$$
\Mm (\Ff /S, G):= \underline{Hom}(\Ff /S, K(G/S,1))\rightarrow S
$$
for principal $G$-bundles over $\Ff$.

\begin{proposition}
\label{itsartin}
The moduli stack $\Mm (\Ff /S, G)$ is an Artin algebraic stack locally of finite
type over $S$.
\end{proposition}
{\em Proof:}
We may choose an embedding $G\subset GL(E)$ with $E=\Oo _S^r$, possibly by
localizing in the etale topology of $S$. We get a morphism
$$
\Mm (\Ff /S, G)\rightarrow
\Mm (\Ff /S, GL(E)).
$$

\newparag{itsartin1}
As stated in \ref{smooth5} and
\ref{lambda}, there is a split almost-polynomial sheaf of rings of differential
operators $\Lambda$ on $X/S$, such that local systems on $\Ff$ are the same
thing as $\Lambda$-modules. (For example if $\Ff = X_{DR}$ then $\Lambda = {\cal
D}_X$ is the usual sheaf of rings of differential operators on $X$.) In
particular, $\Mm (\Ff /S, GL(E))$ is the moduli stack of rank $\Lambda$-modules
on $X/S$ which are locally free of rank $r$ over $\Oo _X$.

The same type of Hilbert scheme construction as used in \cite{Simpson}(iii) I
Theorem 3.8 works for any $\Lambda$ modules rather than just semistable ones;
for any open subfunctor of the moduli functor which is bounded, one obtains a
parametrizing Hilbert scheme; and the quotient by the relevant group action is
an Artin algebraic stack of finite type. As the moduli functor can be
exhausted by
open subfunctors, this shows that $\Mm (\Ff /S, GL(E))$ is an Artin algebraic
stack locally of finite type.

\newparag{itsartin2}
We now have to go from $GL(E)$ to $G$. Note that we have a fibration sequence
relative to $S$,
$$
(GL(E)/G)\rightarrow K(G/S, 1) \rightarrow K(GL(E)/S, 1).
$$
This allows us to identify the fiber of
$$
\underline{Hom}(\Ff /S, K(G/S, 1))\rightarrow
\underline{Hom}(\Ff /S, K(GL(E)/S, 1)).
$$
To be precise, proceed as follows: there is a universal principal $GL(E)$-bundle
$P$ on
$$
\Ff \times _S \Mm (\Ff /S, GL(E))\rightarrow \Mm (\Ff /S, GL(E)).
$$
The bundle associated to $P$ and the action of $GL(E)/G$ is a
morphism
$$
Red(G/GL(E); P):= (GL(E)/G)\times ^{GL(E)}P\rightarrow
\Ff \times _S \Mm (\Ff /S, GL(E)).
$$
This morphism is representable and smooth in the sense that if
$Z$ is a scheme mapping to  the $1$-stack
$\Ff \times _S \Mm (\Ff /S, GL(E))$ then the morphism
$$
Red(G/GL(E); P)\times _{\Ff \times _S \Mm (\Ff /S, GL(E))}Z
\rightarrow Z
$$
is a smooth morphism between algebraic spaces (indeed it is just the bundle
associated to the principal bundle $P$ pulled back to $Z$).

\newparag{itsartin3}
We have the formula
$$
\Mm (\Ff /S, G) = \underline{\Gamma}(
\Ff \times _S \Mm (\Ff /S, GL(E))/\Mm (\Ff /S, GL(E)), Red(G/GL(E); P)),
$$
as a $1$-stack lying over $\Mm (\Ff /S, GL(E))$.

\newparag{itsartin4}
Suppose $Y$ is a scheme with a smooth surjective map
$$
a: Y\rightarrow \Mm (\Ff /S, GL(E)).
$$
Then the map
$$
Y\times _{\Mm (\Ff /S, GL(E))}\Mm (\Ff /S, G)\rightarrow
\Mm (\Ff /S, G)
$$
is again smooth and surjective. Thus to prove that the target
$\Mm (\Ff /S, G)$ is an Artin algebraic stack locally of finite type, it
suffices to prove that for a map $a$, the pullback
$$
Y\times _{\Mm (\Ff /S, GL(E))}\Mm (\Ff /S, G)
$$
is an algebraic space locally of
finite type. The formula \ref{itsartin3} gives
$$
Y\times _{\Mm (\Ff /S, GL(E))}\Mm (\Ff /S, G)
= \underline{\Gamma}(\Ff \times _S Y; Red(G/GL(E); P_Y))
$$
where $P_Y$ is the principal bundle $P$ pulled back to $\Ff \times _S Y$.

The lemma which follows (which is basically the ``Douady space'')  will
therefore
complete the proof of the proposition. \eop
\end{parag}

\begin{lemma}
\label{douady}
Suppose $\Ff \rightarrow S$ is a projective morphism of smooth type from a
formal
category of smooth type to a base scheme $S$. Suppose $W\rightarrow \Ff$ is a
morphism of $1$-stacks such that $W\times _{\Ff} X \rightarrow X$ is a smooth
morphism of algebraic spaces of finite type. Then
$$
\underline{\Gamma} (\Ff /S, W)\rightarrow S
$$
is an algebraic space locally of finite type over $S$.
\end{lemma}
{\em Proof:}
This is a variant of the ``Douady space'' of maps. The algebraic-space version
of the Douady space gives an Artin algebraic space locally of finite type
$$
\underline{\Gamma}(X/S, W) \rightarrow S.
$$
Furthermore, if $N^{(1)}$ denotes the  first infinitesimal neighborhood of
$e(X)$ in the morphism object $N$ of the formal category $\Ff$, then the two
projections $s,t: N^{(1)}\rightarrow X$ induce two morphisms
$$
\underline{\Gamma}(X/S, W)\tworightarrows
\underline{\Gamma}(N^{(1)}/S, W).
$$
The targets of these morphisms are again Artin algebraic spaces locally of
finite type (note that $N^{(1)}$ is flat over $S$). Finally, our
$\underline{\Gamma} (\Ff /S, W)$ is just the equalizer of these two morphisms,
so it is an algebraic space locally of finite type.
\eop

\subnumero{Results for coefficient stacks  with flat linear group schemes
in degree $1$}

Recall from \ref{index1} the definition of the realm $\underline{FV}$.
An $\ngr$-stack $T$ on $\Gg /X$ is in $\underline{FV}(X)$ if and only if
$\pi _0(T)=\ast$, $\pi _1(T,t)$ is a flat linear group scheme (over an etale
covering of $X$ where the basepoint $t$ exists), and the $\pi _i(T)$ (which are
etale-locally defined as sheaves on $X$) are vector sheaves on $X$ for
$i\geq 2$.

\begin{theorem}
\label{fv}
Suppose $(X,\Ff )\rightarrow (S, \Ee )$ is a morphism of smooth type
between formal categories of smooth type, such that $X\rightarrow S$ is
projective. Then the $\underline{FV}$-shape of $\Ff /\Ee$ takes values in
$\underline{VP}^{\rm loc}$, i.e.
$$
{\bf Shape}_{\underline{FV}}(\Ff /\Ee ): \Ee \rightarrow
\underline{Hom}(\underline{FV},\underline{VP}^{\rm loc}).
$$
\end{theorem}
{\em Proof:}
We follow the same outline of proof as \ref{av}.

First apply Proposition \ref{meaningprop}.
Suppose $Z\rightarrow \Ee$ is a morphism from a scheme. Then $\Ff \times _{\Ee}
Z\rightarrow Z$ is a formal category of smooth type over the base scheme $Z$,
itself satisfying the hypotheses \ref{situation}.
Suppose $T\in \underline{FV}(Z)$. By localizing over the base we can assume that
there is a basepoint section $t: Z\rightarrow T$. Let $G= \pi _1(T,t)$ which is
a flat linear group scheme over $Z$. We have
$$
\tau _{\leq 1}(T) = K(G/Z,1).
$$
On the other hand note (by\ref{usualsmooth}) that $\Ff \times _{\Ee}
Z\rightarrow Z$ is a morphism from a formal category of smooth type, to the base
scheme $Z$. By Lemma \ref{douady},
$$
\underline{Hom}(\Ff \times _{\Ee} Z/Z, K(G/Z,1)/Z) = \Mm (\Ff \times _{\Ee}
Z/Z,G)
$$
is an algebraic stack locally of finite type. In particular it is in
$\underline{GE}^{\rm loc}(Z)\subset \underline{PE}^{\rm loc}(Z)$. The
automorphism group
schemes of this algebraic $1$-stack are affine. Therefore it is in
$\underline{VG}^{\rm loc}(Z)\subset \underline{VP}^{\rm loc}(Z)$. In view of
the closure of  $\underline{VP}^{\rm loc}$ under extensions, in order to finish
the proof it suffices to show that the morphism
$$
\underline{Hom}(\Ff \times _{\Ee} Z/Z, T/Z)\rightarrow
\underline{Hom}(\Ff \times _{\Ee} Z/Z, K(G/Z,1)/Z)
$$
is of type $\underline{VP}$.

For this we can apply Corollary \ref{postnicor}
with
$$
A=\Ff \times _{\Ee}Z
$$
and $B = Z$. Differently from the case of the proof of \ref{av}, here we put
$$
R= A\times _Z K(G/Z, 1).
$$
Put the $T$ of Corollary \ref{postnicor} equal to our $T\times _ZA$.
Then the morphism $\underline{\Gamma}(A/B,f)$ of \ref{postnicor}
is the same as our morphism above. The local systems
$\eta ^{\ast}L^k$ which appear are again locally (in the etale topology of the
scheme mapping to $A$) pullbacks of vector sheaves from $Z$. The same
application
as in the proof of \ref{av}, using Theorem \ref{mainab}, completes the
proof that
$\underline{\Gamma}(A/B,f)$ is of type $\underline{VP}$.
\eop

Recall from \ref{index1} the definition of the realm $\underline{FL}$ which is
similar to $\underline{FV}$ but with vector sheaves in degrees $\geq 2$ replaced
by locally free sheaves. An $\ngr$-stack $T$ on $\Gg /X$ is in
$\underline{FL}(X)$ if and only if $\pi _0(T)=\ast$, $\pi _1(T,t)$ is a flat
linear group scheme locally over $X$ where the basepoint is defined, and
the $\pi
_i(T,t)$ are locally free sheaves on $X$ (again, locally in the etale
topology where the basepoint is defined) for $i\geq 2$.

\begin{theorem}
\label{fl}
Suppose $(X,\Ff )\rightarrow (S, \Ee )$ is a morphism of smooth type
between formal categories of smooth type, such that $X\rightarrow S$ is
projective. Then the $\underline{FL}$-shape of $\Ff /\Ee$ takes values in
the realm of locally geometric very presentable $n$-stacks $\underline{VG}^{\rm
loc}$, i.e.
$$
{\bf Shape}_{\underline{FL}}(\Ff /\Ee ): \Ee \rightarrow
\underline{Hom}(\underline{FL},\underline{VG}^{\rm loc}).
$$
\end{theorem}
{\em Proof:}
Combine the proofs of \ref{al} and \ref{fv}.
\eop

\subnumero{A simpler statement for nonabelian cohomology}

The above statements being somewhat complex, we extract here the simpler
statement for nonabelian cohomology with coefficients in a connected very
presentable $n$-stack $T$ (cf \ref{degree3} \ref{index4}).

\begin{theorem}
\label{simpler}
Suppose $(X,\Ff )\rightarrow (S, \Ee )$ is a morphism of smooth type
between formal categories of smooth type, such that $X\rightarrow S$ is
projective. Suppose $T$ is a connected very presentable $\ngr$-stack over
$Spec(k)$, i.e. an $n$-stack with $\pi _0=\ast$, $\pi _1$ an affine algebraic
group scheme of finite type, and the $\pi _i$ being vector spaces for $i\geq
2$.  Then the nonabelian cohomology $\ngr$-stack
$$
\underline{Hom}(\Ff /\Ee , T) \rightarrow \Ee
$$
is relatively a locally geometric very presentable $\ngr$-stack over $\Ee$, i.e.
the corresponding cartesian family is a morphism
$$
\Ee \rightarrow \underline{VG}^{\rm loc}.
$$
\end{theorem}
{\em Proof:}
This is just a restatement of what was said above. The reader
may follow out the proofs of \ref{av}, \ref{al}, \ref{fv} restricting to this
case to simplify.
\eop

The above theorems apply directly to all of the morphisms of formal categories
$\Ff \rightarrow \Ee$ related to Hodge theory, which were constructed in \S 9.

\begin{remark}
One can define a realm $\Ff \Pp \subset 1^{\rm gr}\underline{STACK}$ as
follows: for a scheme $X$, set $\Ff \Pp (X)$ equal to the $2$-category of
projective formal schemes of smooth type over $X$. The above results can then be
interpreted as saying that certain triples $(\Ff \Pp , \underline{R},
\underline{A})$ are well-chosen in the sense of \ref{zoo2}. Notably, we
have that
the following triples are well-chosen:
$$
(\Ff \Pp , \underline{AV} , \underline{VP}),
$$
$$
(\Ff \Pp , \underline{AL} , \underline{VG}),
$$
$$
(\Ff \Pp , \underline{FV} , \underline{VP}^{\rm loc}),
$$
$$
(\Ff \Pp , \underline{FL} , \underline{VG}^{\rm loc}).
$$
\end{remark}

\subnumero{Semistability}

Ideally, one should integrate into the above discussion a notion of
``semistability'' for degree $1$ cohomology. We haven't done this, because in
the present paper we are mostly looking at the effect of higher-degree
cohomology. Very briefly, what needs to be done is as follows. Suppose $\Ff
\rightarrow Z$ is a morphism of smooth type from a formal category of
smooth type
to a base scheme $Z$. If $G\rightarrow Z$ is a flat linear group scheme,
then the
nonabelian cohomology $1$-stack
$$
\underline{Hom}(\Ff /Z, K(G/Z,1))
$$
contains an open substack which could be denoted
$\underline{Hom}^{\rm se}(\Ff /Z, K(G/Z,1))$, consisting of principal
$G$-bundles over $\Ff$ which are semistable and---say---with vanishing Chern
classes.

This semistable cohomology stack  has better properties, for example it
is much closer to being separated. Recall from \cite{Simpson}(i), (vi) that
in the case $Z=\af$ and $\Ff = X_{Hod}:={\bf Hodge}(X_{DR})$, this open subset
of semistable objects with vanishing Chern classes was used as the ``Hodge
filtration'' on the cohomology of $X_{DR}$.

Now if $T\rightarrow Z$ is a coefficient $n$-stack in one of the realms
$\underline{FV}$ or $\underline{FL}$, let $G= \pi _1(T/Z)$. We have a morphism
$$
\underline{Hom}(\Ff /Z, T)\rightarrow
\underline{Hom}(\Ff /Z, K(G/Z,1)).
$$
Define $\underline{Hom}^{\rm se}(\Ff /Z, T)$ to be the inverse image under this
morphism of the open substack
$$
\underline{Hom}^{\rm se}(\Ff /Z, K(G/Z,1)).
$$
If $\Ff \rightarrow \Ee$ is a morphism of smooth type between formal categories
of smooth type, we obtain similarly an open substack
$$
\underline{Hom}^{\rm se}(\Ff /\Ee , T)\subset
\underline{Hom}(\Ff /\Ee , T).
$$
As $T$ varies, we obtain a ``semistable shape'' which could be denoted
$$
{\bf Shape}_{\underline{FV}}(\Ff /\Ee , {\rm se}):\Ee \rightarrow
\underline{Hom}(\underline{FV}, \underline{VP}).
$$
Note that the semistable open substack of the degree $1$ cohomology stack is
of finite type (since we included vanishing of Chern classes in the
definition) so
the answer realm here is $\underline{VP}$ rather than just $\underline{VP}^{\rm
loc}$.

It is the semistable shape which should really be used in defining the Hodge
filtration, regular singularity of the Gauss-Manin connection, etc. One can note
that for $\Ff = X_{DR}$ (and even in the relative ``Gauss-Manin'' situation
where $\Ee = S_{DR}$) the semistability and vanishing of Chern classes are
automatic, so the semistable shape is the same as the whole shape. However, for
Dolbeault cohomology (i.e. over $0\in \af$) they are no longer the same.

Developing the above remarks in a precise way would require a thorough
understanding of semistability for principal $G$-bundles where $G$ is a flat
linear (but not necessarily constant) group scheme. This should be possible
using
the techniques of \cite{Simpson}(iii) and other authors, but that would get
beyond
the scope of the present paper.

\subnumero{Results on the representing stacks}

Here we translate the above results into results about the $n$-stacks ${\bf
rep}(\Ff /S)$ representing the $\underline{AV}$-shape of $\Ff /S$ (see \S 6.4).
The representing $n$-stacks exist essentially only in the cases where classical
rational homotopy theory can do the job, for example in the simply connected
case. In this section we restrict to the hypothesis $H^1(\Ff /S, \Oo )=0$.
The newness of the present treatment resides in the fact that it avoids a whole
lot of calculations with d.g.a.'s and in so doing, allows us to go to the full
level of homotopy-coherence for the Gauss-Manin connection of
Navarro-Aznar \cite{Navarro-Aznar} for example.
Apart from this advantage, everything we say here is essentially already known
e.g. by Navarro-Aznar's paper.

Refer to \S 6.4 for notations and results concerning ${\bf rep}(\Ff /S)$.

\begin{hypothesis}
\label{situationA}
Suppose that $X\rightarrow \Ff \stackrel{p}{\rightarrow} S$ is a
morphism  from a formal category of smooth type to a base scheme $S$, and
suppose
that $X\rightarrow S$ is flat and projective, i.e. $f$ is a projective morphism
of smooth type of formal categories.
\end{hypothesis}

\begin{lemma}
\label{uc}
In the situation \ref{situationA}, let
$$
{\bf R} p_{\ast} (\Omega ^{\cdot}_{X/\Ff })
$$
denote the higher direct image of the de Rham complex. It is a perfect complex
on $S$, and if $V$ is a vector sheaf on $S$ then
$$
H^i(\Ff /S, V) = H^i (V\otimes _{\Oo _S}
{\bf R} p_{\ast} (\Omega ^{\cdot}_{X/\Ff })).
$$
\end{lemma}
{\em Proof:}
It is well-known that the higher direct image of a perfect complex by a flat
morphism is again perfect, see \cite{Mumford} for example. This works even for
complexes where the differentials are not $\Oo$-linear, because the spectral
sequence for hypercohomology expresses the higher direct image complex as a
successive extension.

For the second statement, we use the argument that was sketched at the end
of \S 8.3 and which was made precise at the start of the proof of Theorem
\ref{mainab}. Namely, we get functorially in $Z\rightarrow S$, the formula
$$
H^i(\Ff /S, V)(Z) = {\bf H}^i(X_Z^{\rm zar}, L_Z \otimes _{\Oo _{X_Z}}\Omega
^{\cdot} _{X_Z/\Ff _Z}),
$$
for $L$ the local system on $\Ff$ pull-back of the local system $V$ on $S$.
(This formula is the same as the first displayed formula in the proof of
\ref{mainab}.)

The hypercohomology can be calculated by a \v{C}ech complex for a
fixed affine open covering of $X$, and from this formula (and the fact that the
components of the de Rham complex are localy free) it is easy to see that
tensoring  before taking the hypercohomology is the same as first taking the
higher-derived direct image of the de Rham complex and then tensoring with $V$
and then taking the cohomology. Thus we obtain, again functorially in $Z$,
$$
H^i(\Ff /S, V)(Z)= H^i (V_Z\otimes _{\Oo _Z}
{\bf R} p_{Z,\ast} (\Omega ^{\cdot}_{X_Z/\Ff _Z}))
$$
where $p_Z: X_Z \rightarrow Z$ is the projection (base-change of $p$). This
formula functorially in $Z$, is equivalent to the formula in the statement of
the lemma.
\eop

\begin{corollary}
\label{itsanchored}
In the situation \ref{situationA}, for any vector sheaf $V$ on $S$ the
cohomology
$H^i(\Ff /S, V)$ is a vector sheaf, and the
the endofunctor $V\mapsto H^i(\Ff /S, V)$ on the category of vector sheaves over
$S$, is anchored (cf \ref{defanchored}).
\end{corollary}
{\em Proof:}
If $L^{\cdot}$ is a perfect complex, then the
$H^i(V\otimes _{\Oo _S}L^{\cdot})$ are vector sheaves (this by the closure of
the category of vector sheaves under kernel, cokernel and extension).
Furthermore, it follows from \ref{anchored1} that the functor
$$
V\mapsto H^i(V\otimes _{\Oo _S}L^{\cdot})
$$
is anchored.
\eop

\begin{corollary}
\label{zeroone}
In the situation \ref{situationA}, suppose that for every point $z:Spec
(k(z))\rightarrow S$, $H^0(\Ff _z, \Oo )= \Oo_{Spec(k(z))}$ and $H^1(\Ff _z, \Oo
)= 0$. Then for any vector sheaf $V$ on $S$, we have that $H^0(\Ff /S, V )=
V$ and
$H^1(\Ff /S, V )= 0$.
\end{corollary}
{\em Proof:}
Let $L^{\cdot}$ denote the higher direct image of the de Rham complex appearing
in \ref{uc}. Let $\overline{L}^{\cdot}$ denote the cone of the map
$$
\Oo _S \rightarrow L^{\cdot}
$$
(the bar notation is because this calculates reduced cohomology). The hypothesis
of the current corollary implies that $\overline{L}^{\cdot}$ is exact over every
point of $S$, in degrees $0$ and $1$. Note that  $\overline{L}^{\cdot}$ is a
perfect complex bounded below. Realize it as a complex of vector bundles (we
might have to localize on $S$ for this). Suppose that the first nonzero term is
$\overline{L}^i$. If $i\leq 1$ then the morphism
$$
\overline{L}^i\rightarrow \overline{L}^{i+1}
$$
has the property that it is injective over each fiber. Thus it is a strict
morphism of vector bundles, so the quotient
$$
\overline{L}^{i+1}/\overline{L}^{i}
$$
is again a vector bundle. In this way we can replace $\overline{L}^{\cdot}$ by a
quasiisomorphic complex starting in degree $i+1$. Applying this argument
inductively, we can replace $\overline{L}^{\cdot}$  by a quasiisomorphic
perfect complex $M^{\cdot}$ starting in degree $2$. Now
$$
\overline{H}^i(\Ff /S, V) = H^i(V\otimes _{\Oo _S}M^{\cdot})
$$
where $\overline{H}^i$ denotes the reduced cohomology (dividing by $V$ in degree
$0$). Thus the reduced cohomology vanishes in degrees $\leq 1$ as required.
\eop

\begin{corollary}
\label{last1}
Suppose that $\Ff \rightarrow \Ee$ is a projective morphism of smooth type
between formal categories of smooth type.  Suppose that for every point $z:Spec
(k(z))\rightarrow \Ee$, $H^0(\Ff _z, \Oo )= \Oo_{Spec(k(z))}$ and $H^1(\Ff
_z, \Oo
)= 0$. Then the representing object for the $\underline{AV}$-shape of $\Ff /S$
exists,
$$
\Ff \rightarrow {\bf rep}(\Ff /\Ee )\rightarrow \Ee .
$$
Furthermore ${\bf rep}(\Ff /\Ee )\rightarrow \Ee $ is a $1$-connected geometric
very presentable morphism, corresponding to a morphism
$$
\Ee\rightarrow \underline{AG}.
$$
\end{corollary}
{\em Proof:}
Note that for any scheme $Z\rightarrow \Ee$, the morphism
$$
\Ff _Z:= \Ff \times _{\Ee} Z\rightarrow Z
$$
satisfies the hypothesis \ref{situationA}.  In particular, \ref{itsanchored}
applies. Furthermore, the hypotheses of the present corollary imply those of the
previous corollary for $\Ff _Z\rightarrow Z$, since base-change by a point $z:
Spec (k(z))\rightarrow Z$ is the same thing as base-change of $\Ff \rightarrow
\Ee$ by the composed point $Spec(k(z))\rightarrow \Ee$.  By
\ref{itsanchored} and
\ref{zeroone}, we obtain the hypotheses necessary to apply Corollary
\ref{representable2}, giving the existence of ${\bf rep}(\Ff /\Ee )$.

We now prove that the morphism
${\bf rep}(\Ff /\Ee )\rightarrow \Ee $ is a $1$-connected geometric
very presentable morphism. This statement means that for a  map from a scheme
$Z\rightarrow \Ee$, the object
$$
{\bf rep}(\Ff _Z/Z )=
{\bf rep}(\Ff /\Ee )\times _{\Ee} Z
$$
should be geometric over $Z$. For this, apply Corollary \ref{geometricity}
(which is basically just an application of the criterion of Theorem
\ref{criterion}). Note that the higher direct image complex ${\bf R} p_{Z,\ast}
(\Oo )$ (where $p_Z:\Ff _Z\rightarrow Z$ is the projection) is the same as the
higher direct image of the de Rham complex of $\Ff _Z$. This direct image is a
perfect complex by \ref{uc}. Thus Corollary \ref{geometricity} applies to show
that ${\bf rep}(\Ff _Z/Z )$ is geometric over $Z$. The representing stack
corresponds to a morphism $$
{\bf rep}(\Ff /\Ee ): \Ee \rightarrow \underline{AG}.
$$
\eop

\begin{corollary}
\label{last2}
Suppose that $\Ff \rightarrow \Ee$ is a projective morphism of smooth type
between formal categories of smooth type.
Suppose that the scheme $S$ underlying the formal category $\Ee$, is integral.
Suppose that for every point $z:Spec (k(z))\rightarrow \Ee$, $H^0(\Ff _z, \Oo )=
\Oo_{Spec(k(z))}$ and $H^1(\Ff _z, \Oo )= 0$. Suppose furthermore that for every
$i\geq 2$, the dimensions of the homotopy group vector-spaces
$$
z\mapsto {\rm dim} _{k(z)}\pi _i({\bf rep}(\Ff _z/Spec(k(z))))
$$
are constant as  functions of $z$. Then the homotopy group local systems
$$
\pi _i({\bf rep}(\Ff /\Ee )/ \Ee )
$$
are local systems of vector bundles over $\Ee$. Thus the representing object
corresponds to a morphism
$$
{\bf rep}(\Ff /\Ee ): \Ee \rightarrow \underline{AL}.
$$
\end{corollary}
{\em Proof:}
Apply \ref{locfreeness}
\eop


\begin{thebibliography}{MM2}

\bibitem{Adams}
J. F. Adams. {\em Infinite Loop Spaces.} Princeton University Press {\em Annals
of Math. Studies} {\bf 90} (1978).

\bibitem{Artin}
M. Artin. Versal deformations and algebraic stacks. {\em Inventiones Math.}
{\bf 27} (1974), 165-189.

\bibitem{ArtinMazur}
M. Artin, B. Mazur. Etale Homotopy, {\em Lecture Notes in Mathematics}

\bibitem{SGA4}
M. Artin, J. Verdier, A. Grothendieck. {\em SGA 4: Th\'eorie des Topos et
Cohomologie Etale des Sch\'emas}.
\newline
(i)\, {\sc Lecture Notes in Math.} {\bf 269} Springer, Heidelberg (1972-73).
\newline
(ii)\, {\sc Lecture Notes in Math.} {\bf 270} Springer, Heidelberg (1972-73).
\newline
(iii)\, {\sc Lecture Notes in Math.} {\bf 305} Springer, Heidelberg (1972-73).

\bibitem{Auslander}
M. Auslander. Coherent functors. {\em Proc. Conf. Categorical Algebra (La Jolla,
1965)} Springer (1966), 189-231.

\bibitem{Baez}
J. Baez. An introduction to $n$-categories. {\em Category theory and computer
science (Santa Margherita Ligure 1997)} {\sc Lect.  Notes in Comp. Sci.}
{\bf 1290}, Springer (1997), 1-33.

\bibitem{BaezDolan}
J. Baez, J. Dolan.
\newline
(i)\, Higher dimensional algebra and topological quantum field
theory. {\em J. Math. Phys.} {\bf 36} (1995), 6073-6105 (preprint q-alg/9503).
\newline
(ii)\,  Higher-dimensional algebra III: $n$-categories and the
algebra of opetopes. {\em Adv. Math.} {\bf 135} (1998), 145-206 (preprint
q-alg/9702014).
\newline
(iii)\, Categorification. Preprint math/9802029.


\bibitem{Batanin}
M. Batanin.
\newline
(i)\, Coherent categories with respect to monads and coherent prohomotopy
theory. {\em Cahiers Top. G\'eom. Diff.} {\bf 34} (1993), 279-304.
\newline
(ii)\, Homotopy coherent category theory and $A_{\infty}$-structures in monoidal
categories. {\em J. Pure Appl. Alg.} {\bf 123} (1998), 67-103.
\newline
(iii)\, Categorical strong shape theory. {\em Cahiers Top. G\'eom. Diff.}, to
appear.
\newline
(iv)\, Monoidal globular categories as a natural environment for the
theory of weak $n$-categories. {\em Adv. Math.} {\bf 136} (1998), 39-103.



\bibitem{Berthelot}
P. Berthelot. {\em Cohomologie cristalline des sch\'emas en caract\'eristique
$p>0$.} {\sc Lecture Notes in Math.} {\bf 407} Springer (1974).


\bibitem{BierstoneMilman}
E. Bierstone, P. Milman. Canonical desingularization in characteristic zero by
blowing up the maximum strata of a local invariant. {\em Inventiones} {\bf 128}
(1997), 207-302.



\bibitem{Borsuk}
K. Borsuk. Concerning homotopy properties of compacta. {\em Fund. Math.}
{\bf 62}
(1968), 223-254. (See {\sc MathSci} for many other papers on shape theory.)


\bibitem{Bousfield}
A. Bousfield.
\newline
(i)\, Homotopy spectral sequences and obstructions.
{\em Isr. J. Math} {\bf 66} (1989).
\newline
(ii)\, On the $p$-adic completions of nonnilpotent spaces. {\em Trans. Amer.
Math. Soc.} {\bf 331} (1992), 335-359.

\bibitem{BousfieldGugenheim}
A. Bousfield, V. K. A. M. Gugenheim. {\em On PL de Rham theory and Rational
Homotopy Type} {\sc Memoirs A.M.S.} {\bf 179} (1976).



\bibitem{BousfieldKan}
A. Bousfield, D. Kan. {\em Homotopy limits, completions and localisations}.
Lecture Notes in Math. {\bf 304}, Springer-Verlag (1972).


\bibitem{Breen}
L. Breen.
\newline
(i)\, Extensions du groupe
additif. {\em Publ. IHES} {\bf 48} (1978).
\newline
(ii)\,  Extensions du groue additif sur le site parfait.
{\em Algebraic Surfaces (Orsay, 1976-1978)}
Springer Lecture Notes in Math. {\bf 868} (1981), 238-262.
\newline
(iii)\, On the classification of $2$-gerbs and $2$-stacks. {\em Ast\'erisque}
{\bf 225}, Soc. Math. de France (1994).

\bibitem{EBrown}
E. Brown. Abstract homotopy theory. {\em Trans. Amer. Math. Soc.} {\bf 119}
(1965), 79-85.

\bibitem{BrownSzczarba}
E. Brown, R. Szczarba.
\newline
(i)\, On real homotopy theory. {\em Algebraic Topology, Arcata (1986)} {\sc
Lecture Notes in Math.} {\bf 1370} Springer (1989), 103-116.
\newline
(ii)\, Continuous cohomology and real homotopy type. {\em Trans. Amer. Math.
Soc.} {\bf 311} (1989), 57-106.
\newline
(iii)\, Continuous cohomology and real homotopy type II. International
conference
on homotopy theory (Marseille-Luminy 1988). {\em Ast\'erisque} {\bf 191} (1990),
45-70.

\bibitem{KBrown}
K. Brown. Abstract homotopy theory and generalized sheaf cohomology. {\em
Trans. A.M.S.} {\bf 186} (1973), 419-458.

\bibitem{BrownGersten}
K. Brown, S. Gersten. {\em Algebraic $K$-theory as generalized sheaf cohomology}
Springer {\em Lecture Notes in Math.} {\bf 341} (1973), 266-292.


\bibitem{Burgos}
J. I. Burgos. A $C^{\infty}$ logarithmic Dolbeault complex. {\em Comp. Math.}
{\bf 92} (1994), 61-68.

\bibitem{CKVW}
A. Carboni, G. Kelly, D. Verity, R. Wood. A $2$-categorical approach to change
of base and geometric morphisms II. {\em Theory and Applications of Categories}
{\bf 4} (1998), 82-136.

\bibitem{Cooke}
G. Cooke. Replacing homotopy actions by topological actions. {\em Trans. Amer.
Math. Soc.} {\bf 237} (1978), 391-406.

\bibitem{Cordier}
J.-M. Cordier.
\newline
(i)\, Sur la notion de diagramme homotopiquement coh\'erent. Third
colloquium on categories, Part VI (Amiens, 1980). {\em Cahiers Top. Geom. Diff.
Cat.} {\bf 23} (1982), 93-112.
\newline
(ii)\,  Comparaison de deux cat\'egories d'homotopie de morphismes
coh\'erents. {\em Cahiers Top. Geom. Diff. Cat.} {\bf 30} (1989), 257-275.
\newline
(iii)\,  Sur les limites homotopiques de diagrammes homotopiquement
coh\'erents. {\em Compositio Math.} {\bf 62} (1987), 367-388.

\bibitem{CordierPorter}
J.-M. Cordier, T. Porter.
\newline
(i) \, Functors between shape categories. {\em J. Pure Appl. Algebra} {\bf 27}
(1983), 1-13.
\newline
(ii)\,  Vogt's theorem on categories of homotopy coherent
diagrams. {\em Math. Proc. Camb. Phil. Soc.} {\bf 100} (1986), 65-90.
\newline
(iii) \, {\em Shape Theory. Categorical methods of
approximation.} Halsted Press, New York (1989).
\newline
(iv) \, Fibrant diagrams, rectifications and a construction
of Loday. {\em J. Pure Appl. Alg.} {\bf 67} (1990), 111-124.
\newline
(v)\, Categorical aspects of equivariant homotopy.
The European Colloquium of Category Theory (Tours, 1994). {\em  Appl. Categ.
Structures} {\bf 4} (1996), 195-212.
\newline
(vi)\, Homotopy-coherent category theory. {\em Trans. Amer.
Math. Soc.} {\bf 349} (1997), 1-54.

\bibitem{DeleanuHilton}
A. Deleanu, P. Hilton. On Kan extensions of cohomology theories and Serre
classes of groups. {\em Fund. Math.} {\bf 97} (1977), 157-176.

\bibitem{DeligneMumford}
P. Deligne, D. Mumford. The irreducibility of the space of curves of a given
genus. {\em Publ. Math. I.H.E.S.} {\bf 36} (1969), 75-110.

\bibitem{Deninger}
C. Deninger. On the $\Gamma$-factors attached to motives. {\em Invent. Math.}
{\bf 104} (1991), 245-263.

\bibitem{Dunn}
G. Dunn. Uniqueness of $n$-fold delooping machines. {\em J. Pure and Appl. Alg.}
{\bf 113} (1996), 159-193.

\bibitem{DHK}
W. Dwyer, P. Hirschhorn, D. Kan. Model categories and more general abstract
homotopy theory: a work in what we like to think of as progress. Preprint.

\bibitem{DwyerKan} W. Dwyer, D. Kan.
\newline
(i)\, Simplicial localizations of categories. {\em J. Pure and Appl.
Algebra} {\bf 17} (1980), 267-284.
\newline
(ii)\, Calculating simplicial localizations. {\em J. Pure and
Appl. Algebra} {\bf 18} (1980), 17-35.
\newline
(iii)\, Function complexes in homotopical algebra. {\em Topology}
{\bf 19} (1980), 427-440.
\newline
(iv)\, Equivalences between homotopy theories of diagrams.
{\em Algebraic Topology and Algebraic $K$-theory}, {\em Annals of Math. Studies}
{\bf 113}, Princeton University Press (1987), 180-205.

\bibitem{DwyerKanSmith}
W. Dwyer, D. Kan, J. Smith. Homotopy commutative diagrams and their
realizations.
{\em J. Pure Appl. Algebra} {\bf 57} (1989), 5-24.


\bibitem{FelixThomas}
Y. F\'elix, J.-C. Thomas. Le $Tor$ diff\'erentiel d'une fibration non
nilpotente. {\em J.P.A.A.} {\bf 38} (1985), 217-233.

\bibitem{Frei}
A. Frei. On categorical shape theory. {\em Cahiers Top. Geom. Diff.} {\bf 17}
(1976), 261-294.


\bibitem{GabrielZisman}
P. Gabriel, M. Zisman. {\em Calculus of fractions and homotopy theory.}
Ergebnisse der Math. {\bf 35}, Springer-Verlag, New York (1967).


\bibitem{Giraud} J. Giraud.
\newline
(i)\, {\em Cohomologie Non Ab\'elienne}, Springer-Verlag (1971).
\newline
(ii)\, {\em M\'ethode de la descente}. Memoire {\bf 2} de la S.M.F. (1964).

\bibitem{Gomez-Tato}
A. G\'omez-Tato.
\newline
(i)\, Th\'eorie de Sullivan pour la cohomologie \`a coefficients
locaux. {\em Trans. Amer. Math. Soc.} {\bf 330} (1992), 295-305.
\newline
(ii)\, Mod\`eles minimaux r\'esolubles. {\em J. P. A. A.} {\bf 85} (1993),
43-56.

\bibitem{GriffithsMorgan}
P. Griffiths, J. Morgan. {\em Rational Homotopy Theory and Differential Forms.}
{\sc Progress in Math.} Birkhauser (1981).


\bibitem{Grivel}
P. P. Grivel. Formes diff\'erentielles et suites spectrales. {\em Ann. Inst.
Fourier} {\bf 29} (1979), 17-37.

\bibitem{Grothendieck} A. Grothendieck.
\newline
(i)\,  {\em Rev\^etements \'etales et groupe fondamental.}
Springer Lecture Notes in Mathematics {\bf 224}, Heidelberg (1971).
\newline
(ii)\, On the de Rham cohomology of algebraic varieties. {\em Publ. Math.
I.H.E.S.} {\bf 29} (1966), 96-103.
\newline
(iii)\, {\em Dix expos\'es sur la cohomologie des sch\'emas}, North-Holland.
\newline
(iv)\, Pursuing stacks. Preprint, Montpellier.

\bibitem{GuillenNavarro}
F. Guillen, V. Navarro-Aznar.  Un crit\`ere d'extension d'un foncteur
d\'efini sur les sch\'emas lisses, preprint.

\bibitem{Hain} R. Hain.
\newline
(i)\, Mixed Hodge structures on homotopy groups. {\em Bull. Amer. Math. Soc.}
{\bf 14} (1986), 111-114.
\newline
(ii)\, The de Rham homotopy theory of complex algebraic varieties, I. {\em
K-theory} {\bf 1} (1987), 271-324.
\newline
(iii)\, The de Rham homotopy theory of complex algebraic varieties, II. {\em
K-theory} {\bf 1} (1987), 481-497.
\newline
(iv)\, The geometry of the mixed Hodge structure on the fundamental group. {\em
Algebraic Geometry, Bowdoin} {\sc AMS Proceedings of Symposia in Pure Math.}
{\bf 46} (1987), 247-282.
\newline (v)\, The Hodge de Rham theory of relative Malcev completion. {\em
Ann. Sci. E.N.S.} {\bf 31} (1998), 47-92.

\bibitem{Halperin}
S. Halperin. Lectures on minimal models. {\em M\'em. Soc. Math. France}
{\bf 9-10} (1983), 261pp.

\bibitem{HalperinTanre}
S. Halperin, D. Tanr\'e. Homotopie filtr\'ee et fibr\'es $C^{\infty}$.
{\em Illinois Math. J.} {\bf 34} (1990), 284-324.


\bibitem{Hartshorne}
R. Hartshorne. On the de Rham cohomology of algebraic varieties. {\em Publ.
Math. I.H.E.S.} {\bf 45} (1975), 5-99.


\bibitem{Heller}
A. Heller. Homotopy in functor categories. {\em Trans. Amer. Math. Soc.}
{\bf 272} (1982), 185-202.

\bibitem{Hirschhorn}
P. Hirschhorn. {\em Localization of model categories.} Book-preprint, available
at \newline {\tt http://www-math.mit.edu}

\bibitem{Hirschowitz}
A. Hirschowitz. Coh\'erence et dualit\'e sur le gros site de Zariski. {\em
Algebraic curves and projective geometry} (Ballico, Ciliberto eds.),
{\sc Lecture Notes in Math.} {\bf 1389}, (1989), 91-102.

\bibitem{HirschowitzSimpson}
A. Hirschowitz, C. Simpson. Descente pour les $n$-champs. Preprint
math/9807049

\bibitem{Illusie} L. Illusie.
\newline
(i)\,  Le complexe cotangent I, {\sc Lecture Notes in Mathematics}
{\bf 239}, Springer (1971).
\newline
(ii)\,  Le complexe cotangent II, {\sc Lecture Notes in Mathematics}
{\bf 283}, Springer (1972).

\bibitem{Jaffe}
D. Jaffe. Coherent functors, with application to torsion in the Picard group.
{\em Trans. Amer. Math. Soc.} {\bf 349} (1997), 481-527.

\bibitem{Handbook}
I. M. James, ed. {\em Handbook of Algebraic Topology} Elsevier, Amsterdam
(1995).


\bibitem{Jardine} J. F. Jardine.
\newline
(i)\, Simplicial presheaves {\em J.Pure. Appl.Alg.} {\bf 47}
(1987) 35-87.
\newline
(ii)\,
Stable homotopy theory of simplicial presheaves, {\em Can. J. Math.} {\bf 39}
(1987) 733-747.
\newline
(iii)\, Steenrod operations
in the cohomology of simplicial presheaves.
{\em Algebraic K-theory: connections with geometry and topology}
Jardine-Snaith,
eds, Kluwer Academic Publishers (1989), 103-116.
\newline
(iv)\,
J. F. Jardine. Supercoherence.
{\em J. Pure Appl. Algebra} {\bf 75} (1991), 103-194.
\newline
(v)\, Boolean localization, in practice. {\em Doc. Math.} {\bf 1} (1996),
245-275.

\bibitem{Joshua}
R. Joshua. Generalised Verdier duality for presheaves of spectra--I. {\em J. P.
A. A.} {\bf 70} (1991), 273-289.

\bibitem{Joyal}
A. Joyal, Lettre \`a A. Grothendieck (cf \cite{Jardine}).

\bibitem{Kaenders}
R. Kaenders. On de Rham homotopy theory for plane algebraic curves and their
singularities. Thesis, Catholic University of Nijmegen, October 1997.

\bibitem{Kahn}
B. Kahn. La conjecture de Milnor (d'apr\`es Voevodsky). {\sc Seminaire Bourbaki}
{\bf 834}, juin 1997; {\em Ast\'erisque} {\bf 245} (1997), 379-418.

\bibitem{Katz}
N. Katz. Nilpotent connections and the monodromy theorem. Applications of a
result of Turritin. {\em Publ. Math. I.H.E.S.} {\bf 39} (1971), 175-232.

\bibitem{KatzOda}
N. Katz, T. Oda. On the differentiation of de Rham cohomology classes with
respect to parameters. {\em J. Math. Kyoto Univ.} {\bf 8} (1968), 199-213.


\bibitem{Kelly}
G. M. Kelly, {\em Basic Concepts of Enriched Category Theory},
{\em Springer Lecture Notes in Math.} {\bf 64}, Springer-Verlag.

\bibitem{LaumonMoretBailly}
G. Laumon, L. Moret-Bailly. Champs alg\'ebriques. Preprint {\bf 42}, Orsay
(1992).

\bibitem{Leitch}
R. Leitch. The homotopy-commutative cube. {\em J. London Math. Soc.} {\bf 9}
(1974/75), 23-29.

\bibitem{LemaireSigrist}
J.-M. Lemaire, F. Sigrist. D\'enombrement de types d'homotopie rationnelle. {\em
C.R.A.S.} {\bf 287} (1978), A109--A112.

\bibitem{Lewis}
L. G. Lewis. Is there a convenient category of spectra?
{\em J. Pure Appl. Algebra} {\bf 73} (1991), 233-246.

\bibitem{MardesicSegal}
S. Mardesic, J. Segal. Shapes of compacta and ANR-systems. {\em Fund.
Math.} {\bf
72} (1971), 41-59. (See {\sc MathSci} for many other papers on shape theory.)

\bibitem{Mather}
M. Mather. Pull-backs in homotopy theory. {\em Canadian J. Math.}
{\bf 28} (1976), 225-263.


\bibitem{MorelVoevodsky}
F. Morel, V. Voevodsky. ${\bf A}^1$-homotopy theory of schemes. K-theory
preprint
archives, Oct. 1998 (\cite{Voevodsky} refers to a version dating from 1996 but
\cite{Kahn} cites this as ``in preparation'' in 1997).

\bibitem{Morgan}
J. Morgan. The algebraic topology of smooth algebraic varieties. {\em Publ.
Math. I.H.E.S.} {\bf 48} (1978), 137-204.

\bibitem{Mumford}
D. Mumford. {\em Abelian Varieties} Tata Institute of Fundamental Research,
Studies in Mathematics {\bf 5}, Oxford University Press, London (1970).


\bibitem{Navarro-Aznar} V. Navarro-Aznar.
\newline
(i)\, Sur la th\'eorie de Hodge-Deligne, {\em Inventiones Math.}
{\bf 90}, (1987), 11-76.
\newline
(ii)\, Sur la connexion de Gauss-Manin en homotopie rationnelle. {\em Ann. Sci.
E.N.S. (4)} {\bf 26} (1993), 99-148.

\bibitem{NeisendorferTaylor}
J. Neisendorfer, L. Taylor. Dolbeault homotopy theory. {\em Trans. Amer. Math.
Soc.} {\bf 245} (1978), 183-210.

\bibitem{Porter} T. Porter.
\newline
(i)\, Coherent prohomotopical algebra.
{\em Cahiers Top. Geom. Diff. Cat.} {\bf 18} (1977), 139-179.
\newline
(ii)\,  Coherent prohomotopy theory.
{\em Cahiers Top. Geom. Diff. Cat.} {\bf 19} (1978), 3-46.

\bibitem{Quillen}
D. Quillen.
\newline
(i)\, {\em Homotopical algebra}. Springer Lecture Notes in Mathematics
{\bf 43} (1967).
\newline
(ii)\, Rational Homotopy Theory. {\em Ann. Math.} {\bf 90} (1969), 205-295.

\bibitem{Reedy}
C. Reedy. Homotopy theory of model categories. Preprint (1973) available from P.
Hirschhorn.


\bibitem{Segal}
G. Segal. Categories and cohomology theories. {\em Topology} {\bf 13}
(1974), 293-312.

\bibitem{SchwanzlVogt}
R. Schw\"{a}nzl, R. Vogt. Coherence in homotopy group actions. {\em
Transformation Groups (Poznan, 1985)}, {\em Springer L.N.M.} {\bf 1217} (1986),
364-390.

\bibitem{Siegel}
J. Siegel. Cech extensions and localization of homotopy functors. {\em Fund.
Math.} {\bf 108} (1980), 159-170.

\bibitem{Simpson}
C. Simpson.
\newline
(i)\, Nonabelian Hodge theory. {\em Proceedings of ICM-90}, Springer (1991),
746-756.
\newline
(ii)\, Subspaces of moduli spaces of rank one local systems. {\em Ann. E.N.S.}
{\bf 26} (1993), 361-401.
\newline
(iii)\, Moduli of representations of the fundamental group of a smooth
projective
variety, I: {\em Publ. Math. I.H.E.S.} {\bf 79} (1994), 47-129;
II: {\em Publ. Math. I.H.E.S.} {\bf 80} (1994), 5-79.
\newline
(iv)\, Flexible sheaves. Preprint (1993), later posted as alg-geom 9608025.
\newline
(v)\, Homotopy over the complex numbers and generalized de Rham
cohomology. {\em Moduli of Vector Bundles (Taniguchi symposium December
1994)}, M.
Maruyama, ed., Dekker (1996), 229-263. Toulouse preprint no. 50, April
1995.
\newline
(vi)\, The Hodge filtration on nonabelian cohomology. {\em Algebraic
Geometry (Santa Cruz, 1995)}. A.M.S. Proceedings of Symposia in Pure
Mathematics {\bf 62}, Part 2 (1997), 217-281.
\newline
(vii)\, On a relative notion of algebraic Lie group and applications
to $n$-stacks, alg-geom 9607002.
\newline
(viii)\,  The topological realization of a simplicial presheaf,  q-alg 9609004.
\newline
(ix)\, Algebraic (geometric) $n$-stacks, alg-geom 9609014.
\newline
(x)\, A closed model structure for $n$-categories, internal $Hom$,
$n$-stacks and generalized Seifert-Van Kampen, alg-geom 9704006.
\newline
(xi)\, Limits in $n$-categories, alg-geom 9708010.
\newline
(xii)\, Secondary Kodaira-Spencer classes and nonabelian Dolbeault cohomology,
alg-geom 9712020.
\newline
(xiii)\, On the Breen-Baez-Dolan stabilization hypothesis for Tamsamani's weak
$n$-categories, math.CT/9810058.


\bibitem{Stasheff}
J. Stasheff.
\newline
(i)\, Homotopy associativity of $H$-spaces. {\em Trans. Amer. Math.
Soc.} {\bf 108} (1963), 275-312.
\newline
(ii)\, Grafting Borardman's cherry trees to quantum field theory.
Preprint math/9803156.

\bibitem{Street}
R. Street. Notes on descent theory. Oberwolfach, October 1995.

\bibitem{Sullivan}
D. Sullivan. Infinitesimal computations in topology. {\em Publ. Math. I.H.E.S.}
{\bf 47} (1977), 269-331.

\bibitem{SuslinVoevodsky}
A. Suslin, V. Voevodsky. Singular homology of abstract algebraic varieties.
K-theory preprint archives Nov. 1994.



\bibitem{Tamsamani}
Z. Tamsamani.  Sur des notions de $n$-cat\'egorie et $n$-groupoide
non-strictes via des ensembles multi-simpliciaux. Thesis,
Universit\'e Paul Sabatier (Toulouse). Available as alg-geom 9512006
and alg-geom 9607010. To appear, {\em $K$-Theory}.

\bibitem{Tanre}
D. Tanr\'e.
\newline
(i)\, {\em Homotopie rationnelle: mod\`eles de Chen, Quillen, Sullivan.}
{\sc Lecture Notes in Math.} {\bf 1025}, Springer (1983).
\newline
(ii)\, Mod\`ele de Dolbeault et fibr\'e holomorphe. {\em J.P.A.A.}
{\bf 91} (1994), 333-345.

\bibitem{CTeleman}
C. Teleman. Borel-Weil-Bott theory on the moduli stack of $G$-bundles over a
curve. {\em Inv. Math.} {\bf 134} (1998), 1-57.

\bibitem{Thomason}
R. Thomason.
\newline
(i) \, Algebraic K-theory and etale cohomology. {\em Ann. Sci. E.N.S.} {\bf
13} (1985), 437-552.
\newline
(ii)\, Uniqueness of delooping machines.
{\em Duke Math. J.} {\bf 46} (1979), 217-252.

\bibitem{Voevodsky}
V. Voevodsky. The Milnor conjecture. K-theory preprint archives, Dec. 1996.

\bibitem{Vogt}
R. Vogt.
\newline
(i)\, Homotopy limits and colimits. {\em Math. Z.} {\bf 134} (1973), 11-52.
\newline
(ii)\, Commuting homotopy limits. {\em Math. Z.} {\bf 153} (1977), 59-82.

\bibitem{Wiebel}
C. Wiebel. The mathematical enterprises of Robert Thomason. {\em Bull. A.M.S.}
{\bf 34} (1997), 1-13.

\bibitem{Wojtkowiak}
Z. Wojtkowiak. Cosimplicial objects in algebraic geometry. {\em Algebraic
$K$-theory and algebraic topology} {\sc NATO Adv. Sci. Inst. Ser. C. Math. Phys.
Sci.} {\bf 407}, Kluwer (1993), 287-327.



\end{thebibliography}
\end{document}